\documentclass[a4paper, 12pt]{article2}

\usepackage{amsmath} 
\usepackage{theorem}
\usepackage{amssymb}
\usepackage{amsfonts}
\usepackage{latexsym}
\usepackage{amscd}
\usepackage{diagramb}
\input GrCalc3.sty
\usepackage{graphicx}
\usepackage{xspace,colortbl}
\usepackage{color}
\usepackage{indentfirst}

\newtheorem{theorem}{Theorem}[section]
\newtheorem{lemma}[theorem]{Lemma}
\newtheorem{proposition}[theorem]{Proposition}

\newtheorem{definition}[theorem]{Definition}
\newtheorem{corollary}[theorem]{Corollary}
\newtheorem{remark}[theorem]{Remark}
\newtheorem{problem}[theorem]{Problem}

\newcommand{\qed}{\hfill\quad\fbox{\rule[0mm]{0,0cm}{0,0mm}}\par\bigskip}

\newcommand{\eval}{\operatorname
	{
	\gcn{1}{1}{1}{3}\gcn{1}{1}{3}{1} \gnl
	\gvac{1} \gmp{\ev} \gnl
	}
	}

\newcommand{\crtaeval}{\operatorname
	{
	\gcn{1}{1}{1}{3}\gcn{1}{1}{3}{1} \gnl
	\gvac{1} \gmp{\tilde\ev} \gnl
	}
	}	

\newcommand{\C}{{\mathcal C}}

\newcommand{\F}{{\mathcal F}}
\newcommand{\G}{{\mathcal G}}

\newcommand{\U}{{\mathcal U}}
\newcommand{\R}{{\mathcal R}}
\newcommand{\M}{{\mathcal M}}
\newcommand{\YD}{{\mathcal YD}}

\def\Z{{\mathbb Z}}

\newcommand{\crta}{\overline}
\newcommand{\ev}{\it ev}
\newcommand{\db}{\it db}
\newcommand{\can}{\it can}
\newcommand{\ot}{\otimes}
\newcommand{\Fi}{\varphi}
\newcommand{\iso}{\cong}
\newcommand{\x}{\mbox{-}}
\newcommand{\id}{\it id}
\newcommand{\Id}{\operatorname {Id}}
\newcommand{\Bcgo}{\operatorname {B}}
\newcommand{\Br}{\operatorname {Br}}
\newcommand{\BM}{\operatorname {BM}}
\newcommand{\BW}{\operatorname {BW}}
\newcommand{\Gal}{\operatorname{Gal}}

\newcommand{\Ker}{\operatorname{Ker}}
\newcommand{\Bi}{\operatorname{Im}}

\newcommand{\pH}{\hspace{-0,06cm}\cdot_{\hspace{-0,1cm}_H}\hspace{-0,1cm}}
\newcommand{\pB}{\hspace{-0,06cm}\cdot_{\hspace{-0,1cm}_B}\hspace{-0,1cm}}

\newcommand{\tr}{\hspace{-0,08cm}\triangleright}
\newcommand{\tm}{\times}

\newcommand{\Hc}{\operatorname{H^2}}

\newcommand{\teta}{\theta}
\newcommand{\Teta}{\Theta}
\newcommand{\Epsilon}{\varepsilon}
\newcommand{\Ent}{\mathbb{Z}}

\definecolor{verde}{rgb}{0.,0.7,0.}
\definecolor{turq}{rgb}{.804, .804, .196}
\definecolor{caqui}{rgb}{.741, .718, .42}
\definecolor{rojoo}{rgb}{.980,.502,.447}
\definecolor{rojo}{rgb}{1,0,0}
\definecolor{negro}{rgb}{0,0,0}
\definecolor{section1}{rgb}{.957,.643,.376}

\begin{document}

\title{{\bf A sequence to compute the Brauer group of certain quasi-triangular Hopf algebras}}
\author{{\large J. Cuadra}\thanks{Corresponding author. Email addresses: jcdiaz@ual.es and bfemic@mi.sanu.ac.rs.} \ and {\large B. Femi\'c} \vspace{-2pt} \\
{\normalsize Universidad de Almer\'{\i}a} \vspace{-2pt}\\
{\normalsize Dpto. \'{A}lgebra y An\'{a}lisis Matem\'atico} \vspace{-2pt}\\
{\normalsize E-04120 Almer\'{\i}a, Spain\vspace{-0.2cm}}}
\date{{\small {\it To Freddy Van Oystaeyen on the occasion of his 65th birthday}}}
\maketitle

\vspace{-0.5cm}
\indent \hspace{0.3cm} {\small 2000 Mathematics Subject Classification: 18D10, 18D35, 16W30.}

\vspace{-0.0cm}
\begin{abstract}
A deeper understanding of recent computations of the Brauer group of Hopf algebras is attained by explaining why a direct product decomposition for this group holds and describing the non-interpreted factor occurring in it. For a Hopf algebra $B$ in a braided monoidal category $\C$, and under certain assumptions on the braiding (fulfilled if $\C$ is symmetric), we construct a sequence for the Brauer group $\BM(\C;B)$ of $B$-module algebras, generalizing Beattie's one. It allows one to prove that $\BM(\C;B) \cong \Br(\C) \times \Gal(\C;B),$ where $\Br(\C)$ is the Brauer group of $\C$ and $\Gal(\C;B)$ the group of $B$-Galois objects. We also show that $\BM(\C;B)$ contains a subgroup isomorphic to $\Br(\C) \times \Hc(\C;B,I),$ where $\Hc(\C;B,I)$ is the second Sweedler cohomology group of $B$ with values in the unit object $I$ of $\C$. These results are applied to the Brauer group of a quasi-triangular Hopf algebra that is a Radford biproduct $B \times H$, where $H$ is a usual Hopf algebra over a field $K$, the Hopf subalgebra generated by the quasi-triangular structure $\R$ is contained in $H$ and $B$ is a Hopf algebra in the category ${}_H\M$ of left $H$-modules.  The Hopf algebras whose Brauer group was recently computed fit this framework. We finally show that $\BM(K,H,\R) \times \Hc({}_H\M;B,K)$ is a subgroup of the Brauer group $\BM(K,B \times H,\R),$ confirming the suspicion that a certain cohomology group of $B \times H$ (second lazy cohomology group was conjectured) embeds into $\BM(K,B \times H,\R).$ New examples of Brauer groups of quasi-triangular Hopf algebras are computed using this sequence.
\end{abstract}

\setcounter{tocdepth}{1}

\renewcommand{\contentsname}{\centerline{{\small Contents }}}

\vspace{-0.9cm}
\tableofcontents
\newpage

\renewcommand{\theequation}{\thesection.\arabic{equation}}

\section{Introduction}

 {\bf 1.1} The Brauer group of a field was introduced in 1929 in order to classify finite dimensional division algebras. Its elements are equivalence classes of central simple algebras with the equivalence relation defined in such a way that each equivalence class corresponds to a unique central division algebra. Brauer group theory has strong connections to Galois theory, Homology theory, Group theory, Algebraic Geometry and K-Theory.  \par \medskip

{\bf 1.2} Several versions and generalizations of the Brauer group of a field were proposed through the last century, ending, from a formal point of view, in the construction of the Brauer group of a braided monoidal category by Van Oystaeyen and Zhang \cite{VZ1}. Their construction extends Pareigis' one for a symmetric monoidal category \cite{P4}. An exposition on the way leading to this general construction is offered in \cite{JCD}. A startling use of this Brauer group in Mathematical Physics is pointed out in \cite[Page 408]{FRS}, where it is proposed as a tool to classify full conformal field theories based on a modular tensor category. We mention two of the most important generalizations of the Brauer group of a field. \par \smallskip

(1) It was extended by Auslander and Goldman in \cite{AG} to commutative rings. Central simple algebras are replaced by
central separable algebras, usually called Azumaya algebras. This group plays a relevant role in Algebraic Geometry. \par \smallskip

(2) Wall proposed in \cite{W} the Brauer group of $\Ent_2$-graded algebras, known today as the Brauer-Wall group. The motivation for this generalization is the observation that the Clifford algebra associated to quadratic vector space $V$ is central simple if and only $V$ is even dimensional. However, considering the natural $\Ent_2$-gradation on them, they are always $\Ent_2$-central simple. The Brauer-Wall group is related to Clifford algebras, quadratic extensions and $K$-Theory, see the excellent monograph \cite{Lam}  for a comprehensive account. \par \smallskip

For a quasi-triangular Hopf algebra $(H,\R)$ over a field $K$ we denote by $\BM(K,H,R)$ the Brauer group of the (braided monoidal) category ${}_H\M$ of left $H$-modules, where the braiding stems from the quasi-triangular structure $\R$. This Brauer group, that generalizes Wall's one and introduced in \cite{CVZ2}, is the main object of study of this paper. \par \medskip

{\bf 1.3} Let $K$ be a field with $char(K)\not=2$. Recall that Sweedler Hopf algebra is $H_4=K\langle g,x\vert g^2=1, x^2=0, gx=-xg\rangle$ as an algebra, where $g$ is group-like and $x$ is $(g, 1)$-primitive, with antipode $S(g)=g$ and $S(x)=gx$. In \cite{VZ2} the group $\BM(K,H_4,\R_0)$ corresponding to the quasi-triangular structure $\R_0=\frac{1}{2}(1\ot 1+g\ot 1+1\ot g-g\ot g)$ was computed. It was proved that $$\BM(K,H_4, \R_0)\iso \BW(K) \times (K,+)$$  where $\BW(K)$ denotes the Brauer-Wall group of $K$ and $(K, +)$ the additive group of $K$. Other computations generalizing this one were done later: for Radford Hopf algebra \cite{CC1}, Nichols Hopf algebra \cite{CC2}, and a modified supergroup algebra \cite{Car2}. In all these computations a direct product decomposition for the Brauer group holds and there was no interpretation for one of the factors ($(K,+)$ in the above case) appearing in it. \par \medskip

{\bf 1.4} In this paper we get an insight into these calculations by explaining why a direct product decomposition holds and describing the non-interpreted factor occurring in it. All these Hopf algebras are examples of Radford biproduct (indeed bosonization) Hopf  algebras $B \times H,$ where $(H, \R)$ is a quasi-triangular Hopf algebra, $B \in {}_H\M$ is a braided Hopf algebra and $\R$ is also a quasi-triangular structure for $B \times H$. Here $H$ is viewed as a Hopf subalgebra of $B \times H$ through the canonical inclusion map $\iota: H \rightarrow B \times H$. Radford biproduct construction is an important algebraic link between ordinary Hopf algebras and braided Hopf algebras. We will construct a sequence in a braided monoidal category that will allow us to relate the Brauer group of the braided Hopf algebra $B$ and that of the base category ${}_H\M$ and show that this sequence underlies behind the above-mentioned computations. To explain this we analyze in detail the preceding example. \par \smallskip

Let $\C$ be the category of $K\Z_{2}$-modules with braiding stemming from $\R_0$. It may be identified with the category of $\Z_2$-graded vector spaces with its non-trivial braided structure. Sweedler Hopf algebra is a Radford biproduct, more precisely, a bosonization of $H=K\Z_{2}$ and the braided Hopf algebra $B=K[x]/(x^2)$ in $\C$. The quasi-triangular structure $\R_0$ of $K\Z_2$ is also a quasi-triangular structure for $H_4$. This means that the braiding of $\C$ is $B$-linear and hence the category $_B\C$ of left $B$-modules in $\C$ is also a braided monoidal category with the same braiding (Proposition \ref{Radf extends}). The Brauer group of $\C$ is $\BW(K)$ and the Brauer group of $_B\C$ is  $\BM(K,H_4,\R_0)$ because $_B\C$ and $_{H_4}\M$ are isomorphic as braided monoidal categories (Corollary \ref{braid monoid iso Radf}). By forgetting the $B$-module structure on a $B$-Azumaya algebra, we obtain an Azumaya algebra in $\C$, since the braiding of both categories is the same. Therefore we may consider the following forgetting map $p$ in (\ref{parmor}), that splits by $q$ -- any Azumaya algebra in $\C$ can be equipped with the trivial $B$-module
structure obtaining a $B$-Azumaya algebra:
\begin{equation}\label{parmor}
\bfig
\putmorphism(0, 0)(1, 0)[\Br(\C)`\Br({}_B\C).`q]{600}1a
\putmorphism(0, -50)(-1, 0)[\phantom{\Br(\C)}`\phantom{\Br({}_B\C).}`p]{600}{-1}b
\efig
\end{equation}
\par \smallskip

{\bf 1.5} This situation led us to Beattie's exact sequence \cite{B} where a similar couple of morphisms is considered in the specific case of a finitely generated and projective commutative and cocommutative Hopf algebra $H$ over a commutative ring $R$. Beattie proved that there is a split exact sequence
$$
\bfig
\putmorphism(0, 0)(1, 0)[1`\Br(R)`]{400}1a
\putmorphism(400, 0)(1, 0)[\phantom{\Br(R)}`\BM(R;H)`q]{600}1a
\putmorphism(325, -35)(-1, 0)[\phantom{\BM(R;H)}`\phantom{\Br(R)}`p]{550}{-1}b
\putmorphism(1000, 0)(1, 0)[\phantom{\BM(R;H)}`\Gal(R;H)`\Pi]{700}1a
\putmorphism(1700, 0)(1, 0)[\phantom{\Gal(R;H)}`1.`]{400}1a
\efig
$$
Here $\BM(R;H)$ is the Brauer group of $H$-module algebras, $\Br(R)$ the Brauer group of $R$
and $\Gal(R;H)$ denotes the group of $H$-Galois objects. In our categorical context $\C$ is now the (symmetric) category of $R$-modules, $H$ is a finite Hopf algebra in $\C$ and $_H\C$ is the (symmetric) category of left $H$-modules. The (abelian) Brauer groups $\Br(\C)$ and $\Br({}_H\C)$ are $\Br(R)$ and $\BM(R;H)$ respectively. Our idea was to extend Beattie's exact sequence to any braided monoidal category and to compute $Coker(q)$ in (\ref{parmor}) using this new sequence, explaining so the above computations. The weird factor appearing in the direct product decomposition is the group of Galois objects of the braided Hopf algebra. In our example $(K,+)$ should be the group of $B$-Galois objects. Beattie's sequence was constructed in Fern\'andez Vilaboa's Ph.D. thesis \cite{FV} for a closed symmetric monoidal category with equalizers and coequalizers.  This construction also appears in \cite{AV2}, although proofs are only sketched and for complete proofs one is referred to \cite{FV}. For computations the author uses here tapestry diagrams which are less intuitive and spread than the standard graphical calculus. The case of Radford Hopf algebra escapes from this construction because the base category is not symmetric (see examples in Section 6).\par \medskip

{\bf 1.6} Motivated by this case  we generalize this construction to a braided monoidal category under some assumptions on the braiding, fulfilled in any symmetric monoidal category. Our first main result is:

\begin{theorem}\label{introth1}
Let $(\C,\otimes,\Phi)$ be a closed braided monoidal category with equalizers and coequalizers.
Let $H$ be a finite and commutative Hopf algebra in $\C$. Suppose that the braiding $\Phi$ is $H$-linear and that $\Phi_{T,X}=\Phi^{-1}_{X,T}$ for any $H$-Galois object $T$ and any $X \in \C$. Then there is a split exact sequence
\begin{equation}\label{introbeattie1}
\bfig
\putmorphism(0, 0)(1, 0)[1`\Br(\C)`]{400}1a
\putmorphism(400, 0)(1, 0)[\phantom{\Br(\C)}`\BM(\C;H)`q]{600}1a
\putmorphism(410, -50)(1, 0)[\phantom{\Br(\C)}`\phantom{\Br(\C)}`p]{470}{-1}b
\putmorphism(1000, 0)(1, 0)[\phantom{\BM(\C;H)}`\Gal(\C;H)`\Upsilon]{700}1a
\putmorphism(1700, 0)(1, 0)[\phantom{\Gal(\C;H)}`1.`]{400}1a
\efig
\end{equation}
Here $\BM(\C;H)$ denotes the Brauer group of $H$-module algebras, $\Br(\C)$ the Brauer group of $\C$ and $\Gal(\C;H)$ the group of $H$-Galois objects. Moreover, $\BM(\C;H)\iso \Br(\C)\times \Gal(\C;H).$
\end{theorem}

Some words must be said about the hypotheses. That $\Phi$ is $H$-linear ensures that the category $_H\C$ of left $H$-modules is braided and hence we can consider its Brauer group $\BM(\C;H)$.  It also implies that $H$ is cocommutative (Proposition \ref{braid-lin}). This is why cocommutativity does not explicitly appear in our list of assumptions, unlike in \cite{AV2}. The symmetricity condition $\Phi_{T, X}=\Phi_{X,T}^{-1}$ for any $H$-Galois object $T$ and any $X\in\C$ (i.e. $T$ belongs to the M\"uger's center of $\C$) is needed in order to construct the morphism $\Upsilon$ in (\ref{introbeattie1}) and the group of $H$-Galois objects $\Gal(\C;H)$. This group is constructed in \cite{CF} based on the construction of the group of biGalois objects due to Schauenburg \cite{Sch1}. \par \medskip

{\bf 1.7} The symmetricity condition on the braiding, although satisfied by our examples, seems rare.  If we drop it, we still have a similar sequence for the subgroup $\BM_{inn}(\C; H)$ of $\BM(\C; H)$ consisting of $H$-Azumaya algebras with inner actions. Our second main result asserts:
\begin{theorem}\label{introth2}
Let $\C$ be a closed braided monoidal category with equalizers and coequalizers. Let $H$ be a finite and commutative Hopf algebra in $\C$. Assume that the braiding is $H$-linear. Then there is a split exact sequence
\begin{equation}\label{introbeattie2}
\bfig
\putmorphism(0, 0)(1, 0)[1`\Br(\C)`]{400}1a
\putmorphism(400, 10)(1, 0)[\phantom{\Br(\C)}`\BM_{inn}(\C;H)`q]{680}1a
\putmorphism(305, -35)(-1, 0)[\phantom{\BM(R;H)}`\phantom{\Br(R)}`p]{600}{-1}b
\putmorphism(1080, 0)(1, 0)[\phantom{\BM_{inn}(\C;H)}`\Hc(\C;H,I)`\Upsilon']{800}1a
\putmorphism(1900, 0)(1, 0)[\phantom{\Hc(\C;H,I)}`1`]{500}1a
\efig
\end{equation}
 where $\Hc(\C;H,I)$ is the second Sweedler cohomology group of $H$ with values in the unit object $I$. Furthermore, $\BM_{inn}(\C;H)\iso \Br(\C)\times \Hc(\C;H,I).$
\end{theorem}
This sequence generalizes to {\it any} braided monoidal category the one obtained by Alonso \'Alvarez and Fern\'andez Vilaboa in \cite[Theorem 11]{AV1} and \cite[Proposition 0.3]{AV3} for symmetric monoidal categories. \par \smallskip

We remark that our proof of these two theorems is different in several aspects from the one of these authors for the symmetric case: instead of tapestry notation we use braided notation, more intuitive and widely accepted; we use, as Schauenburg does in \cite{Sch1}, the notion of flatness and faithful flatness in a monoidal category, that allows to simplify the proof of many intermediate results; finally we use Schauenburg's observation from \cite{Sch} that $\Phi_{H,H}=\Phi_{H,H}^{-1}$ for a cocommutative Hopf algebra $H$, allowing one to construct in any braided monoidal category the group $\Gal_{nb}(\C;H)$ of Galois objects with a normal basis, that turns out to be isomorphic to $\Hc(\C;H,I).$ \par \medskip

{\bf 1.8} Our two main results are applied to a certain family of finite dimensional quasi-triangular Hopf algebras that we next describe. Let $(H, \R)$ be a quasi-triangular Hopf algebra over a field $K$ and $\C={}_H\M$ the braided monoidal category of left $H$-modules. We denote by $\Phi_\R$ the braiding of $\C$ stemming from $\R$. Let $B\in \C$ be a Hopf algebra and consider the Radford biproduct Hopf algebra $B\times H$. Denote by $\iota:H\to B\times H$ the canonical inclusion and set $\overline{\R}=(\iota \otimes \iota)(\R)$. Then, $\overline{\R}$ is a quasi-triangular structure for $B\times H$ if and only if the braiding $\Phi_\R$ is $B$-linear (Proposition \ref{Radf extends}). The category of left $B\times H$-modules ${}_{B\times H}\M$ is isomorphic, as a braided monoidal category, to ${}_B\C$ (Corollary \ref{braid monoid iso Radf}). Theorem \ref{introth1} applied to this case states:

\begin{corollary}\label{introcor1}
With notation as above, if $\Phi_\R$ satisfies the symmetricity condition for any $B$-Galois object $T\in {}_H\M$ and any $X\in {}_H\M$, we have
$$\BM(K, B\times H, \overline{\R})\iso\BM(K, H, \R)\times\Gal({}_H\M; B).$$
\end{corollary}
We show that this result underlies in the computation of the Brauer group of Sweedler Hopf algebra, Radford Hopf algebra, Nichols Hopf algebra and modified supergroup algebras (Theorem \ref{Beattieroot}). They are better understood in our categorical framework by means of Beattie's exact sequence. Applying Theorem \ref{introth2} to the same type of Radford biproducts one obtains:

\begin{corollary}\label{introcor2}
$\BM(K, H, \R)\times\Hc({}_H\M; B, K)$ is a subgroup of $\BM(K, B\times H, \overline{\R})$.
\end{corollary}
For this result the symmetricity assumption on the braiding is not needed. It was suspected that the group of lazy 2-cocycles would embed in the Brauer group, \cite[Introduction]{BiCar}. The cohomology group $\Hc({}_H\M; B, K)$ embeds in (and coincides in some cases with) the second lazy cohomology group of $B\times H$ for the examples analyzed, as checked in \cite{CF}. Thus our result clarifies this suspicion. \par \medskip

{\bf 1.9} This paper is organized as follows. In the second section we fix notation and present some preliminaries including notions like (faithful) flatness, closed monoidal categories, inner hom objects, finite and dual objects and algebraic structures in braided monoidal categories. The Morita Theorems for categories of modules constructed over a monoidal category, due to Pareigis, are recalled. \par \medskip

The third section deals with the first and second term of our short exact sequence (\ref{introbeattie1}). We recall from \cite{VZ1} the construction of the Brauer group of a braided monoidal category $\C$. We give an alternative description of one of the functors associated to an Azumaya algebra (Proposition \ref{Az-right-adj}). When the braiding is $H$-linear the Brauer group $\BM(\C;H)$ is defined as the Brauer group of the category ${}_H\C$ of left $H$-modules. Finally, we construct the subgroup $\BM_{inn}(\C;H)$ of $\BM(\C;H)$ consisting on those classes containing a representative which has inner $H$-action. \par \smallskip

The third term of our sequence, the group $\Gal(\C;H)$ of $H$-Galois objects, together with some results on braided Hopf-Galois theory, is presented in the fourth section.  This group is constructed in \cite{CF} taking advantage of the group of biGalois objects due to Schauenburg \cite{Sch1}. To construct it we assume that the braiding satisfies $\Phi_{A, B}=\Phi_{B,A}^{-1}$ for any two $H$-Galois objects $A$ and $B$. Galois objects with a normal basis are defined and it is pointed out that the above-mentioned condition on the braiding is satisfied by them. They form a group $\Gal_{nb}(\C;H)$, that exists in any braided monoidal category, and it is isomorphic to $\Hc(\C;H,I)$. \par \smallskip

Theorems \ref{introth1} and \ref{introth2} are stated and proved in several steps in the fifth section. We first define the map $\Upsilon: \BM(\C; H) \rightarrow \Gal(\C;H)$ in (\ref{introbeattie1}) and check that it is a group morphism. Secondly, we assign to any $H$-Galois object an Azumaya algebra, indeed an $H$-Azumaya algebra, showing that $\Upsilon$ is surjective. We finally prove that the sequence is exact and get the direct product decomposition. Theorem \ref{introth2} is obtained as a byproduct by proving that $\Upsilon(\BM_{inn}(\C; H))=\Gal_{nb}(\C;H).$ \par \smallskip

The sixth section contains the applications, Corollaries \ref{introcor1} and \ref{introcor2}, of our two main results to the class of Radford biproducts previously detailed. We show that our sequence lies behind the computations of the Brauer group of Sweedler Hopf algebra, Radford Hopf algebra, Nichols Hopf algebra and modified supergroup algebras. We will go further by computing the Brauer group of a new family of quasi-triangular Hopf algebras (Subsection 6.3, Theorem \ref{bmnd}). They will provide us with an example of a Hopf algebra in a braided (non-symmetric) monoidal category such that every Galois object has a normal basis, and with an example for which this is not true.  \par \bigskip

{\bf Acknowledgment.} The research of the first named author was supported by projects MTM2008-03339 from MEC and FEDER and P07-FQM03128 from Junta de Andaluc\'{\i}a. The second named author was supported by a predoctoral fellowship from the European Marie Curie project 'LIEGRITS', MRTN-CT 2003-505078. She would like to thank the Department of Algebra and Mathematical Analysis of the University of Almer\'{\i}a for its warm hospitality and the facilities put to her disposal. The authors are grateful to professor Pareigis for providing them his LaTex package to draw braided diagrams.

\section{Preliminaries}
\setcounter{equation}{0}

\renewcommand{\theequation}{\thesection.\arabic{equation}}

\newcounter{preli}[section]
\renewcommand{\thepreli}{{\bf \thesubsection.\arabic{preli}}}

\subsection{Braided monoidal categories}
\setcounter{preli}{0}

We assume that the reader is familiar with general category theory. We recommend the monographs \cite{McL2} and \cite{P} as references. In this paper we will deal with braided monoidal categories. The reader is referred to \cite{JS1}, \cite{JS2}, \cite{K}, \cite{KT}, \cite{Maj}, \cite{McL2} and \cite{T} for basic notions and results on this important sort of categories. The $7$-tuple $(\C,\otimes,a,I,l,r,\Phi)$ will stand for a braided monoidal category, where $\C$ is a category, $\otimes:\C \times \C \rightarrow \C$ is the tensor product functor, $a_{X,Y,Z}:(X \otimes Y) \otimes Z \rightarrow X \otimes (Y \otimes Z)$ is the associativity constraint, that satisfies Mac Lane's pentagonal axiom, $I$ denotes the unit object, $l_X:I \otimes X \rightarrow X$ and $r_X:X \otimes I \rightarrow X$ are the left and right unity constraint respectively and $\Phi_{X,Y}:X \otimes Y \rightarrow Y \otimes X$ denotes the braiding. In view of Mac Lane's Coherence Theorem we may (and we will do) assume that our braided monoidal category is strict, i.e., the associativity and unity constraints are the identity in $\C$. We will use the standard graphical calculus to work in braided monoidal categories. For two objects $V,W$ in $\C$ we denote the braiding between them and its inverse  by
$$\Phi_{V, W}=
\gbeg{1}{3}
\got{1}{V}\got{1}{W} \gnl
\gbr \gnl
\gob{1}{W} \gob{1}{V}
\gend\qquad\textnormal{and}\quad
\Phi_{V, W}^{-1}=
\gbeg{1}{3}
\got{1}{W}\got{1}{V} \gnl
\gibr \gnl
\gob{1}{V} \gob{1}{W}
\gend
$$
respectively. \par \medskip

\addtocounter{preli}{1}\thepreli\label{pr11} {\bf Flatness:} An object $A$ in $\C$ is called {\em flat} if the functor $A\ot -: \C\to\C$ preserves equalizers. If, in addition, it reflects isomorphisms, then $A$ is called {\em faithfully flat}. By naturality of the braiding the functor $A\ot -: \C\to\C$ preserves equalizers (resp. reflects isomorphisms) if and only if $-\ot A: \C\to\C$ does it. The following statements for objects $A,B \in \C$ are easy to prove: \vspace{-5pt}
\begin{enumerate}
\itemsep -2pt
\item[(i)] {\it If $A$ and $B$ are flat, then so is $A\ot B$.}
\item[(ii)] {\it If $A$ and $B$ are faithfully flat, then so is $A\ot B$.}
\item[(iii)] {\it If the functor $A\ot -$ reflects equalizers and $A\ot B$ is faithfully flat, then $B$ is faithfully flat.}
\end{enumerate}

\addtocounter{preli}{1}\thepreli\label{pr12} {\bf Closed categories:} A braided monoidal category $\C$ is called {\em closed} if the functor $-\otimes M:\C \rightarrow \C$ has a right adjoint for all $M \in \C$. The right adjoint, called the {\em inner hom functor}, will be denoted by $[M,-]:\C \rightarrow \C$. For $M,N \in \C$ the counit of the adjunction evaluated at $N$ is denoted by $\ev_{M,N}:[M,N]\otimes M \rightarrow N$. It satisfies the following universal property: {\it for any morphism $f:T \otimes M \rightarrow N$ there is a unique morphism $g:T \rightarrow [M,N]$ such that $f=\ev_{M,N}(g \otimes M)$}. The functor $[M, -]$ is defined on morphisms as follows: Let $f: N\to N'$ be a morphism in $\C$. Then $[M, f]:[M, N]\to [M, N']$ is the unique morphism such that the following diagram commutes:
\begin{eqnarray}\label{rightAdj}
\scalebox{0.8}[0.9]{\bfig
\putmorphism(0, 400)(1, 0)[[M, N]\ot M`\phantom{N}`\ev_{M,N}]{900}1a
\putmorphism(60, 400)(0, -1)[\phantom{[M, N]\ot M`}`\phantom{[M, N']\ot M}`{[M, f]}\ot M]{400}1l
\putmorphism(900, 400)(0, -1)[N`\phantom{N'}`f]{400}1r
\putmorphism(0, 0)(1, 0)[[M, N']\ot M`N'`\ev_{M,N'}]{900}1b
\efig}
\end{eqnarray}
The unit of the adjunction $\alpha_{M,N}: N\to [M, N\ot M]$ satisfies:
\begin{eqnarray} \label{adjunit}
\gbeg{3}{5}
\got{1}{N} \got{3}{M} \gnl
\gbmp{\alpha} \gvac{1} \gcl{1} \gnl
\hspace{-0,26cm} $\eval$ \gnl
\gvac{1} \gob{1}{N\ot M}
\gend=
\gbeg{1}{4}
\got{1}{N\ot M} \gnl
\gcl{2} \gnl
\gob{1}{N \ot M}
\gend
\end{eqnarray}
For the rest of this paper, and unless otherwise stated, $\C$ will stand for a {\it closed braided monoidal category with equalizers and coequalizers}. Some of the results presented in this section hold under weaker assumptions on $\C$. However, our assumption will avoid some technical difficulties and will allow to simplify the statement of some of them. \par \medskip

\addtocounter{preli}{1}\thepreli\label{pr13} {\bf Inner hom algebra:} There is an associative pre-multiplication $\Fi_{M, Y, Z}: [Y, Z]\ot [M, Y]\to [M, Z]$ and a unital morphism $\eta_Y: I\to [Y, Y]$ such that $\Fi_{M, Y, Y}(\eta_Y\ot [M, Y])=\Id_{[M, Y]}$ and $\Fi_{M, M, Z}([M, Z]\ot \eta_M)=\Id_{[M, Z]}.$ The morphisms $\Fi_{M, Y, Z}$ and
$\eta_Y$ are given via the universal properties of $[M, Z]$ and $[Y, Y]$, respectively,
by the diagrams:
$$\begin{array}{cc}
\scalebox{0.75}[0.8]{\bfig
\putmorphism(0, 500)(1, 0)[\phantom{[Y, Z] \ot [M, Y] \ot M}`
\phantom{[Y, Z] \ot Y}`{[Y, Z]}\ot \ev_{M,Y}]{1400}1a
\putmorphism(1400, 500)(0, -1)[[Y, Z] \ot Y` Z`\ev_{Y,Z}]{500}1r
\putmorphism(0, 500)(0, -1)[[Y, Z] \ot [M, Y] \ot M`
\phantom{[M, Z] \ot M}` \Fi_{M,Y,Z} \ot M]{480}1l
\putmorphism(0,0)(1, 0)[[M, Z] \ot M` \phantom{Z}`\ev_{M,Z}]{1400}1b
\efig} & \hspace{2cm}\scalebox{0.8}[0.9]{
\bfig
\putmorphism(0, 0)(1, 0)[[Y,Y] \ot Y`Y`ev_{Y,Y}]{800}1b
\putmorphism(0, 400)(2, -1)[``]{800}1r
\putmorphism(0, 430)(2, -1)[``l_Y]{800}0r
\putmorphism(0, 400)(0, -1)[I \ot Y``\eta_Y \ot Y]{380}1l
\efig}
\end{array}$$
By $\varphi_{M,Y,Z}$ being associative we mean that $\varphi_{M,Y,V}(\varphi_{Y,Z,V} \otimes [M,Y])=\varphi_{M,Z,V}([Z,V] \otimes \varphi_{M,Y,Z})$. For each $M\in\C$, {\it the object $[M,M]$ is equipped with an algebra structure via
$\Fi_{M,M,M}$ and $\eta_M$.} \par \medskip

\addtocounter{preli}{1}\thepreli\label{pr14} Using the braiding, $[M,-]:\C \to \C$ becomes a right adjoint of $M\otimes -:\C \to \C$. The counit of this adjunction will be denoted by $\tilde{ev}_{M,N}:M \otimes [M,N] \rightarrow N$ for any $N \in \C$. Then $\tilde{\ev}_{M,N}=\ev_{M,N}\Phi^{-1}_{[M,N],M}$. The universal property for $\tilde{\ev}_{M,N}$ reads as follows: for any $f:M \otimes T \to N$ there is a unique $g:T\to [M,N]$ such that
\begin{eqnarray} \label{basic-adj-l}
\gbeg{3}{5}
\got{1}{M}\got{3}{T} \gnl
\gcl{1} \gvac{1} \gbmp{g} \gnl
\gcn{1}{1}{1}{3}\gcn{1}{1}{3}{1} \gnl
\gvac{1} \gmp{\tilde{\ev}} \gnl
\gvac{1} \gob{1}{N}
\gend=
\gbeg{2}{5}
\got{1}{M\ot T}\gnl
\gcl{1} \gnl
\gbmp{f} \gnl
\gcl{1} \gnl
\gob{1}{N}
\gend
\end{eqnarray}
The unit of the adjunction $\tilde{\alpha}: N\to [M, M\ot N]$ obeys
\begin{eqnarray} \label{adjunit-op}
\gbeg{3}{5}
\got{1}{M} \got{3}{N} \gnl
\gcl{1} \gvac{1} \gbmp{\tilde\alpha} \gnl
\hspace{-0,28cm} $\crtaeval$ \gnl
\gvac{1} \gob{1}{M\ot N}
\gend=
\gbeg{1}{4}
\got{1}{M\ot N} \gnl
\gcl{2} \gnl
\gob{1}{M\ot N}
\gend
\end{eqnarray}
The relation between $\alpha$ and $\tilde{\alpha}$ is $\tilde{\alpha}=[M,\Phi_{N,M}]\alpha.$ When there is no danger of confusion we will simply write $\alpha$ and $\ev$ instead of $\tilde{\alpha}$ and $\tilde{\ev}.$ \par \medskip

\addtocounter{preli}{1}\thepreli\label{pr15} {\bf Dual objects:} Given $P \in \C$, an object $P^*\in\C$ together with a morphism $e_P:P^*\ot P \to I$ is called a {\em dual object} for $P$ if there exists a morphism $d_P:I\to P\ot P^*$ such that $(P \otimes e_P)(d_P \otimes P)=id_P$ and $(e_P \otimes P^*)(P^* \otimes d_P)=id_{P^*}$. The morphisms $e_P$ and $d_P$ are called {\em evaluation} and {\em dual basis} respectively. In braided diagrams $e_P$ and $d_P$ are denoted by
$$e_P=
\gbeg{3}{1}
\got{2}{\hspace{-0,3cm}P^*} \got{1}{\hspace{-0,5cm}P} \gnl
\gev \gnl
\gob{1}{}
\gend
\qquad\textnormal{and}\qquad
d_P=
\gbeg{3}{3}
\got{1}{} \gnl
\gdb \gnl
\gob{1}{P} \gob{2}{\hspace{-0,2cm}P^*}
\gend$$
Then the conditions in the definition take the form:
\begin{center}
\begin{tabular}{p{6cm}p{1cm}p{6cm}}
\begin{eqnarray} \label{dbev=id}
\scalebox{0.9}[0.9]{
\gbeg{3}{4}
\got{1}{} \got{3}{P} \gnl
\gdb  \gcl{1} \gnl
\gcl{1} \gev \gnl
\gob{1}{P}
\gend=\gbeg{1}{4}
\got{1}{P} \gnl
\gcl{2} \gnl
\gob{1}{P}
\gend}
\end{eqnarray} &  &
\begin{eqnarray} \label{evdb=id}
\scalebox{0.9}[0.9]{
\gbeg{3}{4}
\got{1}{P^*} \got{1}{} \gnl
\gcl{1} \gdb \gnl
\gev \gcl{1} \gnl
\gob{1}{} \gob{4}{\hspace{-0,1cm}P^*}
\gend=\gbeg{1}{4}
\got{1}{P^*} \gnl
\gcl{2} \gnl
\gob{1}{P^*}
\gend}
\end{eqnarray}
\end{tabular}
\end{center}
A dual object is unique up to isomorphism. For a dual object $(P^*, e_P)$ for $P$ the functor $-\ot P^*:\C\to \C$ is a right adjoint of $-\ot P:\C\to \C$. Hence $e_P:P^*\ot P \to I$ satisfies the following universal property: for any object $X \in \C$ and any morphism $f:X \otimes P \rightarrow I$ there is a unique morphism $g:X \rightarrow P^*$ such that $f=e_P(g \otimes P)$. Moreover, $f:P \rightarrow Q$ induces $f^*:Q^* \rightarrow P^*$ such that $e_P(f^* \otimes P)=e_Q(Q^* \otimes f)$. Another consequence of the above adjunction is that the functors $[P,-]$ and $-\otimes P^*$ are isomorphic. If $P$ has a dual object, then $P^*$ has it as well and there is a natural isomorphism $P \cong P^{**}$. \par \medskip

\addtocounter{preli}{1}\thepreli\label{pr16} {\bf Tensor product of dual objects:} If $P,Q \in \C$ have dual objects, then $Q^* \otimes P^*$ is a dual object for $P \otimes Q$ with evaluation morphism $e_{P\otimes Q}=e_Q(Q^* \otimes e_P \otimes Q)$. We may also pair $P^* \otimes Q^*$ with $P \otimes Q$ via $e'_{P\otimes Q}=(e_P \otimes e_Q)(P^*\otimes \Phi_{Q^*,P} \otimes Q)$. There is an isomorphism $\theta_{P,Q}:P^* \otimes Q^* \rightarrow (P \otimes Q)^*$, where $\theta_{P,Q}$ is the unique morphism such that $e'_{P\otimes Q}=e_{P\otimes Q}(\theta_{P,Q} \otimes P \otimes Q)$. In the sequel, we will use $e'_{P\otimes Q}$ as evaluation morphism for $P \otimes Q$ and simply denote it by $e_{P\otimes Q}$.\par \smallskip

Using the universal property of $e_{P \otimes Q}$ one may show that  $\Phi_{P^*,Q^*}=\theta^{-1}_{Q,P}\Phi_{Q,P}^*\theta_{P,Q}$. Hence, {\it if $\Phi_{P,Q}=\Phi_{Q,P}^{-1}$, then
$\Phi_{P^*,Q^*}=\Phi_{Q^*,P^*}^{-1}$.} \par \medskip

\addtocounter{preli}{1}\thepreli\label{pr17} {\bf Finite objects:} The morphism $\db: P \ot [P,I] \to [P,P]$ defined via the universal property of $([P,P], \ev_{P,P}:[P,P]\ot P\to P)$ by
\begin{eqnarray} \label{db-def}
\scalebox{0.9}[0.9]{
\bfig
\putmorphism(130, 440)(2, -1)[\phantom{P\ot [P,I]\ot P}`\phantom{P\iso P \ot I}`]{850}1r
\putmorphism(90, 500)(2, -1)[\phantom{P\ot [P,I]\ot P}`\phantom{P\iso P \ot I}`P\ot \ev_{P,I}]{850}0r
\putmorphism(80, 400)(0, -1)[`\phantom{[P,P] \ot P}`\db\ot P]{400}1l
\putmorphism(-50, 400)(0, -1)[P\ot [P,I]\ot P``]{400}0l
\putmorphism(0, 0)(1, 0)[[P,P] \ot P`P\iso P\ot I`]{900}1b
\putmorphism(0, -25)(1, 0)[\phantom{[P,P] \ot P}`\phantom{P\iso P \ot I}`\ev_{P,P}]{900}0b
\efig}
\end{eqnarray}
is also called the {\em dual basis morphism}. The object $P$ is {\em finite} if $\db$ is an isomorphism. In braided diagrams (\ref{db-def}) reads as:
$$\gbeg{3}{5}
\got{1}{P \ot [P,I]}\got{3}{P} \gnl
\gbmp{db} \gvac{1} \gcl{1} \gnl
\hspace{-0,135cm}$\eval$
\gvac{1} \gob{1}{P}
\gend =
\gbeg{3}{4}
\got{1}{P} \got{2}{\hspace{-0,14cm}[P,I]} \got{1}{P} \gnl
\gcl{2} \gcn{1}{1}{1}{3} \gcn{1}{1}{3}{1} \gnl
\gvac{2} \gmp{\ev} \gnl
\gob{1}{P} \gvac{2} \gob{1}{}
\gend
$$
One may easily prove that {\it if $P$ is finite, then $([P,I], \ev_{P,I})$ is its left dual}. The dual basis morphism is $d_P=db^{-1}\eta_{[P,P]}$.
\par \medskip

\addtocounter{preli}{1}\thepreli\label{pr18} {\it In a closed braided monoidal category $\C$ a finite object $P$ is flat.} \par \smallskip

\begin{proof} Since $P$ is finite, $P^*:=[P,I]$ is a left dual for $P$ and $P\iso P^{**}$. Hence $P$ is a left dual to $P^*$. Then $-\ot P^{**}\iso -\ot P$ is a right adjoint functor to $-\ot P^*$ and, as such, $-\ot P$ preserves equalizers.
\qed\end{proof}

\addtocounter{preli}{1}\thepreli\label{product of duals} {\it Let $M$ and $N$ be finite objects in $\C$. There is a natural isomorphism} $$[M\ot N, I]\iso [N, I]\ot [M, I]= N^*\ot M^*.$$

\begin{proof}
There is a natural isomorphism between  $[M\otimes N,-]$ and $[M,[N,-]]$. Combining 2.1.5 and 2.1.7, we find that there is also a natural isomorphism between $[M,[N,-]]$ and $[N,-] \otimes [M,I].$ Hence we have a natural isomorphism between $[M\otimes N,-]$ and $[N,-] \otimes [M,I].$ Apply now this to $I$ and take into account that $[M,I]$ and $[N,I]$ are the left duals of $M$ and $N$ respectively in view of 2.1.7.\qed
\end{proof}

For practical purposes we will involve the braiding in the above isomorphism and consider a new isomorphism
$$\delta_{M, N}:M^*\ot N^* \to [M\ot N, I]$$
defined as the unique one such that
\begin{eqnarray} \label{duals isom delta}
\scalebox{0.9}[0.9]{
\gbeg{4}{5}
\got{2}{M^*\ot N^*} \gvac{1} \got{1}{M\ot N} \gnl
\gvac{1} \glmptb \gnot{\hspace{-0,4cm}\delta_{M, N}} \grmp \gcl{1} \gnl
\gvac{1} \gcn{1}{1}{1}{3}\gcn{1}{1}{3}{1} \gnl
\gvac{2}  \gmp{\ev} \gnl
\gob{1}{}
\gend} =  \scalebox{0.9}[0.9]{
\gbeg{5}{3}
\got{1}{M^*} \got{1}{N^*}  \got{1}{M} \got{1}{N}\gnl
\gcl{1} \gbr \gcl{1} \gnl
\gev \gev \gnl
\gvac{3} \gob{2}{}
\gend}
\end{eqnarray}

\subsection{Structure transmission lemmas}
\setcounter{preli}{0}

We refer to \cite{Besp} for the definition of algebraic structures in braided monoidal categories like algebras, modules, coalgebras, comodules, Hopf algebras, module algebras and comodule algebras. The list of axioms for each of these structures, expressed in graphical calculus, may be found in \cite[Page 159]{Besp}. We set some notation for the different structures: $A$ is an algebra, $C$ a coalgebra, $H$ a Hopf algebra with antipode $S$, $M$ an $A$-module and $N$ a $C$-comodule.
\begin{center}
\begin{tabular}{ccp{0.5cm}cccp{0.5cm}c}
\multicolumn{2}{c}{{\footnotesize Algebra}} & & \multicolumn{2}{c}{{\footnotesize Coalgebra}} & & {\footnotesize Antipode} \\
{\footnotesize Unit} & {\footnotesize Multiplication} & & {\footnotesize Counit} & {\footnotesize Comultiplication} & & \\
$\eta_A=
\gbeg{1}{3}
\got{1}{} \gnl
\gu{1} \gnl
\gob{1}{A}
\gend$ & $\nabla_A=\gbeg{2}{3}
\got{1}{A} \got{1}{A} \gnl
\gmu \gnl
\gob{2}{A}
\gend$ & & $\varepsilon_C=
\gbeg{1}{3}
\got{1}{C} \gnl
\gcu{1} \gnl
\gob{1}{}
\gend$ &  $\Delta_C=
\gbeg{2}{3}
\got{2}{C} \gnl
\gcmu \gnl
\gob{1}{C} \gob{1}{C}
\gend$ & & $\gbeg{2}{5}
\got{1}{H} \gnl
\gcl{1} \gnl
\gmp{S} \gnl
\gcl{1} \gnl
\gob{1}{H}
\gend$ \vspace{5pt} \\
\multicolumn{2}{c}{{\footnotesize Module}} & & \multicolumn{2}{c}{\footnotesize Comodule} &  \\
{\footnotesize Left} & {\footnotesize Right} & & {\footnotesize Left} & {\footnotesize Right} & \\
$\nu_M= \gbeg{3}{3}
\got{1}{A} \got{1}{M} \gnl
\glm \gnl
\gob{3}{M}
\gend$ & $\mu_M=\gbeg{2}{3}
\got{1}{M} \got{1}{A} \gnl
\grm \gnl
\gob{1}{M}
\gend$ & &  $\lambda_N=\gbeg{1}{3}
\got{3}{N} \gnl
\glcm \gnl
\gob{1}{C} \gob{1}{N}
\gend$ &  $\rho_N=
\gbeg{2}{3}
\got{1}{N} \gnl
\grcm \gnl
\gob{1}{N} \gob{1}{C}
\gend$ & &
\end{tabular}
\end{center}
\par \smallskip

In this subsection we will record in several results how the equalizer and coequalizer of two morphisms with additional structures also inherits these structures. The proofs of these results are standard. \par \medskip

\addtocounter{preli}{1}\thepreli\label{pr21}
Let
$\scalebox{0.9}[0.9]{
\bfig
\putmorphism(0,0)(1,0)[E`A`e]{400}1a
\putmorphism(430,23)(1,0)[``f]{370}1a
\putmorphism(430,0)(1,0)[`B`]{425}0a
\putmorphism(430,-23)(1,0)[``g]{370}1b
\efig}
$
be an equalizer in $\C$. \vspace{-5pt}
\begin{enumerate}
\itemsep -2pt
\item[(i)] If $f$ and $g$ are algebra morphisms, then $E$ is an algebra and $e$ is an algebra morphism. We call $(E,e)$ an algebra pair.
\item[(ii)] If $f$ and $g$ are left (resp. right) $H$-module morphisms, then $E$ is a left (resp. right)
$H$-module and $e$ is a left (resp. right) $H$-module morphism. The pair $(E,e)$ is called an $H$-module pair.
The dual statement follows for $H$-comodules {\it provided that $H$ is flat}. Similarly, $(E,e)$ is said to be an $H$-comodule pair.
\end{enumerate}
Assume that $A,B$ are $H$-comodule algebras, where $H$ is flat, and let $(E, e)$ be the equalizer of $H$-comodule algebra morphisms $f, g:A\to B$. Then $E$ inherits an $H$-comodule algebra structure such that $e$ is a morphism of
$H$-comodule algebras. It is important to stress that once the corresponding structure morphisms are induced on $E$, one does not need anymore the hypothesis on $f$ and $g$ to prove the compatibility conditions in each structure (the reader is encouraged to check it). \par \medskip

\addtocounter{preli}{1}\thepreli\label{pr22} Let $A \in \C$ be an algebra and ${}_A\C$ the category of left $A$-modules. Given two morphisms $f,g:M \rightarrow N$ in ${}_A\C$ consider its equalizer in $\C$:
$$\scalebox{0.9}[0.9]{
\bfig
\putmorphism(0,0)(1,0)[E`M`e]{400}1a
\putmorphism(430,23)(1,0)[``f]{370}1a
\putmorphism(430,0)(1,0)[`N`]{425}0a
\putmorphism(430,-23)(1,0)[``g]{370}1b
\efig}
$$
Then $E$ admits a structure of left $A$-module such that $e$ is a morphism of $A$-modules and $(E,e)$ is the equalizer in $_A\C$ of $f$ and $g$. Hence ${}_A\C$ has equalizers. Moreover, the forgetful functor $\U: {}_A\C\to \C$ preserves equalizers. Both statements hold true for the category $\C^C$ of right $C$-comodules {\it if $C$ is a flat coalgebra}. \par \medskip

\addtocounter{preli}{1}\thepreli\label{pr23} Consider a commutative diagram of morphisms in $\C$
$$\scalebox{0.9}[0.9]{
\square[E_1`A_1`E_2`A_2;e_1`\overline{f}`f`e_2]}
$$
and assume that $e_2$ is a monomorphism.
\begin{enumerate}
\itemsep -2pt
\item[(i)] If $e_1, e_2$ and $f$ are right $H$-comodule morphisms and $e_2\ot H$ is a monomorphism, then $\overline{f}$ is an $H$-comodule morphism.
\item[(ii)] If $e_1, e_2$ and $f$ are right $H$-module morphisms, then $\overline{f}$ is
an $H$-module morphism.
\item[(iii)] If $e_1, e_2$ and $f$ are algebra morphisms, then $\overline{f}$ is an algebra morphism.
\end{enumerate}
As a consequence we have:
\begin{enumerate}
\itemsep -2pt
\item[(a)] Let $f:A\to R$ and $g:B\to R$ be algebra (resp. $H$-comodule) morphisms and assume
that $g$ (resp. $g\ot H$) is a monomorphism. Let $h:A\to B$ be such that $gh=f$. Then $h$
is an algebra (resp. $H$-comodule) morphism.
\end{enumerate}
\par \medskip

\addtocounter{preli}{1}\thepreli\label{pr24} Given two morphisms $f,g:N \rightarrow M$ in ${}_A\C$ consider its coequalizer in $\C$:
$$\scalebox{0.9}[0.9]{
\bfig
\putmorphism(400,0)(1,0)[N`Q`q]{380}1a
\putmorphism(0,20)(1,0)[``f]{350}1a
\putmorphism(-30,0)(1,0)[M``]{425}0a
\putmorphism(0,-20)(1,0)[``g]{350}1b
\efig}$$
Then $Q$ admits a structure of left $A$-module such that $q$ is a morphism of $A$-modules and $(Q,q)$ is the equalizer in $_A\C$ of $f$ and $g$. Hence, ${}_A\C$ has coequalizers. The closedness of $\C$ guarantees the existence of the structure morphism $A\ot Q\to Q$. This is because,  as a left adjoint functor, $A\ot -$ preserves coequalizers. The forgetful functor $\U: {}_A\C\to \C$ preserves coequalizers.

\subsection{Tensor product of modules over an algebra}
\setcounter{preli}{0}

\addtocounter{preli}{1}\thepreli\label{pr31} Given an algebra $A \in \C$, the {\em tensor product over $A$} of a right $A$-module $M$ and a left $A$-module $N$ is the coequalizer in $\C$
$$\scalebox{0.9}[0.9]{\bfig
 \putmorphism(930,0)(1,0)[`M\ot_A N.`\Pi_{M, N}]{550}1a
 \putmorphism(-150,25)(1,0)[\phantom{M\ot A\ot N}`\phantom{M\ot N}`\mu_M\ot N]{950}1a
 \putmorphism(-150,0)(1,0)[M\ot A\ot N`M\ot N`]{950}0a
 \putmorphism(-150,-25)(1,0)[\phantom{M\ot A\ot N}`\phantom{M\ot N}` M\ot \nu_N]{950}1b
\efig}$$ Consider morphisms $f:M\to M'$ in $\C_A$ and $g:N\to N'$ in $_A\C$. Then $f\ot g:M\ot N\to M'\ot N'$ induces a morphism $f\ot_A g:M\ot_A N\to M'\ot_A N'$. Thus we have functors $M \otimes_A -: {}_A\C \rightarrow \C$ and $-\otimes_A N: \C_A \rightarrow \C$. \par \medskip

\addtocounter{preli}{1}\thepreli\label{pr32} Let $A, B, R\in\C$ be algebras and consider the categories of $A\x B$-bimodules ${}_A\C_B$ and $B\x R$-bimodules ${}_B\C_R$. For $M\in {}_A\C_B$ and $N\in{}_B\C_R$, the object $M \otimes N$ admits a structure of left $A$-module and right $R$-module inherited from $M$ and $N$ respectively. Denote the structure morphisms by $\nu_{M\ot N}$ and $\mu_{M\ot N}$ respectively. The object $M\ot_B N$ is an $A\x R$-bimodule with structure morphisms $\crta\nu:A\ot(M\ot_B N)\to M\ot_B N$ and $\crta\mu:(M\ot_B N)\ot R\to M\ot_B N$ defined via the universal properties of the coequalizers $(A\ot(M\ot_B N), A\ot\Pi_{M, N})$ and $((M\ot_B N)\ot R, \Pi_{M, N}\ot R)$ respectively by $\crta\nu(A\ot\Pi_{M, N})=\Pi_{M, N}\nu_{M\ot N}$ and $\crta\mu(\Pi_{M, N}\ot R)=\Pi_{M, N}\mu_{M\ot N}$. Furthermore, the coequalizer morphism $\Pi_{M, N}$ is an $A\x R$-bimodule morphism. Moreover, if $f:M \rightarrow M'$ and $g:N \rightarrow N'$ are morphisms in $_A\C_B$ and $_B\C_R$ respectively, then $f\otimes_B g:M \otimes_B N \rightarrow M' \otimes_B N'$ is a morphism in $_A\C_R.$ In particular, for $M \in {}_A\C_B$ and $N \in {}_B\C_R$ we have functors $M \otimes_B -:{}_B\C \rightarrow {}_A\C$ and $- \otimes_B N:\C_B \rightarrow \C_R.$ \par \medskip

\addtocounter{preli}{1}\thepreli\label{pr33} For $M \in {}_A\C$ there is a natural isomorphism $A\ot_A M \iso M$ of left $A$-modules induced by the structure morphism of $M$. If additionally $M\in{}_A\C_B$, then this is an isomorphism of $A\x B$-bimodules. \par \medskip

\addtocounter{preli}{1}\thepreli\label{pr34} Let $M\in{}_A\C_A$ and $N\in\C$. The object $M \otimes N$ becomes a left $A$-module from the structure of left $A$-module of $M$. Then $M\ot_A (A\ot N)\iso M\ot N$ in ${}_A\C_A$ by
$$\gamma: M\ot_A (A\ot N)\to M\ot N$$ induced on the coequalizer by the commuting
diagram:
$$\scalebox{0.85}[0.85]{
\bfig
\putmorphism(-1050,20)(1,0)[``\mu_M\ot (A\ot N)]{1130}1a
\putmorphism(-1470,0)(1,0)[M\ot A\ot (A\ot N)``]{1130}0a
\putmorphism(-1050,-20)(1,0)[``M\ot {[(\nabla_A\ot N)\alpha_{A, A, N}^{-1}]}]{1130}1b
\putmorphism(360,0)(0,-1)[`(M\ot A)\ot N`\alpha_{M, A, N}^{-1}]{600}1l
\putmorphism(380,0)(1,0)[M\ot (A\ot N)`M\ot_A (A\ot N)`\Pi_{M, A\ot N}]{1050}1a
\putmorphism(1300,-30)(0,-1)[`M\ot N`\gamma]{600}1r
\putmorphism(330,-600)(1,0)[\phantom{(M\ot A)\ot N}``\mu_M\ot N]{800}1r
\efig}
$$
Subsequently, we have an isomorphism in ${}_A\C_A$
$$\omega:=\gamma^{-1}(\delta\ot N):(M\ot_A A)\ot N\to M\ot_A (A\ot N)$$
where $\delta$ is the right version of the morphism from 2.3.3. \par \medskip

\addtocounter{preli}{1}\thepreli\label{pr35} {\bf Associativity of the tensor product:} Since $\C$ is closed, the functor $M \otimes-$ preserves coequalizers. Then, for every $X \in \C, M \in \C_A$ and $N \in {}_A\C$ there is a canonical isomorphism $X \otimes (M \otimes_A N) \cong (X \otimes M) \otimes_A N$ induced by the associativity of the tensor product $\otimes$. Using this isomorphism, the coequalizer universal property and $3\times 3$ Lemma, one
may show that for $M\in\C_A, N\in{}_A\C_B$ and $L\in{}_B\C$ the coequalizers $M\ot_A (N \ot_B L)$ and $(M\ot_A N) \ot_B L$ are isomorphic in $\C$. This isomorphism becomes an isomorphism in ${}_R\C_T$ if $M \in {}_R\C_A$ and $L\in {}_B\C_T$.
\par \smallskip

For general monoidal categories, Pareigis used in \cite[page 202]{P1} the notion of coflatness to assure the associativity of the tensor product. Let $A$ and $B$ be algebras in $\C$. An object $M\in {}_B\C_A$ is called {\em $A$-coflat} if for all algebras $R, T\in\C$ and objects $L\in {}_A\C_R$ the coequalizer $M\ot_A L \in {}_B\C_R$ exists and for every $P\in\C_T$ the natural morphism $M\ot_A(L\ot P) \to (M\ot_A L)\ot P$ in ${}_B\C_T$, induced by the associativity of the tensor product, is an isomorphism.  Symetrically,  one may define the notion of $B$-coflat for $M$. For $M,N$ and $L$ as in the preceding paragraph, there is an isomorphism between $M\ot_A (N \ot_B L)$ and $(M\ot_A N) \ot_B L$ if $M$ is $A$-coflat and $L$ is $B$-coflat. Notice that {\it in a closed braided monoidal category with coequalizers every bimodule $M\in {}_B\C_A$ is coflat} on both sides since the functors $X \otimes -$ and $- \otimes X$, as left adjoint functors, preserve coequalizers for all $X \in \C$. \par \medskip

\addtocounter{preli}{1}\thepreli\label{pr36} {\bf The functor $_A[-,-]:$} For an algebra $A\in\C$ and $M\in {}_A\C$ we can consider the functor $M\ot -:\C\to {}_A\C$. For each $X \in \C$ the object $M \otimes X$ is a left $A$-module in the natural way. We recall from
the discussion preceding \cite[Theorem 3.2]{P1} that this functor has a right adjoint ${}_A[M, -]$. For $N \in{}_A\C$, the object ${}_A[M, N]$ is the following equalizer:
$$\scalebox{0.9}[0.9]{
\bfig
\putmorphism(0,0)(1,0)[{}_A[M, N]`[M, N]`\iota_N]{600}1a
\putmorphism(630,20)(1,0)[\phantom{[M, N]}`\phantom{[A\ot M, N],}`u_N]{720}1a
\putmorphism(630,0)(1,0)[`[A\ot M, N],`]{725}0a
\putmorphism(630,-20)(1,0)[\phantom{[M, N]}`\phantom{[M, N]}`v_N]{600}1b
\efig}
$$
where $u$ and $v$ are given via the commutative diagrams
\begin{eqnarray}\label{Ae-u}
\hspace{-1cm}
\scalebox{0.9}[0.9]{\bfig
\putmorphism(0, 400)(1, 0)[A\ot M\ot [M, N]`\phantom{M\ot [M, N]}`\nu_M\ot {[M, N]}]{1300}1a
\putmorphism(0, 400)(0, -1)[\phantom{A\ot M\ot [M, N]}`\phantom{A\ot M\ot [A\ot M, N]}`A\ot M\ot u_N]{400}1l
\putmorphism(1300, 400)(0, -1)[M\ot [M, N]`\phantom{N}`\ev]{400}1r
\putmorphism(0, 0)(1, 0)[A\ot M\ot [A\ot M, N]`N`\ev]{1300}1b
\efig} & \quad
\scalebox{0.9}[0.9]{\bfig
\putmorphism(0, 400)(1, 0)[A\ot M\ot [M, N]`\phantom{A\ot N}`A\ot \ev]{1100}1a
\putmorphism(0, 400)(0, -1)[\phantom{A\ot M\ot [M, N]}`\phantom{A\ot M\ot [A\ot M, N]}`A\ot M\ot v_N]{400}1l
\putmorphism(1100, 400)(0, -1)[A\ot N`\phantom{N}`\nu_N]{400}1r
\putmorphism(0, 0)(1, 0)[A\ot M\ot [A\ot M, N]`N`\ev]{1100}1b
\efig}
\end{eqnarray}
If $f:N\to N'$ is a morphism in ${}_A\C$, then ${}_A[M, f]: {}_A[M, N]\to {}_A[M, N']$ is
induced as a morphism on an equalizer via the diagram
$$\scalebox{0.9}[0.9]{
\bfig
\putmorphism(0,425)(1,0)[{}_A[M, N]`[M, N]`\iota_N]{600}1a
\putmorphism(0,0)(1,0)[{}_A[M, N']`[M, N']`\iota_{N'}]{600}1a
\putmorphism(0,420)(0,-1)[\phantom{{}_A[M, N]}`\phantom{{}_A[M, N']}`{}_A{[M, f]}]{400}1l
\putmorphism(605,420)(0,-1)[\phantom{[M, N]}`\phantom{[M, N']}`{{[M, f]}}]{400}1r
\putmorphism(725,20)(1,0)[``u_{N'}]{450}1a
\putmorphism(1000,0)(1,0)[`[A\ot M, N']`]{450}0a
\putmorphism(725,-20)(1,0)[``v_{N'}]{450}1b
\efig}
$$
That $[M, f]$ induces ${}_A[M, f]$ is provided by $A$-linearity of $f$. \par \smallskip

The unit $\alpha: N\to [M, M\ot N]$ and the counit $\ev:M\ot [M, N]\to N$ of the adjunction $(M\ot -, [M,-])$ between $\C$ and $\C$, on the one hand, and the unit ${\sf a}:N\to {}_A[M, M\ot N]$ and the counit ${\sf ev}: M\ot {}_A[M, N]\to N$ of the adjunction $(M\ot -, {}_A[M,-])$ between $\C$ and ${}_A\C$, on the other hand, are related by $\alpha=\iota{\sf a}$ and ${\sf ev}=\ev(A\ot\iota)$. \par \smallskip

Similarly to 2.1.3, $\nabla: {}_A[M,M] \otimes {}_A[M,M] \rightarrow {}_A[M,M]$ and $\eta: I \rightarrow {}_A[M,M]$ (given by the universal property of ${\sf ev}$ as the unique morphisms such that ${\sf ev}(M \otimes \nabla)={\sf ev}({\sf ev} \otimes {}_A[M,M])$ and ${\sf ev}(M \otimes \eta)=Id_M$) endow $_A[M,M]$ with a structure of algebra. Moreover, $M$ is an $A\x {}_A[M,M]$-bimodule with right structure morphism ${\sf ev}:M \otimes {}_A[M,M] \rightarrow M.$ \par \smallskip

If $M \in {}_A\C_B$ and $N \in {}_A\C$, then $_A[M,N]$ is a left $B$-module with left structure morphism $\nu:B \otimes {}_A[M,N] \rightarrow {}_A[M,N]$ defined as the unique one such that ${\sf ev}(\mu_M \otimes {}_A[M,N])={\sf ev}(M \otimes \nu)$. Moreover, if $f:N \rightarrow N'$ is a morphism in ${}_A\C$, then $_A[M,f]:{}_A[M,N] \rightarrow {}_A[M,N']$ becomes a morphism in $_B\C$. Thus we have a functor $_A[M,-]:{}_A\C \rightarrow {}_B\C$. For $M \in {}_A\C_B$ the adjunction $(M\ot -, {}_A[M,-])$ between $\C$ and ${}_A\C$ induces an adjunction $(M\ot_B -, {}_A[M,-])$ between $_B\C$ and ${}_A\C$ through the universal property of the coequalizer $M\ot_B -,$ \cite[Proposition 3.10]{P1}. The counit of the adjunction $\overline{{\sf ev}}:M \otimes_B {}_A[M,N] \rightarrow N$ is given as the unique morphism such that ${\sf ev}=\overline{{\sf ev}}\Pi_{M,{}_A[M,N]}.$\par \smallskip

Finally, if $N \in {}_A\C_R$, then $_A[M,N]$ becomes a $B\x R$-bimodule with right structure morphism $\zeta:{}_A[M,N] \otimes R \rightarrow {}_A[M,N]$ satisfying ${\sf ev}(M \otimes \zeta)=\mu_N({\sf ev} \otimes R)$. If $f:N \rightarrow N'$ is a morphism in ${}_A\C_R,$ then $_A[M,f]:[M,N] \rightarrow [M,N']$ is a morphism in $_B\C_R$.

\subsection{Morita theory}
\setcounter{preli}{0}

In this subsection we recall the version for monoidal categories of Morita Theorems due to Paregis \cite[Theorems 5.1 and 5.3]{P3}. Since we are working in the framework of closed braided monoidal categories with equalizers and coequalizers we will adapt these theorems to our setting, losing in generality but avoiding some technical difficulties like the use of coflatness.\par \medskip

\addtocounter{preli}{1}\thepreli\label{pr41} {\bf $\C$-functors:} Let $A,B$ be algebras in $\C$. A functor $\F:{}_A\C\to{}_B\C$ is called a {\em $\C$-functor} if for all $M\in{}_A\C$ there is a natural isomorphism between the functors  $\F(M\ot -)$ and $\F(M)\ot -$ from $\C$ to $_B\C$. It is proved in \cite[Theorems 4.2 and 4.3]{P2} that for $Q\in {}_B\C_A$, the functors $Q\ot_A-:{}_A\C\to{}_B\C$ and ${}_B[Q,-]:{}_B\C\to{}_A\C$ are $\C$-functors. \par \medskip

\addtocounter{preli}{1}\thepreli\label{pr42} {\bf Finite and faithful projective objects:} An object $P\in {}_A\C$ is called {\em finite} if the morphism $\overline{\sf db}: {}_A[P,A] \ot_A P \to {}_A[P,P]$  defined via the universal property of $({}_A[P,P],
{\sf ev}_{P,P}:P\ot {}_A[P,P]\to P)$ by
\begin{eqnarray}
\scalebox{0.9}[0.9]{
\bfig
\putmorphism(130, 440)(2, -1)[\phantom{P\ot {}_A[P,A]\ot P}`\phantom{P\iso P \ot_A A}`]{850}1r
\putmorphism(90, 500)(2, -1)[\phantom{P\ot {}_A[P,A]\ot P}`\phantom{P\iso P \ot_A A}`{{\sf ev}}_{P,A} \ot_A P]{850}0r
\putmorphism(80, 400)(0, -1)[`\phantom{{}_A[P,P] \ot P}`P\ot {\overline {\sf db}}]{400}1l
\putmorphism(-50, 400)(0, -1)[P\ot {}_A[P,A]\ot_A P``]{400}0l
\putmorphism(0, 0)(1, 0)[P \ot {}_A[P,P]`P\iso A\ot_A P`]{900}1b
\putmorphism(0, -25)(1, 0)[\phantom{{}_A[P,P] \ot P}`\phantom{P\iso A \ot_A P}`{{\sf ev}_{P,P}}]{900}0b
\efig}
\end{eqnarray}
is an isomorphism. An object $P\in {}_A\C$ is {\em faithfully projective} if it is finite in ${}_A\C$ and  $\overline{{\sf ev}}_{P,A}: P \ot_{{}_A[P,P]}{}_A[P,A] \to A$ is an isomorphism. Similarly, one may define finiteness and faithful projectiveness for a right module. \par \medskip

\addtocounter{preli}{1}\thepreli\label{pr43} {\bf Morita context:} A {\em Morita context in $\C$} is a sextuple $(A,B,P,Q,f,g)$ consisting of algebras $A, B\in\C,$ bimodules $P \in {_A}\C_B, Q \in {_B}\C_A$ and morphisms $f:P\ot_B Q\to A$ in ${}_A\C_A$ and $g:Q\ot_A P\to B$ in ${}_B\C_B$ such that the diagrams
\begin{flushleft}
\noindent
\begin{tabular}{p{6cm}p{1.3cm}p{6cm}}
$\begin{array}{l}
\scalebox{0.8}[0.8]{\bfig
\putmorphism(0,500)(1,0)[\phantom{P\ot_B(Q\ot_A P) \iso (P\ot_B Q)\ot_A P}`\phantom{A \ot_A P}`f \ot_A P]{1300}1a
\putmorphism(1300,500)(0,-1)[A \ot_A P`P`\nu_P]{500}1r
\putmorphism(0,500)(0,-1)[P\ot_B(Q\ot_A P) \iso (P\ot_B Q)\ot_A P`P\ot_B B`P\ot_B g]{500}1l
\putmorphism(0,0)(1,0)[\phantom{P\ot_B B}`\phantom{P}`\mu_P]{1300}1b
\efig}
\end{array}$ & &
$\begin{array}{l}
\scalebox{0.8}[0.8]{\bfig
\putmorphism(0,500)(1,0)[\phantom{Q\ot_A(P\ot_B Q) \iso (Q\ot_A P)\ot_B Q}`\phantom{B \ot_B Q}`g \ot_B Q]{1300}1a
\putmorphism(1300,500)(0,-1)[B \ot_B Q`Q`\nu_Q]{500}1r
\putmorphism(0,500)(0,-1)[Q\ot_A(P\ot_B Q) \iso (Q\ot_A P)\ot_B Q`Q\ot_A A`Q\ot_A f]{500}1l
\putmorphism(0,0)(1,0)[\phantom{Q\ot_A A}`\phantom{Q}`\mu_Q]{1300}1b
\efig}
\end{array}$
\end{tabular}
\end{flushleft}
are commutative. Here, abusing of notation, $\mu_P,\nu_P$ and $\mu_Q,\nu_Q$ are the isomorphisms induced by the
module structures of $P$ and $Q$ respectively. If moreover there exist morphisms
$\zeta\in\C(I, P\ot_B Q)$ and $\xi\in\C(I, Q\ot_A P)$ such that $f\zeta=\eta_A$ and $g\xi=\eta_B,$ the context is called {\em strict}.\par \medskip

\addtocounter{preli}{1}\thepreli\label{pr44} {\bf Morita Theorem I:} Let $(A, B, P, Q, f, g)$ be a strict Morita context in $\C$. Then $f$ and $g$ are isomorphisms. Moreover:
\begin{enumerate}
\item[(i)] The functors $P\ot_B -:{}_B\C\to{}_A\C$ and $Q\ot_A -:{}_A\C\to{}_B\C$ are inverse $\C$-equivalences. In particular, $P\in{}_A\C$ and $Q\in{}_B\C$ are faithfully projective.
\item[(ii)] The functors $-\ot_B Q: \C_B \to \C_A$ and $-\ot_A P: \C_A \to \C_B$ are inverse $\C$-equivalences. In particular, $Q\in \C_A$ and $P\in \C_B$ are faithfully projective.
\end{enumerate}

\addtocounter{preli}{1}\thepreli\label{pr45} {\bf Morita Theorem II:} Let $\F:{}_A\C\to{}_B\C$ and $\G:{}_B\C\to{}_A\C$ be inverse $\C$-equivalences. Then there exist objects $P\in {}_A\C_B$ and $Q\in {}_B\C_A$ such that:
\begin{enumerate}
\item[(i)]  There are natural isomorphisms $\F(-)\iso Q\ot_A -\iso {}_A[P, -]$ and $\G(-)\iso P\ot_B -\iso {}_B[Q,-].$

\item[(ii)] There are isomorphisms $A\iso P\ot_B Q$ in $_A\C_A$ and $B\iso Q\ot_A P$ in $_B\C_B$ so that the diagrams in 2.4.3 commute.

\item[(iii)] There are isomorphisms ${}_B[Q, B]\iso P$ in ${}_A\C_B$ and ${}_A[P, A]\iso Q$ in ${}_B\C_A.$
	
\item[(iv)] There are isomorphisms ${}_B[Q, Q]\iso A$ in ${}_A\C_A$ and ${}_A[P, P]\iso B$ in ${}_B\C_B$ that are also isomorphisms of algebras.
\end{enumerate}

\addtocounter{preli}{1}\thepreli\label{pr46} {\bf The derived Morita context:} Let $P \in {}_A\C$ and set ${}^\circledast P={}_A[P,A].$ Recall that $P$ is an $A\x {}_A[P,P]$-bimodule and that it induces a structure of ${}_A[P,P]\x A$-bimodule
on ${}^\circledast P.$ Consider now the morphisms $f=\overline{{\sf ev}}: P \ot_{{}_A[P,P]} {}^\circledast P \to A$ and $g=\overline{{\sf db}}: {}^\circledast P \ot_A P \to {}_A[P,P]$  of 2.4.2. Then $(A,{}_A[P,P],P,{}^\circledast P,f,g)$ is a Morita context. If $_AP$ is faithfully projective, then it is a strict Morita context. In particular, {\it if $P \in \C$ is faithfully projective,  then $P\ot-: \C\to {}_{[P,P]}\C$ is a $\C$-equivalence.}

\subsection{Faithful projectiveness and faithful flatness}
\setcounter{preli}{0}

\addtocounter{preli}{1}\thepreli\label{pr51} {\it A faithfully projective object is faithfully flat.}\par \smallskip

\begin{proof}
Let $P$ be faithfully projective in $\C$ and set $P^*=[P, I]$. By 2.4.6,  $(P\ot-, P^*\ot_{[P,P]} -)$ defines a $\C$-equivalence between $\C$ and ${}_{[P,P]}\C$. To prove that $P$ is flat, consider an equalizer in $\C$:
$$\bfig
\putmorphism(0,0)(1,0)[E`M`e]{400}1a
\putmorphism(430,25)(1,0)[``f]{350}1a
\putmorphism(430,0)(1,0)[`N`]{425}0a
\putmorphism(430,-25)(1,0)[``g]{350}1b
\efig$$
Through the equivalence $P\ot -,$
$$\bfig
\putmorphism(-30,0)(1,0)[P\ot E`P\ot M`P\ot e]{620}1a
\putmorphism(730,25)(1,0)[``P\ot f]{470}1a
\putmorphism(730,0)(1,0)[`P\ot N`]{620}0a
\putmorphism(730,-25)(1,0)[``P\ot g]{470}1b
\efig$$
is an equalizer in ${}_{[P,P]}\C$.  By 2.2.2, $P\ot E$ is the equalizer in $\C$ of $P\ot f$ and $P\ot g$. To see that $P\ot -$ reflects isomorphisms, assume $P\ot h$ is an isomorphism in $\C$ for a certain $h:M \rightarrow N$. Considering $P\ot h$ as a morphism in ${}_{[P,P]}\C$, it is then an isomorphism in there. By the inverse equivalence $P^*\ot_{[P,P]}(P\ot h)\iso(P^*\ot_{[P,P]}P)\ot h \iso h$ will then be an isomorphism in $\C$.
\qed\end{proof}

\addtocounter{preli}{1}\thepreli\label{pr52} {\it If $P \in \C$ is finite and $P \otimes -$ reflects isomorphisms, then $P$ is faithfully projective.} \par \smallskip

\begin{proof}
To prove that $P$ is faithfully projective we need to show that $\overline{{\sf ev}}:P^*\ot_{[P,P]}P\to I$ is an isomorphism. Consider the isomorphism $\beta:P\ot(P^*\ot_{[P,P]}P)\to (P\ot P^*)\ot_{[P,P]}P$ induced by the associativity of the tensor product. If we show that the composition
$$P\ot(P^*\ot_{[P,P]}P)\stackrel{\beta}{\iso}(P\ot P^*)\ot_{[P,P]}P \stackrel{\db\ot_{[P,P]}P}{\iso}[P,P]\ot_{[P,P]}P\stackrel{\delta}{\iso}P$$
equals $P\ot \overline{{\sf ev}}$, we would obtain that $\overline{{\sf ev}}$ is an isomorphism because $P\ot -$ reflects isomorphisms and we would be done. \par \smallskip

That $P\ot \overline{{\sf ev}}=\delta(\db\ot_{[P,P]}P)\beta$ we will conclude from the following diagram:
$$\scalebox{0.9}[0.9]{\bfig
\putmorphism(-100,570)(1,0)[\phantom{P\ot(P^*\ot_{[P,P]} P)}`\phantom{P\ot(P^*\ot_{[P,P]}P)}` P\ot \Pi_{P^*,P}]{1400}1b
\putmorphism(1300,570)(1,0)[\phantom{P\ot(P^*\ot_{[P,P]} P)}`P\ot I\iso P`P\ot\overline{{\sf ev}}]{1300}1b
\putmorphism(0,625)(1,0)[\phantom{P\ot(P^*\ot P)}`\phantom{P\ot I\iso P} `P\ot\ev]{2590}1a
\putmorphism(0,600)(0,-1)[P\ot(P^*\ot P)`(P\ot P^*)\ot P`\iso]{600}1l
\putmorphism(1300,580)(0,-1)[\phantom{P\ot I\iso P}``\beta]{580}1l
\putmorphism(1300,550)(0,-1)[P\ot(P^*\ot_{[P,P]}P)``]{600}0l
\putmorphism(1300,620)(0,-1)[`(P\ot P^*)\ot_{[P,P]} P`]{600}0l
\putmorphism(20,0)(1,0)[\phantom{(P\ot P^*)\ot P}`\phantom{(P\ot P^*)\ot_{[P,P]} P}`\Pi_{P\ot P^*,P}]{1300}1a
\putmorphism(0,0)(0,-1)[\phantom{(P\ot P^*)\ot P}`[P, P]\ot P`db \ot P]{600}1l
\putmorphism(20,-600)(1,0)[\phantom{[P, P]\ot P}`[P, P]\ot_{[P,P]} P`\Pi_{[P,P],P}]{1300}1a
\putmorphism(1300,-20)(0,-1)[\phantom{P\ot I\iso P}``\db\ot_{[P,P]} P]{580}1l
\putmorphism(1520,-470)(1,1)[\phantom{P\ot I\iso P}` \phantom{P\ot I\iso P}`\delta]{900}1r
\put(520,250){\fbox{2}}
\put(520,-300){\fbox{3}}
\put(1640,200){\fbox{1}}
\efig}$$
Our question is formulated in the diagram $\langle 1\rangle$. It suffices to prove that it commutes
when composed to $P\ot \Pi_{P^*,P}$, as the latter is an epimorphism. In the top row of the diagram we have $P\ot \ev=(P\ot\overline{{\sf ev}}) (P\ot \Pi_{P^*,P})$ by the definition of $\overline{{\sf ev}}$. Diagrams $\langle 2 \rangle$ and
$\langle3\rangle$ commute by the definition of $\beta$ and $\db\ot_{[P,P]} P$ respectively. A diagram chasing argument will yield our claim when we prove that the outer diagram commutes. \par \smallskip

Note that $\delta: [P,P]\ot_{[P,P]} P\to P$ is induced by $\ev_{P,P}: [P,P]\ot P\to P$ so
that $\ev_{P,P}=\delta\Pi_{[P,P],P}$. On the other hand, by Diagram \ref{db-def}, we have $\ev_{P,P}(\db\ot P)=P\ot\ev$. Combining these two we get $\delta \Pi_{[P,P], P}(\db\ot P)=P\ot\ev.$
\qed\end{proof}

\section{The Brauer group of a braided monoidal category}
\setcounter{equation}{0}

In \cite{VZ1} Van Oystaeyen and Zhang constructed the Brauer group of a braided monoidal category extending
a previous construction of Pareigis for symmetric monoidal categories \cite{P4}. This categorical construction collects practically all examples of Brauer groups ocurring in the literature. In the first part of this section we recall the construction of this group, although we will use braided diagrams instead of applying the Yoneda Lemma, as done in the two cited articles. We will also give an alternative description of one of the functors associated to an Azumaya algebra that will be crucial for future results. In the second part, we will introduce the Brauer group of module algebras and present the left hand side part of the sequence we aim to construct. \par\medskip

As in the previous section, {\em $\C$ will denote a closed braided monoidal category with braiding $\Phi$ and possessing equalizers and coequalizers} (although we will not use coequalizers here).

\subsection{Azumaya algebras}
For an algebra $A \in \C$ its {\em opposite algebra}, denoted by $\crta A$, is defined as follows: As an object
$\crta A=A$, the unit is the same as in $A$ and the multiplication is given by $\nabla_{\crta A}=
\nabla_{\hspace{-3pt}A}\hspace{1pt} \Phi_{A,A}$. Let $\alpha_{A,-}: \Id_{\C} \to [A, -\ot A]$ denote the unit of the adjunction $(-\ot A, [A,-])$. Consider the morphism $f=\nabla_{\hspace{-3pt} A}(\nabla_A\ot A)(A\ot \Phi):(A\otimes \crta A)\otimes A \rightarrow A$. We define $F:A\ot \crta A\to[A, A]$ as the composition
$$\scalebox{0.95}[0.95]{
\bfig
\putmorphism(0, 0)(1, 0)[A\ot \crta A`{[A, (A\ot \crta A)\ot A]}`\alpha_{A,A\ot \crta A}]{1150}1a
\putmorphism(1150, 0)(1, 0)[\phantom{[A, (A\ot \crta A)\ot A]}`[A,A].`]{1150}1a
\putmorphism(1700, 20)(1, 0)[``{{[A, f]}}]{300}0a
\efig}
$$
The morphism $F$ is an algebra morphism. Using the universal property, it can also be described as follows:
\begin{eqnarray} \label{AzDefMap}
\scalebox{1}[1]{
\gbeg{3}{5}
\got{1}{A\ot\crta A} \got{3}{A} \gnl
\gbmp{F}\gvac{1} \gcl{1}\gnl
$\eval$
\gvac{1} \gob{1}{A}
\gend}=\scalebox{1}[1]{
\gbeg{4}{6}
\gvac{1}\got{1}{A\ot\crta A} \got{3}{A} \gnl
\gvac{1}\gbmp{\alpha}\gvac{1} \gcl{2}\gnl
\glmp\gnot{[A,f]}\gcmptb\grmp \gnl
\gvac{1} \gcn{1}{1}{1}{3}\gcn{1}{1}{3}{1} \gnl
\gvac{2} \gmp{\ev} \gnl
\gvac{2} \gob{1}{A}
\gend}\stackrel{(\ref{rightAdj})}{=}
\scalebox{1}[1]{\gbeg{3}{6}
\got{1}{A\ot\crta A} \got{3}{A} \gnl
\gbmp{\alpha}\gvac{1} \gcl{1}\gnl
$\eval$ \gnl
\gnl
\gvac{1} \gbmp{f} \gnl
\gvac{1} \gob{1}{A}
\gend}\stackrel{(\ref{adjunit})}{=}
\scalebox{1}[1]{\gbeg{3}{5}
\got{3}{(A\ot\crta A)\ot A} \gnl
\gvac{1}\gcl{1}\gvac{1}\gnl
\gvac{1}\gbmp{f}\gvac{1}\gnl
\gvac{1}\gcl{1}\gvac{1}\gnl
\gvac{1}\gob{1}{A}\gvac{1}
\gend}=
\scalebox{0.95}[0.95]{\gbeg{4}{5}
\got{1}{A}\got{3}{\crta A}\got{1}{\hspace{-0,8cm}A} \gnl
\gcl{1} \gvac{1}\gbr \gnl
\gwmu{3} \gcl{1} \gnl
\gvac{1}\gwmu{3} \gnl
\gvac{2}\gob{1}{A}
\gend}
\end{eqnarray}
Let $\tilde{\alpha}_{A,-}: \Id_{\C} \to [A, A\ot -]$ denote the unit of the adjunction $(A\ot -, [A,-])$.
Recall from 2.1.4 how it is related to the unit $\alpha_{A,-}: \Id_{\C} \to [A, -\ot A]$. Let $g=\nabla_{\hspace{-3pt}A}(A\otimes\nabla_A)(\Phi \ot A)$. We define $G:\crta A\ot A\to\crta{[A,A]}$ as the composition
$$\bfig
\putmorphism(0, 0)(1, 0)[\crta A\ot A`{[A, A\ot (\crta A\ot A)]}`\tilde\alpha_{A,\crta A\ot A}]{1200}1a
\putmorphism(1200, 0)(1, 0)[\phantom{[A, A\ot (\crta A\ot A)]}`\crta{[A,A]}.`]{1200}1a
\putmorphism(1800, 15)(1, 0)[``{{[A, g]}}]{300}0a
\efig$$
It can be described as follows:
\begin{eqnarray}\label{AzDefMap-G}
\gbeg{3}{5}
\got{1}{A} \got{3}{\crta A \ot A} \gnl
\gcl{1} \gvac{1} \gbmp{G} \gnl
\hspace{-0,14cm}$\crtaeval$
\gvac{1} \gob{1}{A}
\gend=
\gbeg{5}{6}
\got{1}{A} \got{3}{\crta A \ot A} \gnl
\gcl{2}\gvac{1}\gbmp{\tilde\alpha}\gvac{1} \gnl
\gvac{1}\glmp\gnot{[A,g]}\gcmptb\grmp \gnl
\gcn{1}{1}{1}{3}\gcn{1}{1}{3}{1} \gvac{1} \gnl
\gvac{1} \gmp{\tilde\ev} \gvac{1} \gnl
\gvac{1} \gob{1}{A}
\gend\stackrel{(\ref{rightAdj})}{=}
\gbeg{4}{6}
\got{1}{A} \got{3}{\crta A \ot A} \gnl
\gcl{1} \gvac{1} \gbmp{\tilde\alpha} \gnl
\hspace{-0,26cm} $\crtaeval$ \gnl
\gnl
\gvac{1} \gbmp{g}  \gnl
\gvac{1} \gob{1}{A}
\gend\stackrel{(\ref{adjunit-op})}{=}
\gbeg{3}{5}
\got{3}{A\ot(\crta A\ot A)} \gnl
\gvac{1}\gcl{1}\gvac{1}\gnl
\gvac{1}\gbmp{g}\gvac{1}\gnl
\gvac{1}\gcl{1}\gvac{1}\gnl
\gvac{1}\gob{1}{A}
\gend=
\gbeg{4}{5}
\got{1}{A}\got{1}{\crta A} \gvac{1} \got{1}{A} \gnl
\gbr \gvac{1}\gcl{1} \gnl
\gcl{1} \gwmu{3} \gnl
\gwmu{3} \gvac{1}\gnl
\gvac{1} \gob{1}{A}
\gend
\end{eqnarray}
A faithfully projective algebra $A$ in $\C$ is called an {\em Azumaya algebra} if $F$ and $G$ are isomorphisms. We call $F$ and $G$ {\em the Azumaya defining morphisms}. In a symmetric monoidal category $F$ and $G$ coincide. \par \medskip

We recall a characterization of Azumaya algebras which appears in \cite[Theorem 3.1]{VZ1} and in \cite[Proposition 1]{P4} for a symmetric monoidal category. First we observe that for two algebras $A,B \in \C$ the category of $A\x B$-bimodules is isomorphic to the category of left $A\ot\crta B$-modules. An $A\x B$-bimodule $M$ is equipped with a structure of left $A\ot\crta B$-module by
\begin{eqnarray} \label{A-B-bimod}
\gbeg{3}{4}
\got{1}{A\ot\crta B} \got{1}{} \got{1}{M} \gnl
\gcn{2}{1}{1}{3} \gcl{1} \gnl
\gvac{1} \glm \gvac{1} \gnl
\gvac{1} \gob{3}{M}
\gend=
\gbeg{3}{5}
\got{1}{A} \got{1}{B} \got{1}{M} \gnl
\gcl{1} \gbr \gnl
\glm \gcl{1} \gnl
\gvac{1} \grm \gnl
\gob{3}{M}
\gend
\end{eqnarray}
A left $A\ot\crta B$-module $N$ is made into an $A\x B$-bimodule by
\begin{eqnarray}\label{A-B-bimod2}
\gbeg{2}{3}
\got{1}{A} \got{1}{N} \gnl
\glm \gnl
\gvac{1} \gob{1}{N}
\gend =
\gbeg{3}{5}
\got{1}{A} \got{1}{} \got{1}{N} \gnl
\gcl{1} \gu{1} \gcl{1} \gnl
\glmpt \gnot{\hspace{-0,4cm}A\ot\crta B} \grmptb \gcl{1} \gnl
\gvac{1} \glm \gvac{1} \gnl
\gvac{1} \gob{3}{N}
\gend \qquad\textnormal{and}\qquad
\gbeg{2}{3}
\got{1}{N} \got{1}{B} \gnl
\grm \gnl
\gob{1}{N}
\gend =
\gbeg{3}{6}
\got{1}{} \got{1}{N} \got{1}{B} \gnl
\gvac{1} \gibr \gnl
\gu{1} \gcl{1} \gcl{1} \gnl
\glmpt \gnot{\hspace{-0,4cm}A\ot\crta B} \grmptb \gcl{1} \gnl
\gvac{1} \glm \gvac{1} \gnl
\gvac{1} \gob{3}{N}
\gend
\end{eqnarray}
Analogously, the category of $A\x B$-bimodules is also isomorphic to the category of right $\crta A \otimes B$-modules. \par \smallskip

An algebra $A \in \C$ is an $A\x A$-bimodule through the multiplication and then it becomes a left $A \otimes \crta A$-bimodule and a right $\crta A \otimes A$-module. For any $X \in \C$ the object $A \otimes X$ is a left $A \otimes \crta A$-module with the structure of left $A \otimes \crta A$-module inherited from $A$. Thus we can define the functor
$A \otimes -: \C \rightarrow {}_{A\ot\crta A}\C$. Similarly, we can define the functor $-\ot A:\C\to \C_{\crta A\ot A}$. Both are $\C$-funtors and allow to characterize Azumaya algebras.

\begin{theorem}
An algebra $A\in\C$ is Azumaya if and only if $A\ot -:\C\to {}_{A\ot\crta A}\C$ and $-\ot A:\C\to \C_{\crta A\ot A}$ establish $\C$-equivalences of categories.
\end{theorem}

The inverse functor of $A\ot -$ is its right adjoint ${}_{A\ot \crta A} [A, -]:{}_{A\ot\crta A}\C\to\C$. The algebra $A\ot\crta A$ is called {\em the enveloping algebra of $A$} and it is usually denoted by $A^e$. \par \bigskip

Two Azumaya algebras $A$ and $B$ in $\C$ are called {\em Brauer equivalent}, denoted by
$A\sim B$, if there exist faithfully projective objects $P$ and $Q$ in $\C$ such that
$$A\ot [P, P]\iso B\ot [Q, Q]$$
as algebras. This defines an equivalence relation in the (wrongly called) set $B(\C)$ of isomorphism classes of Azumaya algebras. The quotient set $\Br(\C)=B(\C)/\sim$ is a group, called the {\em Brauer group} of $\C$, with product induced by $\ot$, that is, $[A][B]=[A\otimes B]$, identity element $[I]$ or class of $[P,P]$ for a faithfully projective object $P$, and the inverse of $[A]$ is $[\crta A]$. \par \medskip

The rest of this subsection is devoted to provide an alternative description of the functor ${}_{A^e}[A,-].$

\begin{definition} \label{M^A}
For $M \in {}_{A^e}\C$, let $M^A$ denote the equalizer:
$$\scalebox{0.9}[0.9]{
\bfig
\putmorphism(-40,0)(1,0)[M^A`M`j]{340}1a
\putmorphism(330,0)(1,0)[`[A, A\ot M]`\tilde{\alpha}_{A,M}]{490}1a
\putmorphism(820,25)(1,0)[\phantom{[A, A\ot M]}`\phantom{[A,A\ot \crta A\ot M]}`{{[A,A\ot \eta_{\crta A}\ot M]}}]{1400}1a
\putmorphism(820,0)(1,0)[`[A,A\ot \crta A\ot M]`]{1400}0a
\putmorphism(820,-25)(1,0)[\phantom{[A, A\ot M]}` \phantom{[A,A\ot \crta A\ot M]}`{{[A,\eta_{A}\ot \crta A\ot M]}}]{1400}1b
\putmorphism(2220,0)(1,0)[\phantom{[A,A\ot \crta A\ot M]}`[A,M],`{{[A,\nu ]}}]{860}1a
\efig}
$$
where $\nu:A \otimes \crta A \otimes M \rightarrow M$ is the structure morphism.
\end{definition}

For simplicity set $\theta_1:= \nu(A\ot\eta_{\crta A}\ot M)$ and $\theta_2:= \nu(\eta_A\ot \crta A\ot M)$.
We have:
\begin{center}
\begin{tabular}{p{14cm}}
\begin{eqnarray} \label{M^A-left}
\scalebox{0.9}[0.9]{
\gbeg{4}{7}
\got{1}{A} \got{3}{M} \gnl
\gcl{3} \glmp\gnot{\tilde{\alpha}}\gcmptb\grmp \gvac{1} \gnl
\gvac{2}\gcl{1} \gnl
\gvac{1}\glmp\gnot{[A,\hspace{0,1cm}\theta_1]}\gcmptb\grmp \gnl
\gcn{1}{1}{1}{3}\gcn{1}{1}{3}{1} \gvac{1} \gnl
\gvac{1} \gmp{\tilde\ev} \gvac{1} \gnl
\gvac{1} \gob{1}{M} \gvac{1}
\gend} \stackrel{(\ref{rightAdj})}{=}
\scalebox{0.9}[0.9]{
\gbeg{4}{6}
\got{1}{A} \got{3}{M} \gnl
\gcl{1} \gvac{1} \gbmp{\tilde{\alpha}} \gnl
\hspace{-0,26cm} $\crtaeval$ \gnl
\gnl
\gvac{1} \gbmp{\theta_1}  \gnl
\gvac{1} \gob{1}{M}
\gend} \stackrel{(\ref{adjunit-op})}{=}\scalebox{0.9}[0.9]{
\gbeg{2}{5}
\got{1}{A\ot M} \gnl
\gcl{1}\gnl
\gbmp{\hspace{0,1cm}\theta_1}\gnl
\gcl{1}\gnl
\gob{1}{\hspace{0,12cm}M}
\gend} = \scalebox{0.9}[0.9]{
\gbeg{4}{6}
\got{1}{A} \got{4}{M} \gnl
\gcl{1} \gu{1} \gvac{1} \hspace{-0,32cm} \gcl{3} \gnl
\gvac{1} \hspace{-0,37cm} \glmpt \gnot{\hspace{-0,4cm}A\#\crta A} \grmpt  \gnl
\gvac{2} \hspace{-0,18cm} \gcn{2}{1}{1}{3} \gnl
\gvac{3} \hspace{-0,14cm} \glm \gnl
\gvac{4} \gob{1}{M}
\gend} = \scalebox{0.9}[0.9]{
\gbeg{2}{3}
\got{1}{A} \got{1}{M} \gnl
\glm \gnl
\gvac{1} \gob{1}{M}
\gend}
\end{eqnarray} \vspace{-1cm} \\
\begin{eqnarray}\label{M^A-right}
\scalebox{0.9}[0.9]{
\gbeg{4}{7}
\got{1}{A} \got{3}{M} \gnl
\gcl{3} \glmp\gnot{\tilde{\alpha}}\gcmptb\grmp \gvac{1} \gnl
\gvac{2}\gcl{1} \gnl
\gvac{1}\glmp\gnot{[A, \hspace{0,1cm}\theta_2]}\gcmptb\grmp \gnl
\gcn{1}{1}{1}{3}\gcn{1}{1}{3}{1} \gvac{1} \gnl
\gvac{1} \gmp{\tilde\ev} \gvac{1} \gnl
\gvac{1} \gob{1}{M} \gvac{1}
\gend} = \scalebox{0.9}[0.9]{
\gbeg{4}{6}
\got{1}{A} \got{3}{M} \gnl
\gcl{1} \gvac{1} \gbmp{\tilde{\alpha}} \gnl
\hspace{-0,26cm} $\crtaeval$ \gnl
\gnl
\gvac{1} \gbmp{\theta_2}  \gnl
\gvac{1} \gob{1}{M}
\gend} = \scalebox{0.9}[0.9]{
\gbeg{2}{5}
\got{1}{A\ot M} \gnl
\gcl{1}\gnl
\gbmp{\hspace{0,1cm}\theta_2}\gnl
\gcl{1}\gnl
\gob{1}{\hspace{0,12cm}M}
\gend} = \scalebox{0.9}[0.9]{
\gbeg{4}{6}
\got{3}{A} \got{1}{\hspace{-0,22cm}M} \gnl
\gu{1} \gcl{1} \gvac{1} \hspace{-0,32cm} \gcl{3} \gnl
\gvac{1} \hspace{-0,37cm} \glmpt \gnot{\hspace{-0,4cm}A\#\crta A} \grmpt  \gnl
\gvac{2} \hspace{-0,18cm} \gcn{2}{1}{1}{3} \gnl
\gvac{3} \hspace{-0,14cm} \glm \gnl
\gvac{4} \gob{1}{M}
\gend} = \scalebox{0.9}[0.9]{
\gbeg{2}{4}
\got{1}{A} \got{1}{M} \gnl
\gbr \gnl
\grm \gnl
\gob{1}{M}
\gend}
\end{eqnarray}
\end{tabular}
\end{center}

Thus we obtain a short description of the property satisfied by the equalizer $(M^A, j)$:
\begin{equation} \label{M^A-handy}
\gbeg{2}{4}
\got{1}{A} \got{2}{M^A} \gnl
\gcl{1} \gbmp{j} \gnl
\glm \gnl
\gvac{1} \gob{1}{M}
\gend=
\gbeg{2}{5}
\got{1}{A} \got{1}{\hspace{0,3cm}M^A} \gnl
\gcl{1} \gbmp{j} \gnl
\gbr \gnl
\grm \gnl
\gob{1}{M}
\gend
\end{equation}

If $f:M\to N$ is a morphism in ${}_{A^e}\C$, then $f^A: M^A\to N^A$
is induced as a morphism on an equalizer via the diagram
\begin{eqnarray}\label{f^A}
\scalebox{0.9}[0.9]{\bfig
\putmorphism(0,425)(1,0)[M^A`M`j_M]{600}1a
\putmorphism(0,0)(1,0)[N^A`N`j_N]{600}1a
\putmorphism(0,420)(0,-1)[\phantom{M^A}`\phantom{N^A}`f^A]{400}1l
\putmorphism(605,420)(0,-1)[\phantom{M}`\phantom{N}`f]{400}1r
\putmorphism(625,20)(1,0)[``{[A, \theta_1]}\tilde\alpha_{A,N}]{660}1a
\putmorphism(600,0)(1,0)[`[A, N].`]{860}0a
\putmorphism(625,-25)(1,0)[``{[A, \theta_2]}\tilde\alpha_{A,N}]{660}1b
\efig}
\end{eqnarray}
In the sequel we will adopt the following notation for the equalizer $(M^A, j_M)$:
\begin{eqnarray} \label{M^A - Psi's}
\bfig
\putmorphism(0,0)(1,0)[M^A`M`j_M]{400}1a
\putmorphism(430,25)(1,0)[``\Psi_1]{350}1a
\putmorphism(430,0)(1,0)[`[A, M],`]{520}0a
\putmorphism(430,-25)(1,0)[``\Psi_2]{350}1b
\efig
\end{eqnarray}
where $\Psi_1:= [A, \theta_1]\tilde{\alpha}_{A,N}$ and $\Psi_2:= [A, \theta_2]\tilde{\alpha}_{A,N}$.

\begin{remark}\label{morma}
To define a morphism with codomain $M^A$, for example $\hat f: Q\to M^A$,
we will have to give first a morphism $f:Q\to M$ and then check that $\Psi_1 f=\Psi_2 f$. By the equalizer property of $(M^A, j_M)$ we will always do this using the universal property of $([A, M], \tilde{\ev}:A\ot [A, M]\to M)$. We will prove that $\tilde{\ev}(A\ot\Psi_1 f)=\tilde\ev(A\ot\Psi_2 f)$, then it will follow $\Psi_1 f=\Psi_2 f$. Now by Diagrams \ref{M^A-left} and \ref{M^A-right} it is to check
\begin{equation}\label{make morf. on M^A}
\scalebox{0.9}[0.9]{
\gbeg{2}{4}
\got{1}{A} \got{1}{Q} \gnl
\gcl{1} \gbmp{f} \gnl
\glm \gnl
\gvac{1} \gob{1}{M}
\gend} = \scalebox{0.9}[0.9]{
\gbeg{2}{5}
\got{1}{A} \got{1}{Q} \gnl
\gcl{1} \gbmp{f} \gnl
\gbr \gnl
\grm \gnl
\gob{1}{M}
\gend}
\end{equation}
\end{remark}

The functor $(-)^A$ gives another description of the functor ${}_{A^e}[A, -]$ for an algebra $A$. This description is one of the key pieces to prove our main theorem (Theorem \ref{Beattie}) . The same can be said of
Propositions \ref{C-equiv-Az-t} and \ref{Azcounitmult}.

\begin{proposition}\label{Az-right-adj}
Let $A$ be an algebra in $\C$. Then, the functors ${}_{A^e}[A, -]$ and $(-)^A$ are isomorphic. In particular, if $A$ is Azumaya, the pair of functors $(A\ot -, (-)^A)$ establishes an equivalence of categories between $\C$ and ${}_{A^e}\C$.
\end{proposition}

\begin{proof}
Take $M\in{}_{A^e}\C$. We define the morphisms $g:[A, M]\to M$ and $h:M\to [A, M]$ in the following way:
$$\scalebox{0.9}[0.9]{
\gbeg{2}{5}
\got{1}{[A, M]} \gnl
\gcl{1}\gnl
\gbmp{g}\gnl
\gcl{1}\gnl
\gob{1}{M}
\gend}:= \scalebox{0.9}[0.9]{
\gbeg{4}{5}
\got{1}{} \got{3}{[A, M]} \gnl
\gu{1} \gvac{1} \gcl{1} \gnl
\gcn{1}{1}{1}{3}\gcn{1}{1}{3}{1} \gnl
\gvac{1} \gmp{\tilde\ev}\gnl
\gvac{1} \gob{1}{M}
\gend} \qquad \textrm{and} \qquad
\scalebox{0.9}[0.9]{
\gbeg{3}{5}
\got{1}{A} \got{3}{M} \gnl
\gcl{1} \gvac{1} \gbmp{h} \gnl
\hspace{-0,28cm} $\crtaeval$ \gnl
\gnl
\gvac{1} \gob{1}{M}
\gend}:=\scalebox{0.9}[0.9]{
\gbeg{3}{4}
\got{1}{A} \got{1}{M} \gnl
\glm \gnl
\gvac{1} \gcl{1} \gnl
\gvac{1} \gob{1}{M}
\gend}
$$
From the equalizer property of ${}_{A^e}[A, M]$ and the definitions of $u$ and $v$ from (\ref{Ae-u}) we
obtain the equality
\begin{eqnarray} \label{dva right adj-pom}
\gbeg{5}{5}
\got{1}{\hspace{-0,28cm}A\ot\crta A}\got{1}{A}\got{3}{{}_{A^e}[A, M]} \gnl
\glm \gvac{1} \gbmp{\iota} \gnl
\gvac{1}	\gcn{1}{1}{1}{3}\gcn{1}{1}{3}{1} \gnl
\gvac{2} \gmp{\tilde{\ev}}\gnl
\gvac{2} \gob{1}{M}
\gend=
\gbeg{5}{6}
\gvac{1} \got{1}{\hspace{-0,28cm}A\ot\crta A}\got{1}{A}\got{3}{{}_{A^e}[A, M]} \gnl
\gcn{1}{1}{2}{2} \gvac{1} \gcl{1} \gvac{1} \gbmp{\iota} \gnl
\gcn{1}{2}{2}{5} \gvac{1}\gcn{1}{1}{1}{3}\gcn{1}{1}{3}{1} \gnl
\gvac{3} \gmp{\tilde{\ev}}\gnl
\gvac{2} \glm \gnl
\gvac{3} \gob{1}{M}
\gend
\end{eqnarray}
Consider the diagram
$$\scalebox{0.9}[0.9]{
\bfig
\putmorphism(-50,420)(0,-1)[``\crta g]{400}1l
\putmorphism(0,420)(0,-1)[``\crta h]{400}{-1}r
\putmorphism(-30,425)(1,0)[{}_{A^e}[A, M]`[A, M]`\iota]{710}1a
\putmorphism(0,0)(1,0)[ M^A` M` j]{670}1b
\putmorphism(640,420)(0,-1)[`` g]{400}1l
\putmorphism(690,420)(0,-1)[`` h]{400}{-1}r
\putmorphism(680,450)(1,0)[\phantom{[A, M]}``u]{600}1a
\putmorphism(680,425)(1,0)[\phantom{[A, M]}`[(A\ot\crta A)\ot A, M]`]{1000}0a
\putmorphism(680,400)(1,0)[\phantom{[A, M]}``v]{600}1b
\putmorphism(580,25)(1,0)[\phantom{[A, M]}``\Psi_1]{900}1a
\putmorphism(680,0)(1,0)[\phantom{[A, M]}`[A, M]`]{950}0a
\putmorphism(580,-25)(1,0)[\phantom{[A, M]}``\Psi_2]{900}1b
\put(280,180){\fbox{1}}
\efig}
$$
We are going to prove that the composition $g\iota$ induces $\crta g$ and $hj$
induces $\crta h$ so that the square $\langle1\rangle$ commutes. We proceed with $\crta g$:
$$\begin{array}{ll}
\scalebox{0.9}[0.9]{
\gbeg{4}{5}
\got{1}{A} \got{3}{{}_{A^e}[A, M]} \gnl
\gcn{1}{1}{1}{3} \gvac{1} \gbmp{\iota} \gnl
\gvac{1} \gcl{1} \gbmp{g} \gnl
\gvac{1} \glm \gnl
\gvac{2} \gob{1}{M}
\gend} & = \scalebox{0.9}[0.9]{
\gbeg{5}{6}
\got{1}{A} \got{1}{} \got{3}{{}_{A^e}[A, M]} \gnl
\gcl{2} \gu{1} \gvac{1} \gbmp{\iota} \gnl
\gvac{1} \gcn{1}{1}{1}{3}\gcn{1}{1}{3}{1} \gnl
\gcn{1}{1}{1}{3} \gvac{1} \gmp{\tilde{\ev}}\gnl
\gvac{1} \glm \gnl
\gvac{2} \gob{1}{M}
\gend} \stackrel{(\ref{A-B-bimod})}{=}
\scalebox{0.9}[0.9]{
\gbeg{6}{6}
\got{1}{A} \got{2}{} \got{3}{{}_{A^e}[A, M]} \gnl
\gcl{1} \gu{1} \gu{1} \gvac{1} \gbmp{\iota} \gnl
\glmpt\gnot{\hspace{-0,34cm}A\ot\crta A}\grmptb  \gcn{1}{1}{1}{3}\gcn{1}{1}{3}{1} \gnl
\gvac{1} \gcn{1}{1}{1}{3} \gvac{1} \gmp{\tilde{\ev}}\gnl
\gvac{2} \glm \gnl
\gvac{3} \gob{1}{M}
\gend} \stackrel{(\ref{dva right adj-pom})}{=}
\scalebox{0.9}[0.9]{
\gbeg{6}{7}
\got{1}{A} \got{2}{} \got{3}{{}_{A^e}[A, M]} \gnl
\gcl{1} \gu{1} \gu{1} \gvac{1} \gcl{1} \gnl
\glmpt\gnot{\hspace{-0,36cm}A\ot\crta A}\grmptb \gcl{1} \gvac{1} \gbmp{\iota} \gnl
\gvac{1} \glm \gvac{1} \gcl{1} \gnl
\gvac{2} \gcn{1}{1}{1}{3}\gcn{1}{1}{3}{1} \gnl
\gvac{3} \gmp{\tilde{\ev}} \gnl
\gvac{3} \gob{1}{M}
\gend} \stackrel{(\ref{A-B-bimod})}{=}
\scalebox{0.9}[0.9]{
\gbeg{6}{8}
\got{1}{A} \got{2}{} \got{2}{{}_{A^e}[A, M]} \gnl
\gcl{2} \gu{1} \gu{1} \gcn{1}{1}{3}{1} \gnl
\gvac{1} \gbr \gbmp{\iota} \gnl
\gcn{1}{1}{1}{2} \gmu \gcl{2} \gnl
\gvac{1} \hspace{-0,34cm} \gmu \gnl
\gvac{2} \gcn{1}{1}{0}{2}\gcn{1}{1}{2}{0} \gnl
\gvac{3} \hspace{-0,34cm} \gmp{\tilde{\ev}} \gnl
\gvac{3} \gob{1}{M}
\gend} \vspace{3pt} \\
 & = \scalebox{0.9}[0.9]{\gbeg{3}{5}
\got{1}{A} \got{3}{{}_{A^e}[A, M]} \gnl
\gcl{1} \gvac{1} \gbmp{\iota} \gnl
\gcn{1}{1}{1}{3}\gcn{1}{1}{3}{1} \gnl
\gvac{1} \gmp{\tilde{\ev}} \gnl
\gvac{1} \gob{1}{M}
\gend}= \scalebox{0.9}[0.9]{\gbeg{5}{8}
\got{1}{} \got{1}{A} \got{4}{{}_{A^e}[A, M]} \gnl
\gu{1} \gcl{1} \gu{1} \gcn{1}{1}{3}{1} \gnl
\gcl{1} \gbr \gbmp{\iota} \gnl
\gcn{1}{1}{1}{2} \gmu \gcl{2} \gnl
\gvac{1} \hspace{-0,34cm} \gmu \gnl
\gvac{2} \gcn{1}{1}{0}{2}\gcn{1}{1}{2}{0} \gnl
\gvac{3} \hspace{-0,34cm} \gmp{\tilde{\ev}} \gnl
\gvac{3} \gob{1}{M}
\gend} \stackrel{(\ref{A-B-bimod})}{=} \scalebox{0.9}[0.9]{
\gbeg{6}{7}
\got{1}{} \got{1}{A} \got{4}{{}_{A^e}[A, M]} \gnl
\gu{1} \gcl{1} \gu{1} \gvac{1} \gcl{1} \gnl
\glmpt\gnot{\hspace{-0,36cm}A\ot\crta A}\grmptb \gcl{1} \gvac{1} \gbmp{\iota} \gnl
\gvac{1} \glm \gvac{1} \gcl{1} \gnl
\gvac{2} \gcn{1}{1}{1}{3}\gcn{1}{1}{3}{1} \gnl
\gvac{3} \gmp{\tilde{\ev}} \gnl
\gvac{3} \gob{1}{M}
\gend} \stackrel{(\ref{dva right adj-pom})}{=} \scalebox{0.9}[0.9]{
\gbeg{5}{7}
\got{1}{} \got{1}{A} \got{1}{} \got{3}{{}_{A^e}[A, M]} \gnl
\gu{1} \gcl{1} \gu{1} \gvac{1} \gbmp{\iota} \gnl
\glmpt\gnot{\hspace{-0,34cm}A\ot\crta A}\grmptb \gcn{1}{1}{1}{3}\gcn{1}{1}{3}{1} \gnl
\gvac{1} \gcn{1}{1}{1}{3} \gvac{1} \gmp{\tilde{\ev}}\gnl
\gvac{2} \glm \gnl
\gvac{3} \gob{1}{M}
\gend} \vspace{5pt} \\
 & \stackrel{(\ref{A-B-bimod})}{=}
\scalebox{0.9}[0.9]{
\gbeg{4}{7}
\got{1}{A} \got{1}{} \got{3}{{}_{A^e}[A, M]} \gnl
\gcl{2} \gu{1} \gvac{1} \gbmp{\iota} \gnl
\gvac{1} \gcn{1}{1}{1}{3}\gcn{1}{1}{3}{1} \gnl
\gcn{1}{1}{1}{3} \gvac{1} \gmp{\tilde{\ev}}\gnl
\gvac{1} \gbr \gnl
\gvac{1} \grm \gnl
\gvac{1} \gob{1}{M}
\gend} =
\scalebox{0.9}[0.9]{
\gbeg{4}{6}
\got{1}{A} \got{3}{{}_{A^e}[A, M]} \gnl
\gcn{1}{1}{1}{3} \gvac{1} \gbmp{\iota} \gnl
\gvac{1} \gcl{1} \gbmp{g} \gnl
\gvac{1} \gbr \gnl
\gvac{1} \grm \gnl
\gvac{1} \gob{1}{M}
\gend}
\end{array}$$
By Remark \ref{morma} and Diagram \ref{make morf. on M^A}, $g\iota$ induces $\crta g$ so that $j\crta g=g\iota$. Let us now prove that $\crta h$ is well defined. We compute:
$$\begin{array}{ll}
\scalebox{0.9}[0.9]{
\gbeg{5}{6}
\got{1}{\hspace{-0,28cm}A\ot\crta A}\got{1}{A}\got{3}{M^A} \gnl
\gcl{1} \gcl{1} \gvac{1} \gbmp{j} \gnl
\glm \gvac{1} \gbmp{h} \gnl
\gvac{1}\gcn{1}{1}{1}{3}\gcn{1}{1}{3}{1} \gnl
\gvac{2} \gmp{\tilde{\ev}}\gnl
\gvac{2} \gob{1}{M}
\gend} & = \scalebox{0.9}[0.9]{\gbeg{5}{5}
\got{1}{\hspace{-0,28cm}A\ot\crta A}\got{1}{A}\got{3}{M^A} \gnl
\glm \gvac{1} \gbmp{j} \gnl
\gvac{1} \gcn{1}{1}{1}{3} \gvac{1} \gcl{1} \gnl
\gvac{2} \glm \gnl
\gvac{3} \gob{1}{M}
\gend} \stackrel{(\ref{A-B-bimod})}{=}
\scalebox{0.9}[0.9]{
\gbeg{5}{7}
\got{1}{A} \got{1}{A} \got{1}{A} \got{2}{M^A} \gnl
\gcl{1} \gbr \gcl{1} \gnl
\gcn{1}{1}{1}{2} \gmu \gbmp{j} \gnl
\gvac{1} \hspace{-0,36cm} \gmu \gcn{1}{2}{2}{2} \gnl
\gvac{1} \gcn{1}{1}{2}{4} \gnl
\gvac{3} \hspace{-0,34cm} \glm \gnl
\gvac{4} \gob{1}{M}
\gend} \stackrel{ass.}{=}
\scalebox{0.9}[0.9]{
\gbeg{5}{7}
\got{1}{A} \got{1}{A} \got{1}{A} \got{2}{M^A} \gnl
\gcl{1} \gbr \gcl{1} \gnl
\gmu \gcn{1}{1}{1}{0} \gbmp{j} \gnl
\gvac{1} \hspace{-0,36cm} \gmu \gcn{1}{2}{2}{2} \gnl
\gvac{1} \gcn{1}{1}{2}{4} \gnl
\gvac{3} \hspace{-0,34cm} \glm \gnl
\gvac{4} \gob{1}{M}
\gend} \stackrel{mod.}{=} \scalebox{0.9}[0.9]{
\gbeg{5}{6}
\got{1}{A} \got{1}{A} \got{1}{A} \got{2}{M^A} \gnl
\gcl{1} \gbr \gbmp{j} \gnl
\gmu \glm \gnl
\gcn{1}{1}{2}{5} \gvac{2} \gcl{1} \gnl
\gvac{2} \hspace{-0,14cm} \glm \gnl
\gvac{3} \gob{1}{M}
\gend} \vspace{7pt} \\
 & \stackrel{(\ref{M^A-handy})}{\stackrel{mod.}{=}}
\scalebox{0.9}[0.9]{
\gbeg{5}{7}
\got{1}{A} \got{1}{A} \got{1}{A} \got{2}{M^A} \gnl
\gcl{1} \gbr \gbmp{j} \gnl
\gcl{2} \gcl{2} \gbr \gnl
\gvac{2} \grm \gnl
\gcn{1}{1}{1}{3} \glm \gnl
\gvac{1} \glm \gnl
\gvac{2} \gob{1}{M}
\gend} \stackrel{bimod.}{=} \scalebox{0.9}[0.9]{
\gbeg{5}{7}
\got{1}{A} \got{1}{A} \got{1}{A} \got{2}{M^A} \gnl
\gcl{1} \gbr \gbmp{j} \gnl
\gcl{2} \gcl{1} \gbr \gnl
\gvac{1} \glm \gcl{1} \gnl
\gcn{1}{1}{1}{3} \gvac{1} \grm \gnl
\gvac{1} \glm \gnl
\gvac{2} \gob{1}{M}
\gend} \stackrel{nat.}{=} \scalebox{0.9}[0.9]{
\gbeg{5}{7}
\got{1}{A} \got{1}{A} \got{1}{A} \got{2}{M^A} \gnl
\gcl{3} \gcl{1} \gcl{1} \gbmp{j} \gnl
\gvac{1} \gcn{1}{1}{1}{3} \glm \gnl
\gvac{2} \gbr \gnl
\gcn{1}{1}{1}{3} \gvac{1} \grm \gnl
\gvac{1} \glm \gnl
\gvac{2} \gob{1}{M}
\gend} \stackrel{(\ref{A-B-bimod})}{=} \scalebox{0.9}[0.9]{
\gbeg{5}{5}
\got{1}{} \got{1}{\hspace{-0,2cm}A\ot\crta A} \got{1}{A} \got{2}{M^A} \gnl
\gvac{1} \gcl{1} \gcl{1} \gbmp{j} \gnl
\gvac{1} \gcn{1}{1}{1}{3} \glm \gnl
\gvac{2} \glm \gnl
\gvac{3} \gob{1}{M}
\gend } \vspace{7pt} \\
 & = \scalebox{0.9}[0.9]{\gbeg{5}{7}
\gvac{1} \got{1}{\hspace{-0,28cm}A\ot\crta A}\got{1}{A}\got{3}{M^A} \gnl
\gvac{1} \gcn{1}{2}{0}{0} \gcl{1} \gvac{1} \gbmp{j} \gnl
\gvac{2} \gcl{1} \gvac{1} \gbmp{h} \gnl
\gcn{1}{2}{2}{5} \gvac{1}\gcn{1}{1}{1}{3}\gcn{1}{1}{3}{1} \gnl
\gvac{3} \gmp{\tilde{\ev}}\gnl
\gvac{2} \glm \gnl
\gvac{3} \gob{1}{M}
\gend}
\end{array}$$
Thus $hj$ induces $\crta h$ so that $\iota\crta h=hj$. We now prove that $\crta g$ and $\crta h$ are inverse to each other. We have:
$$\scalebox{0.9}[0.9]{
\gbeg{4}{7}
\got{1}{A} \got{3}{{}_{A^e}[A, M]} \gnl
\gcl{3} \gvac{1} \gbmp{\iota} \gnl
\gvac{2} \gbmp{g} \gnl
\gvac{2} \gbmp{h} \gnl
\gcn{1}{1}{1}{3}\gcn{1}{1}{3}{1} \gnl
\gvac{1} \gmp{\tilde{\ev}} \gnl
\gvac{1} \gob{1}{M}
\gend} = \scalebox{0.9}[0.9]{
\gbeg{4}{5}
\got{1}{A} \got{3}{{}_{A^e}[A, M]} \gnl
\gcl{1} \gvac{1} \gbmp{\iota} \gnl
\gcn{1}{1}{1}{3} \gvac{1} \gbmp{g} \gnl
\gvac{1} \glm \gnl
\gvac{2} \gob{1}{M}
\gend} =
\scalebox{0.9}[0.9]{
\gbeg{3}{5}
\got{1}{A} \got{3}{{}_{A^e}[A, M]} \gnl
\gcl{1} \gvac{1} \gbmp{\iota} \gnl
\gcn{1}{1}{1}{3}\gcn{1}{1}{3}{1} \gnl
\gvac{1} \gmp{\tilde{\ev}} \gnl
\gvac{1} \gob{1}{M}
\gend}$$
where the last equality was shown when discussing the definition of $\overline{g}.$ We conclude from here $\iota=hg\iota=hj\crta g=\iota\crta h\crta g.$ Since $\iota$ is a monomorphism, $\crta h\crta g=\id_{{}_{A^e}[A, M]}$. On the other hand, it is $gh=\tilde{\ev}(\eta_A\ot h)=\nu(\eta_A\ot M)=\id_M$. Composing this from the right with $j$ and using the fact that it is a monomorphism, similarly as above we obtain $\crta g\crta h=\id_{M^A}$. \par \medskip

We finally prove that the isomorphism ${}_{A^e}[A, M]\iso M^A$ is natural. Let $N\in{}_{A^e}\C$
and $f:M\to N$ be a morphism in ${}_{A^e}\C$. Observe the following diagram:
$$\scalebox{0.85}[0.85]{
\bfig
\putmorphism(-30,425)(1,0)[{}_{A^e}[A, M]`[A, M]`\iota_M]{1320}1a
\putmorphism(0,420)(0,-1)[``{{ {}_{A^e}[A, f]} }]{400}1l
\putmorphism(-30,420)(0,-1)[`{{ {}_{A^e}[A, N] }}`]{400}0l
\putmorphism(1290,420)(0,-1)[`M`g_M]{560}1r
\putmorphism(1290,-150)(0,-1)[`N`f]{420}1r
\putmorphism(0,0)(0,-1)[`[A, N]`\iota_N]{560}1l
\putmorphism(130,-565)(1,0)[``g_N]{1120}1b
\putmorphism(50,450)(2,-1)[\phantom{{}_{A^e}[A, N]} ` ` ]{700}1l
\putmorphism(50,60)(2,-1)[\phantom{{}_{A^e}[A, N]} ` ` ]{700}1l
\putmorphism(120,400)(2,-1)[\phantom{{}_{A^e}[A, N]} ` ` \crta g_M]{670}0l
\putmorphism(120,10)(2,-1)[\phantom{{}_{A^e}[A, N]} ` ` \crta g_N]{670}0l
\putmorphism(50,500)(2,-1)[\phantom{{}_{A^e}[A, N]} ` M^A` ]{670}0l
\putmorphism(50,100)(2,-1)[\phantom{{}_{A^e}[A, N]} ` N^A` ]{670}0l
\putmorphism(690,160)(0,-1)[`` f^A]{400}1l
\putmorphism(640,170)(2,-1)[\phantom{[A, M]}``]{700}1l
\putmorphism(640,-225)(2,-1)[\phantom{[A, M]}``]{700}1l
\putmorphism(720,130)(2,-1)[\phantom{[A, M]}``j_M]{700}0l
\putmorphism(700,-265)(2,-1)[\phantom{[A, M]}``j_N]{700}0l
\put(280,60){\fbox{1}}
\efig}
$$
The upper and lower triangles in this picture commute by the definitions of $\crta g_M$
and $\crta g_N$, respectively. The right inner trapeze commutes by the definition of
$f^A$. The outer diagram commutes as well, it can be seen as a juxtaposition
$$\scalebox{0.85}[0.85]{
\bfig
\putmorphism(0,420)(0,-1)[``{{ {}_{A^e}[A, f] }}]{400}1l
\putmorphism(-30,425)(1,0)[{}_{A^e}[A, M]`[A, M]`\iota_M]{710}1a
\putmorphism(-30,0)(1,0)[ {{ {}_{A^e}[A, N] }}` [A, N]` \iota_N]{710}1a
\putmorphism(665,420)(0,-1)[`` {[A, f]}]{400}1l
\putmorphism(810,425)(1,0)[`M`g_M]{600}1a
\putmorphism(1400,420)(0,-1)[`` f]{400}1r
\putmorphism(660,0)(1,0)[\phantom{[A, M]}`N`g_N]{740}1a
\efig}
$$
where the left inner parallelogram commutes by the definition of the morphism
${}_{A^e}[A, f]$, and the right one since $g:[A, -]\to \Id_{\C}$, due to the definition,
is a natural transformation. Now a diagram chasing argument provides
$j_N\crta g_N{}_{A^e}[A, f]=j_N f^A\crta g_M$, which yields that $\crta g$, and hence also $\crta h$, is a
natural transformation.
\qed\end{proof}

\begin{remark}\label{Az-counit}
The counit of the adjunction $(A\ot -, {}_{A^e}[A, -])$ is $\tilde{{\sf ev}}=\tilde{\ev}(A\ot\iota)$ and
hence the counit of the adjunction $(A\ot -, (-)^A)$ will be $\beta:=\nu(A\ot j)$.
\end{remark}
\par \smallskip

\begin{lemma}\label{Az-counit-inverse}
The unit of the adjunction $(A\ot -, (-)^A)$ evaluated at $M\in\C$ is the
morphism $\zeta_M:M\to (A\ot M)^A$ induced by $\eta_A \otimes M$, i.e., the following diagram is commutative:
\begin{eqnarray} \label{counitazu}
\scalebox{0.8}[0.8]{
\bfig
\putmorphism(100,0)(1,0)[(A\ot M)^A`A\ot M`j_{A\ot M}]{800}1b
\putmorphism(515,470)(1,-1)[M``\eta_A\ot M]{440}1r
\putmorphism(515,465)(-1,-1)[``\zeta_M]{460}1l
\efig}
\end{eqnarray}
\end{lemma}

\begin{proof}
The morphism $\eta_A\ot M$ factors through $(A\ot M)^A$ because
$$\scalebox{0.9}[0.9]{
\gbeg{4}{4}
\got{1}{A} \got{1}{M} \gnl
\gcl{1} \glmptb \gnot{\eta_A\ot M} \gcmp \grmp \gnl
\glm \gnl
\gvac{1} \gob{1}{A\ot M}
\gend} = \scalebox{0.9}[0.9]{
\gbeg{3}{4}
\got{1}{A} \got{1}{} \got{1}{M} \gnl
\gcl{1} \gu{1} \gcl{2} \gnl
\gmu \gnl
\gob{2}{A} \gob{1}{M}
\gend} = \scalebox{0.9}[0.9]{
\gbeg{3}{4}
\got{1}{} \got{1}{A} \got{1}{M} \gnl
\gu{1} \gcl{1} \gcl{2} \gnl
\gmu  \gnl
\gob{2}{A} \gob{1}{M}
\gend} \stackrel{nat.}{=} \scalebox{0.9}[0.9]{
\gbeg{4}{7}
\got{1}{A} \got{1}{} \got{1}{M} \gnl
\gcl{1} \gu{1} \gcl{2} \gnl
\gbr \gnl
\gcl{1} \gbr \gnl
\gcl{1} \gibr  \gnl
\gmu \gcl{1} \gnl
\gob{2}{A} \gob{1}{M}
\gend} = \scalebox{0.9}[0.9]{
\gbeg{4}{5}
\got{1}{A} \got{1}{M} \gnl
\gcl{1} \glmptb \gnot{\eta_A\ot M} \gcmp \grmp \gnl
\gbr \gnl
\grm \gnl
\gob{1}{A\ot M}
\gend}$$
where in the last equality we applied the right $A$-module structure of $A\ot M$. Thus $\eta_A\ot M$ induces a morphism $\zeta_M:M\to (A\ot M)^A$ so that the Diagram \ref{counitazu} commutes. We are going to prove that $\zeta_M$ is the unit of the adjunction $(A\ot -, (-)^A)$.

Let $\tilde{{\sf a}}:M\to {}_{A^e}[A, A\ot M]$ denote the unit of the adjunction $(A \ot -,  {}_{A^e}[A,-])$. Then $\iota\tilde{{\sf a}}=\tilde{\alpha}$ and $\tilde{{\sf ev}}(A \otimes \tilde{{\sf a}})=id_{A \otimes M}$. We compute:
$$\begin{array}{ll}
\tilde{{\sf ev}}(A \otimes \tilde{{\sf a}}) & = (\nabla_A\ot M)(A\ot\eta_A\ot M) \vspace{3pt} \\
 & =  \tilde{\ev}(A \otimes h_{A \otimes M}) (A\ot j_{A\ot M})(A\ot \zeta_M) \vspace{3pt} \\
&  =  \tilde{\ev}(A\ot \iota)(A \ot \overline{h}_{A \otimes M})(A\ot \zeta_M) \vspace{3pt} \\
& = \tilde{{\sf ev}}(A \otimes (\overline{h}_{A \otimes M}\zeta_M))
\end{array}$$
This implies $\tilde{{\sf a}}=\crta h\zeta_M$, proving that $\zeta_M$
is the unit of the adjunction $(A\ot -, (-)^A)$.
\qed\end{proof}
\par \medskip

In the sequel we are going to deduce some properties from the adjunction $(A\ot -, (-)^A)$.
It is clear that $A\ot -: \C\to {}_{A^e}\C$ is a $\C$-functor. Suppose that $A$ is an
Azumaya algebra. Then ${}_{A^e}[A, -]$ is as well a $\C$-functor and by Proposition
\ref{Az-right-adj} so is $(-)^A$. Therefore for every $M\in {}_{A^e}\C$
and $V\in\C$ we know that $M^A\ot V\iso(M\ot V)^A$. We explicitly construct this isomorphism.
Recall that $\beta_M=\nu_M(A\ot j): A\ot M^A\to M$ is the counit of the adjunction $(A\ot -, (-)^A)$ and
let $t_{M, V}:M^A\ot V\to (M\ot V)^A$ be the unique morphism that makes the triangle
\begin{eqnarray}\label{t_{M,V}}
\scalebox{0.9}[0.9]{\bfig
\putmorphism(130, 440)(2, -1)[\phantom{A\ot M^A\ot V}`\phantom{M\ot V}`]{850}1r
\putmorphism(135, 450)(2, -1)[\phantom{A\ot M^A\ot V}`\phantom{M\ot V}`\beta_M\ot V]{850}0r
\putmorphism(0, 400)(0, -1)[A\ot M^A\ot V`\phantom{A \ot (M\ot V)^A}`A\ot t_{M, V}]{400}1l
\putmorphism(0, 0)(1, 0)[A \ot (M\ot V)^A`M\ot V`]{900}1b
\putmorphism(0, -20)(1, 0)[\phantom{A \ot (M\ot V)^A}`\phantom{M\ot V}`\beta_{M\ot V}]{900}0b
\efig}
\end{eqnarray}
commutative. Taking into account that $\beta$'s are isomorphisms, $A\ot t_{M, V}$ is an
isomorphism as well. Further, as an Azumaya algebra, $A$ is faithfully projective, hence faithfully
flat, so $t_{M, V}$ is an isomorphism.

Using another approach, we will now find an alternative description of $t_{M, V}$. We claim that
$j_{M} \ot V$ induces a morphism $\chi_{M, V}:M^A\ot V\to (M\ot V)^A$. To this end we check:
$$
\gbeg{3}{5}
\got{1}{A} \got{2}{\hspace{-0,2cm}M^A} \got{1}{\hspace{-0,6cm}V} \gnl
\gcl{1} \gbmp{j} \gcl{1} \gnl
\hspace{-0,14cm} \gcn{2}{1}{1}{2} \hspace{-0,41cm} \glmpt \gnot{\hspace{-0,36cm}M\ot V} \grmpt  \gnl
\gvac{2} \hspace{-0,2cm} \glm \gnl
\gvac{3} \gob{1}{M\ot V}
\gend=
\gbeg{3}{4}
\got{1}{A} \got{2}{\hspace{-0,2cm}M^A} \got{1}{\hspace{-0,6cm}V} \gnl
\gcl{1} \gbmp{j} \gcl{2} \gnl
\glm \gnl
\gvac{1} \gob{1}{M} \gob{1}{V}
\gend\stackrel{j}{=}
\gbeg{3}{5}
\got{1}{A} \got{2}{\hspace{-0,2cm}M^A} \got{1}{\hspace{-0,6cm}V} \gnl
\gcl{1} \gbmp{j} \gcl{3} \gnl
\gbr \gnl
\grm \gnl
\gob{1}{M} \gvac{1} \gob{1}{V}
\gend\stackrel{nat.}{=}
\gbeg{3}{7}
\got{1}{A} \got{2}{\hspace{-0,2cm}M^A} \got{1}{\hspace{-0,6cm}V} \gnl
\gcl{1} \gbmp{j} \gcl{2} \gnl
\gbr \gnl
\gcl{1} \gbr \gnl
\gcl{1} \gibr \gnl
\grm \gcl{1} \gnl
\gob{1}{M} \gvac{1} \gob{1}{V}
\gend=
\gbeg{3}{6}
\got{1}{A} \got{2}{\hspace{-0,2cm}M^A} \got{1}{\hspace{-0,6cm}V} \gnl
\gcl{1} \gbmp{j} \gcl{1} \gnl
\hspace{-0,13cm} \gcn{2}{1}{1}{2} \hspace{-0,42cm} \glmpt \gnot{\hspace{-0,36cm}M\ot V} \grmpt  \gnl
\gvac{2} \hspace{-0,22cm} \gbr \gnl
\gvac{2} \grm \gnl
\gvac{2} \gob{1}{M\ot V}
\gend
$$
where in the first and last equalities we applied the left and right $A$-module structures of $M\ot V$.
Thus there is $\chi_{M, V}:M^A\ot V\to (M\ot V)^A$ such that $j_{M\ot V}\chi_{M, V}=j_M\ot V$. This
makes the left triangle in the following picture commutative:
$$\scalebox{0.9}[0.9]{
\bfig
\putmorphism(40, 680)(4, -1)[\phantom{A\ot M^A\ot V}`\phantom{M\ot V}` ]{1960}1r
\putmorphism(40, 730)(4, -1)[\phantom{A\ot M^A\ot V}`\phantom{M\ot V}` \beta_M\ot V]{1970}0r 
\putmorphism(-70, 730)(4, -1)[\phantom{A\ot M^A\ot V}`M\ot V` ]{1985}0a 
\putmorphism(-35, 600)(3, -2)[\phantom{A\ot M^A\ot V}`A\ot M\ot V`]{600}1r
\putmorphism(550, 220)(1, 0)[\phantom{A\ot M\ot V}`\phantom{M\ot V}`\nu_M\ot V]{1300}1a
\putmorphism(-160, 625)(2, -1)[\phantom{A\ot M^A\ot V}``A\ot j_M\ot V]{850}0r
\putmorphism(0, 600)(0, -1)[A\ot M^A\ot V`A \ot (M\ot V)^A`A\ot \chi_{M, V}]{750}1l
\putmorphism(220, 0)(2, 1)[\phantom{A\ot M^A\ot V}`\phantom{A\ot M\ot V}`]{100}1r
\putmorphism(220, -10)(2, 1)[\phantom{A\ot M^A\ot V}`\phantom{A\ot M\ot V}`A\ot j_{M\ot V}]{100}0r
\putmorphism(610, -75)(4, 1)[\phantom{A\ot M^A\ot V}`\phantom{A\ot M\ot V}` ]{830}1r
\putmorphism(610, -120)(4, 1)[\phantom{A\ot M^A\ot V}`\phantom{A\ot M\ot V}` \beta_{M\ot V}]{830}0r
\efig}
$$
The other two triangles in the picture also commute by definition of $\beta_M$ and $\beta_{M\ot V}$, bearing in mind that $\nu_{M \ot V}=\nu_M\ot V$. Now from the commutativity of the outer diagram we deduce that $\chi_{M, V}$
satisfies the same property as $t_{M, V}$. Then $\chi_{M, V}=t_{M, V}$. Thus
$j_{M\ot V}t_{M, V}=j_M\ot V$ and we can state the following:

\begin{proposition} \label{C-equiv-Az-t}
Let $A$ be an Azumaya algebra in $\C$. For every $M\in {}_{A^e}\C$ and $V\in\C$ we have a natural
isomorphism
$$M^A\ot V\iso(M\ot V)^A$$
given by $t_{M, V}:M^A\ot V\to (M\ot V)^A$ from Diagram \ref{t_{M,V}}. This isomorphism
is such that the diagram
$$\bfig
\putmorphism(0,0)(1,0)[(M\ot V)^A`M\ot V`j_{{M\ot V}}]{700}1a
\putmorphism(0,300)(0,-1)[M^A\ot V``t_{M, V}]{300}1l
\putmorphism(50,360)(2,-1)[\phantom{M^A\ot V}`\phantom{M\ot V}`]{750}1r
\putmorphism(20,420)(2,-1)[\phantom{M^A\ot V}`\phantom{M\ot V}`j_M\ot V]{780}0r
\efig
$$
commutes.
\end{proposition}

\vspace{0,3cm}

We next obtain another result that will be necessary in a subsequent section. Let $A$ and $B$ be
two Azumaya algebras and take $M\in {}_{A^e}\C$ and $N\in{}_{B^e}\C$. Then $M\ot N\in {}_{(A\ot B)^e}\C$
with the structures given by
\begin{eqnarray} \label{(AtensB)^e-mod}
\scalebox{0.9}[0.9]{
\gbeg{4}{4}
\got{1}{A\ot B} \got{3}{M\ot N} \gnl
\gcn{2}{1}{1}{3} \gcl{1} \gnl
\gvac{1} \glm \gvac{1} \gnl
\gvac{1} \gob{3}{M\ot N}
\gend} = \scalebox{0.9}[0.9]{
\gbeg{4}{4}
\got{1}{A} \got{1}{B} \got{1}{M} \got{1}{N} \gnl
\gcl{1} \gbr \gcl{1} \gnl
\glm \glm \gnl
\gvac{1}\gob{1}{M}\gvac{1}\gob{1}{N}
\gend} \qquad\textnormal{and}\qquad
\scalebox{0.9}[0.9]{
\gbeg{4}{4}
\got{1}{M\ot N} \got{3}{A\ot B} \gnl
\gcl{1} \gcn{2}{1}{3}{1} \gnl
\grm \gvac{1} \gnl
\gob{1}{M\ot N}
\gend} = \scalebox{0.9}[0.9]{
\gbeg{4}{4}
\got{1}{M} \got{1}{N} \got{1}{A} \got{1}{B} \gnl
\gcl{1} \gbr \gcl{1} \gnl
\grm \grm \gnl
\gob{1}{M}\gvac{1}\gob{1}{N}
\gend}
\end{eqnarray}
Employing the same strategy as for the proof of Proposition \ref{C-equiv-Az-t}, we prove that the
bifunctors
$$(-)^A\ot (-)^B: {}_{A^e}\C \ot {}_{B^e}\C \to \C$$
and
$$(-\ot -)^{A\ot B}: {}_{A^e}\C \ot {}_{B^e}\C \to \C$$
are isomorphic if the braiding fulfills $\Phi_{N^B, A}=\Phi^{-1}_{A,N^B}$ for every
$N\in{}_{B^e}\C$. \par \medskip

In the adjunction isomorphism ${}_{(A\ot B)^e}\C(A\ot B\ot X, Y)\iso \C(X, Y^{A\ot B})$ put
$X=M^A\ot N^B$ and $Y=M\ot N$. Let $\zeta_{M, N}:M^A\ot N^B\to (M\ot N)^{A\ot B}$ denote the
image of $(\beta_M\ot\beta_N)(A\ot \Phi_{B, M^A}\ot N^B)$ through this isomorphism. This morphism is the unique one
making the diagram
\begin{eqnarray}\label{zeta_{M,N}}
\scalebox{0.9}[0.9]{\bfig
\putmorphism(0, 400)(1, 0)[A\ot B\ot M^A\ot N^B`\phantom{A\ot M^A\ot B\ot N^A}`
  A\ot \Phi_{B, M^A}\ot N^B]{1800}1a
\putmorphism(0, 400)(0, -1)[\phantom{A\ot B\ot M^A\ot N^B}`\phantom{A\ot B\ot (M\ot N)^{A\ot B}}
  `A\ot B\ot \zeta_{M, N}]{400}1l
\putmorphism(1800, 400)(0, -1)[A\ot M^A\ot B\ot N^B` \phantom{M\ot N}`\beta_M\ot\beta_N]{400}1r
\putmorphism(0, 0)(1, 0)[A\ot B\ot (M\ot N)^{A\ot B}`M\ot N`\beta_{M\ot N}]{1800}1b
\efig}
\end{eqnarray}
commutative. Observe that, as a tensor product of Azumaya algebras, $A\ot B$ is an
Azumaya algebra. Then we have that all $\beta$'s in the diagram are isomorphisms, so
$A\ot B\ot \zeta_{M, N}$ turns out to be one as well. Further, as an Azumaya algebra,
$A\ot B$ is faithfully projective, hence faithfully flat, so $\zeta_{M, N}$ is an
isomorphism.

Analogously as in Proposition \ref{C-equiv-Az-t}, we will find another property of $\zeta_{M, N}$
similar to the one of $t_{M, V}$. We claim that $j_{M} \ot j_{N}$ induces a morphism
$\chi_{M, N}:M^A\ot N^B\to (M\ot N)^{A\ot B}$. We compute:
$$\scalebox{0.9}[0.9]{
\gbeg{4}{5}
\got{1}{\hspace{-0,3cm}A\ot B} \got{4}{\hspace{-0,2cm}M^A\ot N^B} \gnl
\gcl{1} \glmp \gcmptb \gnot{\hspace{-0,8cm} j_M\ot j_N} \grmp  \gnl
 \gcn{1}{1}{1}{3} \gvac{1} \gcl{1}  \gnl
\gvac{1} \glm \gnl
\gvac{2} \gob{1}{M\ot N}
\gend} = \scalebox{0.9}[0.9]{
\gbeg{4}{5}
\got{1}{A} \got{1}{B} \got{2}{\hspace{-0,2cm}M^A} \got{1}{\hspace{-0,4cm}N^B} \gnl
\gcl{1} \gcl{1} \gbmp{\hspace{0,14cm}j_M} \gbmp{\hspace{0,14cm}j_N} \gnl
\gcl{1} \gbr \gcl{1} \gnl
\glm \glm \gnl
\gvac{1} \gob{1}{M} \gvac{1} \gob{1}{N}
\gend} = \scalebox{0.9}[0.9]{
\gbeg{4}{6}
\got{1}{A} \got{1}{B} \got{2}{\hspace{-0,2cm}M^A} \got{1}{\hspace{-0,4cm}N^B} \gnl
\gcl{1} \gcl{1} \gbmp{\hspace{0,13cm}j_M} \gbmp{\hspace{0,14cm}j_N} \gnl
\gcl{1} \gbr \gcl{1} \gnl
\gbr \gbr \gnl
\grm \grm \gnl
\gob{1}{M} \gvac{1} \gob{1}{N}
\gend} = \scalebox{0.9}[0.9]{
\gbeg{4}{8}
\got{1}{A} \got{1}{B} \got{2}{\hspace{-0,2cm}M^A} \got{1}{\hspace{-0,4cm}N^B} \gnl
\gcl{1} \gcl{1} \gbmp{\hspace{0,14cm}j_M} \gbmp{\hspace{0,14cm}j_N} \gnl
\gcl{1} \gbr \gcl{1} \gnl
\gbr \gbr \gnl
\gcl{1} \gbr \gcl{1} \gnl
\gcl{1} \gibr \gcl{1} \gnl
\grm \grm \gnl
\gob{1}{M} \gvac{1} \gob{1}{N}
\gend} \hspace{2pt} = \hspace{2pt} \scalebox{0.9}[0.9]{
\gbeg{4}{8}
\got{1}{\hspace{-0,2cm}A\ot B} \got{2}{\hspace{-0,2cm}M^A} \got{1}{\hspace{-0,4cm}N^B} \gnl
\gcl{1} \gbmp{\hspace{0,148cm}j_M} \gbmp{\hspace{0,14cm}j_N} \gnl
\gbr \gcl{1} \gnl
\gcl{1} \gbr \gnl
\gcl{1} \gcl{1} \glmp \gnot{\hspace{-0,36cm}A\ot B} \grmp \gnl
\gcl{1} \gbr \gcl{1} \gnl
\grm \grm \gnl
\gob{1}{M} \gvac{1} \gob{1}{N}
\gend} \hspace{3pt} = \scalebox{0.9}[0.9]{
\gbeg{3}{6}
\got{1}{\hspace{-0,3cm}A\ot B} \got{4}{\hspace{-0,2cm}M^A\ot N^B} \gnl
\gcl{1} \glmp \gcmptb \gnot{\hspace{-0,8cm} j_M\ot j_N} \grmp  \gnl
 \gcn{1}{1}{1}{3} \gvac{1} \gcl{1}  \gnl
\gvac{1} \gbr \gnl
\gvac{1} \grm \gnl
\gvac{1} \gob{1}{M\ot N}
\gend}
$$
The middle three equalities are respectively due to: equalizer property of $M^A$
and $N^B$, naturality and because of the assumption $\Phi_{N^B, A}=\Phi^{-1}_{A,N^B}$
-- applying naturality and the structure from Diagram \ref{(AtensB)^e-mod}. Then there
exists $\chi_{M, N}: M^A\ot N^B\to (M\ot N)^{A\ot B}$ such that $j_{M\ot N}\chi_{M, N}=j_M\ot j_N$.
With this we have:
\begin{flushleft}
$\begin{array}{ll} \vspace{3pt}
(A\ot \Phi_{B, M}\ot N)(A\ot B\ot (j_{M\ot N}\chi_{M, N})) & =(A\ot \Phi_{B, M}\ot N^B)(A\ot B\ot j_M\ot j_N)
\vspace{3pt} \\
 & =(A\ot j_M\ot B\ot j_N)(A\ot \Phi_{B, M^A}\ot N^B)
\end{array}$
\end{flushleft}
because of naturality. This makes the left pentagram in the following picture commutative:
$$\scalebox{0.85}[0.85]{
\bfig
\putmorphism(-600, 600)(1, 0)[\phantom{A\ot B\ot M^A\ot N^B}`A \ot M^A\ot B\ot N^B`
   A\ot \Phi_{B, M^A}\ot N^B]{1700}1a
\putmorphism(1140, 600)(2, -1)[\phantom{A\ot M^A\ot V}`\phantom{M\ot N}` ]{800}1r
\putmorphism(1130, 640)(2, -1)[\phantom{A\ot M^A\ot V}`\phantom{M\ot N}` \beta_M\ot \beta_N]{800}0r
\putmorphism(500, 220)(1, 0)[\phantom{A\ot M\ot B\ot N}`M\ot N`\nu_M\ot \nu_N]{1360}1a
\putmorphism(1150, 620)(-2, -1)[\phantom{A\ot M^A\ot V}``A\ot j_M\ot B\ot j_N]{800}1l
\putmorphism(900, 625)(-1, -1)[\phantom{A\ot M^A\ot V}`A\ot M\ot B\ot N`]{400}0l 
\putmorphism(-600, 600)(0, -1)[A\ot B\ot M^A\ot N^B`A \ot B\ot (M\ot N)^{A\ot B}`A\ot B\ot \chi_{M, N}]{750}1l
\putmorphism(1340, -60)(2, 1)[\phantom{A \ot B\ot M\ot N}`\phantom{M\ot N}` ]{320}1r
\putmorphism(1370, -60)(2, 1)[\phantom{A \ot B\ot M\ot N}`\phantom{M\ot N}`
   \nu_{M\ot N}]{320}0r
\putmorphism(-700, -200)(1, 0)[\phantom{A\ot B\ot (M\ot N)^{A\ot B}}`A \ot B\ot M\ot N`
   A\ot B\ot j_{M\ot N}]{1800}1b
\putmorphism(340, 220)(2, -1)[\phantom{A\ot M\ot B\ot N}``]{830}{-1}l
\putmorphism(310, 210)(2, -1)[\phantom{A\ot M\ot B\ot N}``
   A\ot\Phi_{B, M}\ot N]{830}0l
\efig}
$$
The other two triangles in the picture also commute by the definitions of $\beta_M$ and
$\beta_N$ and the left $A\ot B$-module structure of $M\ot N$. From the commutativity
of the outer diagram, bearing in mind that by definition $\beta_{M\ot N}=\nu_{M\ot N}
(A\ot B\ot j_{M\ot N})$, we deduce that $\chi_{M, N}$ satisfies the same property as
$\zeta_{M, N}$. Then $\chi_{M, N}=\zeta_{M, N}$. Thus
$j_{M\ot N}\zeta_{M, N}=j_M\ot j_N$. We now check the naturality. \par \smallskip

Let $f:M\to M'$ and $g:N\to N'$ be morphisms in ${}_{A^e}\C$ and ${}_{B^e}\C$ respectively.
Because of the definition of $\zeta$ the diagrams $\langle1\rangle$ and $\langle4\rangle$ in
the following picture commute:
$$\scalebox{0.85}[0.85]{\bfig
\putmorphism(770, 1500)(1, 0)[M^A\ot N^B` (M\ot N)^{A\ot B}`\zeta_{M, N}]{1350}1a
\putmorphism(1400, 1100)(0, -1)[M\ot N`M'\ot N'`f\ot g]{400}1r
\putmorphism(770,290)(1, 0)[M^A\ot N^B`(M'\ot N')^{A\ot B}`\zeta_{M', N'}]{1350}1b
\putmorphism(1300, 700)(2, -1)[``]{800}1r
\putmorphism(1240, 760)(2, -1)[``j_{{M'\ot N'}}]{800}0r
\putmorphism(1560, 1240)(2, 1)[``]{230}{-1}r
\putmorphism(1530, 1160)(2, 1)[``j_{{M\ot N}}]{230}0r
\putmorphism(670, 1500)(2, -1)[``]{800}1l
\putmorphism(825, 1400)(2, -1)[``j_M\ot j_N]{800}0l
\putmorphism(970,450)(2, 1)[``]{230}1l
\putmorphism(1075,520)(2, 1)[``j_M'\ot j_N']{230}0l
 \putmorphism(2110, 1500)(0, -1)[``(f\ot g)^{A\ot B}]{1200}1r
\putmorphism(750, 1500)(0, -1)[``f^A\ot g^B]{1200}1l
\put(1350,1300){\fbox{1}}
\put(950,880){\fbox{2}}
\put(1800,880){\fbox{3}}
\put(1350,400){\fbox{4}}
\efig}
$$
Diagrams $\langle2\rangle$ and $\langle3\rangle$ commute by the definitions
of the induced morphisms $f^A, g^B$ and $(f\ot g)^{A\ot B},$ see Diagram \ref{f^A}. From the commutativity of the outer diagram we obtain that $\zeta: (-)^A\ot (-)^B\to (-\ot -)^{A\ot B}$ is a natural transformation. We can finally claim:

\begin{proposition} \label{Azcounitmult}
Let $A$ and $B$ be Azumaya algebras in $\C$. For $M\in {}_{A^e}\C$ and $N\in{}_{B^e}\C,$ if
$\Phi_{N^B, A}=\Phi^{-1}_{A,N^B}$, then there is a natural isomorphism
$\zeta_{M, N}:M^A\ot N^B\to (M\ot N)^{A\ot B}$ such that the diagram
$$\bfig
\putmorphism(0,0)(1,0)[(M\ot N)^{A\ot B}`M\ot N`j_{{M\ot N}}]{700}1a
\putmorphism(0,300)(0,-1)[M^A\ot N^B``\zeta_{M, N}]{300}1l
\putmorphism(50,360)(2,-1)[\phantom{M^A\ot N^B}`\phantom{M\ot N}`]{750}1r
\putmorphism(20,420)(2,-1)[\phantom{M^A\ot N^B}`\phantom{M\ot N}`j_M\ot j_N]{780}0r
\efig
$$
commutes.
\end{proposition}

\subsection{$H$-Azumaya algebras}\label{H-Azu}

In this subsection we present the main protagonist of our study, that is, the Brauer group of $H$-module algebras for a Hopf algebra $H \in \C$. We will start by showing how the category of left $H$-modules inherits the structure of closed monoidal category from $\C$. The tensor product of two $H$-modules $M,N\in \C$ is again an $H$-module via
\begin{eqnarray}\label{tensorH}
\gbeg{4}{4}
\got{1}{H} \got{3}{M\ot N} \gnl
\gcn{2}{1}{1}{3} \gcl{1} \gnl
\gvac{1} \glm \gvac{1} \gnl
\gvac{1} \gob{3}{M\ot N}
\gend :=
\gbeg{4}{5}
\got{2}{H} \got{1}{M} \got{1}{N} \gnl
\gcmu \gcl{1} \gcl{1} \gnl
\gcl{1} \gbr \gcl{1} \gnl
\glm \glm \gnl
\gvac{1}\gob{1}{M}\gvac{1}\gob{1}{N}
\gend
\end{eqnarray}
The following two results are not difficult to prove:
\begin{lemma}\label{mod vs algEnd}
Let $H$ be a Hopf algebra in $\C$.
\begin{enumerate}
\item[(i)] An object $M$ in $\C$ is a left $H$-module if and only if there is an algebra morphism  $\teta: H \to [M, M]$ in $\C$. If $\nu:H \otimes M \rightarrow M$ is the structure morphism, then $\theta$ is the unique morphism such that $ev(\theta \otimes M)=\nu$.

\item[(ii)] If $M$ and $N$ are left $H$-modules, then so is
$[M, N]$ with the action given by:
$$
\gbeg{3}{5}
\got{1}{H} \got{3}{[M, N]} \gnl
\gcn{1}{1}{1}{3} \gvac{1} \gcl{1} \gnl
\gvac{1} \glm \gnl
\gvac{2}\gcl{1}\gnl
\gvac{2}\gob{1}{[M, N]}
\gend=
\gbeg{6}{8}
\got{2}{H}\got{1}{[M, N]} \gnl
\gcmu  \gcl{2} \gvac{1} \gnl
\gcl{1} \gmp{S} \gnl
\gcl{1} \gbr  \gnl
\gbmp{\teta'} \gcl{1} \gbmp{\teta}  \gnl
\gcn{1}{1}{1}{0} \gmu \gnl
\hspace{-0,22cm}\gwmu{3} \gnl
\gvac{1} \gob{1}{[M, N]}
\gend
$$
where $\teta: H\to [M, M]$ and $\teta': H\to [N, N]$ are the algebra morphisms from
(i). We abused of notation in the diagram and wrote as a multiplication what is indeed the premultiplication from 2.1.3.
\end{enumerate}
\end{lemma}
Evaluating on $M$ we obtain another description of the action of $H$ on $[M,N]$:
\begin{eqnarray} \label{inner H-mod eval}
\gbeg{5}{6}
\got{1}{H} \got{3}{[M, N]} \got{1}{M} \gnl
\gcn{1}{1}{1}{3} \gvac{1} \gcl{1} \gvac{1} \gcl{2} \gnl
\gvac{1} \glm \gnl
\gvac{2} \gcn{1}{1}{1}{3}\gcn{1}{1}{3}{1} \gnl
\gvac{3} \gmp{\ev} \gnl
\gvac{3} \gob{1}{N}
\gend=
\gbeg{6}{9}
\got{2}{H}\got{1}{[M, N]}\got{3}{M} \gnl
\gcmu  \gcl{2}\gvac{1} \gcl{2}\gnl
\gcl{1}\gmp{S} \gvac{2}  \gnl
\gcl{1} \gbr \gcn{1}{1}{3}{1} \gnl
\gcl{1} \gcl{1} \glm  \gnl
\gcl{1} \gcn{1}{1}{1}{3}\gcn{1}{1}{3}{1} \gnl
\gcn{1}{1}{1}{3} \gvac{1} \gmp{\ev} \gnl
\gvac{1} \glm \gvac{1} \gnl
\gvac{2} \gob{1}{N}
\gend
\end{eqnarray}

\begin{proposition} \label{H-inner adjunction}
Let $H \in \C$ be a Hopf algebra and $M,P,Q \in \C$ left $H$-modules. Consider the adjunction isomorphism $\Teta: \C(M\ot P, Q)\cong\C(M,[P,Q])$. By restriction $\Teta$ induces an isomorphism ${}_H\C(M\ot P, Q)\cong {}_H \C(M,[P,Q])$. Analogously, for the functor $P\ot -$ we have ${}_H \C(P\ot M, Q) \cong {}_H \C(M,[P,Q])$.
\end{proposition}

Majid pointed out in \cite[Proposition 2.5]{Maj2} that for a Hopf algebra $H \in \C,$ the category of left $H$-modules ${}_{H}\C$ is monoidal, where the action on the tensor product of two $H$-modules is given as in (\ref{tensorH}) and the monoidal structure is inherited from that of $\C$. This together with the preceding proposition allows us to state:

\begin{proposition} \label{H-C closed}
Let $H \in \C$ be a Hopf algebra and assume that the braiding is $H$-linear. Then ${}_H\C$
is a closed braided monoidal category.
\end{proposition}

\begin{remark} \label{Hmod-innerhom}
Applying 2.1.3 to the closed monoidal category ${}_H\C$, we obtain that for $M\in {}_H\C$ the
inner $\hom$-object $[M, M]$ is an algebra in ${}_H\C$. In particular, it is an $H$-module
algebra in $\C$, with the $H$-module structure given by Lemma \ref{mod vs algEnd} (ii).
\end{remark}

The following proposition gives necessary and sufficient conditions for the braiding to be left $H$-linear. They will be of fundamental importance through this paper, since they are one of the pillars which hold up Theorem \ref{innerbeat}, one of our two main theorems.

\begin{proposition} \label{braid-lin}
Let $H\in\C$ be a Hopf algebra.
\begin{enumerate}
\item[(i)] The braiding $\Phi$ of $\C$ is left $H$-linear if and only if $\Phi_{H, X}=\Phi^{-1}_{X, H}$ for
any $X\in\C$ and $H$ is cocommutative. When the above condition on $\Phi$ is satisfied, we will say that the braiding is symmetric on $H\ot X$.

\item[(ii)] The braiding $\Phi$ of $\C$ is right $H$-colinear if and only if $\Phi_{H, X}=\Phi^{-1}_{X, H}$ for any $X\in\C$ and $H$ is commutative.

\item[(iii)] $H$ is commutative and $\Phi$ is left $H$-linear if and only if
$H$ is cocommutative and $\Phi$ is right $H$-colinear.
\end{enumerate}
\end{proposition}

\begin{proof}
(i) Suppose that $\Phi$ is $H$-linear. Then we have:
$$\scalebox{0.9}[0.9]{
\gbeg{2}{4} \got{2}{H} \gnl
\gcmu \gnl
\gcl{1} \gcl{1} \gnl
\gob{1}{H} \gob{1}{H}
\gend} \stackrel{unit}{=}
\scalebox{0.9}[0.9]{
\gbeg{3}{6} \got{2}{H} \gnl
\gcmu \gu{1} \gnl
\gcl{1} \gmu \gnl
\gcl{1} \gu{1} \gcn{1}{1}{0}{1} \gnl
\gmu \gcl{1} \gnl
\gob{2}{H} \gob{1}{H}
\gend} \stackrel{nat.}{=}
\scalebox{0.9}[0.9]{
\gbeg{4}{7} \got{2}{H} \gnl
\gcmu \gu{1} \gu{1} \gnl
\gcl{1} \gmu \gcl{1} \gnl
\gcl{1} \gcn{1}{1}{2}{1} \gcn{1}{1}{3}{1} \gnl
\gcl{1} \gbr \gnl \gmu \gcl{1} \gnl
\gob{2}{H} \gob{1}{H}
\gend} \stackrel{nat.}{=}
\scalebox{0.9}[0.9]{
\gbeg{4}{6} \got{2}{H} \gnl
\gcn{1}{1}{2}{2} \gvac{1} \gu{1} \gu{1} \gnl
\gcmu \gbr \gnl
\gcl{1} \gbr \gcl{1} \gnl
\gmu \gmu \gnl
\gob{2}{H} \gob{2}{H}
\gend} \stackrel{\Phi}{\stackrel{H\x lin.}{=}}
\scalebox{0.9}[0.9]{
\gbeg{4}{7} \got{2}{H} \gnl
\gcmu \gu{1} \gu{1} \gnl
\gcl{1} \gbr \gcl{1} \gnl
\gmu \gmu \gnl
\gcn{1}{1}{2}{3} \gcn{1}{1}{4}{3} \gnl
\gvac{1} \gbr \gnl
\gvac{1} \gob{1}{H} \gob{1}{H}
\gend}=
\scalebox{0.9}[0.9]{
\gbeg{2}{4} \got{2}{H} \gnl
\gcmu \gnl
\gbr \gnl
\gob{1}{H} \gob{1}{H}
\gend}
$$
i.e., $H$ is cocommutative. \par \smallskip

Since $\Phi$ is $H$-linear, so is $\Phi^{-1}$. For $X\in\C$ the object $H\ot X$ is a left $H$-module via $\nabla_H\ot X$. Now, $H$-linearity applied to $\Phi^{-1}_{H\ot X,H}$ means:
$$\scalebox{0.9}[0.9]{
\gbeg{5}{8}
\got{2}{H} \got{1}{H} \got{1}{H} \got{1}{X} \gnl
\gcmu \gcl{1} \gcl{1} \gcl{5} \gnl
\gcl{1} \gbr \gcl{1} \gnl
\gmu \gmu \gnl
\gcn{1}{1}{2}{5} \gvac{1} \gcn{1}{1}{2}{3}  \gnl
\gvac{2} \gibr \gnl
\gvac{2} \gcl{1} \gibr \gnl
\gvac{2} \gob{1}{H} \gob{1}{X} \gob{1}{H}
\gend} = \scalebox{0.9}[0.9]{
\gbeg{5}{6}
\got{2}{H} \got{1}{H} \got{1}{H} \got{1}{X} \gnl
\gcmu \gibr \gcl{1} \gnl
\gcl{1} \gbr \gibr \gnl
\gmu \gbr \gcl{1} \gnl
\gcn{1}{1}{2}{3} \gvac{1} \gcl{1} \gmu \gnl
\gob{3}{H} \gob{1}{\hspace{-0,7cm}X} \gob{1}{\hspace{-0,4cm}H}
\gend}$$
Composing this with $H\ot\eta_H\ot\eta_H\ot X$ we get:
$$\scalebox{0.9}[0.9]{
\gbeg{5}{8}
\got{2}{H} \got{1}{} \got{1}{} \got{1}{X} \gnl
\gcmu \gu{1} \gu{1} \gcl{5} \gnl
\gcl{1} \gbr \gcl{1} \gnl
\gmu \gmu \gnl
\gcn{1}{1}{2}{5} \gvac{1} \gcn{1}{1}{2}{3}  \gnl
\gvac{2} \gibr \gnl
\gvac{2} \gcl{1} \gibr \gnl
\gvac{2} \gob{1}{H} \gob{1}{X} \gob{1}{H}
\gend} = \scalebox{0.9}[0.9]{
\gbeg{5}{7}
\got{3}{H} \got{1}{} \got{1}{X} \gnl
\gcn{1}{1}{3}{2} \gvac{1} \gu{1} \gu{1} \gcl{1} \gnl
\gcmu \gibr \gcl{1} \gnl
\gcl{1} \gbr \gibr \gnl
\gmu \gbr \gcl{1} \gnl
\gcn{1}{1}{2}{3} \gvac{1} \gcl{1} \gmu \gnl
\gob{3}{H} \gob{1}{\hspace{-0,7cm}X} \gob{1}{\hspace{-0,4cm}H,}
\gend} \qquad\textnormal{implying}\qquad
\scalebox{0.9}[0.9]{
\gbeg{3}{5}
\got{2}{H} \got{1}{X} \gnl
\gcmu \gcl{2} \gnl
\gibr \gnl
\gcl{1} \gibr\gnl
\gob{1}{H} \gob{1}{X}  \gob{1}{H}
\gend} = \scalebox{0.9}[0.9]{
\gbeg{2}{4}
\got{2}{H} \got{1}{X} \gnl
\gcmu \gcl{1} \gnl
\gcl{1} \gbr \gnl
\gob{1}{H} \gob{1}{X} \gob{1}{H}
\gend}$$
Finally, composing with $\Epsilon_H\ot X\ot H$ we obtain the desired claim. \par \medskip

Conversely, suppose that $H$ is cocommutative and that $\Phi_{H, X}=\Phi^{-1}_{X, H}$
for any $X\in\C$. Pick $X, Y\in {}_H\C$. Applying the condition on $\Phi$ in the first equality we have:
$$\scalebox{0.9}[0.9]{
\gbeg{4}{7}
\got{2}{H} \got{1}{X} \got{1}{Y} \gnl
\gcmu \gcl{1} \gcl{1} \gnl
\gcl{1} \gbr \gcl{1} \gnl
\glm \glm \gnl
\gvac{1} \gcn{1}{1}{1}{3} \gvac{1} \gcl{1} \gnl
\gvac{2} \gbr \gnl
\gvac{2} \gob{1}{Y}  \gob{1}{X}
\gend} = \scalebox{0.9}[0.9]{
\gbeg{4}{7}
\got{2}{H} \got{1}{X} \got{1}{Y} \gnl
\gcmu \gcl{1} \gcl{1} \gnl
\gcl{1} \gibr \gcl{1} \gnl
\glm \glm \gnl
\gvac{1} \gcn{1}{1}{1}{3} \gvac{1} \gcl{1} \gnl
\gvac{2} \gbr \gnl
\gvac{2} \gob{1}{Y}  \gob{1}{X}
\gend} \stackrel{nat.}{=} \scalebox{0.9}[0.9]{
\gbeg{4}{7}
\got{2}{H} \got{1}{X} \got{1}{Y} \gnl
\gcmu \gcl{1} \gcl{1} \gnl
\gbr \gcl{1} \gcl{1} \gnl
\gcl{1} \glm \gcl{1} \gnl
\gcn{1}{1}{1}{3} \gvac{1} \gbr \gnl
\gvac{1} \glm \gcl{1} \gnl
\gvac{2} \gob{1}{Y} \gob{1}{X}
\gend} \stackrel{coc.}{=} \scalebox{0.9}[0.9]{
\gbeg{4}{6}
\got{2}{H} \got{1}{X} \got{1}{Y} \gnl
\gcmu \gcl{1} \gcl{1} \gnl
\gcl{1} \glm \gcl{1} \gnl
\gcn{1}{1}{1}{3} \gvac{1} \gbr \gnl
\gvac{1} \glm \gcl{1} \gnl
\gvac{2} \gob{1}{Y} \gob{1}{X}
\gend} \stackrel{nat.}{=}
\scalebox{0.9}[0.9]{
\gbeg{4}{5}
\got{2}{H} \got{1}{X} \got{1}{Y} \gnl
\gcmu \gbr \gnl
\gcl{1} \gbr \gcl{1} \gnl
\glm \glm \gnl
\gvac{1} \gob{1}{Y} \gvac{1} \gob{1}{X}
\gend}$$
i.e., $\Phi$ is left $H$-linear. \par \medskip

(ii) Dual to (i). Turn the page 180 degrees. \par \medskip

(iii) Consequence of (i) and (ii).
\qed\end{proof}

As a consequence of (i), {\it in a symmetric monoidal category the braiding is $H$-linear if and only if $H$ is cocommutative}. In the last section we will give more examples of $H$-linear braidings for certain Hopf algebras $H$ in a non-symmetric braided monoidal category. 

\begin{remark} {\em The symmetricity condition for the braiding on $H \otimes X$ occurring in (i), (ii) was already considered in the literature and gives rise to the important notion of M\"uger's center \cite[Definition 2.9]{Mu}. The monodromy of two objects
$X,Y \in \C$ is the morphism $M_{X,Y}=\Phi_{Y,X}\Phi_{X,Y}: X \otimes Y \rightarrow X \otimes Y$. The {\em M\"uger's center} of $\C$ is defined to be $Z(\C)=\{X \in \C: M_{X,Y}=\Id_{X,Y}\ \forall Y \in \C\}$. Then (i) would read as $\Phi$ is left $H$-linear if and only if $H \in Z(\C)$ and $H$ is cocommutative.}
\end{remark}

Let $H\in\C$ be a Hopf algebra and suppose that the braiding is $H$-linear. Since the category of left $H$-modules is braided, we define the {\it Brauer group of $H$-module algebras}, denoted by $\BM(\C;H)$, as $\Br({}_H \C)$. Azumaya algebras in ${}_H\C$ will be called {\em $H$-Azumaya algebras}. Notice that an object $A \in \C$ is an algebra in ${}_H\C$ if and only if $A$ is an $H$-module algebra in $\C$. For an algebra $A\in{}_H\C$, the $H$-Azumaya defining morphisms
$$F_H:A\ot\crta A\to [A, A]\quad\textnormal{and}\quad G_H:\crta A\ot A\to \crta{[A, A]}$$
in ${}_H\C$ are given respectively by
$$\scalebox{0.9}[0.9]{
\bfig
\putmorphism(0, 0)(1, 0)[A\ot \crta A`{[A, (A\ot \crta A)\ot A]}`\alpha_{A,A\ot \crta A}]{1150}1a
\putmorphism(1150, 0)(1, 0)[\phantom{[A, (A\ot \crta A)\ot A]}`[A,A]`]{1930}1a
\putmorphism(1930, 25)(1, 0)[``{{[A, \nabla_A (\nabla_A\ot A)(A\ot \Phi)]}}]{600}0a
\efig}
$$
and
$$\scalebox{0.9}[0.9]{
\bfig
\putmorphism(0, 0)(1, 0)[\crta A\ot A`{[A, A\ot (\crta A\ot A)]}`\tilde\alpha_{A,\crta A\ot A}]{1200}1a
\putmorphism(1200, 0)(1, 0)[\phantom{[A, A\ot (\crta A\ot A)]}`\crta{[A,A]}`]{1900}1a
\putmorphism(1960, 25)(1, 0)[``{{[A,\nabla_A(A\ot \nabla_A)(\Phi \ot A)]}}]{600}0a
\efig}
$$
where $\alpha_{A,-}: \Id_{\C} \to [A, -\ot A]$ and $\tilde\alpha_{A,-}: \Id_{\C} \to [A, A\ot -]$ are the units of the adjunctions $(-\ot A, [A,-])$ and $(A\ot -, [A,-])$, respectively, between ${}_H\C$ and ${}_H\C$ (the adjunction we saw in Proposition \ref{H-inner adjunction}). Note that $\alpha_{A,-}$ and $\tilde\alpha_{A,-}$ are precisely the units of the respective adjunctions between $\C$ and $\C$, only here they are evaluated on $H$-modules. Thus the morphisms $F_H$ and $G_H$ in ${}_H\C$ are indeed the Azumaya defining mor\-phisms
$F:A\ot\crta A\to [A, A]$ and $G:\crta A\ot A\to \crta{[A, A]}$ of $A$ in $\C$. Summing up, we have obtained that an $H$-module algebra $A$ is {\em $H$-Azumaya if and only if it is Azumaya in $\C$.} Now, it is easy to see that:

\begin{proposition}
The forgetful functor induces a group morphism $p:\BM(\C;H)\to \Br(\C), [A]\mapsto[A]$ by forgetting the $H$-module structure of an $H$-Azumaya algebra. This morphism splits by the group morphism $q:\Br(\C)\to\BM(\C;H)$ induced by  assigning to any Azumaya algebra the same algebra equipped with the trivial $H$-module structure.
\end{proposition}

The first aim of this paper is to compute the cokernel of the morphism $q$. We will prove that it is isomorphic to the group of $H$-Galois objects. We will recall them in the next section. Before we introduce a subgroup of $\BM(\C;H)$, whose study will be our second purpose.

\begin{definition}
The action by $H$ on an $H$-Azumaya algebra $A$ is said to be {\em inner} if there exists
a convolution invertible morphism $f:H\to A$ satisfying
$$\scalebox{0.9}[0.9]{
\gbeg{2}{4}
\got{1}{H} \got{1}{A} \gnl
\glm \gnl
\gvac{1} \gcl{1} \gnl
\gvac{1} \gob{1}{A}
\gend} = \scalebox{0.9}[0.9]{
\gbeg{4}{7}
\got{2}{H} \gvac{1} \got{1}{A} \gnl
\gcmu \gvac{1} \gcl{2} \gnl
\gbmp{f} \glmpt \gnot{\hspace{-0,2cm} f^{-1}}\grmpb \gnl
\gcl{1} \gvac{1} \gbr \gnl
\gcn{1}{1}{1}{2}  \gvac{1} \gmu \gnl
\gvac{1} \hspace{-0,2cm} \gwmu{3} \gnl
\gvac{2} \gob{1}{A}
\gend}
$$
\end{definition}

\begin{lemma}
Assume that the braiding of $\C$ is $H$-linear. The subset $$\BM_{inn}(\C; H)=\{[A] \in \BM(\C; H): A \ \textrm{has inner action}\}$$ is a subgroup of $\BM(\C; H)$.
\end{lemma}

\begin{proof}
Let $P \in {}_H\C$ be faithfully projective. Taking into account the $H$-module structure of $[P, P]$ (Lemma \ref{mod vs algEnd}(i)) and that the morphism $\teta: H \to [P, P]$ from Lemma \ref{mod vs algEnd}(ii), has $\teta S$ as convolution inverse, it is clear that the $H$-action on $[P,P]$ is inner. Let $A$ and $B$ be $H$-Azumaya algebras having inner actions with corresponding morphisms $f:H\to A$ and $g:H\to B$. We will prove that the morphism $h:=(f\ot g)\Delta_H:H\to A\ot B$ makes the action of $A\ot B$ inner. Notice first that applying in an appropriate
way four times coassociativity and once cocommutativity of $H$ we have:
\begin{eqnarray} \label{coc-coass}
\scalebox{0.9}[0.9]{
\gbeg{4}{5}
\gvac{1} \got{1}{H} \gnl
\gwcm{3}  \gnl
\hspace{-0,36cm} \gcmu \gcmu \gnl
\gcl{1} \gibr \gcl{1} \gnl
\gob{1}{H} \gob{1}{H} \gob{1}{H} \gob{1}{H}
\gend} = \scalebox{0.9}[0.9]{
\gbeg{5}{4}
\gvac{2} \got{1}{H} \gnl
\gvac{1} \gwcm{3} \gnl
\gvac{1} \hspace{-0,34cm} \gcmu \gcmu  \gnl
\gvac{1} \gob{1}{H} \gob{1}{H} \gob{1}{H} \gob{1}{H}
\gend}
\end{eqnarray}
Using this equality it is easy to prove that the convolution inverse
for $h$ is $(f^{-1}\ot g^{-1})\Delta_H$. We now compute:
$$\begin{array}{ll}
\scalebox{0.85}[0.85]{
\gbeg{3}{5}
\got{1}{H} \gvac{1} \got{1}{A\ot B} \gnl
\gcn{1}{1}{1}{3} \gvac{1} \gcl{1} \gnl
\gvac{1} \glm \gnl
\gvac{2} \gcl{1} \gnl
\gvac{2} \gob{1}{A\ot B}
\gend} & = \scalebox{0.85}[0.85]{
\gbeg{4}{5}
\got{2}{H} \got{1}{A} \got{1}{B} \gnl
\gcmu \gcl{1} \gcl{1} \gnl
\gcl{1} \gbr \gcl{1} \gnl
\glm \glm \gnl
\gvac{1}\gob{1}{A}\gvac{1}\gob{1}{B}
\gend} = \scalebox{0.85}[0.85]{
\gbeg{8}{10}
\gvac{2} \got{1}{H} \gvac{1} \got{1}{A} \gvac{2} \got{1}{B} \gnl
\gvac{1} \gwcm{3} \gcl{1} \gvac{2} \gcl{5} \gnl
\gwcm{3} \gbr \gnl
\gbmp{f} \glmp \gnot{\hspace{-0,2cm} f^{-1}}\grmptb \gcl{1} \gcn{1}{1}{1}{3} \gnl
\gcl{1} \gvac{1} \gbr \gwcm{3} \gnl
\gcn{1}{1}{1}{2} \gvac{1} \gmu \gbmp{g} \glmp \gnot{\hspace{-0,2cm} g^{-1}}\grmptb \gnl
\gvac{1} \hspace{-0,2cm} \gwmu{3} \gvac{1} \hspace{-0,36cm} \gcl{1} \gvac{1} \gbr \gnl
\gvac{2} \gcn{1}{2}{2}{2} \gvac{2} \gcn{1}{1}{1}{2}  \gvac{1} \gmu \gnl
\gvac{6} \hspace{-0,2cm} \gwmu{3} \gnl
\gvac{3} \gob{1}{A} \gvac{3} \gob{1}{B}
\gend} \stackrel{2\times nat.}{=} \scalebox{0.85}[0.85]{
\gbeg{8}{11}
\gvac{2} \got{1}{H} \gvac{3} \got{1}{A} \got{1}{B} \gnl
\gvac{1} \gwcm{3} \gvac{2} \gcl{2} \gcl{3} \gnl
\gvac{1} \hspace{-0,36cm} \gcmu \gcmu \gnl
\gvac{1} \gbmp{f} \gcl{2} \gbmp{g} \glmp \gnot{\hspace{-0,2cm} g^{-1}}\grmpb \gcn{1}{1}{2}{1} \gnl
\gvac{1} \gcl{4} \gvac{1} \gcn{1}{1}{1}{3} \gvac{1} \gbr \gcn{1}{1}{2}{1} \gnl
\gvac{2} \gcn{1}{1}{1}{3} \gvac{1} \gbr \gbr \gnl
\gvac{3} \gbr \gcn{1}{1}{1}{2} \gmu \gnl
\gvac{2} \gcn{1}{1}{3}{2} \gvac{1} \hspace{-0,28cm}\glmp \gnot{\hspace{-0,2cm} f^{-1}}\grmp \hspace{-0,08cm} \gmu \gnl
\gvac{1} \gcn{1}{1}{2}{3} \gvac{1} \gwmu{3} \gcn{1}{2}{2}{2} \gnl
\gvac{2} \gwmu{3} \gnl
\gvac{3} \gob{1}{A} \gvac{2} \gob{2}{B}
\gend} \vspace{3mm} \\
 & \stackrel{nat.}{=} \scalebox{0.85}[0.85]{
\gbeg{9}{11}
\gvac{2} \got{1}{H} \gvac{3} \got{1}{A} \got{2}{B} \gnl
\gvac{1} \gwcm{3} \gvac{2} \gcl{2} \gcn{1}{5}{2}{2} \gnl
\gvac{1} \hspace{-0,36cm} \gcmu \gcmu \gnl
\gvac{1} \gbmp{f} \gibr \glmp \gnot{\hspace{-0,2cm} g^{-1}}\grmpb \gcn{1}{1}{2}{1} \gnl
\gvac{1} \gcl{4} \gbmp{g} \gcn{1}{1}{1}{3} \gvac{1} \gbr \gnl
\gvac{2} \gcl{2} \gvac{1} \gbr \gcn{1}{1}{1}{3} \gnl
\gvac{4} \gcl{1} \glmp \gnot{\hspace{-0,2cm} f^{-1}}\grmp \gbr \gnl
\gvac{2} \gcn{1}{1}{1}{4} \gvac{1} \gmu \gvac{1} \gmu \gnl
\gvac{1} \gcn{1}{1}{1}{2} \gvac{2} \hspace{-0,34cm} \gbr \gvac{2} \gcl{1} \gnl
\gvac{2} \gwmu{3} \gwmu{4} \gnl
\gvac{3} \gob{1}{A} \gvac{2} \gob{2}{B}
\gend} \stackrel{(\ref{coc-coass})}{\stackrel{nat.}{=}} \scalebox{0.85}[0.85]{
\gbeg{9}{13}
\gvac{2} \got{1}{H} \gvac{3} \got{1}{A} \got{1}{B} \gnl
\gvac{1} \gwcm{3} \gvac{2} \gcl{2} \gcl{3} \gnl
\gvac{1} \hspace{-0,36cm} \gcmu \gcmu \gnl
\gvac{1} \gbmp{f} \gbmp{g} \gcl{1} \glmpt \gnot{\hspace{-0,2cm} g^{-1}}\grmpb \gcn{1}{1}{2}{1} \gnl
\gvac{1} \gcl{5} \gcl{5} \gcn{1}{1}{1}{3} \gvac{1} \gbr \gcn{1}{1}{2}{1} \gnl
\gvac{4} \gbr \gbr \gnl
\gvac{3} \gcn{1}{1}{3}{1} \gvac{1} \gbr \gcl{2} \gnl
\gvac{3} \gcl{2} \gvac{1} \gibr \gnl
\gvac{2} \gvac{2} \glmpb \gnot{\hspace{-0,2cm} f^{-1}}\grmp \gmu \gnl
\gvac{1} \gcl{1} \gcn{1}{1}{1}{2} \gmu \gvac{1} \gcn{1}{1}{2}{2} \gnl
\gvac{1} \gcn{1}{1}{1}{2} \gvac{1} \hspace{-0,34cm}\gbr \gvac{1} \gcn{1}{1}{3}{1} \gnl
\gvac{2} \gmu \gwmu{3} \gnl
\gvac{2} \gob{2}{A} \gvac{1} \gob{1}{B}
\gend} \stackrel{Prop. \ \ref{braid-lin}(i)}{=} \scalebox{0.85}[0.85]{
\gbeg{9}{13}
\gvac{2} \got{1}{H} \gvac{3} \got{1}{A} \got{1}{B} \gnl
\gvac{1} \gwcm{3} \gvac{2} \gcl{2} \gcl{3} \gnl
\gvac{1} \hspace{-0,36cm} \gcmu \gcmu \gnl
\gvac{1} \gbmp{f} \gbmp{g} \gcl{1} \glmpt \gnot{\hspace{-0,2cm} g^{-1}}\grmpb \gcn{1}{1}{2}{1} \gnl
\gvac{1} \gcl{5} \gcl{5} \gcn{1}{1}{1}{3} \gvac{1} \gbr \gcn{1}{1}{2}{1} \gnl
\gvac{4} \gbr \gbr \gnl
\gvac{3} \gcn{1}{1}{3}{1} \gvac{1} \gbr \gcl{2} \gnl
\gvac{3} \gcl{2} \gvac{1} \gbr \gnl
\gvac{2} \gvac{2} \glmpb \gnot{\hspace{-0,2cm} f^{-1}}\grmp \gmu \gnl
\gvac{1} \gcl{1} \gcn{1}{1}{1}{2} \gmu \gvac{1} \gcn{1}{1}{2}{2} \gnl
\gvac{1} \gcn{1}{1}{1}{2} \gvac{1} \hspace{-0,34cm}\gbr \gvac{1} \gcn{1}{1}{3}{1} \gnl
\gvac{2} \gmu \gwmu{3} \gnl
\gvac{2} \gob{2}{A} \gvac{1} \gob{1}{B}
\gend} \vspace{3mm} \\
 & \stackrel{2\times nat.}{=} \scalebox{0.85}[0.85]{
\gbeg{8}{13}
\gvac{2} \got{1}{H} \gvac{3} \got{1}{A} \got{1}{B} \gnl
\gvac{1} \gwcm{3} \gvac{2} \gcl{2} \gcl{3} \gnl
\gvac{1} \hspace{-0,36cm} \gcmu \gcmu \gnl
\gvac{1} \gbmp{f} \gbmp{g} \gcl{1} \glmpt \gnot{\hspace{-0,2cm} g^{-1}}\grmpb \gcn{1}{1}{2}{1} \gnl
\gvac{1} \gcl{5} \gcl{5} \glmpt \gnot{\hspace{-0,2cm} f^{-1}}\grmpb \gbr \gcn{1}{1}{2}{1} \gnl
\gvac{4} \gbr \gbr \gnl
\gvac{4} \gcl{2} \gbr \gcl{2} \gnl
\gvac{5} \gbr \gnl
\gvac{2} \gvac{2} \gmu \gmu \gnl
\gvac{1} \gcl{1} \gcn{1}{1}{1}{2} \gvac{1} \gcn{1}{1}{2}{0} \gvac{1} \gcn{1}{1}{2}{2} \gnl
\gvac{1} \gcn{1}{1}{1}{2} \gvac{1} \hspace{-0,34cm}\gbr \gvac{1} \gcn{1}{1}{3}{1} \gnl
\gvac{2} \gmu \gwmu{3} \gnl
\gvac{2} \gob{2}{A} \gvac{1} \gob{1}{B}
\gend} = \scalebox{0.85}[0.85]{
\gbeg{4}{7}
\gvac{1} \got{1}{H} \gvac{1} \got{1}{A\ot B} \gnl
\gwcm{3} \gcl{2} \gnl
\gbmp{h} \glmp \gnot{\hspace{-0,2cm}h^{-1}}\grmptb \gnl
\gcl{1} \gvac{1} \gbr \gnl
\gcn{1}{1}{1}{2} \gvac{1} \gmu \gnl
\gvac{1} \hspace{-0,2cm} \gwmu{3} \gnl
\gvac{1} \gob{3}{A\ot B}
\gend}
\end{array}$$
Finally, using the cocommutativity of $H$ and that $\Phi_{H, A}$ is symmetric
by Proposition \ref{braid-lin}, it is easy to show that $f^{-1}$ makes the
$H$-action of $\crta A$ inner.
\qed\end{proof}

\section{The group of Galois objects}
\setcounter{equation}{0}

In this section we introduce the third term of our sequence, the group of Galois objects. We introduced it in \cite{CF} taking advantage of the construction of the group of biGalois objects due to Schauenburg, \cite{Sch1}. We will present here the main results on Galois objects and this group without proofs to keep this paper in a moderate length. The reader is referred to \cite{CF} for them.

\subsection{Relative Hopf modules}

The category of relative Hopf modules and the functor of coinvariant objects is defined in this subsection. \par \medskip

Let $H$ be a Hopf algebra and $A$ a right $H$-comodule algebra in $\C$. If $M$ is a right $A$-module, then it is not difficult to check that
\begin{eqnarray} \label{A-mod str on MxH}
\scalebox{0.9}[0.9]{\gbeg{3}{4}
\got{1}{M\ot H} \got{3}{A} \gnl
\gcl{1} \gcn{2}{1}{3}{1} \gnl
\grm \gnl
\gob{1}{M \ot H}
\gend:=
\gbeg{3}{5}
\got{1}{M} \got{1}{H} \got{1}{A} \gnl
\gcl{1} \gcl{1} \grcm \gnl
\gcl{1} \gbr \gcl{1} \gnl
\grm \gmu \gnl
\gob{1}{M} \gvac{2} \gob{1}{\hspace{-0,4cm}H}
\gend}
\end{eqnarray}
equips $M\ot H$ with a structure of right $A$-module, usually referred to as the
{\em codiagonal} one. A {\em right relative Hopf module} (or an {\em $(A,H)$-Hopf module}) is a right
$H$-comodule and a right $A$-module $M\in\C$ such that the $H$-comodule structure morphism of $M$ is right
$A$-linear, with the  $A$-module structure on $M\ot H$ as in (\ref{A-mod str on MxH}). The
compatibility condition takes the form:
$$\scalebox{0.9}[0.9]{\gbeg{2}{4}
\got{1}{M} \got{1}{A} \gnl
\grm \gnl
\grcm \gnl
\gob{1}{M} \gob{1}{H}
\gend=
\gbeg{3}{5}
\got{1}{M} \got{3}{A} \gnl
\grcm \grcm \gnl
\gcl{1} \gbr \gcl{1} \gnl
\grm \gmu \gnl
\gob{1}{M} \gvac{2} \gob{1}{\hspace{-0,4cm}H}
\gend}$$
We will denote by $\C^H_A$ the category whose objects are right relative Hopf modules
and whose morphisms are $A$-linear $H$-colinear morphisms. One may prove using 2.2.1 (ii) and 2.2.2 that if $H$ is flat, then an equalizer in $\C$ of two morphisms in $\C^H_A$ is an equalizer in $\C^H_A$. Moreover, the forgetful functor $\U: \C^H_A\to \C$ preserves equalizers.
\par \medskip

Let $M$ be a right $H$-comodule in $\C$. The {\em object of $H$-coinvariants of $M$} is the equalizer:
$$\scalebox{0.9}[0.9]{ \bfig
 \putmorphism(0,0)(1,0)[M^{coH}`M`i]{500}1a
 \putmorphism(500,20)(1,0)[\phantom{M}``\rho_M]{600}1a
 \putmorphism(500,0)(1,0)[`M\ot H.`]{800}0a
 \putmorphism(500,-20)(1,0)[\phantom{M}`\phantom{M\ot H.}`
 M\ot \eta_H]{750}1b
 \efig}$$

Actually, $(-)^{coH}$ defines a functor from $\C^H$ to $\C$. If $f:M\to N$ is a morphism
in $\C^H$, then $f^{coH}:M^{coH}\to N^{coH}$ is defined, by the universal property of the equalizer $(N^{coH},i_N)$, as the unique morphism such that $i_Nf^{coH}=fi_M$. Clearly, the functor $(-)^{coH}$ also acts on $\C^H_A$. Indeed, it is part of an adjoint pair:

\begin{proposition}\cite[Proposition 3.3]{CF} \label{adj-Hopf mod}
With $A, H\in \C$ as above, $\F: \C \to \C^H_A, N\mapsto N\ot A$ is a
left adjoint to $\G:\C^H_A \to \C, M\mapsto M^{coH}$.
\end{proposition}

Given $N\in\C$ we view $N\ot A$ as a right $H$-comodule by the coaction $\rho_{N\ot A}=N\ot\rho_A$ and a right $A$-module by the action $\mu_{N\ot A}= N\ot\nabla_A$. Then the compatibility condition of $\C^H_A$ holds for $N\ot A$, since
$A$ is an $H$-comodule algebra. The definition of $\F$ on morphisms is clear. We next describe the isomorphism of the adjunction $\Teta$ and its inverse $\Psi:$
$$\scalebox{0.9}[0.9]{\bfig
\putmorphism(0,0)(1,0)[\C^H_A(N\ot A, M)`\C(N, M^{coH})`]{1000}0a
\putmorphism(0,25)(1,0)[\phantom{\C^H_A(N\ot A, M)}`\phantom{\C(N, M^{coH})}`\Teta]{1000}1a
\putmorphism(330,-20)(1,0)[``\Psi]{400}{-1}b
\efig}$$

Let $N\in\C$ and $M\in\C^H_A$. For $f\in\C^H_A(N\ot A, M)$ the image $\Teta(f)\in\C(N, M^{coH})$ is defined, by the universal property of $(M^{coH},i)$, as the unique morphism such that
$$\scalebox{0.9}[0.9]{
\gbeg{3}{4}
\got{1}{N} \gnl
\hspace{-0,6cm} \glmp \gnot{\Teta(f)}\gcmptb \grmp \gnl
\gvac{1} \gbmp{i} \gnl
\gvac{1} \gob{1}{M}
\gend} = \scalebox{0.9}[0.9]{
\gbeg{3}{5}
\got{1}{N} \gnl
\gcl{1} \gu{1} \gnl
\glmpt \gnot{\hspace{-0,4cm}f} \grmpt \gnl
\gvac{1} \hspace{-0,2cm} \gcl{1} \gnl
\gvac{1} \gob{1}{M}
\gend}
$$
whereas for $g\in\C(N, M^{coH})$ the image $\Psi(g)\in\C^H_A(N\ot A, M)$ is given by
$$\scalebox{0.9}[0.9]{
\gbeg{2}{5}
\got{1}{N\ot A} \gnl
\gcl{1} \gnl
\hspace{-0,56cm} \glmp \gnot{\Psi(g)}\gcmptb  \grmp \gnl
\gvac{1} \gcl{1} \gnl
\gvac{1} \gob{1}{M}
\gend} = \scalebox{0.9}[0.9]{
\gbeg{2}{5}
\got{1}{N} \got{1}{A} \gnl
\gbmp{g} \gcl{1} \gnl
\gbmp{i}  \gcl{1} \gnl
\grm \gnl
\gob{1}{M}
\gend}
$$

The unit of the adjunction $\alpha_N=\Theta_{N, N\ot A}(\id_{N\ot A}):N \rightarrow (N \otimes A)^{coH}$ is the unique morphism such that
\begin{eqnarray} \label{unit-coH}
i_{N\ot A}\alpha_N=N\ot\eta_A.
\end{eqnarray}
The counit of the adjunction $\beta_M=\Psi_{M^{coH}, M}(\id_{M^{coH}}):M^{coH} \otimes A \rightarrow M$ is given via
\begin{eqnarray} \label{counit-coH}
\beta_M=\gbeg{4}{4}
\got{2}{} \got{1}{\hspace{-0,3cm}M^{coH}} \got{1}{\hspace{-0,1cm}A} \gnl
\gvac{2} \hspace{-0,8cm} \glmp \gnot{\hspace{-0,4cm}i} \grmptb \gcl{1} \gnl
\gvac{2} \gvac{1} \grm \gnl
\gvac{2} \gvac{1} \gob{1}{M}
\gend
\end{eqnarray}

\subsection{Galois objects}

The $H$-comodule algebras $A$ for which the above pair of functors establishes an equivalence between the categories $\C$ and $\C^H_A$ are precisely the $H$-Galois objects.

\begin{definition}\label{rightGalois}
A right $H$-comodule algebra $A$ in $\C$ is called an
 $H$-Galois object if the following two conditions are satisfied:\vspace{-5pt}
\begin{enumerate}
\itemsep 0pt
\item[(i)] $A$ is faithfully flat.
\item[(ii)] The  canonical morphism $can=(\nabla_A\ot H)(A\ot \rho_A):A\ot A \rightarrow  A \ot H$ is an isomorphism.
\end{enumerate}
\end{definition}

Considering $A\ot A$ and $A\ot H$ as right $H$-comodules by the structure morphisms $A\ot\rho_A$ and $A\ot\Delta_H$ respectively, $can$ is right $H$-colinear. If we view $A\ot A$ as a right $A$-module by the structure morphism $A\ot\nabla_A$ and equip $A\ot H$ with the codiagonal structure, $can$ is also right $A$-linear. For $H$-Galois objects the object of $H$-coinvariants turns out to be trivial.

\begin{proposition}\cite[Proposition 3.9]{CF} \label{Ulbrich lemma}
Let $A \in \C$ be a right $H$-Galois object. There is an isomorphism $\crta\eta:I\to A^{coH}$
such that $i_A\crta\eta=\eta_A$. In particular,
\begin{equation}\label{acoh}
\scalebox{0.9}[0.9]{
 \bfig
 \putmorphism(0,0)(1,0)[I`A`\eta_A]{500}1a
 \putmorphism(500,20)(1,0)[\phantom{A}``\rho_A]{600}1a
 \putmorphism(500,0)(1,0)[`A\ot H`]{800}0a
 \putmorphism(500,-20)(1,0)[\phantom{A}`\phantom{A\ot H}`
 A\ot \eta_H]{750}1b
 \efig}
\end{equation}
is an equalizer.
\end{proposition}

The above definition of $H$-Galois object is stronger than Schauenburg's one \cite[Definition 3.1]{Sch1}, where an $H$-Galois object is defined to be an $H$-comodule algebra $A$ such that (\ref{acoh}) is an equalizer and $can:A \otimes A \rightarrow A \otimes H$ is an isomorphism. However, to define the group of $H$-biGalois objects faithfully flat $H$-Galois objects are considered. This is one reason to make our definition. Another one is the following result:

\begin{theorem}\cite[Theorem 3.11]{CF}\label{HopfThm}
Let $A\in\C$ be a right $H$-comodule algebra and suppose that $H$ is flat. The
following statements are equivalent:\vspace{-5pt}
\begin{enumerate}
\itemsep -3pt
\item[(i)] $A$ is a right $H$-Galois object.
\item[(ii)] The functors $\bfig
\putmorphism(0,30)(1,0)[-\ot A: \C` \C^H_A: (-)^{coH}`]{800}1a
\putmorphism(-40,-10)(1,0)[\phantom{\C^H: (-)^{coH}}`\phantom{-\ot A: \C}` ]{760}{-1}b
\efig$ establish an equivalence of categories.
\end{enumerate}
\end{theorem}

Clearly $H$ is a comodule algebra over itself by the bialgebra property so we have that
$\bfig
\putmorphism(200,30)(1,0)[-\ot H: \C` \C^H_H: (-)^{coH}`]{700}0a
\putmorphism(170,50)(1,0)[\phantom{\C^H: (-)^{coH}}`\phantom{-\ot H: \C}` ]{700}1a
\putmorphism(160,10)(1,0)[\phantom{\C^H: (-)^{coH}}`\phantom{-\ot H: \C}` ]{680}{-1}b
\efig$ is an adjoint pair of functors and the Fundamental Theorem of Hopf modules \cite[Theorem 1.1]{Lyu} states that these functors establish an equivalence if $H$ is flat. As a consequence we have:

\begin{corollary} \label{H:fl<->f.f.}
Let $H$ be a Hopf algebra in $\C$. \vspace{-5pt}
\begin{enumerate}
\itemsep 0pt
\item[(i)] If $H$ is flat, then it is faithfully flat.
\item[(ii)] If $\C$ is closed and $H$ is finite, then it is faithfully flat and faithfully projective.
\end{enumerate}
\end{corollary}

\begin{proof}
(i) It follows from the Fundamental Theorem of Hopf modules and Theorem \ref{HopfThm}. \par \medskip

(ii) From 2.1.8 and (i), $H$ is faithfully flat. By 2.5.2, $H$ is faithfully projective.
\qed\end{proof}

The next proposition will be essential in proving future results.

\begin{proposition} \label{comod-alg-iso}
Let $H\in\C$ be a flat Hopf algebra. An $H$-comodule algebra morphism $f:A\to B$
between two $H$-Galois objects $A$ and $B$ is an isomorphism.
\end{proposition}

\begin{proof}
As an $H$-Galois object, $B$ is a right $H$-comodule. Equip it with the right
$A$-module structure given by $\mu_B:=\nabla_B(B\ot f)$. With these structures, $B$ lies in $\C^H_A$,
$$\scalebox{0.9}[0.9]{
\gbeg{2}{4}
\got{1}{B} \got{1}{A} \gnl
\grm \gnl
\grcm \gnl
\gob{1}{B} \gob{1}{H}
\gend} = \scalebox{0.9}[0.9]{
\gbeg{3}{5}
\got{1}{B} \got{1}{A} \gnl
\gcl{1} \gbmp{f} \gnl
\gmu \gnl
\gvac{1} \hspace{-0,2cm} \grcm \gnl
\gvac{1} \gob{1}{B} \gob{1}{H}
\gend} \stackrel{comod.}{\stackrel{alg.}{=}}
\scalebox{0.9}[0.9]{\gbeg{4}{6}
\got{1}{B} \got{1}{} \got{1}{A} \gnl
\gcl{1} \gvac{1} \gbmp{f} \gnl
\grcm \grcm \gnl
\gcl{1} \gbr \gcl{1} \gnl
\grm \gmu \gnl
\gob{1}{B} \gvac{2} \gob{1}{\hspace{-0,4cm}H}
\gend} \stackrel{f}{\stackrel{H\x colin.}{=}}
\scalebox{0.9}[0.9]{\gbeg{4}{6}
\got{1}{B} \got{1}{} \got{1}{A} \gnl
\grcm \grcm \gnl
\gcl{1} \gcl{1} \gbmp{f} \gcl{1} \gnl
\gcl{1} \gbr \gcl{1} \gnl
\gmu \gmu \gnl
\gob{2}{B} \gvac{1} \gob{1}{\hspace{-0,4cm}H}
\gend} \stackrel{nat.}{=}
\scalebox{0.9}[0.9]{\gbeg{4}{6}
\got{1}{B} \got{1}{} \got{1}{A} \gnl
\grcm \grcm \gnl
\gcl{1} \gbr \gcl{1} \gnl
\gcl{1} \gbmp{f} \gcl{1} \gcl{1} \gnl
\gmu \gmu \gnl
\gob{2}{B} \gvac{1} \gob{1}{\hspace{-0,4cm}H}
\gend} \stackrel{A\x act.}{=}
\scalebox{0.9}[0.9]{\gbeg{3}{5}
\got{1}{B} \got{1}{} \got{1}{A} \gnl
\grcm \grcm \gnl
\gcl{1} \gbr \gcl{1} \gnl
\grm \gmu \gnl
\gob{1}{B} \gvac{1} \gob{2}{H}
\gend}
$$
Having that $A$ is an $H$-Galois object, by the preceding theorem the counit $\beta$ of the adjunction
$\bfig
\putmorphism(0,30)(1,0)[-\ot A: \C` \C^H_A: (-)^{coH}`]{800}1a
\putmorphism(-40,-10)(1,0)[\phantom{\C^H: (-)^{coH}}`\phantom{-\ot A: \C}` ]{760}{-1}b
\efig$, given in (\ref{counit-coH}), is an isomorphism. Let $\crta\eta: I\to B^{coH}$ be the isomorphism from Proposition \ref{Ulbrich lemma}. Then we obtain that
$$\beta_B(\crta\eta\ot A)=
\scalebox{0.9}[0.9]{\gbeg{2}{5}
\got{1}{I} \got{1}{A} \gnl
\gbmp{\crta\eta} \gcl{1} \gnl
\gbmp{i}  \gcl{1} \gnl
\grm \gnl
\gob{1}{B}
\gend} =
\scalebox{0.9}[0.9]{\gbeg{2}{4}
\got{1}{} \got{1}{A} \gnl
\gu{1} \gcl{1} \gnl
\grm \gnl
\gob{1}{B}
\gend} = \scalebox{0.9}[0.9]{
\gbeg{2}{4}
\got{1}{} \got{1}{A} \gnl
\gu{1} \gbmp{f}\gnl
\gmu \gnl
\gob{2}{B}
\gend} = \scalebox{0.9}[0.9]{
\gbeg{2}{5}
\got{1}{A} \gnl
\gcl{1} \gnl
\gbmp{f}\gnl
\gcl{1} \gnl
\gob{1}{B}
\gend}
$$
is an isomorphism.
\qed\end{proof}

For a cocommutative Hopf algebra $H$ the group of $H$-Galois objects will be the set of isomorphism classes of $H$-Galois objects. The multiplication in this group will be induced by the cotensor product. The {\em cotensor product over $H$} of a right $H$-comodule $M$ and a left $H$-comodule $N$ is the equalizer
\begin{eqnarray} \label{e-def}
\bfig
\putmorphism(0,0)(1,0)[M\Box_H N`M\ot N`e]{600}1a
\putmorphism(600,25)(1,0)[\phantom{M\ot N}``\rho_M\ot N]{950}1a
\putmorphism(600,0)(1,0)[`M\ot H\ot N,`]{1200}0a
\putmorphism(600,-25)(1,0)[\phantom{M\ot N}`\phantom{M\ot H\ot N,}`M\ot \lambda_N]{1200}1b
\efig
\end{eqnarray}
where $\rho_M$ and $\lambda_N$ are the structure morphisms of $M$ and $N$ respectively. Let $H \in \C$ be a cocommutative and flat Hopf algebra and $B$ a right $H$-comodule algebra. Suppose that  $\Phi_{B,H}=\Phi_{H,B}^{-1}$. Then $B$ is a left $H$-comodule algebra (indeed a bicomodule) with structure morphism $\lambda_B=\Phi^{-1}_{H,B}\rho_B$. If $A$ is another $H$-comodule algebra, then $A\Box_H B$ is an $H$-comodule algebra and the equalizer morphism $e_{A, B}: A\Box_H B\to A\ot B$ is an $H$-comodule algebra morphism.

\begin{theorem}\cite[Theorem 4.14]{CF}\label{abelian subgroup}
Let $H \in \C$ be a cocommutative and flat Hopf algebra. Assume that for any pair of $H$-Galois objects $A,B$ the condition
$\Phi_{A,B}=\Phi_{B,A}^{-1}$ holds. Then the set of isomorphism classes of $H$-Galois objects, denoted by $\Gal(\C;H)$, is an abelian  group. The multiplication is given by $[A][B]=[A \Box_H B]$, the identity element is [H], and the inverse of $[A]$ is $[\crta A],$ where  $\crta A$ is the opposite algebra of $A$ endowed with the following right $H$-comodule structure
$$\rho_{\crta A}=
\gbeg{2}{3}
\got{1}{\crta A} \gnl
\grcm  \gnl
\gob{1}{\crta A} \gob{1}{H}
\gend:=
\gbeg{2}{4}
\got{1}{A} \gnl
\grcm \gnl
\gcl{1} \gnot{S} \gmp \gnl
\gvac{1} \gob{1}{\hspace{-2,6cm}A} \gob{1}{\hspace{-2,5cm}H}
\gend
$$
\end{theorem}
\par \medskip

An $H$-Galois object which is isomorphic to $H$ as a right $H$-comodule is called
{\em an $H$-Galois object with a normal basis}. The subset of $\Gal(\C; H)$ consisting of
$H$-Galois objects with a normal basis, denoted by $\Gal_{nb}(\C; H)$, is a subgroup. We must stress
here that the hypothesis required for $H$-Galois objects in the above theorem is automatically fulfilled by this kind
of Galois objects. The reason is Schauenburg's observation (\cite{Sch}) that for a cocommutative Hopf algebra $H$ it is
$\Phi_{H,H}=\Phi_{H,H}^{-1}.$ Hence the group of Galois objects with a normal basis may be defined
in any braided monoidal category with equalizers. \par \medskip

Let $\Hc(\C;H,I)$ be the second Sweedler cohomology group with values in the unit object $I$. We refer to \cite[Section 5]{CF} for the construction of such a group in a braided monoidal category. We must stress that it is possible to construct this group thanks to
Schauenburg's observation. Given a $2$-cocycle $\sigma:H \otimes H \rightarrow I$ we define $H_{\sigma}:=H$ as a right $H$-comodule. Then $H_{\sigma}$ is a right $H$-comodule algebra with multiplication and unit given by \begin{eqnarray} \label{H-sigma structure}
\gbeg{2}{3}
\got{1}{H_{\sigma}} \got{1}{H_{\sigma}} \gnl
\gmu \gnl
\gob{2}{H_{\sigma}}
\gend=
\gbeg{5}{5}
\got{2}{H} \got{2}{H} \gnl
\gcmu \gcmu \gvac{1} \gnl
\gcl{1} \gbr \gcl{1} \gnl
\glmptb \gnot{\hspace{-0,3cm} \sigma} \grmpt \gmu \gnl
\gvac{2} \gob{2}{H}
\gend
\qquad\textnormal{and}\qquad
\gbeg{1}{4}
\got{1}{} \gnl
\gu{1} \gnl
\gcl{1} \gnl
\gob{1}{H_{\sigma}}
\gend=
\gbeg{1}{5}
\got{1}{} \got{1}{} \gnl
\gu{1} \gu{1} \gnl
\glmptb \gnot{\hspace{-0,2cm} \sigma^{-1}}\grmpt \gnl
\gvac{1} \hspace{-0,2cm} \gu{1} \gnl
\gob{3}{H}
\gend
\end{eqnarray}
where $\sigma^{-1}$ is the convolution inverse of $\sigma$. Moreover, if $\sigma$ is normalized, then the unit on $H_{\sigma}$ coincides with $\eta_H$. \par \medskip

\begin{theorem}\cite[Theorem 5.14]{CF}\label{cohpical}
Each $H_{\sigma}$ is an $H$-Galois object and the map $$\zeta:\Hc(\C;H,I) \rightarrow \Gal_{nb}(\C;H), [\sigma] \mapsto [H_{\sigma}]$$ is an isomorphism. 
\end{theorem}

\section{Beattie's sequence in a braided monoidal category}
\setcounter{equation}{0}

In this section we will assume that {\it $\C$ is a closed braided monoidal category with equalizers and coequalizers, $H\in\C$ is a flat Hopf algebra and the braiding of $\C$ is $H$-linear.} Nevertheless, for some of the results not all of the assumptions will be necessary.

\subsection{The map $\Upsilon$ assigning an $H$-Galois object to an $H$-Azumaya algebra}

We start by recalling from \cite[Proposition 2.3]{Maj3} the definition of the smash product algebra. For a left $H$-module algebra $A$ in $\C$ the {\em smash product algebra} $A\# H$ is defined as follows: as an object $A\# H=A\ot H$, the multiplication and unit are given by
$$\gbeg{4}{4}
\got{3}{A\# H} \got{1}{A\# H} \gnl
\gvac{1}\gwmu{3} \gnl
\gvac{2} \gcl{1} \gvac{1} \gnl
\gvac{2} \hspace{-0,5cm} \gob{2}{A\# H}\gvac{1}
\gend :=
\gbeg{4}{6}
\got{1}{A} \got{2}{H} \got{1}{A} \got{1}{H} \gnl
\gcl{1} \gcmu \gcl{1} \gcl{1} \gnl
\gcl{1} \gcl{1} \gbr \gcl{1} \gnl
\gcl{1} \glm \gmu \gnl
\gwmu{3}\gvac{1} \hspace{-0,22cm} \gcl{1} \gnl
\gob{4}{A}  \gob{1}{H}
\gend \qquad \textnormal{and} \qquad
\gbeg{1}{4}
\got{1}{}  \gnl
\gu{1} \gnl
\gcl{1} \gnl
\gob{1}{A\# H}
\gend :=
\gbeg{1}{4}
\got{1}{}  \gnl
\gu{1} \gu{1} \gnl
\gcl{1} \gcl{1} \gnl
\gob{1}{A} \gob{1}{H}
\gend$$
This algebra becomes a right $H$-comodule algebra with the structure of right $H$-comodule given by
$A\ot\Delta_H$. It also admits a structure of $A\x A$-bimodule as stated in the following result.
Its proof is easy.

\begin{lemma}
The object $A\# H$ has a structure of a left $A\ot \crta A$-module via:
\begin{equation} \label{A-Aop-AsmashH}
\gbeg{4}{4}
\got{1}{A\ot \crta A} \got{3}{A\# H} \gnl
\gcn{2}{1}{1}{3} \gcl{1} \gnl
\gvac{1} \glm \gvac{1} \gnl
\gvac{1} \gob{3}{A\# H}
\gend :=
\gbeg{7}{7}
\got{1}{A} \got{3}{\crta A} \got{1}{\hspace{-0,8cm}A} \got{1}{\hspace{-0,8cm}H} \gnl
\gcl{1} \gu{1}  \gbr \gcl{1} \gu{1} \gnl
\gcl{1} \gcl{1}  \gcl{1} \gbr  \gcl{1} \gnl\glmpt \gnot{\hspace{-0,4cm}A\# H} \grmpt
   \glmpt \gnot{\hspace{-0,4cm}A\# H} \grmpt
  \glmpt \gnot{\hspace{-0,4cm}A\# H} \grmpt  \gnl
\gvac{1} \hspace{-0,2cm}\gcl{1} \gvac{1} \gwmu{3} \gnl
\gvac{1} \gwmu{4} \gnl
\gvac{2} \hspace{-0,22cm} \gob{2}{A\# H}
\gend = \gbeg{7}{10}
\got{1}{A} \got{1}{\crta A} \got{1}{A} \got{1}{H} \gnl
\gcl{7} \gbr \gcl{1} \gnl
\gvac{1} \gcl{4} \gbr  \gnl
\gvac{2} \hspace{-0,2cm}\gcmu \gcn{1}{1}{0}{1} \gnl
\gvac{2} \gcn{1}{1}{1}{1} \gbr \gnl
\gvac{2} \glm \gcl{4} \gnl
\gvac{2} \gcn{1}{1}{0}{1} \gcl{1}   \gnl
\gvac{2} \gmu \gnl
\gvac{1} \hspace{-0,2cm} \gwmu{3}  \gnl
\gvac{2} \gob{1}{A} \gvac{1} \gob{2}{H}
\gend
\end{equation}
\end{lemma}
As a consequence, according to Diagram \ref{A-B-bimod2}, the left and right $A$-module structures of $A\# H$
are given respectively  by
\begin{equation}\label{left-right mod Asm.H}
\nu_{A\# H}^A=\gbeg{3}{5}
\got{1}{A} \got{1}{} \got{1}{A\# H} \gnl
\gcl{1} \gu{1} \gcl{1} \gnl
\glmpt \gnot{\hspace{-0,4cm}A\ot\crta A} \grmptb \gcl{1} \gnl
\gvac{1} \glm \gvac{1} \gnl
\gvac{1} \gob{3}{A\# H}
\gend =
\gbeg{5}{3}
\got{1}{A} \got{4}{A\# H} \gnl
\gwmu{3} \gcl{1} \gnl
\gvac{1} \gob{1}{A} \gob{3}{H}
\gend
\qquad\textnormal{and}\qquad
\mu_{A\# H}^A=\gbeg{4}{6}
\gvac{1} \got{3}{\hspace{-0,2cm} A\#H} \got{1}{\hspace{-0,6cm}A} \gnl
\gvac{2} \gibr  \gnl
\gvac{1} \gu{1} \gcl{1} \gcl{1} \gnl
\gvac{1} \glmpt \gnot{\hspace{-0,4cm}A\ot\crta A} \grmptb \gcl{1} \gnl
\gvac{2} \glm \gnl
\gvac{3} \gob{1}{A\# H}
\gend=
\gbeg{5}{7}
\got{4}{A\# H} \got{1}{\hspace{-0,4cm}A} \gnl
\gvac{1} \gcl{3} \hspace{-0,2cm}\gcmu  \gcl{1} \gnl
\gvac{2} \gcl{1} \gbr \gnl
\gvac{2} \glm \gcl{3} \gnl
\gvac{1} \gcn{1}{1}{2}{3} \gvac{1} \gcl{1}  \gnl
\gvac{2} \gmu \gnl
\gvac{2} \gob{2}{A} \gob{1}{H}
\gend
\end{equation}

The first goal of this subsection is to prove that if $A$ is an $H$-Azumaya algebra in $\C$, then $(A\# H)^A$
is an $H$-Galois object ($H$ flat). The proof will be completed once we prove:
\begin{enumerate}
\item[1)] $(A\# H)^A$ is an $H$-comodule algebra;
\item[2)] $(A\# H)^A$ is faithfully flat;
\item[3)] $can_{(A\# H)^A}$ is an isomorphism.
\end{enumerate}
\vspace{0.4cm}

We proceed to prove 1). For this part the algebra $A$ does not have to be Azumaya.

\begin{lemma} \label{(AsmashH)A-alg}
Let $A$ be an $H$-module algebra in $\C$. Then $(A\# H)^A$ has an algebra structure such that
the equalizer morphism $j:(A\# H)^A \to A\# H$ is an algebra morphism.
\end{lemma}

\begin{proof}
We will show that the multiplication $\nabla_{A\# H}$ on $A\# H$ induces a
multiplication on $(A\# H)^A$. We will prove that
$f:=\nabla_{A\# H}(j\ot j): (A\# H)^A \ot (A\# H)^A\to A\# H$ satisfies the identity from Diagram
\ref{make morf. on M^A}, with $Q=(A\# H)^A \ot (A\# H)^A$, $M=A\# H$ and $j=j_M$.

Bear in mind the $A\ot\crta A$-module structure of $A\# H$ (Diagram \ref{A-Aop-AsmashH}) and consider it as an $A$-bimodule, Diagram \ref{left-right mod Asm.H}. According to (\ref{M^A-handy}), $(A\# H)^A$ then satisfies
\begin{eqnarray}\label{Asmash etaH-prop.}
\scalebox{0.9}[0.9]{
\gbeg{4}{5}
\got{1}{A} \got{1}{} \got{1}{(A\# H)^A} \gnl
\gcl{1} \gu{1} \gbmp{j} \gvac{1} \gnl
\glmptb \gnot{\hspace{-0,4cm}A\# H} \grmpt \gcl{1} \gnl
\gwmu{3} \gnl
\gvac{1} \gob{1}{A\# H}
\gend} = \scalebox{0.9}[0.9]{
\gbeg{5}{6}
\got{1}{\hspace{-0,2cm}A} \got{1}{\hspace{0.5cm}(A\# H)^A} \gnl
\gcl{1} \gbmp{j} \gvac{1} \gnl
\gbr \gu{1} \gnl
\gcl{1} \glmpt \gnot{\hspace{-0,4cm}A\# H} \grmptb \gnl
\gwmu{3} \gnl
\gvac{1} \gob{1}{A\# H}
\gend}
\end{eqnarray}
By successive applications of the associativity of $A\# H$ and equality (\ref{Asmash etaH-prop.}) we have:
$$\begin{array}{ll}
\scalebox{0.9}[0.9]{
\gbeg{5}{5}
\got{1}{A} \got{1}{} \got{1}{(\hspace{-2pt}A\# H\hspace{-2pt})^{^A}}
\got{3}{\hspace{5mm}(\hspace{-2pt}A\# H\hspace{-2pt})^{^A}} \gnl
\gcl{1} \gu{1} \gbmp{j} \gvac{1} \gbmp{j} \gnl
\glmpt \gnot{\hspace{-0,4cm}A\# H} \grmptb \gwmu{3} \gnl
\gvac{1} \gwmu{3} \gnl
\gvac{2} \gob{1}{A\# H}
\gend} & \hspace{-2pt}\stackrel{ass.}{=} \scalebox{0.9}[0.9]{
\gbeg{5}{6}
\got{1}{A} \got{1}{} \got{1}{(\hspace{-2pt}A\# H\hspace{-2pt})^{^A}} \got{3}{\hspace{5mm}(\hspace{-2pt}A\# H\hspace{-2pt})^{^A}} \gnl
\gcl{1} \gu{1} \gbmp{j} \gvac{1} \gbmp{j} \gnl
\glmpt \gnot{\hspace{-0,4cm}A\# H} \grmpt \gcl{1} \gvac{1} \gcl{2} \gnl
\gwmu{3} \gnl
\gvac{1} \gwmu{4} \gnl
\gvac{2} \gob{3}{A\# H}
\gend} \hspace{2mm}\stackrel{(\ref{Asmash etaH-prop.})}{=}\hspace{2mm}
\scalebox{0.9}[0.9]{
\gbeg{4}{7}
\got{1}{\hspace{-0,1cm}A} \got{1}{\hspace{4mm}(\hspace{-2pt}A\# H\hspace{-2pt})^{^A}} \got{3}{\hspace{7mm}(\hspace{-2pt}A\# H\hspace{-2pt})^{^A}} \gnl
\gcl{1} \gbmp{j} \gvac{1} \gbmp{j} \gnl
\gbr \gu{1} \gcl{3} \gnl
\gcl{1} \glmpt \gnot{\hspace{-0,4cm}A\# H} \grmptb \gnl
\gwmu{3} \gnl
\gvac{1} \gwmu{3} \gnl
\gvac{2} \gob{1}{A\# H}
\gend} \hspace{2mm}\stackrel{ass.}{=}\hspace{2mm}
\scalebox{0.9}[0.9]{
\gbeg{4}{7}
\got{1}{\hspace{-2mm}A} \got{1}{\hspace{4mm}(\hspace{-2pt}A\# H\hspace{-2pt})^{^A}} \got{3}{\hspace{7mm}(\hspace{-2pt}A\# H\hspace{-2pt})^{^A}} \gnl
\gcl{1} \gbmp{j} \gvac{1} \gbmp{j} \gnl
\gbr \gu{1} \gcl{2} \gnl
\gcl{1} \glmptb \gnot{\hspace{-0,4cm}A\# H} \grmpt \gnl
\gcl{1} \gwmu{3} \gnl
\gwmu{3} \gnl
\gvac{1} \gob{1}{A\# H}
\gend} \hspace{2mm}\stackrel{(\ref{Asmash etaH-prop.})}{=}\hspace{2mm}
\scalebox{0.9}[0.9]{
\gbeg{4}{8}
\got{1}{\hspace{-2mm}A} \got{1}{\hspace{3mm}(\hspace{-2pt}A\# H\hspace{-2pt})^{^A}} \got{3}{\hspace{7mm}(\hspace{-2pt}A\# H\hspace{-2pt})^{^A}} \gnl
\gcl{1} \gbmp{j} \gvac{1} \gbmp{j} \gnl
\gbr \gcn{1}{1}{3}{1} \gnl
\gcl{1} \gbr \gu{1} \gnl
\gcl{1} \gcl{1} \glmpt \gnot{\hspace{-0,4cm}A\# H} \grmpt \gnl
\gcl{1} \gwmu{3} \gnl
\gwmu{3} \gnl
\gvac{1} \gob{1}{A\# H}
\gend} \vspace{1mm} \\
& \hspace{2mm}\stackrel{ass.}{=}\hspace{2mm}
\scalebox{0.9}[0.9]{
\gbeg{4}{8}
\got{1}{\hspace{-2mm}A} \got{1}{\hspace{3mm}(\hspace{-2pt}A\# H\hspace{-2pt})^{^A}} \got{3}{\hspace{7mm}(\hspace{-2pt}A\# H\hspace{-2pt})^{^A}} \gnl
\gcl{1} \gbmp{j} \gvac{1} \gbmp{j} \gnl
\gbr \gcn{1}{1}{3}{1} \gnl
\gcl{1} \gbr \gu{1} \gnl
\gcl{1} \gcl{1} \glmpt \gnot{\hspace{-0,4cm}A\# H} \grmpt \gnl
\gmu \gvac{1}  \hspace{-0,2cm}\gcl{1} \gnl
\gvac{1} \hspace{-0,14cm} \gwmu{3} \gnl
\gvac{2} \gob{1}{A\# H}
\gend} \hspace{2mm}\stackrel{nat.}{=}
\scalebox{0.9}[0.9]{
\gbeg{5}{7}
\got{1}{\hspace{-2mm}A} \got{1}{\hspace{3mm}(\hspace{-2pt}A\# H\hspace{-2pt})^{^A}}
   \got{3}{\hspace{7mm}(\hspace{-2pt}A\# H\hspace{-2pt})^{^A}} \gnl
\gcl{1} \gbmp{j} \gvac{1} \gbmp{j} \gnl
\gcn{1}{1}{1}{3} \gwmu{3} \gnl
\gvac{1} \gbr \gu{1} \gnl
\gvac{1} \gcl{1} \glmpt \gnot{\hspace{-0,4cm}A\# H} \grmptb \gnl
\gvac{1} \gwmu{3} \gnl
\gvac{2} \gob{1}{A\# H}
\gend}
\end{array}$$
This means that Diagram \ref{make morf. on M^A} holds for the upper choice of $Q,M$ and $f$. In view of Remark
\ref{morma}, $\nabla_{A\# H}$ induces then a morphism $\crta\nabla:
(A\# H)^A \ot (A\# H)^A\to (A\# H)^A$ making the following diagram commutative:
$$\scalebox{0.9}[0.9]{
\bfig
\putmorphism(-150,425)(1,0)[(A\# H)^A \ot (A\# H)^A` (A\# H) \ot (A\# H)`j\ot j]{1100}1a
\putmorphism(0,0)(1,0)[(A\# H)^A`A\# H` j]{950}1a
\putmorphism(-150,420)(0,-1)[\phantom{(A\# H)^A \ot (A\# H)^A}`\phantom{(A\# H)^A}`\crta\nabla]{400}1l
\putmorphism(950,420)(0,-1)[\phantom{A\# H \ot A\# H}`\phantom{A\# H.}`\nabla_{A\# H}]{400}1r
\efig}
$$
We prove that $\eta_{A\# H}=\eta_A\#\eta_H:I\to A\# H$ induces a unit on $(A\# H)^A$.
We have:
$$\scalebox{0.9}[0.9]{
\gbeg{4}{5}
\got{1}{A} \got{1}{} \got{1}{} \gnl
\gcl{1} \gu{1} \gu{1}  \gu{1}  \gnl
\glmpt \gnot{\hspace{-0,4cm}A\# H} \grmpt \glmpt \gnot{\hspace{-0,4cm}A\# H} \grmpt \gnl
\gvac{1} \hspace{-0,2cm} \gwmu{3} \gnl
\gvac{2} \gob{1}{A\# H}
\gend}= A\ot\eta_H =\scalebox{0.9}[0.9]{
\gbeg{4}{8}
\got{1}{A}  \gnl
\gcl{1} \gu{1} \gu{1} \gnl
\gcl{1} \glmpt \gnot{\hspace{-0,4cm}A\# H} \grmpt \gnl
\gbr \gcl{1} \gnl
\gcl{1} \gbr \gu{1} \gnl
\glmpt \gnot{\hspace{-0,4cm}A\# H} \grmp \glmpt \gnot{\hspace{-0,4cm}A\# H} \grmp \gnl
\gvac{1} \hspace{-0,2cm} \gwmu{3} \gnl
\gvac{2} \gob{1}{A\# H}
\gend}$$
Then, there is a morphism $\crta\eta:I\to (A\# H)^A$ such that $j\crta{\eta}=\eta_{A\# H}.$ As done for 2.2.1(i), $((A\# H)^A, \crta\nabla, \crta\eta)$ is an algebra and $j:(A\# H)^A\to A\# H$ is an algebra morphism.
\qed\end{proof}

We will next prove that $(A\# H)^A$ is an $H$-subcomodule of $A\# H$.

\begin{lemma} \label{M^A-comod}
Let $H \in \C$ be a flat Hopf algebra, $A\in\C$ an algebra and $M\in \C^H$ with structure morphism $\rho_M:M\to M\ot H$. Assume that $M$ is a left $A \otimes \crta{A}$-module and that $\rho_M$ is $A \otimes \crta{A}$-linear, that is,
$$\gbeg{4}{4}
\got{1}{\hspace{-5pt} A \otimes \crta{A}} \got{1}{\hspace{5pt} M} \gnl
\gcl{1} \grcm \gnl
\glm \gcl{1} \gnl
\gvac{1} \gob{1}{M} \gob{1}{H}
\gend=
\gbeg{4}{4}
\got{1}{\hspace{-5pt} A \otimes \crta{A}} \got{1}{\hspace{5pt} M} \gnl
\glm \gnl
\gvac{1} \grcm \gnl
\gvac{1} \gob{1}{M} \gob{1}{H}
\gend$$
Then $M^A$ is a right $H$-comodule via $\crta\rho:=
t_{M, H}^{-1}\rho_M^A:M^A\to (M\ot H)^A\iso M^A\ot H$ and $j_M:M^A\to M$ is
an $H$-comodule morphism (here $t_{M, H}$ is that from Diagram \ref{t_{M,V}}).
\end{lemma}

\begin{proof}
Since $\rho_M:M\to M\ot H$ is left $A\ot\crta A$-linear, it induces $\rho_M^A:M^A\to (M\ot H)^A$ so that the square in the diagram
$$\scalebox{0.9}[0.85]{
\bfig
\putmorphism(0,425)(1,0)[M^A`M`j_M]{700}1a
\putmorphism(0,0)(1,0)[(M\ot H)^A`M\ot H`j_{{M\ot H}}]{700}1a
\putmorphism(0,420)(0,-1)[\phantom{M^A}`\phantom{(M\ot H)^A}`\rho_M^A]{400}1l
\putmorphism(705,420)(0,-1)[\phantom{M}`\phantom{M\ot H}`\rho_M]{400}1r
\putmorphism(0,0)(0,-1)[`M^A\ot H`t_{M, H}^{-1}]{300}1l
\putmorphism(350,-235)(2,1)[\phantom{M^A\ot H}`\phantom{M\ot H}`j_M\ot H]{200}1r
\efig}
$$
commutes. The triangle below commutes due to Proposition \ref{C-equiv-Az-t}. As in the proof of 2.2.1(ii), comodule version, one may show that $\crta\rho=t_{M, H}^{-1}\rho_M^A$ makes $M^A$ into an $H$-comodule and that
$j_M$ is right $H$-colinear.
\qed\end{proof}

In view of the preceding lemma we will have that $(A\# H)^A$ is an $H$-subcomodule
of $A\# H$ and $j:(A\# H)^A\to A\# H$ an $H$-comodule morphism if $\rho_{A\# H}=A \ot \Delta_H$
is $A\ot\crta A$-linear. This holds if $\Phi_{H, A}=\Phi^{-1}_{A, H}:$
$$\scalebox{0.85}[0.85]{
\gbeg{5}{6}
\got{1}{A\ot\crta A} \gvac{1} \got{1}{A\#H} \gnl
\gcn{1}{1}{1}{3} \gvac{1} \gcl{1} \gnl
\gvac{1} \gcl{1} \grcm \gnl
\gvac{1} \glm \gcl{1} \gnl
\gvac{2} \gcl{1} \gcn{1}{1}{1}{3} \gnl
\gvac{2} \gob{1}{A\#H} \gvac{1} \gob{1}{H}
\gend}\hspace{-3mm}=\hspace{-2mm}\scalebox{0.85}[0.85]{
\gbeg{5}{10}
\got{1}{A} \got{1}{\crta A} \got{1}{A} \got{2}{H} \gnl
\gcl{7} \gbr \gcmu \gnl
\gvac{1} \gcl{4} \gbr \gcl{7} \gnl
\gvac{2} \hspace{-0,34cm} \gcmu \gcn{1}{1}{0}{1} \gnl
\gvac{2} \gcl{1} \gbr \gnl
\gvac{2} \glm \gcl{4} \gnl
\gvac{1} \gcn{1}{1}{2}{3} \gvac{1} \gcl{1}   \gnl
\gvac{2} \gmu \gvac{1} \gnl
\gvac{1} \hspace{-0,36cm} \gwmu{3} \gnl
\gvac{1} \gob{3}{A} \gob{2}{H} \gob{1}{\hspace{-0,7cm}H}
\gend}\hspace{-2mm}\stackrel{nat.}{=}\hspace{-3mm}\scalebox{0.85}[0.85]{
\gbeg{6}{13}
\got{1}{A} \got{1}{\crta A} \got{1}{A} \got{1}{} \got{1}{H} \gnl
\gcl{10} \gbr \gvac{1} \gcl{2} \gnl
\gvac{1} \gcl{7} \gcn{1}{1}{1}{3}  \gnl
\gvac{3} \gbr \gvac{1} \gnl
\gvac{1} \gvac{1} \gwcm{3} \gcn{1}{1}{-1}{1} \gnl
\gvac{2} \hspace{-0,34cm} \gcmu \gvac{1} \hspace{-0,22cm} \gibr \gnl
\gvac{3} \hspace{-0,2cm} \gcl{2} \gcl{1} \gcn{1}{1}{2}{1} \gcn{1}{6}{2}{2} \gnl
\gvac{4} \gbr \gnl
\gvac{3} \glm \gcl{4} \gnl
\gvac{2} \gcn{1}{1}{2}{3} \gvac{1} \gcl{1} \gnl
\gvac{3} \gmu \gvac{2} \gnl
\gvac{2} \hspace{-0,36cm} \gwmu{3} \gnl
\gvac{3} \gob{1}{A} \gvac{1} \gob{2}{H} \gob{1}{H}
\gend}\hspace{-2mm} \stackrel{coass.}{=}\hspace{-3mm} \scalebox{0.85}[0.85]{
\gbeg{7}{13}
\got{1}{A} \got{1}{\crta A} \got{1}{A} \got{1}{} \got{1}{H} \gnl
\gcl{9} \gbr \gvac{1} \gcl{2} \gnl
\gvac{1} \gcl{7} \gcn{1}{1}{1}{3}  \gnl
\gvac{3} \gbr \gvac{1} \gnl
\gvac{2} \gwcm{3} \hspace{-0,42cm} \gcn{1}{1}{1}{4} \gnl
\gvac{3} \gcl{2} \gvac{1} \hspace{-0,34cm} \gcmu \gcl{1} \gnl
\gvac{5} \gcl{1} \gibr \gnl
\gvac{4} \gcn{1}{1}{0}{1} \gbr \gcl{5} \gnl
\gvac{4} \glm \gcl{4} \gnl
\gvac{3} \gcn{1}{1}{0}{1} \gvac{1} \gcl{1}   \gnl
\gvac{2} \gcn{1}{1}{0}{1} \gwmu{3} \gnl
\gvac{2} \gwmu{3}  \gnl
\gvac{2} \gob{3}{A} \gob{3}{H} \gob{1}{\hspace{-0,7cm}H}
\gend}\hspace{-2mm}\stackrel{cond.\Phi}{\stackrel{nat.}{=}}\hspace{-2mm} \scalebox{0.85}[0.85]{
\gbeg{5}{10}
\got{1}{A} \got{1}{\crta A} \got{1}{A} \got{1}{H} \gnl
\gcl{7} \gbr \gcl{1} \gnl
\gvac{1} \gcl{4} \gbr  \gnl
\gvac{2} \hspace{-0,21cm} \gcmu \gcn{1}{1}{0}{1} \gnl
\gvac{2} \hspace{-0,14cm} \gcl{1} \gbr \gnl
\gvac{2} \glm \gcl{3} \gnl
\gvac{1} \gcn{1}{1}{2}{3} \gvac{1} \gcl{1}   \gnl
\gvac{2} \gmu \gvac{1} \gnl
\gvac{1} \hspace{-0,34cm} \gwmu{3} \gcmu \gnl
\gvac{1} \gob{3}{A} \gob{1}{H} \gob{1}{H}
\gend}\hspace{-2mm}=\hspace{-3mm}\scalebox{0.85}[0.85]{
\gbeg{4}{6}
\got{1}{A\ot\crta A} \gvac{1} \got{1}{A\#H} \gnl
\gcn{1}{1}{1}{3} \gvac{1} \gcl{1} \gnl
\gvac{1} \glm \gnl
\gvac{2} \grcm \gnl
\gvac{2} \gcl{1} \gcn{1}{1}{1}{3} \gnl
\gvac{2} \gob{1}{A\#H} \gvac{1} \gob{1}{H}
\gend}
$$
As a subalgebra and a subcomodule of the $H$-comodule algebra $A\# H$, the object
$(A\# H)^A$ is itself an $H$-comodule algebra (2.2.1). We record this fact in the following
result.

\begin{corollary} \label{(Asm.H)A+j com.alg.}
Let $A$ be an $H$-module algebra, where $H$ is flat, and suppose that the braiding satisfies
$\Phi_{H, A}=\Phi^{-1}_{A,H}$. Then $(A\# H)^A$ is a right $H$-comodule algebra and the morphism
$j:(A\# H)^A\to A\# H$ is an $H$-comodule algebra morphism.
\end{corollary}

In our context the above condition on $\Phi$ is satisfied because the braiding is $H$-linear and by Proposition \ref{braid-lin}. This corollary proves 1). We now proceed to prove 2), that is, that $(A\# H)^A$ is faithfully flat. \par \bigskip

Observe that $A\ot (A\# H)^A\iso A\# H =A\ot H$ in $\C$. The first isomorphism holds because $A$ is Azumaya
and hence we have an equivalence of categories  $A\ot -:\C\to {}_{A\ot\crta A}\C:(-)^A$, Proposition \ref{Az-right-adj}. Recall that $A$ is Azumaya in $\C$ because it is so in ${}_H\C.$ As an Azumaya algebra, $A$ is faithfully projective and thus faithfully flat (2.5.1). From Corollary \ref{H:fl<->f.f.}(i) we know that $H$ is faithfully flat. Then we have by 2.1.1(ii) that $A\ot H$,
and hence also $A\ot (A\# H)^A$, is faithfully flat. \par \smallskip

Let $f,g:M\to N$ and $e:E\to M$ be morphisms in $\C$. Assume that
$$\scalebox{0.9}[0.9]{\bfig
\putmorphism(50,0)(1,0)[A\ot E`A\ot M`A\ot e]{700}1a
\putmorphism(880,25)(1,0)[`\phantom{A\ot H}`A\ot f]{720}1a
\putmorphism(880,0)(1,0)[`A\ot N`]{720}0a
\putmorphism(880,-25)(1,0)[`\phantom{A\ot H}`A\ot g]{720}1b
\efig}
$$
is an equalizer in $\C$. It is also an equalizer in ${}_{A\ot\crta A}\C$
due to 2.2.2. Being an equivalence, $A\ot -:\C\to {}_{A\ot\crta A}\C$
reflects equalizers. Thus $(E, e)$ is an equalizer in $\C$ and so $A\ot -$
reflects equalizers in $\C$. Now we may apply 2.1.1(iii) to conclude that $(A\# H)^A$ is
faithfully flat. \par \bigskip

We finally prove 3), that is, that $\can_{(A\#H)^A}$ is an isomorphism. This will not be an easy task. An Azumaya algebra $A$ in ${}_H\C$ is in particular an Azumaya algebra in $\C$. The functor $A\ot -: \C\to {}_{A\ot\crta A}\C$ is then an equivalence of categories. Our strategy will be to prove that $A\ot\can_{(A\#H)^A}$ is an isomorphism, then so will be $\can_{(A\#H)^A}$. \par \smallskip

From Remark \ref{Az-counit} the counit $\beta: A\ot (A\#H)^A\to A\#H$ of the adjunction
$(A\ot -, (-)^A)$ is given by the morphism $\beta=\nu^A_{A\# H}(A\ot j)$. Since $A$ is an
Azumaya algebra, we know that $\beta$ is an isomorphism in ${}_A\C_A$. Let $\delta:(A\#H)\ot_A A\to A\#H$ denote the corresponding isomorphism in ${}_A\C_A$ from 2.3.3 and $\omega: [(A\#H)\ot_A A]\ot (A\#H)^A\to (A\#H)\ot_A [A\ot (A\#H)^A]$ the isomorphism in ${}_A\C_A$ from 2.3.4. Let $\sigma: A\ot (A\#H)^A\ot (A\#H)^A \to (A\#H)\ot_A (A\#H)$
in ${}_A\C_A$ be defined as the following composition of isomorphisms:
$$\begin{array}{rl}
A\ot (A\#H)^A \ot (A\#H)^A \hskip-1em& \stackrel{\beta\ot id}{\iso} (A\#H)\ot (A\#H)^A \\
&\hspace{-0,14cm}\stackrel{\delta^{-1}\ot id}{\iso} [(A\#H)\ot_A A]\ot (A\#H)^A \\
&\hspace{0,2cm}\stackrel{\omega}{\iso}(A\#H)\ot_A [A\ot (A\#H)^A] \\
&\hspace{-0,04cm}\stackrel{id\ot_A\beta}{\iso} (A\#H)\ot_A (A\#H).
\end{array}$$
We further set $\tau:=\beta\ot H: A\ot (A\#H)^A \ot H \to (A\#H) \ot H.$ Now we define a morphism $\xi: (A\#H)\ot_A (A\#H) \to (A\#H)\ot H$ as $\xi=\tau (A \otimes can_{(A\#H)^A})\sigma^{-1}$. If we show that $\xi$
is an isomorphism, then so will be $A\ot \can_{(A\#H)^A}$ and we will be done.

We define the morphism
$$\Lambda:=(\crta\nabla_{A\#H}\ot H)\lambda((A\#H)\ot_A\rho_{A\#H}):
(A\#H)\ot_A (A\#H)\to (A\#H)\ot H,$$
where each morphism in this composition is induced like
the following diagram indicates ($\Lambda$ is the composition on the right-hand side edge):
$$\scalebox{0.85}[0.85]{
\bfig
\putmorphism(0, 500)(1, 0)[\phantom{(A\#H)\ot (A\#H)}`\phantom{(A\#H)\ot_A (A\#H)} `\Pi]{1400}1a
\putmorphism(1400, 500)(0, -1)[(A\#H)\ot_A (A\#H)`(A\#H)\ot_A [(A\#H)\ot H]`(A\#H)\ot_A\rho_{A\#H}]{500}1r
\putmorphism(0, 500)(0, -1)[(A\#H)\ot (A\#H)`\phantom{(A\#H)\ot_A (A\#H)}`(A\#H)\ot\rho_{A\#H}]{480}1l
\putmorphism(0,0)(1, 0)[(A\#H)\ot [(A\#H)\ot H]`\phantom{(A\#H)\ot_A [(A\#H)\ot H]}`\Pi']{1400}1b
\putmorphism(0,0)(0, -1)[\phantom{(A\#H)\ot [(A\#H)\ot H]}`[(A\#H)\ot (A\#H)]\ot H`\iso]{500}1l
\putmorphism(0,-500)(1, 0)[\phantom{[(A\#H)\ot (A\#H)]\ot H}` [(A\#H)\ot_A (A\#H)]\ot H`\Pi\ot H]{1400}1a
\putmorphism(1400,0)(0, -1)[\phantom{(A\#H)\ot_A [(A\#H)\ot H]}`\phantom{[(A\#H)\ot_A (A\#H)]\ot H}`\lambda]{500}1r
\putmorphism(1400,-500)(0, -1)[\phantom{[(A\#H)\ot_A (A\#H)]\ot H}`(A\#H)\ot H`\crta\nabla_{A\#H}\ot H]{500}1r
\putmorphism(-200,-500)(3, -1)[\phantom{[(A\#H)\ot (A\#H)]\ot H}`\phantom{(A\#H)\ot H}`]{1500}1l
\putmorphism(-170,-560)(3, -1)[``\nabla_{A\#H}\ot H]{1500}0l
\efig}
$$
Here $\Pi:=\Pi_{A\#H,A\#H}$ and $\Pi':=\Pi_{A\#H,(A\#H)\ot H}$. Observe that
$\lambda$ is an isomorphism because the third row of the diagram is a part of a
coequalizer, since $\C$ is closed. From the diagram it is clear that
$\Lambda\Pi=can_{A\#H}.$ We are going to prove first that $\xi=\Lambda$ and then we will find the inverse for $\Lambda$. So $\xi$ will be an isomorphism, as desired. For that purpose we consider the following diagram. The composition of morphisms on the right hand-side edge represents $\xi=\tau (A\ot can_{(A\#H)^A})\sigma^{-1}$
whereas on the left-hand side edge are mostly the morphisms that induce the latter ones:
$$\scalebox{0.8}[0.8]{
\bfig
\putmorphism(0, 500)(1, 0)[\phantom{(A\#H)\ot (A\#H)}`
  \phantom{(A\#H)\ot_A (A\#H)} `\Pi]{1400}1a
\putmorphism(1400, 500)(0, -1)[(A\#H)\ot_A (A\#H)` (A\#H)\ot_A [A\ot (A\#H)^A]`
  id\ot_A\beta^{-1}]{500}1r
\putmorphism(-40, 500)(0, -1)[(A\#H)\ot (A\#H)`
  \phantom{(A\#H)\ot_A (A\#H)}` id\ot\beta^{-1}]{480}1l
\putmorphism(-110,0)(1, 0)[(A\#H)\ot [A\ot (A\#H)^A]` \phantom{(A\#H)\ot_A [A\ot (A\#H)^A]}`
  \Pi'']{1470}1a
\putmorphism(-40,0)(0, -1)[\phantom{(A\#H)\ot [A\ot (A\#H)^A]}`\phantom{[(A\#H)\ot A]\ot (A\#H)^A}
  `\iso]{500}1l
\putmorphism(-150,-500)(1, 0)[[(A\#H)\ot A]\ot (A\#H)^A` [(A\#H)\ot_A A]\ot (A\#H)^A`
  \Pi_1\ot id]{1550}1a
\putmorphism(1400,0)(0, -1)[\phantom{(A\#H)\ot_A [(A\#H)\ot H]}`
  \phantom{[(A\#H)\ot_A (A\#H)]\ot H}`\omega^{-1}]{500}1r
\putmorphism(1400,-500)(0, -1)[\phantom{[(A\#H)\ot_A A]\ot (A\#H)^A}` (A\#H)\ot (A\#H)^A`
  \delta\ot id]{500}1r
\putmorphism(-200,-500)(3, -1)[\phantom{(A\#H)\ot (A\#H)^A}` \phantom{(A\#H)\ot H}`]{1500}1l
\putmorphism(-260,-520)(3, -1)[`\phantom{[A\ot (A\#H)^A]\ot A\ot (A\#H)^A} ` \mu^A_{A\#H}\ot id]
  {1500}0l
\putmorphism(-230,-520)(0, 1)[`A\ot (A\#H)^A\ot A\ot (A\#H)^A ` \beta^{-1}\ot A\ot id]{500}1l
\putmorphism(-200,-1030)(3,-1)[` `]{1500}1r
\putmorphism(-300,-940)(3, -1)[`  ` \mu^A_{A\ot (A\#H)^A}\ot id]{1500}0r
\putmorphism(1400,-1000)(0, -1)[\phantom{(A\#H)\ot (A\#H)^A}`  ` \beta^{-1}\ot id]{500}1r
\putmorphism(-80,-1030)(3,-1)[`[A\ot(A\#H)^A]\ot(A\#H)^A`]{1500}0r
\put(650,250){\fbox{1}}
\put(650,-250){\fbox{2}}
\put(1000,-750){\fbox{3}}
\put(590,-1050){\fbox{4}}
\putmorphism(-230,-1020)(0, 1)[`A\ot A\ot(A\#H)^A\ot(A\#H)^A`A\ot\Phi^{-1}\ot id]{500}1l
\putmorphism(-170,-1500)(3, -1)[`` \nabla_A\ot id\ot id]{1500}0r
\putmorphism(1400,-1500)(0, -1)[  ``\iso]{500}1r
\putmorphism(-200,-1530)(3,-1)[``]{1500}1r
\putmorphism(-400,-2030)(1, 0)[A\ot (A\#H)\ot (A\#H)` A\ot [(A\#H)^A\ot (A\#H)^A]`
  A\ot j\ot j]{1800}{-1}a
\putmorphism(-230,-2030)(0, -1)[\phantom{A\ot (A\#H)\ot (A\#H)}` \phantom{A\ot (A\#H)\ot H}`
  A\ot \can_{A\#H}]{500}1r
\putmorphism(-130,-2530)(1, 0)[A\ot (A\#H)\ot H` \phantom{A\ot [(A\#H)^A\ot H]}`
  A\ot j\ot H]{1530}{-1}a
\putmorphism(1400,-2030)(0, -1)[\phantom{A\ot [(A\#H)^A\ot (A\#H)^A]}` A\ot [(A\#H)^A\ot H]`
  A\ot can_{(A\#H)^A}]{500}1r
\putmorphism(1400,-2530)(0, -1)[\phantom{A\ot [(A\#H)^A\ot H]}` [A\ot (A\#H)^A]\ot H]`
  \iso]{400}1r
\putmorphism(1400,-2930)(0, -1)[\phantom{[A\ot (A\#H)^A]\ot H]}` `\beta\ot H]{400}1r
\putmorphism(1400,-2950)(0, -1)[\phantom{[A\ot (A\#H)^A]\ot H]}` (A\#H)\ot H`]{400}0r
\putmorphism(-260,-2560)(2, -1)[\phantom{A\ot (A\#H)\ot H}` \phantom{(A\#H)\ot H}` ]{1550}1l
\putmorphism(-340,-2560)(2, -1)[\phantom{A\ot (A\#H)\ot H}` \phantom{(A\#H)\ot H}`
  \nabla_A\ot H\ot  H]{1550}0l
\putmorphism(-550,-2030)(0, -1)[\phantom{A\ot (A\#H)\ot (A\#H)}` (A\#H)\ot (A\#H)`]{1330}1l
\putmorphism(-550,-2055)(0, -1)[` `\nabla_A\ot H\ot (A\#H)]{1330}0l
\putmorphism(-570,-3360)(1, 0)[\phantom{(A\#H)\ot (A\#H)}` \phantom{(A\#H)\ot H]}`
  \can_{A\#H}]{1950}1b
\put(490,-1550){\fbox{5}}
\put(490,-2290){\fbox{6}}
\put(590,-2790){\fbox{7}}
\put(-440,-3050){\fbox{8}}
\efig}
$$
Here $\Pi'':=\Pi_{A\#H, A\ot (A\#H)^A}$ and $\Pi_1:=\Pi_{A\#H, A}$. Note that the
diagrams $\langle1\rangle$, $\langle2\rangle$, $\langle3\rangle$, $\langle6\rangle$ and
$\langle7\rangle$ commute as they define morphisms on their right hand side edges. Diagram
$\langle4\rangle$ commutes since $\beta$ as counit (and hence also $\beta^{-1}$) is a
morphism in ${}_A\C_A$, in particular it is right $A$-linear. Diagram $\langle5\rangle$ commutes
by the way the right $A$-module structure of $A\ot (A\#H)^A$ is defined. The object $A\ot (A\#H)^A$ has a structure of an $A$-bimodule inherited from the one of its left tensor factor $A$. Finally, $\langle8\rangle$ commutes since $\can$ is left $A$-linear. So the whole big diagram commutes.

Applying associativity in $A$, the equalizer property of $(A\#H)^A$ and naturality respectively,
we have:
$$\begin{array}{l}
(\nabla_A\ot H)(\nabla_A\ot j)(A\ot\Phi^{-1}_{A,(A\#H)^A})= \\
\hspace{1.5cm}=(\nabla_A\ot H)(A\ot\nabla_A\ot H)(A\ot A\ot j)(A\ot\Phi^{-1}_{A,(A\#H)^A}) \\
\hspace{1.5cm}= (\nabla_A\ot H)(A\ot\mu^A_{A\#H})(A\ot\Phi_{A, A\#H})(A\ot A\ot j)(A\ot\Phi^{-1}_{A,(A\#H)^A}) \\
\hspace{1.5cm}=(\nabla_A\ot H)(A\ot\mu^A_{A\#H})(A\ot j\ot A).
\end{array}$$
Tensoring this equality with $j$ on the right, we get that the composition of morphisms
from the edges of diagrams $\langle5\rangle$, $\langle6\rangle$ and $\langle8\rangle$
$$(\nabla_A\ot H\ot (A\#H))(A\ot j\ot j)(\nabla_A\ot\id\ot\id)(A\ot\Phi^{-1}_{A,(A\#H)^A}\ot\id)$$
becomes
$$(\nabla_A\ot H\ot \id)(A\ot\mu^A_{A\#H}\ot \id)(A\ot j\ot A\ot j).$$
We substitute this and, from the outer edges of the big diagram, we get the following situation where the outer diagram commutes:
$$\scalebox{0.8}[0.8]{\bfig \hspace{-0,8cm}
\putmorphism(0, 2000)(0, -1)[(A\#H)\ot (A\#H)` A\ot (A\#H)^A\ot A\ot (A\#H)^A`
  \beta^{-1}\ot\beta^{-1}]{500}1l
 \putmorphism(0, 2000)(1, 0)[\phantom{(A\#H)\ot (A\#H)}` (A\#H)\ot_A (A\#H)`\Pi]{2820}1a
\putmorphism(0, 1500)(0, -1)[\phantom{A\ot (A\#H)^A\ot A\ot (A\#H)^A}` ` A\ot j\ot A\ot j]{480}1l
\putmorphism(-200, 1500)(0, -1)[` A\ot (A\#H)\ot A\ot (A\#H)` ]{480}0l
 \putmorphism(0, 1500)(1, 0)[\phantom{A\ot (A\#H)^A\ot A\ot (A\#H)^A}` (A\#H)\ot (A\#H)`\beta\ot\beta]{2000}1a
 \putmorphism(-200, 1000)(1, 0)[\phantom{A\ot (A\#H)\ot A\ot (A\#H)}` A\ot (A\#H)\ot  (A\#H)`A\ot\id\ot{}\nu_{A \# H}^A]{1800}1a
 \putmorphism(1600, 1000)(0, -1)[\phantom{(A\#H)\ot (A\#H)}`[A\ot (A\#H)]\ot_A (A\#H)`\Pi_2]{600}1r
\putmorphism(0,1000)(0, -1)[`A\ot (A\#H)\ot (A\#H)`]{600}1l
 \putmorphism(0,400)(1, 0)[\phantom{A\ot (A\#H)\ot (A\#H)}`\phantom{[A\ot (A\#H)]\ot_A (A\#H)}`\Pi_2]{1600}1a
\putmorphism(0,1000)(0, -1)[`  `\mu^A_{A\ot (A\#H)}\ot\id]{600}0l
\putmorphism(0,400)(0, -1)[\phantom{A\ot (A\#H)\ot (A\#H)}`
 (A\#H)\ot (A\#H)`\nu_{A \# H}^A\ot\id]{500}1l
\putmorphism(2900, 2000)(0, -1)[\phantom{(A\#H)\ot_A (A\#H)}` (A\#H)\ot H`\xi]{2100}1r
\putmorphism(0,-100)(1, 0)[\phantom{(A\#H)\ot (A\#H)}`(A\#H)\ot_A (A\#H)`\Pi]{2000}1a
\putmorphism(2000,-100)(1, 0)[\phantom{(A\#H)\ot_A (A\#H)}`\phantom{(A\#H)\ot H}`\Lambda]{900}1a
\putmorphism(1600, 400)(1, -1)[``\nu_{A \# H}^A\ot_A\id]{500}1l
\putmorphism(1660, 1140)(1, 1)[``\nu_{A \# H}^A\ot\id]{230}1l
 \putmorphism(2210, 1500)(0, -1)[``\Pi]{1600}1r
\put(2500,1250){\fbox{9}}
\put(650,1250){\fbox{10}}
\put(650,680){\fbox{11}}
\put(1840,680){\fbox{12}}
\put(650,140){\fbox{12}}
\efig}
$$
Here diagram $\langle10\rangle$ commutes by the definition of $\beta$. Diagram
$\langle11\rangle$ commutes by the coequalizer property of $[A\ot (A\#H)]\ot_A (A\#H)$.
Diagram $\langle12\rangle$ commutes by the definition of $\nu_{A \# H}^A\ot_A \id_{A\#H}$ ($\Pi_2$ is the
coequalizer morphism). Since also the outer diagram commutes, diagram $\langle9\rangle$ commutes as well,
yielding $\xi\Pi=\Lambda\Pi$. But $\Pi$ is an epimorphism, so $\xi=\Lambda$. \par \medskip

The next step is to find the inverse of $\Lambda$. For that purpose we define the morphism
$\teta: (A\#H)\ot H\to (A\#H)\ot (A\#H)$ by
$$\teta:=\scalebox{0.9}[0.9]{
\gbeg{4}{5}
\got{1}{\hspace{0,1cm}A} \got{1}{\#H} \got{3}{H} \gnl
\gcl{2} \gcl{2} \gwcm{3} \gvac{1} \gnl
\gvac{2} \gnot{S}\gmp \gnl \gvac{1} \gcl{1} \gnl
\gcl{1} \gmu \gu{1} \gcl{1} \gnl
\gob{1}{A} \gob{2}{\hspace{-0,2cm}\#H} \gob{2}{A\#H}
\gend}
$$
and let $\crta\teta:=\Pi\teta: (A\#H)\ot H\to (A\#H)\ot_A (A\#H)$. Let us prove that $\crta\teta$
is the inverse of $\Lambda$. We have:
$$\begin{array}{ll}
\crta\teta\Lambda\Pi=\crta\teta\can_{A\#H} & =
\scalebox{0.8}[0.8]{\gbeg{7}{10}
\got{1}{A} \got{2}{H} \got{1}{A} \got{2}{H} \gnl
\gcl{3} \gcmu \gcl{1} \gcmu \gvac{1} \gnl
\gvac{1} \gcl{1} \gbr \gcl{1} \gcl{3} \gnl
\gvac{1} \glm \gmu \gnl
\gwmu{3} \gvac{1} \hspace{-0,36cm} \gcl{3} \gnl
\gvac{2} \gcn{1}{2}{0}{4} \gvac{2} \gcmu \gnl
\gvac{5} \gnot{S}\gmp \gnl \gcn{1}{1}{1}{3} \gnl
\gvac{4} \hspace{-0,36cm} \gcl{1} \hspace{-0,2cm} \gmu \gu{1} \gcl{1} \gnl
\gvac{4} \gnot{\hspace{1,7cm}\Pi} \glmp \gcmp \gcmpb \gcmp \grmp \gnl
\gvac{6} \gob{1}{(A\#H)\ot_A (A\#H)}
\gend} \stackrel{coass.}{\stackrel{ass.}{=}}
\scalebox{0.8}[0.8]{\gbeg{7}{8}
\got{1}{A} \got{2}{H} \got{1}{A} \got{2}{H} \gnl
\gcl{3} \gcmu \gcl{1} \gcmu \gvac{1} \gnl
\gvac{1} \gcl{1} \gbr \hspace{-0,22cm} \gcmu \gcn{1}{1}{0}{1} \gnl
\gvac{2} \hspace{-0,34cm} \glm \gcl{2} \hspace{-0,22cm} \gcl{1} \gnot{S}\gmp \gnl \gcl{3} \gnl
\gvac{2} \hspace{-0,34cm} \gwmu{3} \gvac{1} \hspace{-0,36cm} \gmu \gnl
\gvac{4} \hspace{-0,2cm} \gcn{1}{1}{1}{3} \gvac{1} \gmu \hspace{-0,2cm} \gu{1} \gnl
\gvac{6} \hspace{-0,42cm} \gnot{\hspace{1,7cm}\Pi} \glmp \gcmp \gcmpb \gcmp \grmp \gnl
\gvac{8} \gob{1}{(A\#H)\ot_A (A\#H)}
\gend} \stackrel{antip.}{\stackrel{counit}{\stackrel{unit}{=}}}
\scalebox{0.8}[0.8]{\gbeg{6}{8}
\got{1}{A} \got{2}{H} \got{1}{A} \got{1}{H} \gnl
\gcl{1} \gcmu \gcl{1} \gcl{5} \gvac{1} \gnl
\gcl{1} \gcl{1} \gbr \hspace{-0,2cm} \gnl
\gvac{1} \hspace{-0,5cm} \gcl{1} \glm \gcl{1} \gnl
\gvac{1} \gwmu{3} \hspace{-0,41cm} \gcn{1}{2}{3}{1}  \gnl
\gvac{3} \gcl{1} \gvac{1} \gu{1} \gnl
\gvac{3} \hspace{-0,32cm} \gnot{\hspace{1,7cm}\Pi} \glmp \gcmp \gcmpb \gcmp \grmp \gnl
\gvac{5} \gob{1}{(A\#H)\ot_A (A\#H)}
\gend} \vspace{3mm} \\
 & = \scalebox{0.8}[0.8]{
\gbeg{5}{5}
\got{1}{A\#H} \got{3}{A} \got{1}{\hspace{-0,7cm}H} \gnl
\gcl{1} \gcn{2}{1}{3}{1} \gcl{2} \gnl
\grm \gu{1} \gnl
\hspace{-0,32cm} \gnot{\hspace{1,7cm}\Pi} \glmp \gcmp \gcmpb \gcmp \grmp \gnl
\gvac{2} \gob{1}{(A\#H)\ot_A (A\#H)}
\gend} \hspace{2mm} \stackrel{\Pi}{=}\hspace{2mm}
\scalebox{0.8}[0.8]{\gbeg{4}{5}
\got{1}{\hspace{-0,06cm}A\#H} \got{1}{\hspace{0,04cm}A} \got{3}{H} \gnl
\gcl{1} \gcl{1} \gu{1} \gcl{2} \gnl
\gcl{1} \gmu \gnl
\hspace{-0,2cm} \gnot{\hspace{1,7cm}\Pi} \glmp \gcmp \gcmptb \gcmp \grmp \gnl
\gvac{2} \gob{1}{(A\#H)\ot_A (A\#H)}
\gend}=\Pi.
\end{array}$$
From here, we get $\crta\teta\Lambda=\id_{(A\#H)\ot_A (A\#H)}$. Finally:
$$\begin{array}{lll}
\Lambda\crta\teta=\Lambda\Pi\teta=\can_{A\#H}\teta & =
\scalebox{0.8}[0.8]{\gbeg{6}{9}
\got{1}{A} \got{1}{\hspace{0,1cm}H} \got{3}{H} \gnl
\gcl{6} \gcl{2} \gwcm{3} \gnl
\gvac{2} \gnot{S}\gmp \gnl \gvac{1} \gcl{1} \gnl
\gcl{1} \gmu \gvac{1} \gcn{1}{2}{1}{2}  \gnl
\gvac{1} \gcmu \gu{1} \gvac{1} \gnl
\gvac{1} \gcl{1} \gbr \gcmu \gnl
\gvac{1} \glm \gmu \gcl{2} \gnl
\gwmu{3} \gvac{1} \hspace{-0,34cm} \gcl{1} \gnl
\gvac{1} \gob{2}{A} \gob{3}{H} \gob{1}{\hspace{-0,5cm}H}
\gend} \stackrel{mod.\x alg.}{=}
\scalebox{0.8}[0.8]{\gbeg{5}{9}
\got{1}{A} \got{1}{\hspace{0,1cm}H} \got{3}{H} \gnl
\gcl{6} \gcl{2} \gwcm{3} \gnl
\gvac{2} \gnot{S}\gmp \gnl \gvac{1} \gcl{3} \gnl
\gvac{1} \gmu \gnl
\gvac{1} \gcmu \gvac{1} \gnl
\gvac{1} \gcu{1} \gcn{1}{1}{1}{2} \gvac{1} \hspace{-0,35cm} \gcmu \gnl
\gvac{2} \hspace{-0,22cm} \gu{1} \gvac{1} \hspace{-0,34cm} \gmu \gcl{2} \gnl
\gvac{2} \hspace{-0,35cm} \gmu \gvac{1}  \gcl{1} \gnl
\gvac{2} \gob{2}{A} \gob{3}{H} \gob{1}{\hspace{-0,5cm}H}
\gend} & \vspace{3mm} \\
 & \stackrel{counit}{\stackrel{unit}{=}}
\scalebox{0.8}[0.8]{\gbeg{5}{6}
\got{1}{A} \got{1}{\hspace{0,1cm}H} \got{3}{H} \gnl
\gcl{4} \gcl{2} \gwcm{3} \gnl
\gvac{2} \gnot{S}\gmp \gnl \gvac{1} \hspace{-0,22cm} \gcmu \gnl
\gvac{2} \hspace{-0,2cm} \gmu \hspace{-0,2cm} \gvac{1} \hspace{-0,3cm} \gcl{1} \gcl{2} \gnl
\gvac{3} \gwmu{3} \gnl
\gvac{2} \gob{1}{\hspace{-0,4cm}A} \gob{3}{H} \gob{1}{H}
\gend} \stackrel{coass.}{\stackrel{ass.}{=}}
\scalebox{0.8}[0.8]{\gbeg{5}{7}
\got{1}{\hspace{0,6cm}A\#H} \got{6}{H} \gnl
\gcl{5} \gcl{4} \gvac{1} \gcmu \gnl
\gvac{3} \hspace{-0,2cm} \gcmu \hspace{-0,22cm} \gcl{4} \gnl
\gvac{4} \hspace{-0,34cm} \gnot{S}\gmp \gnl \gcl{1} \gnl
\gvac{4} \gmu \gnl
\gvac{3} \hspace{-0,22cm} \gwmu{3} \gnl
\gvac{2} \gob{1}{A} \gob{3}{H} \gob{1}{H}
\gend} & \stackrel{antip.}{\stackrel{counit}{\stackrel{unit}{=}}}
\id_{(A\#H)\ot H}.
\end{array}$$
This finishes the proof of 3). Thus we have established:

\begin{proposition} \label{(Asm.H)A is Gal.}
Let $\C$ be a closed braided monoidal category with equalizers and coequalizers.
Let $H\in \C$ be a flat Hopf algebra and suppose that the braiding is
$H$-linear. If $A \in \C$ is an $H$-Azumaya algebra, then $(A\# H)^A$ is an $H$-Galois
object.
\end{proposition}

We are next going to show that the assignment
$$\Upsilon: \BM(\C; H)\to\Gal(\C; H), \quad [A]\mapsto [(A\# H)^A]$$
just established is a group morphism. The proof will be long and technical.

\begin{proposition} \label{Pi multiplicative}
Assumptions are like in the previous proposition. Suppose in addition that $\Phi_{T,X}=\Phi^{-1}_{X,T}$ for any $H$-Galois object $T$ and $X \in \C$. If $A$ and $B$ are $H$-Azumaya
algebras in $\C$, then there is an isomorphism of $H$-Galois objects
$$[(A\ot B)\#H]^{A\ot B} \iso (A\#H)^A\Box_H (B\#H)^B.$$
\end{proposition}

\begin{proof}
As the braiding is $H$-linear, $H$ is cocommutative by Proposition \ref{braid-lin}. That the above two objects are $H$-Galois we know from Proposition \ref{(Asm.H)A is Gal.} and because the cotensor product of two $H$-Galois objects
is such too. In order to prove that they are isomorphic as
$H$-Galois objects, it suffices to find an $H$-comodule algebra morphism between them, in virtue of Proposition \ref{comod-alg-iso}. Now we explain the strategy of the proof. Observe the following diagram
\begin{eqnarray}\label{diagp69}
\scalebox{0.8}[0.8]{\hspace{-3cm}
\bfig
\putmorphism(-590,20)(0, -1)[`\phantom{[(A\#H)\ot (B\#H)]^{A\ot B}}`]{1010}1r
\putmorphism(-590,-270)(0, -1)[`\phantom{[(A\#H)\ot (B\#H)]^{A\ot B}}`\alpha_2]{1000}0r
\putmorphism(-170,0)(0, -1)[[(A\ot B)\#H]^{A\ot B}`\phantom{[(A\#H)\ot (B\#H)]^{A\ot B}}`
   \alpha_1]{500}1l
\putmorphism(-200,0)(1, 0)[\phantom{[(A\ot B)\#H)]^{A\ot B}}` (A\ot B)\#H `j]{1330}1a
\putmorphism(-50,-500)(1, 0)[[(A\#H)\ot (B\#H)]^{A\ot B}` (A\#H)\ot (B\#H)` \breve j]{1250}1a
\putmorphism(1180,0)(0, -1)[\phantom{(A\#H)\ot_A [(A\#H)\ot H]}`
  \phantom{[(A\#H)\ot_A (A\#H)]\ot H}`\alpha]{500}1r
\putmorphism(1180,-500)(0, -1)[\phantom{[(A\#H)\ot_A A]\ot (A\#H)^A}`
   (A\#H)^A\ot (B\#H)^B` j_{A'}\ot j_{B'}]{500}{-1}r
\putmorphism(-300,-500)(3, -1)[\phantom{(A\#H)^A\ot (B\#H)^B}` \phantom{(A\#H)\ot H}`]{1500}1l
\putmorphism(-340,-520)(3, -1)[`\phantom{[A\ot (A\#H)^A]\ot A\ot (A\#H)^A} ` \zeta^{-1}_{A', B'}]
  {1500}0l
\putmorphism(-120,-980)(1, 0)[(A\#H)^A\Box_H (B\#H)^B` \phantom{(A\#H)^A\ot (B\#H)^B}` e]{1250}1a
\putmorphism(1400,-1000)(1, 0)[\phantom{(A\#H)^A\ot (B\#H)^B}` (A\#H)^A\ot H\ot (B\#H)^B` ]{1250}0a
\putmorphism(1200,-980)(1, 0)[\phantom{(A\#H)^A\ot (B\#H)^B}` \phantom{(A\#H)^A\ot H\ot (B\#H)^B}`
   \rho'\ot\id]{1450}1a
\putmorphism(1200,-1020)(1, 0)[\phantom{(A\#H)^A\ot (B\#H)^B}` \phantom{(A\#H)^A\ot H\ot (B\#H)^B}`
   \id\ot\lambda' ]{1450}1b
\putmorphism(1180,-1000)(0, -1)[\phantom{(A\#H)^A\ot (B\#H)^B}`
  (A\#H)\ot (B\#H)`j_{A'}\ot j_{B'}]{500}1l
\putmorphism(1400,-1500)(1, 0)[\phantom{(A\#H)\ot (B\#H)}` (A\#H)\ot H\ot(B\#H),` ]{1250}0a
\putmorphism(1180,-1480)(1, 0)[\phantom{(A\#H)\ot (B\#H)}` \phantom{(A\#H)\ot (B\#H)}`
   \rho\ot\id_{B'}]{1360}1a
\putmorphism(1180,-1520)(1, 0)[\phantom{(A\#H)\ot (B\#H)}` \phantom{(A\#H)\ot (B\#H)}`
   \id_{A'}\ot\lambda]{1360}1b
\putmorphism(2680,-1000)(0, -1)[\phantom{(A\#H)^A\ot (B\#H)^B}`
  \phantom{(A\#H)\ot (B\#H)}`j_{A'}\ot H\ot j_{B'}]{500}1l
\put(450,-250){\fbox{1}}
\put(900,-750){\fbox{2}}
\put(1580,-1280){\fbox{3}}
\put(-350,-750){\fbox{4}}
\efig}
\end{eqnarray}
where the notation we enlighten here:
$$\begin{array}{lll}
A':=A\# H, & \qquad B':=B\# H, & \qquad e:=e_{(A\#H)^A, (B\#H)^B}, \vspace{2mm} \\
 & \qquad j:=j_{(A\ot B)\#H}, & \qquad \breve j:=j_{(A\#H)\ot (B\#H)} \vspace{2mm} \\
 & \qquad \rho':=\rho_{(A\#H)^A}, & \qquad \lambda':=\lambda_{(B\#H)^B} \vspace{2mm} \\
 & \qquad \rho:=\rho_{A\#H}, & \qquad \lambda:=\lambda_{B\#H}.
\end{array}$$
Since $(B\# H)^B$ is an $H$-Galois object, the condition $\Phi_{(B\# H)^B, A}=\Phi_{A,(B\# H)^B}^{-1}$ is fulfilled by hypothesis. We can use then the natural transformation $\zeta$ from Proposition \ref{Azcounitmult}. We are going to define a morphism $\alpha:(A\ot B)\#H\to (A\#H)\ot (B\#H)$ which
will induce $\alpha_1$. Then $\zeta^{-1}_{A', B'}\alpha_1$ will induce $\alpha_2$
and it will be an $H$-comodule algebra morphism, which would finish the proof.

 Set $\alpha=(A \otimes \Phi_{B,H} \otimes B)(A \otimes B \otimes \Delta_H)$.  We first show that $\alpha$ induces $\alpha_1$ and that they are comodule algebra morphisms.
In view of Remark \ref{morma} we should prove that $f:=\alpha j:
[(A\ot B)\# H]^{A\ot B}\to (A\# H)\ot (B\# H)$ satisfies the equality
$$\scalebox{0.9}[0.9]{
\gbeg{7}{9}
\got{1}{A\ot B} \got{1}{} \got{1}{} \got{3}{[(A\ot B)\# H]^{A\ot B}} \gnl
\gcl{1} \gvac{1} \gnot{\hspace{1,3cm}j} \glmp \gcmptb \gcmp \grmp \gnl
\gcl{1} \gvac{1} \gnot{\hspace{1,3cm}(A\ot B)\# H} \glmpb \gcmptb \gcmp \grmp \gnl
\gcl{1} \gvac{1} \gcl{1} \gcl{1} \gcmu \gnl
\gcl{1} \gvac{1} \gcl{1} \gbr \gcl{1} \gnl
\gcl{1} \gnot{\hspace{2,2cm}(A\# H)\ot (B\# H)} \glmp \gcmpt \gcmpt \gcmpt \gcmpt \grmp \gnl
\gcn{1}{1}{1}{6} \gvac{3} \hspace{-0,34cm}\gcl{1} \gnl
\gvac{3} \glm \gnl
\gvac{4} \gob{1}{(A\#H)\ot (B\#H)}
\gend} = \scalebox{0.9}[0.9]{
\gbeg{7}{10}
\got{1}{A\ot B} \got{1}{} \got{1}{} \got{3}{[(A\ot B)\# H]^{A\ot B}} \gnl
\gcl{1} \gvac{1} \gnot{\hspace{1,3cm}j} \glmp \gcmptb \gcmp \grmp \gnl
\gcl{1} \gvac{1} \gnot{\hspace{1,3cm}(A\ot B)\# H} \glmpb \gcmptb \gcmp \grmp \gnl
\gcl{1} \gvac{1} \gcl{1} \gcl{1} \gcmu \gnl
\gcl{1} \gvac{1} \gcl{1} \gbr \gcl{1} \gnl
\gcl{1} \gnot{\hspace{2,2cm}(A\# H)\ot (B\# H)} \glmp \gcmpt \gcmpt \gcmpt \gcmpt \grmp \gnl
\gcn{1}{1}{1}{6} \gvac{3} \hspace{-0,34cm}\gcl{1} \gnl
\gvac{3} \gbr \gnl
\gvac{3} \grm \gnl
\gvac{3} \gob{1}{(A\#H)\ot (B\#H)}
\gend}$$
Denote the left-hand side by $\Sigma$ and the right one by $\Omega$. Recalling the
$A\ot B$-bimodule structure of $(A\#H)\ot (B\#H)$ from Diagram \ref{(AtensB)^e-mod}, we compute
$$\Omega=\scalebox{0.8}[0.75]{
\gbeg{8}{14}
\got{1}{A} \got{1}{B} \got{1}{} \got{1}{} \got{3}{[(A\ot B)\# H]^{A\ot B}} \gnl
\gcl{1} \gcl{1} \gvac{1} \gnot{\hspace{1,3cm}j} \glmp \gcmptb \gcmp \grmp \gnl
\gcl{1} \gcl{1} \gvac{1} \gnot{\hspace{1,3cm}(A\ot B)\# H} \glmpb \gcmptb \gcmp \grmp \gnl
\gcl{1} \gcl{1} \gvac{1} \gcl{1} \gcl{1} \gcmu \gnl
\gcl{1} \gcl{1} \gvac{1} \gcl{1} \gbr \gcl{1} \gnl
\gcl{1} \gcl{1} \gnot{\hspace{2,2cm}(A\# H)\ot (B\# H)} \glmp \gcmpt \gcmpt \gcmpt \gcmpt \grmp \gnl
\gcn{1}{1}{1}{3} \gcn{1}{1}{1}{3} \gvac{1} \gcl{1} \gvac{2} \hspace{-0,32cm} \gcl{1} \gnl
\gvac{2} \hspace{-0,38cm} \gcl{1} \gbr \gvac{1} \hspace{-0,4cm} \gcn{1}{1}{4}{1} \gnl
\gvac{3} \hspace{-0,16cm} \gbr \gbr \gnl
\gvac{2} \gcn{1}{1}{3}{1} \gvac{1} \gbr \gcl{1} \gnl
\gvac{1} \glmp \gnot{\hspace{-0,4cm}A\#H}  \grmp \gvac{1} \gbr \gcn{1}{2}{1}{3} \gnl
\gvac{2} \gcl{1} \gcn{1}{1}{3}{1} \gvac{1} \glmp \gnot{\hspace{-0,3cm}B\#H}  \grmp \gnl
\gvac{2} \grm \gvac{2} \grm \gnl
\gvac{2} \gob{1}{A\#H} \gvac{3} \gob{1}{B\#H}
\gend} \stackrel{nat.}{\stackrel{(\ref{left-right mod Asm.H})}{=}}
\scalebox{0.8}[0.75]{\gbeg{6}{16}
\got{1}{A} \got{1}{B}  \got{1}{} \got{3}{\hspace{2mm} [(A\ot B)\# H]^{A\ot B}} \gnl
\gcl{1} \gcl{1} \gnot{\hspace{1,3cm}j} \glmp \gcmptb \gcmp \grmp \gnl
\gcl{1} \gcl{1} \gnot{\hspace{1,3cm}(A\ot B)\# H} \glmp \gcmptb \gcmpb \grmpb \gnl
\gcn{1}{1}{1}{3} \gcn{1}{1}{1}{3} \gvac{1} \gcl{1} \gcl{2} \gcl{3} \gnl
\gvac{1} \gcl{1} \gbr  \gnl
\gvac{1} \gbr \gbr  \gnl
\gvac{1} \gcl{1} \gbr \gbr \gnl
\gvac{1} \gcl{1} \gcl{1} \gbr \gcl{1} \gnl
\gcn{1}{1}{3}{1} \gcn{1}{1}{3}{2} \gvac{1} \hspace{-0,36cm} \gcmu \gcn{1}{1}{0}{1}
    \gvac{1} \hspace{-0,56cm} \gcn{1}{2}{0}{2} \gnl
\gvac{2} \hspace{-0,34cm} \gcl{4} \gvac{1} \hspace{-0,36cm} \gbr \gbr \gnl
\gvac{3} \gcn{1}{1}{3}{2} \gvac{1} \gbr \hspace{-0,2cm} \gcmu \gcl{1} \gnl
\gvac{4} \hspace{-0,38cm} \gcmu \gvac{1} \hspace{-0,4cm} \gcl{1} \gcl{1} \hspace{-0,2cm} \gcl{1} \gbr \gnl
\gvac{6} \hspace{-0,24cm} \gcl{1} \gbr \gcn{1}{1}{1}{2} \hspace{-0,2cm} \glm \gcl{3} \gnl
\gvac{6} \hspace{-0,34cm} \gcn{1}{1}{2}{3} \hspace{-0,02cm} \glm \gcl{2} \gvac{1} \hspace{-0,34cm} \gmu \gnl
\gvac{8} \hspace{-0,21cm} \gmu \gvac{2} \gcl{1} \gnl
\gvac{8} \gob{2}{A} \gob{1}{H} \gvac{1} \gob{1}{B} \gob{2}{H}
\gend} \stackrel{2\times nat.}{=}
\scalebox{0.8}[0.75]{\gbeg{7}{18}
\got{1}{A} \got{1}{B}  \got{1}{} \got{3}{\hspace{2mm}[(A\ot B)\# H]^{A\ot B}} \gnl
\gcl{1} \gcl{1} \gnot{\hspace{1,3cm}j} \glmp \gcmptb \gcmp \grmp \gnl
\gcl{1} \gcl{1} \gnot{\hspace{1,3cm}(A\ot B)\# H} \glmp \gcmptb \gcmpb \grmpb \gnl
\gcn{1}{1}{1}{3} \gcn{1}{1}{1}{3} \gvac{1} \gcl{1} \gcl{2} \gcl{3} \gnl
\gvac{1} \gcl{1} \gbr  \gnl
\gvac{1} \gbr \gbr  \gnl
\gvac{1} \gcl{1} \gbr \gbr \gnl
\gvac{1} \gcl{1} \gcl{1} \gbr \gcl{1} \gnl
\gcn{1}{1}{3}{0} \gcn{1}{1}{3}{0} \gvac{1} \hspace{-0,36cm} \gcmu \gcn{1}{2}{0}{2}
    \gvac{1} \hspace{-0,56cm} \gcn{1}{2}{0}{2} \gnl
\gvac{1} \gcl{4} \gcl{4} \gcn{1}{1}{3}{1} \gvac{1} \hspace{-0,34cm} \gcmu  \gnl
\gvac{3} \gcmu \gcl{1} \gbr \gcl{1} \gnl
\gvac{3}  \gcl{2} \gcl{1} \gbr \gbr  \gnl
\gvac{4} \gbr \glm \gcl{5} \gnl
\gvac{1} \gcn{1}{2}{2}{3} \gcn{1}{1}{2}{3} \glm \gcl{2} \gvac{1} \gcl{3} \gnl
\gvac{3} \gbr \gnl
\gvac{2} \gmu \gbr \gnl
\gvac{3} \hspace{-0,22cm} \gcl{1} \gvac{1} \hspace{-0,34cm} \gcl{1} \gwmu{3} \gnl
\gvac{3} \gob{1}{A} \gvac{1} \gob{1}{H} \gvac{1} \gob{1}{B} \gvac{1} \gob{1}{H}
\gend} \stackrel{(\ref{coc-coass})}{\stackrel{nat.}{=}}
\scalebox{0.8}[0.75]{\gbeg{7}{19}
\got{1}{A} \got{1}{B}  \got{1}{} \got{3}{\hspace{2mm} [(A\ot B)\# H]^{A\ot B}} \gnl
\gcl{1} \gcl{1} \gnot{\hspace{1,3cm}j} \glmp \gcmptb \gcmp \grmp \gnl
\gcl{1} \gcl{1} \gnot{\hspace{1,3cm}(A\ot B)\# H} \glmp \gcmptb \gcmp \grmpb \gnl
\gcn{1}{1}{1}{3} \gcn{1}{1}{1}{3} \gvac{1} \gcl{1} \gcl{2} \gcl{3} \gnl
\gvac{1} \gcl{1} \gbr  \gnl
\gvac{1} \gbr \gbr  \gnl
\gvac{1} \gcl{1} \gbr \gbr \gnl
\gvac{1} \gcl{1} \gcl{1} \gbr \gcl{1} \gnl
\gcn{1}{1}{3}{0} \gcn{1}{1}{3}{0} \gvac{1} \hspace{-0,36cm} \gcmu \gcn{1}{1}{0}{1}
    \gvac{1} \hspace{-0,56cm} \gcn{1}{2}{0}{3} \gnl
\gvac{1} \gcl{7} \gcl{5} \gcn{1}{1}{3}{2} \gvac{1} \gcn{1}{1}{1}{2} \gcn{1}{1}{1}{3} \gnl
\gvac{3} \gcmu \gcmu  \gcl{1} \gcl{2} \gnl
\gvac{3}  \gcl{3} \gibr \gbr \gnl
\gvac{4} \gcl{1} \gbr \gbr \gnl
\gvac{4} \gbr \glm \gcl{5} \gnl
\gvac{2} \gcn{1}{1}{1}{3} \glm \gcl{1} \gcn{1}{1}{3}{1} \gnl
\gvac{3} \gbr \gibr \gnl
\gvac{1} \gwmu{3} \gmu \gcn{1}{1}{1}{0} \gnl
\gvac{2} \gcl{1} \gvac{2} \hspace{-0,36cm} \gbr \gnl
\gvac{2} \gob{2}{A} \gvac{1} \gob{1}{H} \gob{1}{B} \gvac{1} \gob{2}{H}
\gend}
$$
$$\stackrel{nat.}{\stackrel{Prop.\ref{braid-lin}}{=}}
\scalebox{0.8}[0.75]{\gbeg{6}{17}
\got{1}{A} \got{1}{B}  \got{1}{} \got{3}{\hspace{2mm}[(A\ot B)\# H]^{A\ot B}} \gnl
\gcl{1} \gcl{1} \gnot{\hspace{1,3cm}j} \glmp \gcmptb \gcmp \grmp \gnl
\gcl{1} \gcl{1} \gnot{\hspace{1,3cm}(A\ot B)\# H} \glmp \gcmptb \gcmpb \grmpb \gnl
\gcn{1}{1}{1}{3} \gcn{1}{1}{1}{3} \gvac{1} \gcl{1} \gcl{2} \gcl{3} \gnl
\gvac{1} \gcl{1} \gbr  \gnl
\gvac{1} \gbr \gbr  \gnl
\gvac{1} \gcl{1} \gbr \gbr \gnl
\gvac{1} \gcl{1} \gcl{1} \gbr \gcl{1} \gnl
\gcn{1}{1}{3}{0} \gcn{1}{1}{3}{0} \gvac{1} \hspace{-0,36cm} \gcmu \gcn{1}{1}{0}{1}
    \gvac{1} \hspace{-0,56cm} \gcn{1}{2}{0}{1} \gnl
\gvac{1} \gcl{5} \gcl{3} \gcn{1}{1}{3}{2} \gvac{1} \gbr \gnl
\gvac{3} \gcmu \gcl{1} \gbr \gnl
\gvac{3}  \gcl{1} \gbr \gcl{1} \gcl{1} \gnl
\gvac{2} \gcn{1}{1}{1}{3} \glm \glm \gcn{1}{1}{1}{2} \gnl
\gvac{3} \gbr \gvac{1} \gcl{1} \gcmu \gnl
\gvac{1} \gwmu{3} \gwmu{3} \hspace{-0,42cm} \gcn{1}{1}{3}{1} \gvac{1} \gcl{2} \gnl
\gvac{3} \gcl{1} \gvac{2} \gbr \gnl
\gvac{3} \gob{1}{A} \gvac{2} \gob{1}{H} \gob{1}{B} \gvac{1} \gob{1}{H}
\gend} \stackrel{j}{=}
\scalebox{0.8}[0.75]{\gbeg{6}{8}
\got{1}{A} \got{1}{B}  \got{1}{} \got{3}{\hspace{2mm}[(A\ot B)\# H]^{A\ot B}} \gnl
\gcl{3} \gcl{2} \gnot{\hspace{1,3cm}j} \glmp \gcmptb \gcmp \grmp \gnl
\gvac{2} \gnot{\hspace{1,3cm}(A\ot B)\# H} \glmpb \gcmptb \gcmp \grmp \gnl
\gcl{1} \gbr \gcl{1} \gcl{2} \gnl
\gmu \gmu \gnl
\gvac{1} \hspace{-0,2cm} \gcl{1} \gvac{1} \gcl{1} \gcmu \gnl
\gvac{1} \gcl{1} \gvac{1} \gbr \gcl{1} \gnl
\gvac{1} \gob{1}{A} \gvac{1} \gob{1}{H} \gob{1}{B} \gob{1}{H}
\gend}\stackrel{nat.}{=}
\scalebox{0.8}[0.75]{\gbeg{6}{9}
\got{1}{A} \got{1}{B}  \got{1}{} \got{3}{\hspace{2mm}[(A\ot B)\# H]^{A\ot B}} \gnl
\gcl{4} \gcl{4} \gnot{\hspace{1,3cm}j} \glmp \gcmptb \gcmp \grmp \gnl
\gvac{2} \gnot{\hspace{1,3cm}(A\ot B)\# H} \glmpb \gcmptb \gcmp \grmp \gnl
\gvac{2} \gcl{1} \gcl{1} \gcmu \gnl
\gvac{2} \gcl{1} \gbr \gcl{4} \gnl
\gcl{1} \gbr \gcl{1} \gcl{2} \gnl
\gmu \gbr \gnl
\gvac{1} \hspace{-0,22cm} \gcl{1} \gvac{1} \hspace{-0,34cm} \gcl{1} \gmu \gnl
\gvac{1} \gob{2}{A} \gob{1}{H} \gob{2}{B} \gob{1}{H}
\gend}\stackrel{(\ref{left-right mod Asm.H})}{=}\Sigma.
$$
We will show that $\alpha$ is an algebra and a right $H$-comodule
morphism. For multiplicativity it is to prove that the diagrams $\Lambda$ and $\Gamma$ below are equal:
$$\Lambda:=
\scalebox{0.8}[0.8]{\gbeg{7}{5}
\gvac{1} \got{1}{(A\ot B)\# H} \gvac{3} \got{1}{(A\ot B)\# H} \gnl
\gvac{1} \hspace{-0,26cm} \gwmu{5}  \gnl
\gvac{3} \gbmp{\alpha} \gnl
\gvac{3} \gcl{1} \gnl
\gvac{3} \gob{1}{(A\# H)\ot (B\# H)}
\gend}=
\scalebox{0.8}[0.75]{\gbeg{8}{10}
\got{1}{\hspace{0,2cm}A} \got{1}{\hspace{0,2cm}B}  \gvac{2} \got{1}{\hspace{-0,2cm}H}
   \got{1}{\hspace{0,2cm}A} \got{1}{\hspace{0,2cm}B} \got{1}{\hspace{0,2cm}H}\gnl
\gvac{1} \hspace{-0,36cm} \gcl{1} \gcl{1} \gvac{1} \gcmu \gcl{1} \gcl{2} \gcl{3} \gnl
\gvac{1} \gcl{5} \gcl{3} \gcn{1}{1}{3}{2} \gvac{1} \gbr \gnl
\gvac{3} \gcmu \gcl{1} \gbr \gnl
\gvac{3}  \gcl{1} \gbr \gcl{1} \gmu \gnl
\gvac{2} \gcn{1}{1}{1}{3} \glm \glm \gvac{1} \hspace{-0,22cm} \gcl{1} \gnl
\gvac{4} \hspace{-0,34cm} \gbr \gvac{1} \gcl{1} \gcmu \gnl
\gvac{2} \gwmu{3} \gwmu{3} \hspace{-0,42cm} \gcn{1}{1}{3}{1} \gvac{1} \gcl{2} \gnl
\gvac{4} \gcl{1} \gvac{2} \gbr \gnl
\gvac{4} \gob{1}{A} \gvac{2} \gob{1}{H} \gob{1}{B} \gvac{1} \gob{1}{H}
\gend} \hspace{1cm}
\Gamma:=\scalebox{0.8}[0.8]{
\gbeg{7}{5}
\got{2}{(A\ot B)\# H} \gvac{3} \got{1}{(A\ot B)\# H} \gnl
\gvac{1} \gbmp{\alpha} \gvac{3} \gbmp{\alpha} \gnl
\gvac{1} \gwmu{5} \gnl
\gvac{3} \gcl{1} \gnl
\gvac{2} \gob{3}{(A\# H)\ot (B\# H)}
\gend}=\scalebox{0.8}[0.75]{
\gbeg{10}{11}
\got{1}{} \got{1}{A} \got{1}{B}  \gvac{1} \got{1}{\hspace{-0,3cm}H} \got{1}{A}
   \got{1}{B} \got{2}{H} \gnl
\gvac{1} \gcl{1} \gcl{1} \gcmu \gcl{1} \gcl{1} \gcmu \gnl
\gvac{1} \gcl{1} \gbr \gcl{1} \gcl{1} \gbr \gcl{2} \gnl
\gvac{1} \gcn{1}{2}{1}{-1} \gcn{1}{1}{1}{0}  \gcl{1} \gbr \gcl{1} \gcl{1} \gnl
\gvac{1} \gcmu \gbr \gbr \gcn{1}{2}{1}{3} \gcn{1}{2}{1}{3} \gnl
\gcl{2} \gcl{1} \gbr \gbr \gcn{1}{1}{1}{2} \gnl
\gvac{1} \glm \gmu \gcl{3} \gcmu \gcl{1} \gcl{2} \gnl
\gwmu{3} \gcn{1}{3}{2}{2} \gvac{2} \gcl{1} \gbr \gnl
\gvac{1} \gcl{2} \gvac{4} \glm \gmu \gnl
\gvac{5} \gwmu{3} \gcn{1}{1}{2}{2} \gnl
\gvac{1} \gob{1}{A} \gvac{1} \gob{2}{H} \gvac{1} \gob{1}{B} \gvac{1} \gob{2}{H}
\gend}
$$
We develop $\Gamma$ as follows:
$$\begin{array}{ll}
\Gamma & \stackrel{3\times nat.}{=}
\scalebox{0.8}[0.75]{\gbeg{10}{10}
\got{1}{A} \got{1}{B}  \gvac{1} \got{2}{H} \got{2}{A} \got{2}{B} \got{1}{\hspace{-0,4cm}H}\gnl
\gcl{7} \gcl{2} \gvac{1} \hspace{-0,36cm} \gwcm{3} \gcn{1}{1}{1}{2} \gvac{1} \gbr \gnl
\gvac{3} \hspace{-0,2cm} \gcmu \gcmu  \gcl{1} \gcmu \gcn{1}{1}{0}{1} \gnl
\gvac{2} \gbr \gcl{1} \gcl{1} \gbr \gcl{1} \gibr \gnl
\gvac{2} \gcl{2} \gbr \gbr \gbr \gcl{1} \gcl{2} \gnl
\gvac{3} \gcl{1} \gbr \gbr \gbr \gnl
\gvac{2} \gcl{1} \gbr \gbr \glm \gmu \gnl
\gvac{2} \glm \gmu \gwmu{3} \gvac{1} \hspace{-0,36cm} \gcl{2} \gnl
\gvac{2} \hspace{-0,34cm} \gwmu{3} \gvac{1} \hspace{-0,34cm} \gcl{1} \gvac{2} \hspace{-0,34cm}\gcl{1} \gnl
\gvac{4} \gob{1}{A} \gvac{1} \gob{2}{H} \gvac{1} \gob{1}{B} \gvac{1} \gob{2}{H}
\gend}\stackrel{coass.}{\stackrel{coc.}{=}}
\scalebox{0.8}[0.75]{\gbeg{10}{11}
\got{1}{A} \got{1}{B}  \gvac{1} \got{1}{H} \gvac{2} \got{1}{A} \got{1}{B} \got{1}{H}\gnl
\gcl{8} \gcl{3} \gwcm{3} \gvac{1} \gcl{3} \gbr \gnl
\gvac{2} \gcl{2}  \gvac{1} \hspace{-0,36cm} \gcmu
   \gvac{2} \hspace{-0,2cm} \gcn{1}{1}{1}{2} \gcn{1}{1}{1}{3} \gnl
\gvac{2} \gcl{1} \gcl{1} \gcmu \gcn{1}{1}{0}{1} \gvac{1} \gcmu \gcl{1} \gnl
\gvac{2} \gbr \gbr \gbr \gcl{1} \gibr \gnl
\gvac{2} \gcl{2} \gbr \gbr \gbr \gcl{1} \gcl{2} \gnl
\gvac{3} \gcl{1} \gbr \gbr \gbr \gnl
\gvac{2} \gcl{1} \gbr \gbr \glm \gmu \gnl
\gvac{2} \glm \gmu \gwmu{3} \gvac{1} \hspace{-0,34cm} \gcl{2} \gnl
\gvac{2} \hspace{-0,34cm} \gwmu{3} \gvac{1} \hspace{-0,36cm} \gcl{1} \gvac{2} \hspace{-0,36cm} \gcl{1} \gnl
\gvac{4} \gob{1}{A} \gvac{1} \gob{2}{H} \gvac{1} \gob{1}{B} \gvac{1} \gob{2}{H}
\gend}\stackrel{nat.}{=}
\scalebox{0.8}[0.75]{\gbeg{10}{12}
\got{1}{A} \got{1}{B}  \gvac{1} \got{1}{H} \gvac{2} \got{1}{A} \got{2}{B} \got{1}{\hspace{-0,4cm}H}\gnl
\gcl{8} \gcl{6} \gwcm{3} \gvac{1} \gcl{3} \gbr \gnl
\gvac{2} \gcl{5}  \gvac{1} \hspace{-0,36cm} \gcmu
   \gvac{2} \hspace{-0,2cm} \gcn{1}{1}{1}{2} \gcn{1}{1}{1}{3} \gnl
\gvac{2} \gcl{1} \gcl{1} \gcmu \gcn{1}{1}{0}{1} \gvac{1} \gcmu \gcl{1} \gnl
\gvac{4} \gbr \gbr \gcl{1} \gibr \gnl
\gvac{4} \gcl{1} \gbr \gbr \gcl{1} \gcl{2} \gnl
\gvac{4} \gbr \gbr \gbr \gnl
\gvac{2} \gcn{1}{1}{1}{3} \glm \gmu \glm \gmu \gnl
\gvac{3} \gbr \gcn{1}{1}{2}{1} \gvac{1} \gcn{1}{1}{3}{1} \gvac{2} \hspace{-0,36cm} \gcl{3} \gnl
\gvac{2} \hspace{-0,34cm} \gwmu{3} \gbr \gvac{1} \gcl{1} \gnl
\gvac{3} \gcl{1} \gvac{1} \gcl{1} \gwmu{3} \gnl
\gvac{3} \gob{1}{A} \gvac{1} \gob{1}{H} \gvac{1} \gob{1}{B} \gvac{2} \gob{2}{H}
\gend} \vspace{3mm} \\
 & \stackrel{nat.}{=}
\scalebox{0.8}[0.75]{\gbeg{9}{12}
\got{1}{\hspace{0,2cm}A} \got{1}{\hspace{0,2cm}B}  \gvac{1} \got{1}{\hspace{-0,2cm}H}
   \got{1}{\hspace{0,2cm}A} \gvac{1} \got{1}{\hspace{0,2cm}B} \got{1}{\hspace{0,2cm}H}\gnl
\gvac{1} \hspace{-0,36cm} \gcl{1} \gcl{1} \gcmu \gcl{1} \gvac{1} \gbr \gnl
\gvac{1} \gcl{7} \gcl{5} \gcl{1} \gbr \gvac{1} \gcl{1} \gcn{1}{1}{1}{2} \gnl
\gvac{3} \glm  \gcl{1} \gvac{1} \hspace{-0,35cm} \gcmu \gcl{1} \gnl
\gvac{4} \gcn{1}{1}{2}{0} \gcmu \gcl{1} \gibr \gnl
\gvac{4} \hspace{-0,22cm} \gcl{2} \gcmu \hspace{-0,2cm} \gbr \gcl{1} \gcl{2} \gnl
\gvac{6} \hspace{-0,36cm} \gbr \gcn{1}{1}{0}{1} \hspace{-0,2cm} \gbr \gnl
\gvac{5} \hspace{-0,22cm} \gbr \gcl{1} \gbr \gcn{1}{1}{0}{1} \hspace{-0,2cm} \gmu \gnl
\gvac{6} \hspace{-0,22cm} \gcl{1} \gcn{1}{1}{1}{2} \gmu \glm \gvac{1} \hspace{-0,4cm} \gcl{3} \gnl
\gvac{6} \hspace{-0,14cm} \gmu \gvac{1} \hspace{-0,24cm} \gbr \gvac{1} \gcn{1}{1}{2}{1} \gnl
\gvac{7} \gcl{1} \gvac{1} \gcl{1} \gwmu{3} \gnl
\gvac{7} \gob{1}{A} \gvac{1} \gob{1}{H} \gvac{1} \gob{1}{B} \gvac{1} \gob{2}{H}
\gend} \stackrel{coass.}{\stackrel{nat.}{=}}
\scalebox{0.8}[0.75]{\gbeg{9}{11}
\got{1}{A} \got{1}{B} \got{2}{H} \got{1}{A} \gvac{2} \got{1}{B} \got{1}{H}\gnl
\gcl{1} \gcl{1} \gcmu \gcl{1} \gvac{2} \gbr \gnl
\gcl{6} \gcl{4} \gcl{1} \gbr \gvac{2} \gcn{1}{1}{1}{0} \gcl{2} \gnl
\gvac{2} \glm \hspace{-0,2cm} \gcmu \gvac{1} \hspace{-0,22cm} \gcmu \gnl
\gvac{4} \gcn{1}{1}{1}{-1} \hspace{-0,2cm} \gcl{3} \hspace{-0,22cm} \gcmu \gcl{1} \gibr \gnl
\gvac{4} \gcl{1} \gvac{1} \gcl{1} \gbr \gcl{1} \gcl{2} \gnl
\gvac{3} \gbr \gvac{1} \gmu \gbr \gnl
\gvac{3} \gcl{1} \gcn{1}{1}{1}{2} \gvac{1} \hspace{-0,35cm} \gbr \gcn{1}{1}{2}{1}
   \gvac{1} \hspace{-0,34cm} \gmu \gnl
\gvac{3} \hspace{-0,02cm} \gmu \gvac{1} \hspace{-0,21cm} \gbr
   \glm \gvac{2} \hspace{-0,4cm} \gcl{2} \gnl
\gvac{5} \hspace{-0,16cm} \gcl{1} \gvac{1} \gcl{1} \gwmu{3} \gnl
\gvac{5} \gob{1}{A} \gvac{1} \gob{1}{H} \gvac{1} \gob{1}{B} \gvac{2} \gob{1}{H}
\gend}\stackrel{Prop. \ref{braid-lin}}{\stackrel{nat.}{=}}
\scalebox{0.8}[0.75]{\gbeg{9}{12}
\got{1}{A} \got{1}{B} \got{2}{H} \got{1}{A} \gvac{2} \got{1}{B} \got{1}{H}\gnl
\gcl{1} \gcl{1} \gcmu \gcl{1} \gvac{2} \gbr \gnl
\gcl{7} \gcl{4} \gcl{1} \gbr \gvac{2} \gcn{1}{1}{1}{0} \gcl{4} \gnl
\gvac{2} \glm \hspace{-0,2cm} \gcmu \gvac{1} \hspace{-0,22cm} \gcmu \gnl
\gvac{4} \gcn{1}{1}{1}{-1} \hspace{-0,2cm} \gcl{3} \hspace{-0,22cm} \gcmu \gcl{1} \gcl{2} \gnl
\gvac{4} \gcl{1} \gvac{1} \gcl{1} \gbr \gvac{1} \gnl
\gvac{3} \gbr \gvac{1} \gmu \gmu \gcn{1}{1}{1}{0} \gnl
\gvac{3} \gcl{2} \gcl{1} \gvac{1} \hspace{-0,35cm} \gbr \gvac{1} \gbr \gnl
\gvac{4} \gcn{1}{1}{2}{3} \gvac{1} \gcl{1} \gcl{1} \gcn{1}{1}{3}{1} \gvac{1} \gcl{3} \gnl
\gvac{3} \hspace{-0,34cm} \gmu \gvac{1} \hspace{-0,21cm} \gbr \glm \gnl
\gvac{4} \gcl{1} \gvac{1} \gcl{1} \gwmu{3} \gnl
\gvac{4} \gob{1}{A} \gvac{1} \gob{1}{H} \gvac{1} \gob{1}{B} \gvac{2} \gob{1}{H}
\gend} \vspace{3mm} \\
 & \stackrel{bialg.}{=}
\scalebox{0.8}[0.75]{\gbeg{9}{11}
\got{1}{A} \got{1}{B} \got{2}{H} \got{1}{A} \gvac{2} \got{1}{B} \got{1}{H}\gnl
\gcl{1} \gcl{1} \gcmu \gcl{1} \gvac{2} \gbr \gnl
\gcl{6} \gcl{4} \gcl{1} \gbr \gvac{2} \gcn{1}{1}{1}{0} \gcl{3} \gnl
\gvac{2} \glm \hspace{-0,2cm} \gcmu \gvac{1} \gcl{1} \gnl
\gvac{3} \gcn{1}{1}{2}{0} \gcl{2} \gwmu{3} \gnl
\gvac{3} \hspace{-0,36cm} \gcl{1} \gvac{2} \hspace{-0,34cm} \gwcm{3} \gcn{1}{1}{2}{1} \gnl
\gvac{3} \hspace{-0,38cm} \gbr \gvac{1} \hspace{-0,2cm} \gbr \gvac{1}  \gbr \gnl
\gvac{4} \hspace{-0,22cm} \gcl{1} \gcn{1}{1}{1}{2} \gvac{1} \hspace{-0,34cm}
   \gcl{1} \gcl{1} \gcn{1}{1}{3}{1} \gvac{1} \gcl{3} \gnl
\gvac{4} \hspace{-0,34cm} \gmu \gvac{1} \hspace{-0,21cm} \gbr \glm \gnl
\gvac{5} \gcl{1} \gvac{1} \gcl{1} \gwmu{3} \gnl
\gvac{5} \gob{1}{A} \gvac{1} \gob{1}{H} \gvac{1} \gob{1}{B} \gvac{2} \gob{1}{H}
\gend}=:\Gamma'
\end{array}$$

Observe that
$$
\gbeg{3}{5}
\got{1}{H} \got{1}{H} \got{1}{B} \gnl
\gmu \gcl{2} \gnl
\gcmu \gnl
\gcl{1} \gbr \gnl
\gob{1}{H} \gob{1}{B} \gob{1}{H}
\gend=
\gbeg{5}{7}
\got{1}{H} \got{1}{H} \got{1}{B} \gnl
\gcl{1} \gbr \gnl
\gbr \gcl{1} \gnl
\gcl{2} \gmu \gnl
\gvac{1} \gcmu \gnl
\gibr \gcl{1} \gnl
\gob{1}{H} \gob{1}{B} \gob{1}{H}
\gend
$$
Note that the left-hand side of this diagram appears in $\Gamma'$. We substitute the right-hand side diagram in $\Gamma'$ and obtain
$$\Gamma'=
\scalebox{0.8}[0.75]{\gbeg{8}{13}
\got{1}{A} \got{1}{B} \got{2}{H} \got{1}{A} \gvac{1} \got{1}{B} \got{1}{H}\gnl
\gcl{1} \gcl{1} \gcmu \gcl{1} \gvac{1} \gbr \gnl
\gcl{5} \gcl{4} \gcl{1} \gbr \gvac{1} \gcn{1}{1}{1}{0} \gcn{1}{1}{1}{0} \gnl
\gvac{2} \glm \hspace{-0,2cm} \gcmu \gcl{1} \gcl{1} \gnl
\gvac{3} \gcn{1}{1}{2}{0} \gcl{5} \gcl{1} \gbr \gnl
\gvac{3} \gcn{1}{1}{0}{0} \gvac{1} \gbr \gcl{1} \gnl
 \gvac{2} \hspace{-0,36cm} \gbr \gvac{2} \hspace{-0,2cm} \gcl{2} \gmu  \gnl
\gvac{2} \hspace{-0,22cm} \gmu \gcn{1}{1}{1}{2} \gvac{3} \hspace{-0,34cm} \gcmu \gnl
\gvac{3} \gcl{3} \gvac{1} \gcl{2} \gvac{1} \gibr \gcl{4} \gnl
\gvac{6} \gbr \gcl{1} \gnl
\gvac{5} \gbr \glm \gnl
\gvac{3} \gcl{1} \gvac{1} \gcl{1} \gwmu{3} \gnl
\gvac{3} \gob{1}{A} \gvac{1} \gob{1}{H} \gvac{1} \gob{1}{B} \gvac{1} \gob{1}{H}
\gend} \hspace{-0,4cm}\stackrel{Prop. \ref{braid-lin}}{\stackrel{nat.}{=}}\hspace{-0,2cm}
\scalebox{0.8}[0.75]{\gbeg{7}{11}
\got{1}{A} \got{1}{B} \got{2}{H} \got{1}{A} \gvac{1} \got{1}{B} \got{1}{H}\gnl
\gcl{1} \gcl{1} \gcmu \gcl{1} \gvac{1} \gbr \gnl
\gcl{5} \gcl{4} \gcl{1} \gbr \gvac{1} \gcn{1}{1}{1}{0} \gcn{1}{1}{1}{0} \gnl
\gvac{2} \glm \hspace{-0,2cm} \gcmu \gcl{1} \gcl{1} \gnl
\gvac{3} \gcn{1}{1}{2}{0} \gcl{3} \gcl{1} \gbr \gnl
\gvac{3} \gcn{1}{1}{0}{0} \gvac{1}  \gbr \gcl{1} \gnl
 \gvac{2} \hspace{-0,36cm} \gbr \gvac{2} \hspace{-0,2cm} \gcl{1} \gmu  \gnl
\gvac{2} \hspace{-0,22cm} \gmu \gcn{1}{1}{1}{2} \gvac{1} \hspace{-0,34cm} \glm \gcmu \gnl
\gvac{3} \gcl{2} \gvac{1} \gwmu{3} \gcn{1}{1}{1}{-1} \gcl{2} \gnl
\gvac{6} \gbr \gnl
\gvac{3} \gob{1}{A} \gvac{1} \gob{1}{H} \gvac{1} \gob{1}{B} \gvac{1} \gob{1}{H}
\gend} \stackrel{nat.}{\stackrel{coass.}{=}}\hspace{-0,2cm}
\scalebox{0.8}[0.75]{\gbeg{8}{10}
\got{1}{\hspace{0,2cm}A} \got{1}{\hspace{0,2cm}B}  \gvac{2} \got{1}{\hspace{-0,2cm}H}
   \got{1}{\hspace{0,2cm}A} \got{1}{\hspace{0,2cm}B} \got{1}{\hspace{0,2cm}H}\gnl
\gvac{1} \hspace{-0,36cm} \gcl{1} \gcl{1} \gvac{1} \gcmu \gcl{1} \gbr \gnl
\gvac{1} \gcl{5} \gcl{3} \gcn{1}{1}{3}{2} \gvac{1} \gbr \gbr \gnl
\gvac{3} \gcmu \gcl{1} \gbr \gcl{1} \gnl
\gvac{3}  \gcl{1} \gbr \gcl{1} \gmu \gnl
\gvac{2} \gcn{1}{1}{1}{3} \glm \glm \gvac{1} \hspace{-0,2cm} \gcl{1} \gnl
\gvac{4} \hspace{-0,34cm} \gbr \gvac{1} \gcl{1} \gcmu \gnl
\gvac{2} \gwmu{3} \gwmu{3} \hspace{-0,42cm} \gcn{1}{1}{3}{1} \gvac{1} \gcl{2} \gnl
\gvac{4} \gcl{1} \gvac{2} \gbr \gnl
\gvac{4} \gob{1}{A} \gvac{2} \gob{1}{H} \gob{1}{B} \gvac{1} \gob{1}{H}
\gend} \hspace{-0,1cm}\stackrel{Prop. \ref{braid-lin}}{=} \Lambda.
$$
This proves that $\alpha$ is compatible with multiplication. The compatibility
with unit is obvious. \par \smallskip

We view $(A\ot B)\# H$ as a right $H$-comodule by the structure morphism
$(A\ot B)\# \Delta_H$ and $(A\# H)\ot (B\# H)$ by $(A\# H)\ot (B\# \Delta_H)$. That
$\alpha$ is right $H$-colinear is true by the coassociativity of $H$. In Corollary
\ref{(Asm.H)A+j com.alg.} we have proved that $([(A\ot B)\# H]^{A\ot B},
j_{(A\ot B)\# H})$ is an $H$-comodule algebra pair. Analogously as in Lemma \ref{(AsmashH)A-alg} one proves that
$(((A\# H)\ot (B\# H))^{A\ot B}, j_{(A\# H)\ot (B\# H)})$ is an algebra pair. On
the other hand, it is immediate that the right $H$-comodule structure morphism $(A\# H)\ot \rho_{B\# H}$ is
$(A \otimes B)^e$-linear. Thus, by Lemma \ref{M^A-comod},
$(((A\# H)\ot (B\# H))^{A\ot B}, j_{(A\# H)\ot (B\# H)})$ is an $H$-comodule pair.
As a subcomodule and a subalgebra of an $H$-comodule algebra, $((A\# H)\ot (B\# H))^{A\ot B}$
is such as well, 2.2.1. Having that $\alpha$ is an $H$-comodule algebra
morphism, we obtain as for 2.2.3 that $\alpha_1$ is such too.

We now prove that $\zeta^{-1}_{A', B'}\alpha_1$ induces $\alpha_2$ and that they are comodule
algebra morphisms. Note that the triangle $\langle2\rangle$ in Diagram \ref{diagp69} is the one from
Proposition \ref{Azcounitmult}. From Proposition \ref{(Asm.H)A is Gal.} we know that
$(B\# H)^B$ is an $H$-Galois object. By hypothesis, the conditions of Proposition \ref{Azcounitmult} are fulfilled and we have that $\langle2\rangle$ commutes. Further, the square $\langle3\rangle$ commutes by the
way $\rho'$ and $\lambda'$ are induced. Observe that
$(A\# H)^A$ and $(B\# H)^B$ are equalizers and that both are (faithfully) flat.
Therefore we have that $((A\# H)^A\ot (B\# H)^B, \id_{(A\# H)^A}\ot j_{B'})$ and
$((A\# H)^A\ot (B\# H)^B, j_{A'}\ot \id_{(B\# H)^B})$ are equalizers. By
flatness of $H$ we get further that $(A\# H)^A\ot (B\# H)^B\ot H$, and consequently
$(A\# H)^A\ot H\ot (B\# H)^B,$ are equalizers with the respective morphisms. This gives us in particular that
$j_{A'}\ot H\ot j_{B'}$ is a monomorphism, as a composition of the latter two
equalizer morphisms.

We furthermore have $((A\# H)\ot\lambda)\alpha=(\rho\ot (B\# H))\alpha$.
Indeed,
$$\scalebox{0.8}[0.8]{
\gbeg{5}{7}
\got{1}{A} \got{1}{B} \got{2}{H} \gnl
\gcl{1} \gcl{1} \gcmu \gnl
\gcl{1} \gbr \gcn{1}{1}{1}{2} \gnl
\gcl{1} \gcl{1} \gcl{1} \gcmu \gnl
\gcl{1} \gcl{1} \gcl{1} \gibr \gnl
\gcl{1} \gcl{1} \gibr \gcl{1} \gnl
\gob{1}{A} \gob{1}{H} \gob{1}{H} \gob{1}{B} \gob{1}{H}
\gend}\stackrel{coc.}{=}
\scalebox{0.8}[0.8]{\gbeg{5}{6}
\got{1}{A} \got{1}{B} \got{2}{H} \gnl
\gcl{1} \gcl{1} \gcmu \gnl
\gcl{1} \gbr \gcn{1}{1}{1}{2} \gnl
\gcl{1} \gcl{1} \gcl{1} \gcmu \gnl
\gcl{1} \gcl{1} \gibr \gcl{1} \gnl
\gob{1}{A} \gob{1}{H} \gob{1}{H} \gob{1}{B} \gob{1}{H}
\gend} \stackrel{coass.}{=}
\scalebox{0.8}[0.8]{\gbeg{5}{7}
\got{1}{A} \got{1}{B} \got{1}{} \got{2}{H} \gnl
\gcl{5} \gcl{3} \gvac{1} \gcmu \gnl
\gvac{3} \hspace{-0,22cm} \gcmu \hspace{-0,2cm} \gcl{4} \gnl
\gvac{3} \hspace{-0,16cm} \gcn{1}{1}{2}{1} \gcn{1}{1}{2}{2} \gnl
\gvac{2} \gbr \gcn{1}{1}{2}{1} \gnl
\gvac{2} \gcl{1} \gibr \gnl
\gvac{1} \gob{1}{A} \gob{1}{H} \gob{1}{H} \gob{1}{B} \gob{1}{H}
\gend} \stackrel{Prop. \ref{braid-lin}}{\stackrel{nat.}{=}}
\scalebox{0.8}[0.8]{\gbeg{5}{6}
\got{1}{A} \got{1}{B} \got{2}{H} \gnl
\gcl{1} \gcl{1} \gcmu \gnl
\gcl{1} \gbr \gcl{2} \gnl
\gcl{1} \gcl{1} \gcl{1} \gnl
\gcl{1} \hspace{-0,2cm} \gcmu \hspace{-0,22cm} \gcl{1} \gcl{1} \gnl
\gvac{1} \gob{1}{A} \gob{1}{H} \gob{1}{B} \gob{1}{H}
\gend}
$$
We now compute looking at the Diagram \ref{diagp69}:
$$\begin{array}{l}
\hspace{-0,2cm}(j_{A'}\ot H\ot j_{B'})(\rho'\ot\id_{(B\# H)^B})\zeta^{-1}_{A', B'}
\alpha_1  = \\
 \hspace{3,5cm} \stackrel{\langle3\rangle}{=} (\rho\ot\id_{B'})(j_{A'}\ot j_{B'})\zeta^{-1}_{A', B'}\alpha_1 \\
 \hspace{3,5cm} \stackrel{\langle2\rangle}{=} (\rho\ot\id_{B'})\breve j\alpha_1 \\
 \hspace{3,5cm}  \stackrel{\langle1\rangle}{=} (\rho\ot\id_{B'})\alpha j \\
 \hspace{3,5cm} = (\id_{A'}\ot\lambda)\alpha j \\
  \hspace{3,5cm} \stackrel{\langle1\rangle}{=} (\id_{A'}\ot\lambda)\breve j\alpha_1\\
   \hspace{3,5cm} \stackrel{\langle2\rangle}{=} (\id_{A'}\ot\lambda)(j_{A'}\ot j_{B'})
\zeta^{-1}_{A', B'}\alpha_1 \\
 \hspace{3,5cm} \stackrel{\langle3\rangle}{=} (j_{A'}\ot H\ot j_{B'})(\id_{(A\# H)^A}
\ot\lambda')\zeta^{-1}_{A', B'}\alpha_1.
\end{array}$$
Since $j_{A'}\ot H\ot j_{B'}$ is a monomorphism, we obtain that $\zeta^{-1}_{A', B'}\alpha_1$
induces $\alpha_2$ so that the diagram $\langle4\rangle$ commutes.

It remains to prove that $\alpha_2$ is an $H$-comodule algebra morphism.
From Corollary \ref{(Asm.H)A+j com.alg.} we know that $((A\# H)^A, j_{A\# H})$ and $((B\# H)^B, j_{B\# H})$
are $H$-comodule algebra pairs. Viewing $(A\# H)\ot (B\# H)$ and
$(A\# H)^A\ot (B\# H)^B$ as $H$-comodules via $(A\# H)\ot \rho_{B\# H}$ and $(A\# H)^A\ot
\rho_{(B\# H)^B}$ respectively, we have that they are $H$-comodule algebras.
We have commented before that $(((A\# H)\ot (B\# H))^{A\ot B}, j_{(A\# H)\ot (B\# H)})$ is an
$H$-comodule algebra pair. Now 2.2.3(a) applies to the triangle $\langle2\rangle$,
giving us that  $\zeta^{-1}_{A', B'}$ is an $H$-comodule algebra morphism. That $\alpha_1$ is
such we have seen above. The equalizer morphism $e: (A\# H)^A\Box_H (B\# H)^B\to (A\# H)^A\ot (B\# H)^B$
is an $H$-comodule algebra one, since $(A\# H)^A$ and $(B\# H)^B$ are $H$-comodule algebras. Recall that $(A\# H)^A\Box_H (B\# H)^B$ has the $H$-comodule structure via $(A\# H)^A\Box_H \rho_{(B\# H)^B}$. This time
2.2.3(a) applies to the diagram $\langle4\rangle$ inside Diagram \ref{diagp69} and we obtain the
statement on $\alpha_2$.
\qed\end{proof}

\begin{proposition}\label{gmorphism}
Let $\C$ be a closed braided monoidal category with equalizers and coequalizers. Let $H$
be a flat and commutative Hopf algebra. Suppose that the braiding is $H$-linear and that $\Phi_{T,X}=\Phi^{-1}_{X,T}$ for any $H$-Galois object $T$ and $X \in \C$. The map
$$\Upsilon: \BM(\C; H) \to \Gal(\C; H), \quad [A]\mapsto[(A\# H)^A]$$
is a group morphism.
\end{proposition}

\begin{proof}
Since the braiding is $H$-linear, $H$ is cocommutative, Proposition \ref{braid-lin}. In virtue of Theorem
\ref{abelian subgroup}, $\Gal(\C; H)$ is a group. By Proposition \ref{(Asm.H)A is Gal.}, $(A\# H)^A$ is an $H$-Galois object for any $H$-Azumaya algebra $A$. We will show that $\Upsilon$ does not depend on the representative of a class in $\BM(\C; H)$ and that it is compatible with product.

We will first prove that $\Upsilon([M, M])=[H]$ in $\Gal(\C; H)$ for any faithfully projective $H$-module $M$ in
$\C$ (abusing of notation we suppressed the brackets to denote the class of $[M,M]$). To this end we will define a morphism $\crta\sigma: H\to ([M, M]\# H)^{[M, M]}$. By Proposition \ref{comod-alg-iso} it will be an isomorphism of $H$-Galois objects provided that it is an $H$-comodule algebra morphism. Let $\teta: H \to [M, M]$ be the algebra morphism induced by the $H$-module structure of $M$ as in Lemma \ref{mod vs algEnd}. The morphism $\crta\sigma$ will be induced by $\sigma=(\theta S \otimes H)\Delta_H:H \rightarrow [M,M] \otimes H.$ We check that it factors through $([M, M]\# H)^{[M, M]}$. Due to Remark \ref{morma} and Diagram \ref{make morf. on M^A} this will hold if we prove that
$$\scalebox{0.9}[0.9]{
\gbeg{4}{4}
\got{1}{[M, M]} \got{1}{} \got{1}{H} \gnl
\gcn{1}{1}{1}{3} \gvac{1} \gbmp{\sigma} \gnl
\gvac{1} \glm \gnl
\gvac{2} \gob{1}{[M, M]\# H}
\gend}
=\scalebox{0.9}[0.9]{
\gbeg{3}{5}
\got{1}{[M, M]} \got{1}{} \got{1}{H} \gnl
\gcn{1}{1}{1}{3} \gvac{1} \gbmp{\sigma} \gnl
\gvac{1} \gbr \gnl
\gvac{1} \grm \gnl
\gvac{1} \gob{1}{[M, M]\# H}
\gend}
$$
Applying the structure of an $[M, M]$-bimodule on $[M, M]\# H$ described in Diagram
\ref{left-right mod Asm.H}, we get that the above equality amounts to
$$\scalebox{0.8}[0.8]{
\gbeg{5}{5}
\got{1}{[M, M]} \got{3}{H} \gnl
\gcl{2} \gvac{1} \gbmp{\sigma} \gnl
\gvac{1} \glmp \gnot{\hspace{0,34cm}[M, M]\# H} \gcmpt \gcmp \grmp \gnl
\gwmu{3} \gcl{1} \gnl
\gvac{1} \gob{1}{[M, M]} \gvac{1} \gob{1}{H}
\gend}=
\scalebox{0.8}[0.75]{\gbeg{5}{12}
\got{1}{[M, M]} \got{3}{H} \gnl
\gcl{2} \gvac{1} \gbmp{\sigma} \gnl
\gvac{1} \glmp \gnot{\hspace{0,34cm}[M, M]\# H} \gcmptb \gcmpb \grmp \gnl
\gcn{1}{1}{1}{3} \gvac{1} \gcl{1} \gcl{2} \gnl
\gvac{1} \gbr \gcl{1} \gnl
\gvac{1} \gcn{1}{1}{1}{0} \gbr  \gnl
\gvac{1} \hspace{-0,2cm} \gcl{3} \gcmu \hspace{-0,2cm}\gcn{1}{1}{1}{2} \gvac{1} \gnl
\gvac{3} \hspace{-0,5cm} \gcl{1} \gbr \gnl
\gvac{4} \hspace{-0,42cm} \glm \gcl{3} \gnl
\gvac{4} \hspace{-0,36cm} \gcn{1}{1}{0}{2} \gvac{1} \hspace{-0,34cm} \gcl{1}   \gnl
\gvac{5} \gmu \gnl
\gvac{5} \gob{2}{[M, M]} \gob{2}{\hspace{-0,2cm}H}
\gend}
$$
The left-hand side diagram we denote by $\Sigma$ and the right one by $\Omega$.
We develop $\Omega$ as below, where in the third equation we apply the $H$-module structure of $[M, M]$ described in Lemma \ref{mod vs algEnd}(ii),
$$\begin{array}{lll}
\Omega & \stackrel{nat.}{=}
\scalebox{0.8}[0.75]{\gbeg{5}{11}
\gvac{1} \got{1}{[M, M]} \gvac{1} \got{1}{H} \gnl
\gvac{1} \gcn{1}{1}{1}{3} \gvac{1} \gcl{1} \gnl
\gvac{2} \gbr \gnl
\gvac{2} \gbmp{\sigma} \gcn{1}{1}{1}{3} \gnl
\glmp \gnot{\hspace{0,34cm}[M, M]\# H} \gcmpb \gcmpt \grmp \gcl{2} \gnl
\gvac{1} \gcl{3} \gcmu \gnl
\gvac{3} \hspace{-0,42cm} \gcl{1} \gbr \gnl
\gvac{4} \hspace{-0,42cm} \glm \gcl{3} \gnl
\gvac{4} \hspace{-0,36cm} \gcn{1}{1}{0}{2} \gvac{1} \hspace{-0,34cm} \gcl{1}   \gnl
\gvac{5} \gmu \gnl
\gvac{5} \gob{2}{[M, M]} \gob{2}{\hspace{-0,2cm}H}
\gend} \stackrel{\sigma}{=}
\scalebox{0.8}[0.75]{\gbeg{5}{13}
\gvac{1} \got{1}{[M, M]} \gvac{1} \got{1}{H} \gnl
\gvac{1} \gcn{1}{1}{1}{3} \gvac{1} \gcl{1} \gnl
\gvac{2} \gbr \gnl
\gvac{2} \hspace{-0,2cm} \gcmu \gcn{1}{1}{0}{3} \gnl
\gvac{2} \gnot{S}\gmp \gnl \gcl{1} \gvac{1} \gcl{3} \gnl
\gvac{2} \gbmp{\teta} \gcl{1} \gnl
\gvac{1} \glmp \gnot{\hspace{0,34cm}[M, M]\# H} \gcmpb \gcmptb \grmp \gnl
\gvac{2} \gcl{3} \hspace{-0,22cm} \gcmu \hspace{-0,2cm} \gcn{1}{1}{3}{2} \gnl
\gvac{4} \hspace{-0,34cm} \gcl{1} \gbr \gnl
\gvac{5} \hspace{-0,42cm} \glm \gcl{3} \gnl
\gvac{5} \hspace{-0,34cm} \gcn{1}{1}{1}{2} \gvac{1} \hspace{-0,34cm} \gcl{1}   \gnl
\gvac{6} \gmu \gnl
\gvac{6} \gob{2}{[M, M]} \gob{2}{\hspace{-0,2cm}H}
\gend} \hspace{2mm} \stackrel{Lem.\ \ref{mod vs algEnd}}{=}
\scalebox{0.8}[0.75]{\gbeg{6}{13}
\gvac{1} \got{1}{[M, M]} \gvac{1} \got{1}{H} \gnl
\gvac{1} \gcn{1}{1}{1}{3} \gvac{1} \gcl{1} \gnl
\gvac{2} \gbr \gnl
\gvac{1} \gwcm{3} \gcn{1}{1}{-1}{3} \gnl
\gcn{1}{1}{3}{2} \gvac{1} \gwcm{3} \gcl{1} \gnl
\gvac{1} \hspace{-0,34cm} \gnot{S}\gmp \gnl \gvac{1} \hspace{-0,44cm} \gcmu
   \gvac{1} \hspace{-0,2cm} \gbr \gnl
\gvac{3} \hspace{-0,2cm} \gbmp{\teta} \gbmp{\teta} \gnot{S}\gmp \gnl
   \gcn{1}{1}{2}{1} \gvac{1} \hspace{-0,36cm} \gcl{6} \gnl
\gvac{3} \gcn{1}{1}{2}{1} \gcn{1}{1}{2}{1} \gvac{1} \hspace{-0,34cm} \gbr \gnl
\gvac{4} \hspace{-0,22cm} \gcl{3} \gcl{2} \gvac{1} \hspace{-0,34cm} \gcl{1} \gbmp{\teta} \gnl
\gvac{7} \gmu \gnl
\gvac{6} \hspace{-0,22cm} \gwmu{3} \gnl
\gvac{5} \gwmu{3} \gnl
\gvac{6} \gob{1}{[M, M]} \gvac{3} \gob{1}{H}
\gend} \stackrel{coass.}{\stackrel{ass.}{=}}
\scalebox{0.8}[0.75]{\gbeg{6}{11}
\gvac{1} \got{1}{[M, M]} \gvac{1} \got{1}{H} \gnl
\gvac{1} \gcn{1}{1}{1}{3} \gvac{1} \gcl{1} \gnl
\gvac{2} \gbr \gnl
\gvac{1} \gwcm{3} \gcn{1}{1}{-1}{2} \gnl
\gvac{1} \hspace{-0,36cm} \gcmu \gcmu \gcl{1} \gnl
\gvac{1} \gnot{S}\gmp \gnl \gcl{1} \gnot{S}\gmp \gnl \gbr \gnl
\gvac{1} \gbmp{\teta} \gbmp{\teta} \gbmp{\teta} \gcl{1} \gcl{4} \gnl
\gvac{1} \gcl{1} \gcl{1} \gbr \gnl
\gvac{1} \gmu \gmu \gnl
\gvac{2} \hspace{-0,2cm} \gwmu{3} \gnl
\gvac{3} \gob{1}{[M, M]} \gvac{1} \gob{2}{H}
\gend} & \vspace{3mm} \\
 & \stackrel{\teta}{\stackrel{alg.m.}{=}}
\scalebox{0.8}[0.75]{\gbeg{6}{11}
\gvac{1} \got{1}{[M, M]} \gvac{1} \got{1}{H} \gnl
\gvac{1} \gcn{1}{1}{1}{3} \gvac{1} \gcl{1} \gnl
\gvac{2} \gbr \gnl
\gvac{1} \gwcm{3} \gcn{1}{1}{-1}{2} \gnl
\gvac{1} \hspace{-0,36cm} \gcmu \gcmu \gcl{1} \gnl
\gvac{1} \gnot{S}\gmp \gnl \gcl{1} \gnot{S}\gmp \gnl \gbr \gnl
\gvac{1} \gmu \gbmp{\teta} \gcl{1} \gcl{4} \gnl
\gvac{2} \hspace{-0,36cm} \gbmp{\teta} \gvac{1} \hspace{-0,36cm} \gbr \gnl
\gvac{3} \hspace{-0,2cm} \gcl{1} \gvac{1} \hspace{-0,36cm} \gmu \gnl
\gvac{4} \hspace{-0,2cm} \gwmu{3} \gnl
\gvac{5} \gob{1}{[M, M]} \gvac{1} \gob{2}{H}
\gend} \stackrel{antip.}{=}
\scalebox{0.8}[0.75]{\gbeg{5}{11}
\got{1}{[M, M]} \gvac{1} \got{1}{H} \gnl
\gcn{1}{1}{1}{3} \gvac{1} \gcl{1} \gnl
\gvac{1} \gbr \gnl
\gwcm{3} \gcn{1}{1}{-1}{2} \gnl
\gcu{1} \gvac{1} \hspace{-0,36cm} \gcmu \gcl{1} \gnl
\gvac{2} \gnot{S}\gmp \gnl \gbr \gnl
\gvac{1} \hspace{-0,2cm} \gu{1} \gvac{1} \hspace{-0,36cm} \gbmp{\teta} \gcl{1} \gcl{4} \gnl
\gvac{2} \hspace{-0,36cm} \gbmp{\teta} \gvac{1} \hspace{-0,36cm} \gbr \gnl
\gvac{3} \hspace{-0,2cm} \gcl{1} \gvac{1} \hspace{-0,36cm} \gmu \gnl
\gvac{4} \hspace{-0,2cm} \gwmu{3} \gnl
\gvac{5} \gob{1}{[M, M]} \gvac{1} \gob{2}{H}
\gend} \stackrel{\teta}{\stackrel{alg.m.}{=}}
\scalebox{0.8}[0.75]{\gbeg{4}{9}
\got{1}{[M, M]} \gvac{1} \got{1}{H} \gnl
\gcn{1}{1}{1}{3} \gvac{1} \gcl{1} \gnl
\gvac{1} \gbr \gnl
\gvac{1} \hspace{-0,22cm} \gcmu \gcn{1}{1}{0}{1} \gnl
\gvac{1} \gnot{S}\gmp \gnl \gbr \gnl
\gvac{1} \gbmp{\teta} \gcl{1} \gcl{3} \gnl
\gvac{1} \gbr \gnl
\gvac{1} \gmu \gnl
\gvac{1} \gob{1}{[M, M]} \gvac{1} \gob{1}{H}
\gend} \stackrel{Prop. \ref{braid-lin}}{\stackrel{nat.}{=}}
\scalebox{0.8}[0.75]{\gbeg{4}{6}
\got{1}{[M, M]} \gvac{1} \got{2}{H} \gnl
\gcl{3} \gvac{1} \gcmu \gnl
\gvac{2} \gnot{S}\gmp \gnl \gcl{3} \gnl
\gvac{2} \gbmp{\teta} \gcl{1} \gnl
\gwmu{3} \gnl
\gvac{1} \gob{1}{[M, M]} \gvac{1} \gob{1}{H}
\gend} & \stackrel{\sigma}{=}\Sigma.
\end{array}$$
This proves that $\sigma$ induces the morphism $\crta\sigma: H\to ([M, M]\# H)^{[M, M]}$ such that
$\sigma=j_{[M, M]\# H}\crta{\sigma}$.

We have that $j_{[M, M]\# H}$ and, since $H$ is flat, also $j_{[M, M]\# H}\ot H$
are monomorphisms. From Corollary \ref{(Asm.H)A+j com.alg.} we know that $j_{[M, M]\# H}$
is an $H$-comodule algebra morphism. Then by 2.2.3(a), $\crta\sigma$ will be an $H$-comodule algebra
morphism if so is $\sigma$.

That $\sigma$ is $H$-colinear is clear by the coassociativity of $H$. The multiplicativity of $\sigma$ follows from:
$$\scalebox{0.8}[0.8]{
\gbeg{3}{6}
\got{1}{H} \got{3}{H} \gnl
\gwmu{3} \gnl
\gwcm{3} \gnl
\gmp{S} \gvac{1} \gcl{2} \gnl
\gbmp{\teta} \gnl
\gob{1}{[M, M]} \gvac{1} \gob{1}{H}
\gend} \stackrel{bialg.}{=}
\scalebox{0.8}[0.8]{\gbeg{4}{7}
\got{2}{H} \got{2}{H} \gnl
\gcmu \gcmu \gnl
\gcl{1} \gbr \gcl{1} \gnl
\gmu \gmu \gnl
\gvac{1} \hspace{-0,2cm} \gmp{S} \gvac{1} 
   \gcl{2} \gnl
\gvac{1} \gbmp{\teta} \gnl
\gvac{1} \gob{1}{[M, M]} \gvac{1} \gob{1}{H}
\gend} \stackrel{S}{\stackrel{antih.}{=}}
\scalebox{0.8}[0.8]{\gbeg{4}{8}
\got{2}{H} \got{2}{H} \gnl
\gcmu \gcmu \gnl
\gcl{1} \gbr \gcl{1} \gnl
\gmp{S} \gmp{S} \gmu \gnl
\gbr \gvac{1} \hspace{-0,2cm} \gcl{3} \gnl
\gvac{1} \hspace{-0,34cm} \gmu \gnl
\gvac{2} \hspace{-0,22cm} \gbmp{\teta} \gnl
\gvac{2} \gob{1}{[M, M]} \gvac{1} \gob{1}{H}
\gend} \stackrel{comm.}{=}
\scalebox{0.8}[0.8]{\gbeg{4}{7}
\got{2}{H} \got{2}{H} \gnl
\gcmu \gcmu \gnl
\gcl{1} \gbr \gcl{1} \gnl
\gmp{S} \gmp{S} \gmu \gnl
\gmu \gvac{1} \hspace{-0,2cm} \gcl{2} \gnl
\gvac{1} \gbmp{\teta} \gnl
\gvac{1} \gob{1}{[M, M]} \gvac{1} \gob{1}{H}
\gend} \stackrel{\teta}{\stackrel{alg.m.}{=}}
\scalebox{0.8}[0.8]{\gbeg{4}{7}
\got{2}{H} \got{2}{H} \gnl
\gcmu \gcmu \gnl
\gcl{1} \gbr \gcl{1} \gnl
\gmp{S} \gmp{S} \gmu \gnl
\gbmp{\teta} \gbmp{\teta} \gvac{1} \hspace{-0,34cm} \gcl{2} \gnl
\gvac{1} \hspace{-0,36cm} \gmu \gnl
\gvac{1} \gob{2}{[M, M]} \gvac{1} \gob{1}{H}
\gend} =
\scalebox{0.8}[0.8]{\gbeg{6}{7}
\got{2}{H} \got{2}{H} \gnl
\gcmu \gcmu \gnl
\gmp{S} \gcl{1} \gmp{S} \gcl{3} \gnl
\gbmp{\teta} \gcl{1} \gbmp{\teta} \gnl
\gcl{1} \gbr \gnl
\gmu \gmu \gnl
\gob{2}{[M, M]} \gob{2}{H}
\gend}
$$
It is clear that $\sigma$ is also compatible with unit. This finishes the proof that $\Upsilon([M, M])=[H]$
in $\Gal(\C;H)$. \par \medskip

We now prove that $\Upsilon$ does not depend on a representative of the class in $\BM(\C; H)$.
Take two $H$-Azumaya algebras $A$ and $B$ such that $[A]=[B]$ in $\BM(\C; H)$. Then
there are faithfully projective $H$-modules $P$ and $Q$ such that $A\ot [P, P]\iso B\ot [Q, Q]$
as $H$-module algebras. Using the result established in the previous paragraph and
Proposition \ref{Pi multiplicative}, we have:
$$\begin{array}{l}
((A\ot [P, P])\# H)^{A\ot [P, P]}\iso(A\# H)^A\Box_H ([P, P]\# H)^{[P, P]} \iso(A\# H)^A\Box_H H \iso(A\# H)^A \vspace{3pt} \\
((B\ot [Q, Q])\# H)^{B\ot [Q, Q]}\iso(B\# H)^B\Box_H ([Q, Q]\# H)^{[Q, Q]} \iso(B\# H)^B\Box_H H \iso(B\# H)^B.
\end{array}$$

The two expressions on the left-hand sides are isomorphic because of the
assumption $A\ot [P, P]\iso B\ot [Q, Q]$. The isomorphism of the expressions
on the right hand-sides then means that $\Upsilon([A])=\Upsilon([B])$ in $\Gal(\C; H)$.
Thus $\Upsilon$ is well defined. From Proposition \ref{Pi multiplicative} it follows that $\Upsilon$ is a group
morphism.
\qed\end{proof}

\begin{proposition} \label{Gal-nb-BM-inn}
Let $\C$ be a closed braided monoidal category with equalizers and coequalizers. Let $H$
be a flat and commutative Hopf algebra. Suppose that the braiding is $H$-linear. Then, the action on an $H$-Azumaya algebra $A$ is inner if and only if $(A\# H)^A$ is a Galois object with a normal basis. Moreover, the map
$$\Upsilon':\BM_{inn}(\C;H) \to \Gal_{nb}(\C;H),[A] \mapsto [(A \# H)^A]$$ is a group morphism.
\end{proposition}

\begin{proof}
Suppose that $A$ is an $H$-Azumaya algebra with inner action and corresponding morphism
$f:H\to A$. Since $A$ is an Azumaya algebra, the adjunction $(A\ot -, (-)^A)$ is an
equivalence of categories. Hence $A^A\iso(A\ot I)^A\iso I$ and the equalizer $(A^A, j_A)$ from
Definition \ref{M^A} is isomorphic to the equalizer
$$\scalebox{0.9}[0.9]{
\bfig
\putmorphism(-40,0)(1,0)[I`A`\eta_A]{340}1a
\putmorphism(330,0)(1,0)[`[A, A\ot A]`\tilde{\alpha}_{A,A}]{530}1a
\putmorphism(850,25)(1,0)[\phantom{[A, A\ot A]}`\phantom{[A,A\ot A]}`{{[A,\nabla_A]}}]{1000}1a
\putmorphism(440,0)(1,0)[`[A,A\ot A].`]{1400}0a
\putmorphism(850,-25)(1,0)[\phantom{[A, A\ot A]}` \phantom{[A,A\ot \crta A\ot A]}` {{[A,\nabla_A\Phi]}}]{1100}1b
\efig}$$
Having that $H$ is flat, we obtain that
$$\scalebox{0.9}[0.9]{
\bfig
\putmorphism(-100,0)(1,0)[H`A\ot H`\eta_A\ot H]{600}1a
\putmorphism(660,0)(1,0)[`[A, A\ot A]\ot H`\tilde{\alpha}_{A,A}\ot H]{800}1a
\putmorphism(1460,25)(1,0)[\phantom{[A, A\ot A]\ot H}`\phantom{[A,A\ot A]\ot H}`{{[A,\nabla_A]\ot H}}]{1400}1a
\putmorphism(1460,0)(1,0)[`[A,A\ot A]\ot H`]{1400}0a
\putmorphism(1460,-25)(1,0)[\phantom{[A, A\ot A]\ot H}` \phantom{[A,A\ot \crta A\ot A]\ot H}`{{[A,\nabla_A\Phi]
\ot H}}]{1500}1b
\efig}$$
is an equalizer too. Define $\delta:A\ot H\to A\ot H$ as $\delta=(\nabla_A \otimes H)(A \otimes f \otimes H)(A \otimes \Delta_H)$. Let us prove that $\delta j_{A\# H}:(A\# H)^A\to A\ot H$ induces
$\crta\delta:(A\# H)^A\to H$ using the equalizer property of $(H, \eta_A\ot H)$. For
this we will need the following identities:
\begin{eqnarray}\label{j-Phi}
\scalebox{0.9}[0.9]{
\gbeg{3}{4}
\got{1}{A} \got{3}{(A\# H)^A} \gnl
\gcl{1} \glmpb \gnot{\hspace{-0,32cm} j_{A\# H}}\grmptb \gnl
\gmu \gcl{1} \gnl
\gob{2}{A} \gob{1}{H}
\gend} \stackrel{j_{A\# H}}{=} \scalebox{0.9}[0.9]{
\gbeg{5}{9}
\got{1}{A} \gvac{1} \got{3}{(A\# H)^A} \gnl
\gcn{1}{1}{1}{3} \gvac{1} \glmpb \gnot{\hspace{-0,32cm} j_{A\# H}}\grmptb \gnl
\gvac{1} \gbr \gcl{1} \gnl
\gcn{1}{1}{3}{2} \gvac{1} \gbr  \gnl
\gcn{1}{3}{2}{2} \gvac{1} \hspace{-0,34cm} \gcmu \gcn{1}{1}{0}{1} \gnl
\gvac{2} \hspace{-0,14cm} \gcl{1} \gbr \gnl
\gvac{2} \glm \gcl{2} \gnl
\gvac{1} \gwmu{3} \gnl
\gvac{2} \gob{1}{A} \gvac{1}\gob{1}{H}
\gend} \stackrel{Prop. \ref{braid-lin}}{=} \scalebox{0.9}[0.9]{
\gbeg{5}{9}
\got{1}{A} \gvac{1} \got{3}{(A\# H)^A} \gnl
\gcn{1}{1}{1}{3} \gvac{1} \glmpb \gnot{\hspace{-0,32cm} j_{A\# H}}\grmptb \gnl
\gvac{1} \gbr \gcl{1} \gnl
\gcn{1}{1}{3}{2} \gvac{1} \gbr  \gnl
\gcn{1}{3}{2}{2} \gvac{1} \hspace{-0,34cm} \gcmu \gcn{1}{1}{0}{1} \gnl
\gvac{2} \hspace{-0,14cm} \gcl{1} \gibr \gnl
\gvac{2} \glm \gcl{2} \gnl
\gvac{1} \gwmu{3} \gnl
\gvac{2} \gob{1}{A} \gvac{1}\gob{1}{H}
\gend} \stackrel{nat.}{=} \scalebox{0.9}[0.9]{
\gbeg{5}{8}
\got{1}{A} \gvac{1} \got{3}{(A\# H)^A} \gnl
\gcl{1} \gvac{1} \glmpb \gnot{\hspace{-0,32cm} j_{A\# H}}\grmptb \gnl
\gcn{1}{1}{1}{2} \gcn{1}{1}{3}{2} \gvac{1} \hspace{-0,34cm} \gcmu \gnl
\gvac{1} \gbr \gcl{1} \gcl{4} \gnl
\gvac{1} \gcl{2} \gbr \gnl
\gvac{2} \glm \gnl
\gvac{1} \gwmu{3} \gnl
\gvac{2} \gob{1}{A} \gvac{1} \gob{1}{H}
\gend}
\end{eqnarray}
\begin{eqnarray}\label{assoc-corr}
\scalebox{0.9}[0.9]{
\gbeg{6}{6}
\got{1}{A} \got{1}{A} \got{1}{A} \got{1}{A} \got{1}{A} \gnl
\gcl{2} \gcn{1}{1}{1}{2} \gmu \gcl{2} \gnl
\gvac{2} \hspace{-0,34cm} \gmu \gnl
\gvac{1} \hspace{-0,22cm}\gwmu{3} \gvac{1} \gcl{1} \gnl
\gvac{2} \gwmu{4} \gnl
\gvac{3} \gob{2}{A}
\gend} = \scalebox{0.9}[0.9]{
\gbeg{6}{5}
\got{1}{A} \got{1}{A} \got{1}{A} \got{1}{A} \got{1}{A} \gnl
\gmu \gcn{1}{1}{1}{2} \gmu \gnl
\gcn{1}{1}{2}{3} \gvac{2} \hspace{-0,34cm}\gmu \gnl
\gvac{2} \hspace{-0,2cm} \gwmu{3} \gnl
\gvac{3} \gob{1}{A}
\gend}
\end{eqnarray}
We compute:
$$\begin{array}{ll}
\scalebox{0.85}[0.85]{
\gbeg{4}{5}
\got{1}{A} \got{3}{(A\# H)^A} \gnl
\gcl{1} \glmpb \gnot{\hspace{-0,32cm} j_{A\# H}}\grmptb \gnl
\gcl{1} \glmptb \gnot{\hspace{-0,2cm} \delta}\grmptb \gnl
\gmu \gcl{1} \gnl
\gob{2}{A} \gob{1}{H}
\gend} & = \scalebox{0.85}[0.85]{
\gbeg{5}{7}
\got{1}{A} \gvac{1} \got{3}{(A\# H)^A} \gnl
\gcl{4} \gvac{1} \glmpb \gnot{\hspace{-0,32cm} j_{A\# H}}\grmptb \gnl
\gvac{1} \gcn{1}{1}{3}{2} \gvac{1} \hspace{-0,34cm} \gcmu \gnl
\gvac{2} \gcl{1} \gbmp{f} \gcl{3} \gnl
\gvac{2} \gmu \gnl
\gvac{1} \hspace{-0,2cm}\gwmu{3} \gnl
\gvac{2} \gob{1}{A} \gvac{1} \gob{2}{H}
\gend} \stackrel{ass.}{=} \scalebox{0.85}[0.85]{
\gbeg{5}{6}
\got{1}{A} \gvac{1} \got{3}{(A\# H)^A} \gnl
\gcl{1} \gvac{1} \glmpb \gnot{\hspace{-0,32cm} j_{A\# H}}\grmptb \gnl
\gcl{1} \gcn{1}{1}{3}{1} \gvac{1} \hspace{-0,34cm} \gcmu \gnl
\gvac{1} \hspace{-0,22cm} \gmu \gvac{1} \hspace{-0,2cm} \gbmp{f} \gcl{2} \gnl
\gvac{2} \hspace{-0,14cm}\gwmu{3} \gnl
\gvac{3} \gob{1}{A} \gvac{1} \gob{1}{H}
\gend} \stackrel{(\ref{j-Phi})}{=} \scalebox{0.85}[0.85]{
\gbeg{5}{9}
\got{1}{A} \gvac{1} \got{3}{(A\# H)^A} \gnl
\gcl{1} \gvac{1} \glmpb \gnot{\hspace{-0,32cm} j_{A\# H}}\grmptb \gnl
\gcn{1}{1}{1}{2} \gcn{1}{1}{3}{2} \gvac{1} \hspace{-0,34cm} \gcmu \gnl
\gvac{1} \gbr \gcl{1} \gcn{1}{1}{1}{2} \gnl
\gvac{1} \gcl{2} \gbr \gcmu \gnl
\gvac{2} \glm \gbmp{f} \gcl{3} \gnl
\gvac{1} \gwmu{3} \gcl{1} \gnl
\gvac{2} \gwmu{3} \gnl
\gvac{3} \gob{1}{A} \gvac{1} \gob{1}{H}
\gend} \stackrel{\mu_A}{\stackrel{inner}{=}} \scalebox{0.85}[0.85]{
\gbeg{8}{14}
\got{1}{A} \gvac{1} \got{3}{(A\# H)^A} \gnl
\gcl{1} \gvac{1} \glmpb \gnot{\hspace{-0,32cm} j_{A\# H}}\grmptb \gnl
\gcn{1}{1}{1}{3} \gcn{1}{1}{3}{3} \gvac{1} \gcn{1}{1}{1}{3} \gnl
\gvac{1} \gbr \gwcm{3} \gnl
\gcn{1}{1}{3}{1} \gvac{1} \gbr \gvac{1} \gcn{1}{1}{1}{2} \gnl
\gcl{5} \gcn{1}{1}{3}{2} \gvac{1} \gcn{1}{1}{1}{3} \gvac{1} \gcmu \gnl
\gvac{1} \gcmu \gvac{1} \gcl{2} \gbmp{f} \gcl{7} \gnl
\gvac{1} \gbmp{f} \glmpt \gnot{\hspace{-0,2cm} f^{-1}}\grmpb \gvac{1} \gcl{4} \gnl
\gvac{1} \gcl{1} \gvac{1} \gbr \gnl
\gvac{1} \gcn{1}{1}{1}{2} \gvac{1} \gmu \gnl
\gcn{1}{1}{1}{2}\gvac{1} \hspace{-0,2cm} \gwmu{3} \gnl
\gvac{1} \gwmu{3} \gvac{1} \gcn{1}{1}{2}{1} \gnl
\gvac{2} \gwmu{4} \gnl
\gvac{3} \gob{2}{A} \gvac{1} \gob{2}{H}
\gend} \\
 & \stackrel{Prop. \ref{braid-lin}}{=}
\scalebox{0.85}[0.85]{
\gbeg{8}{14}
\got{1}{A} \gvac{1} \got{3}{(A\# H)^A} \gnl
\gcl{1} \gvac{1} \glmpb \gnot{\hspace{-0,32cm} j_{A\# H}}\grmptb \gnl
\gcn{1}{1}{1}{3} \gcn{1}{1}{3}{3} \gvac{1} \gcn{1}{1}{1}{3} \gnl
\gvac{1} \gbr \gwcm{3} \gnl
\gcn{1}{1}{3}{1} \gvac{1} \gibr \gvac{1} \gcn{1}{1}{1}{2} \gnl
\gcl{5} \gcn{1}{1}{3}{2} \gvac{1} \gcn{1}{1}{1}{3} \gvac{1} \gcmu \gnl
\gvac{1} \gcmu \gvac{1} \gcl{2} \gbmp{f} \gcl{7} \gnl
\gvac{1} \gbmp{f} \glmpt \gnot{\hspace{-0,2cm} f^{-1}}\grmpb \gvac{1} \gcl{4} \gnl
\gvac{1} \gcl{1} \gvac{1} \gbr \gnl
\gvac{1} \gcn{1}{1}{1}{2} \gvac{1} \gmu \gnl
\gcn{1}{1}{1}{2}\gvac{1} \hspace{-0,2cm} \gwmu{3} \gnl
\gvac{1} \gwmu{3} \gvac{1} \gcn{1}{1}{2}{1} \gnl
\gvac{2} \gwmu{4} \gnl
\gvac{3} \gob{2}{A} \gvac{1} \gob{2}{H}
\gend} \stackrel{nat.}{\stackrel{coass.}{=}} \scalebox{0.85}[0.85]{
\gbeg{7}{11}
\got{1}{A} \gvac{1} \got{3}{(A\# H)^A} \gnl
\gcl{2} \gvac{1} \glmpb \gnot{\hspace{-0,32cm} j_{A\# H}}\grmptb \gnl
\gvac{1} \gcn{1}{1}{3}{1} \gcn{1}{1}{3}{3} \gnl
\gbr \gwcm{3} \gnl
\gcl{4} \gibr \gwcm{3} \gnl
\gvac{1} \gbmp{f} \gcl{1} \glmptb \gnot{\hspace{-0,2cm} f^{-1}}\grmp \gcl{1} \gnl
\gvac{1} \gcn{1}{1}{1}{2} \gmu \gwcm{3} \gnl
\gvac{2} \hspace{-0,34cm} \gmu \gvac{1} \hspace{-0,22cm} \gbmp{f} \gvac{1} \gcl{3} \gnl
\gvac{1} \hspace{-0,14cm}\gwmu{3} \gvac{1} \gcl{1} \gnl
\gvac{2} \gwmu{4} \gnl
\gvac{3} \gob{2}{A} \gvac{2} \gob{1}{H}
\gend} \stackrel{Prop. \ref{braid-lin}}{\stackrel{nat.}{=}} \scalebox{0.85}[0.85]{
\gbeg{7}{11}
\got{1}{A} \gvac{1} \got{3}{(A\# H)^A} \gnl
\gcl{2} \gvac{1} \glmpb \gnot{\hspace{-0,32cm} j_{A\# H}}\grmptb \gnl
\gvac{1} \gcn{1}{1}{3}{1} \gcn{1}{1}{3}{3} \gnl
\gbr \gwcm{3} \gnl
\gcl{4} \gcl{1} \gbmp{f} \gwcm{3} \gnl
\gvac{1} \gbr \glmptb \gnot{\hspace{-0,2cm} f^{-1}}\grmp \gcl{1} \gnl
\gvac{1} \gcn{1}{1}{1}{2} \gmu \gwcm{3} \gnl
\gvac{2} \hspace{-0,34cm} \gmu \gvac{1} \hspace{-0,22cm} \gbmp{f} \gvac{1} \gcl{3} \gnl
\gvac{1} \hspace{-0,14cm}\gwmu{3} \gvac{1} \gcl{1} \gnl
\gvac{2} \gwmu{4} \gnl
\gvac{3} \gob{2}{A} \gvac{2} \gob{1}{H}
\gend} \vspace{3mm} \\
 & \stackrel{(\ref{assoc-corr})}{\stackrel{coass.}{=}}
\scalebox{0.85}[0.85]{
\gbeg{7}{11}
\got{1}{A} \gvac{1} \got{3}{(A\# H)^A} \gnl
\gcl{2} \gvac{1} \glmpb \gnot{\hspace{-0,32cm} j_{A\# H}}\grmptb \gnl
\gvac{1} \gcn{1}{1}{3}{1} \gcn{1}{1}{3}{4} \gnl
\gbr \gwcm{4} \gnl
\gcl{3} \gcl{1} \gbmp{f} \gvac{1} \gwcm{3} \gnl
\gvac{1} \gbr \gwcm{3} \gcl{5} \gnl
\gvac{1} \gcl{1} \gcl{2} \glmptb \gnot{\hspace{-0,2cm} f^{-1}}\grmp \gbmp{f} \gnl
\gmu \gvac{1} \gwmu{3} \gnl
\gcn{1}{1}{2}{3} \gvac{1} \gwmu{3} \gnl
\gvac{1} \hspace{-0,14cm} \gwmu{3} \gnl
\gvac{2} \gob{1}{A} \gvac{3} \gob{1}{H}
\gend} \stackrel{f^{-1}*f}{\stackrel{nat.}{=}}\scalebox{0.85}[0.85]{
\gbeg{5}{8}
\got{2}{A} \gvac{1} \got{2}{(A\# H)^A} \gnl
\gcn{1}{3}{2}{2}\gvac{1} \glmpb \gnot{\hspace{-0,32cm} j_{A\# H}}\grmptb \gnl
\gvac{1} \gcn{1}{1}{3}{2} \gvac{1} \hspace{-0,34cm} \gcmu \gnl
\gvac{2} \gcl{1} \gbmp{f} \gcl{4} \gnl
\gvac{1} \gcn{1}{1}{1}{2} \gmu \gnl
\gvac{2} \hspace{-0,2cm}\gbr \gnl
\gvac{2} \gmu \gnl
\gvac{2} \gob{2}{A} \gob{2}{H}
\gend} = \scalebox{0.85}[0.85]{
\gbeg{4}{6}
\got{1}{A} \got{3}{(A\# H)^A} \gnl
\gcl{1} \glmpb \gnot{\hspace{-0,32cm} j_{A\# H}}\grmptb \gnl
\gcl{1} \glmptb \gnot{\hspace{-0,2cm} \delta}\grmptb \gnl
\gbr \gcl{1} \gnl
\gmu \gcl{1} \gnl
\gob{2}{A} \gob{1}{H}
\gend}
\end{array}
$$
This proves that $\delta j_{A\# H}:(A\# H)^A\to A\ot H$ induces $\crta\delta:(A\# H)^A\to H$ such that
$\delta j_{A\# H}=(\eta_A\ot H)\crta\delta$. Clearly, $\delta$ and $\eta_A\ot H$ are right $H$-colinear and by Proposition \ref{braid-lin} and Corollary \ref{(Asm.H)A+j com.alg.} we know that $j_{A\# H}$ is such as well. Then by 2.2.3(i), $\crta\delta$ is right $H$-colinear too. To find the inverse of $\crta\delta$
we prove that $\zeta:=(f^{-1}\ot H)\Delta_H:H\to A\ot H$ factors through
$(A\# H)^A$. We need the following equalities: \vspace{-3mm}

\begin{center}
\begin{tabular}{p{6.8cm}p{1cm}p{6.8cm}}
\begin{eqnarray} \label{coccan8-1}
\scalebox{0.9}[0.9]{
\gbeg{5}{5}
\gvac{3} \got{1}{H} \gnl
\gvac{2} \gwcm{3} \gnl
\gvac{1} \gwcm{3} \gcl{2} \gnl
\gwcm{3} \gcl{1} \gnl
\gob{1}{H} \gvac{1} \gob{1}{H} \gob{1}{H} \gob{1}{H}
\gend} = \scalebox{0.9}[0.9]{
\gbeg{5}{5}
\gvac{1} \got{1}{H} \gnl
\gwcm{3} \gnl
\gcl{1} \gwcm{3}  \gnl
\gcn{1}{1}{1}{0} \hspace{-0,2cm} \gcmu \gvac{1} \hspace{-0,2cm} \gcl{1} \gnl
\gob{2}{H} \gob{1}{\hspace{-0,2cm}H} \gob{1}{\hspace{-0,2cm}H} \gob{1}{H}
\gend}
\end{eqnarray} & &
\begin{eqnarray}\label{assoc-corr1}
\scalebox{0.9}[0.9]{
\gbeg{4}{5}
\got{1}{A} \got{1}{A} \got{1}{A} \got{1}{A} \gnl
\gcl{2} \gcn{1}{1}{1}{2} \gmu \gnl
\gvac{2} \hspace{-0,2cm} \gmu \gnl
\gvac{1} \hspace{-0,22cm}\gwmu{3} \gnl
\gvac{2} \gob{1}{A}
\gend} = \scalebox{0.9}[0.9]{
\gbeg{4}{4}
\got{1}{A} \got{1}{A} \gvac{1} \got{1}{A} \got{1}{A} \gnl
\gmu \gvac{1} \gmu \gnl
\gvac{1} \hspace{-0,22cm} \gwmu{4} \gnl
\gvac{2} \gob{2}{A}
\gend}
\end{eqnarray}
\end{tabular}
\end{center}
We now have:
$$\begin{array}{ll}
\scalebox{0.85}[0.85]{
\gbeg{4}{5}
\got{1}{A} \got{3}{H} \gnl
\gcl{1} \gwcm{3} \gnl
\gcl{1} \glmptb \gnot{\hspace{-0,2cm} f^{-1}}\grmp \gcl{2} \gnl
\gmu \gnl
\gob{2}{A} \gvac{1} \gob{1}{H}
\gend} & \stackrel{unit}{\stackrel{counit}{\stackrel{nat.}{=}}}
\scalebox{0.85}[0.85]{
\gbeg{5}{10}
\gvac{1} \got{1}{A} \got{3}{H} \gnl
\gvac{1} \gcl{1} \gwcm{3} \gnl
\gvac{1} \gibr \gvac{1} \gcl{7} \gnl
\gcn{1}{1}{3}{2} \gvac{1} \gcl{1} \gnl
\gcmu \gcn{1}{1}{1}{3} \gnl
\gcu{1} \glmpt \gnot{\hspace{-0,2cm} f^{-1}}\grmpb \gcl{1} \gnl
\gu{1} \gvac{1} \gbr \gnl
\gcn{1}{1}{1}{2} \gvac{1} \gmu \gnl
\gvac{1} \hspace{-0,2cm} \gwmu{3} \gnl
\gvac{2} \gob{1}{A} \gvac{1} \gob{2}{H}
\gend} \stackrel{f*f^{-1}}{\stackrel{Prop. \ref{braid-lin}}{=}}
\scalebox{0.85}[0.85]{
\gbeg{7}{11}
\gvac{3} \got{1}{A} \got{3}{H} \gnl
\gvac{3} \gcl{1} \gwcm{3} \gnl
\gvac{3} \gbr \gvac{1} \gcl{8} \gnl
\gvac{2} \gcn{1}{1}{3}{1} \gvac{1} \gcl{1} \gnl
\gvac{1} \gwcm{3} \gcn{1}{1}{1}{3} \gnl
\gwcm{3} \glmpt \gnot{\hspace{-0,2cm} f^{-1}}\grmpb \gcl{1} \gnl
\glmptb \gnot{\hspace{-0,2cm} f^{-1}}\grmp \gbmp{f} \gvac{1} \gbr \gnl
\gwmu{3} \gvac{1} \gmu \gnl
\gvac{1} \gcn{1}{1}{1}{2} \gvac{2} \gcn{1}{1}{2}{2} \gnl
\gvac{2} \hspace{-0,34cm} \gwmu{4} \gnl
\gvac{3} \gob{2}{A} \gvac{1} \gob{2}{H}
\gend} \stackrel{nat.}{\stackrel{(\ref{coccan8-1})}{\stackrel{(\ref{assoc-corr1})}{=}}}
\scalebox{0.85}[0.85]{
\gbeg{8}{11}
\gvac{2} \got{1}{A} \got{3}{H} \gnl
\gvac{2} \gcl{1} \gwcm{3} \gnl
\gvac{2} \gbr \gwcm{3} \gnl
\gvac{1} \gcn{1}{1}{3}{1} \gvac{1} \gbr \gvac{1} \gcl{7} \gnl
\gvac{1} \gcl{3} \gwcm{3} \gcn{1}{1}{-1}{1} \gnl
\gvac{2} \gbmp{f} \glmp \gnot{\hspace{-0,2cm} f^{-1}}\grmptb \gcl{1} \gnl
\gvac{2} \gcl{1} \gvac{1} \gbr \gnl
\glmp \gnot{\hspace{-0,2cm} f^{-1}}\grmptb \gcn{1}{1}{1}{2} \gvac{1} \gmu \gnl
\gvac{1} \gcn{1}{1}{1}{2}\gvac{1} \hspace{-0,2cm} \gwmu{3} \gnl
\gvac{2} \gwmu{3} \gnl
\gvac{3} \gob{1}{A} \gvac{2} \gob{2}{H}
\gend} \vspace{3mm} \\
 & \stackrel{nat.}{\stackrel{\mu_A}{\stackrel{inner}{=}}}
\scalebox{0.85}[0.85]{
\gbeg{6}{8}
\got{1}{A} \gvac{1} \got{3}{H} \gnl
\gcl{2} \gvac{1} \gwcm{3} \gnl
\gvac{1} \glmpb \gnot{\hspace{-0,2cm} f^{-1}}\grmpt \gwcm{3} \gnl
\gbr \gcn{1}{1}{3}{1} \gvac{2} \gcl{4} \gnl
\gcl{2} \gbr \gnl
\gvac{1} \glm \gnl
\gwmu{3} \gnl
\gvac{1} \gob{1}{A} \gvac{3} \gob{1}{H}
\gend} \stackrel{nat.}{=} \scalebox{0.85}[0.85]{
\gbeg{5}{10}
\gvac{1} \got{1}{A} \gvac{1} \got{3}{(A\# H)^A} \gnl
\gvac{1} \gcl{2} \gwcm{3} \gnl
\gvac{2} \glmptb \gnot{\hspace{-0,2cm} f^{-1}}\grmp \gcl{1} \gnl
\gvac{1} \gbr \gcn{1}{1}{3}{1} \gnl
\gcn{1}{1}{3}{2} \gvac{1} \gbr  \gnl
\gcn{1}{3}{2}{2} \gvac{1} \hspace{-0,34cm} \gcmu \gcn{1}{1}{0}{1} \gnl
\gvac{2} \hspace{-0,14cm} \gcl{1} \gibr \gnl
\gvac{2} \glm \gcl{2} \gnl
\gvac{1} \gwmu{3} \gnl
\gvac{2} \gob{1}{A} \gvac{1}\gob{1}{H}
\gend} \stackrel{Prop. \ref{braid-lin}}{=}
\scalebox{0.85}[0.85]{
\gbeg{5}{10}
\gvac{1} \got{1}{A} \gvac{1} \got{3}{(A\# H)^A} \gnl
\gvac{1} \gcl{2} \gwcm{3} \gnl
\gvac{2} \glmptb \gnot{\hspace{-0,2cm} f^{-1}}\grmp \gcl{1} \gnl
\gvac{1} \gbr \gcn{1}{1}{3}{1} \gnl
\gcn{1}{1}{3}{2} \gvac{1} \gbr  \gnl
\gcn{1}{3}{2}{2} \gvac{1} \hspace{-0,34cm} \gcmu \gcn{1}{1}{0}{1} \gnl
\gvac{2} \hspace{-0,14cm} \gcl{1} \gbr \gnl
\gvac{2} \glm \gcl{2} \gnl
\gvac{1} \gwmu{3} \gnl
\gvac{2} \gob{1}{A} \gvac{1}\gob{1}{H}
\gend}
\end{array}$$
So $\zeta:H\to A\ot H$ induces $\crta\zeta:H\to (A\# H)^A$ such that $j_{A\# H}\crta\zeta=\zeta.$
We now prove that $\crta\zeta$ and $\crta\delta$ are inverses of each other. It is easy to see that
$(\nabla_A\ot H)(A\ot\zeta)\delta=id_{A\# H}.$ Then
$$\begin{array}{ll}
j_{A\# H} & =(\nabla_A\ot H)(A\ot\zeta)\delta j_{A\# H} \vspace{2pt} \\
          & =(\nabla_A\ot H)(A\ot\zeta)(\eta_A\ot H)\crta\delta \vspace{2pt} \\
          & =(\nabla_A(\eta_A \ot A)\ot H)\zeta\crta\delta \vspace{2pt} \\
          & =\zeta\crta\delta \vspace{2pt} \\
          & =j_{A\# H}\crta\zeta\hspace{2pt}\crta\delta.
\end{array}$$
Since $j_{A\# H}$ is a monomorphism, we get $\crta\zeta\hspace{2pt}\crta\delta=id_{(A\# H)^A}.$ Similarly, one may easily check that $\delta\zeta=\eta_A \otimes H$. Then, $\eta_A\ot H =\delta\zeta=\delta j_{A\#H}\crta\zeta=(\eta_A\ot H)\crta\delta\hspace{2pt}\crta\zeta.$ Having that $\eta_A\ot H$ is a monomorphism, we obtain $\crta\delta\hspace{2pt}\crta\zeta=id_H$. This proves that $(A\# H)^A$ is a Galois object with normal basis.
\par\medskip

Conversely, suppose that $\zeta:H\to (A\# H)^A$ is a right $H$-comodule isomorphism. The unit $\beta=(\nabla_A\ot H)(A\ot j_{A\#H}):A\ot (A\# H)^A\to A\ot H$ of the adjunction $(A\ot -, (-)^A)$ is an isomorphism because
$A$ is an Azumaya algebra, Remark \ref{Az-counit}. Let us prove that
$$v:=\scalebox{0.9}[0.9]{
\gbeg{3}{6}
\got{1}{} \got{1}{H} \gnl
\gu{1} \gcl{1} \gnl
\glmptb \gnot{\hspace{-0,22cm} \beta^{-1}}\grmptb \gnl
\gcl{2} \glmptb \gnot{\hspace{-0,22cm} \zeta^{-1}}\grmp \gnl
\gvac{1} \gcu{1} \gnl
\gob{1}{A}
\gend} \quad\textnormal{and}\quad
u:=\scalebox{0.9}[0.9]{
\gbeg{3}{5}
\got{1}{H} \gnl
\glmptb \gnot{\hspace{-0,26cm} \zeta}\grmpb \gnl
\glmptb \gnot{\hspace{-0,28cm} j_{A\#H}}\grmptb \gnl
\gcl{1} \gcu{1} \gnl
\gob{1}{A}
\gend}
$$
are convolution inverses of each other. Being $j_{A\#H}$ right $H$-colinear, then so
are clearly $\beta$ and $\beta^{-1}$ Moreover, $\zeta^{-1}$ is such as well, so
$$v*u=\scalebox{0.9}[0.9]{
\gbeg{5}{7}
\got{1}{} \got{3}{H} \gnl
\gu{1} \gwcm{3} \gnl
\glmptb \gnot{\hspace{-0,2cm} \beta^{-1}}\grmptb \gvac{1} \gbmp{\zeta} \gnl
\gcl{1} \glmptb \gnot{\hspace{-0,2cm} \zeta^{-1}}\grmp \glmptb \gnot{\hspace{-0,28cm} j_{A\#H}}\grmp \gnl
\gcl{1} \gcu{1} \gvac{1} \gcl{1} \gcu{1} \gnl
\gwmu{4} \gnl
\gvac{1} \gob{2}{A}
\gend} = \scalebox{0.9}[0.9]{
\gbeg{4}{9}
\got{1}{} \got{1}{H} \gnl
\gu{1} \gcl{1} \gnl
\glmptb \gnot{\hspace{-0,2cm} \beta^{-1}}\grmptb \gnl
\gcl{1} \glmpt \gnot{\hspace{-0,2cm} \zeta^{-1}}\grmpb \gnl
\gcl{1} \gwcm{3} \gnl
\gcl{2} \gcu{1} \gvac{1} \gbmp{\zeta} \gnl
\gvac{2} \glmpb \gnot{\hspace{-0,28cm} j_{A\#H}}\grmptb \gnl
\gwmu{3} \gcu{1} \gnl
\gvac{1} \gob{1}{A}
\gend} = \scalebox{0.9}[0.9]{
\gbeg{4}{6}
\got{1}{} \got{1}{H} \gnl
\gu{1} \gcl{1} \gnl
\glmptb \gnot{\hspace{-0,2cm} \beta^{-1}}\grmptb \gnl
\gcl{1} \glmpt \gnot{\hspace{-0,28cm} j_{A\#H}}\grmpb \gnl
\gmu \gcu{1} \gnl
\gob{2}{A}
\gend} \stackrel{\beta}{=} \scalebox{0.9}[0.9]{
\gbeg{3}{6}
\got{1}{H} \gnl
\gcl{1} \gnl
\gcu{1} \gnl
\gu{1} \gnl
\gcl{1} \gnl
\gob{1}{A}
\gend}
$$
Since $\beta$ is left $A$-linear, then so is $\beta^{-1}$, hence
$$\begin{array}{l}
u*v=\scalebox{0.9}[0.9]{
\gbeg{5}{7}
\got{4}{H} \gnl
\gwcm{4} \gnl
\gbmp{\zeta} \gvac{1} \gu{1} \gcl{1} \gnl
\glmptb \gnot{\hspace{-0,28cm} j_{A\#H}}\grmp \glmptb \gnot{\hspace{-0,2cm} \beta^{-1}}\grmptb \gnl
\gcl{1} \gcu{1} \gcl{1} \glmptb \gnot{\hspace{-0,2cm} \zeta^{-1}}\grmp \gnl
\gwmu{3} \gcu{1} \gnl
\gvac{1} \gob{1}{A}
\gend} = \scalebox{0.9}[0.9]{
\gbeg{5}{10}
\got{4}{H} \gnl
\gwcm{4} \gnl
\gbmp{\zeta} \gvac{2} \gcl{3} \gnl
\glmptb \gnot{\hspace{-0,28cm} j_{A\#H}}\grmp \gu{1} \gnl
\gcl{1} \gcu{1} \gcl{1} \gnl
\gwmu{3} \gcn{1}{1}{1}{-1} \gnl
\gvac{1} \glmptb \gnot{\hspace{-0,2cm} \beta^{-1}}\grmptb \gnl
\gvac{1} \gcl{1} \glmptb \gnot{\hspace{-0,2cm} \zeta^{-1}}\grmp \gnl
\gvac{1} \gcl{1} \gcu{1} \gnl
\gvac{1} \gob{1}{A}
\gend} =
\scalebox{0.9}[0.9]{
\gbeg{5}{10}
\got{5}{H} \gnl
\gvac{1} \gwcm{3} \gnl
\gvac{1} \gbmp{\zeta} \gvac{1} \gcl{3} \gnl
\gu{1} \glmptb \gnot{\hspace{-0,28cm} j_{A\#H}}\grmpb \gnl
\gmu \gcu{1} \gnl
\gvac{1} \gcn{1}{1}{0}{2} \gcn{1}{1}{3}{2} \gnl
\gvac{2} \hspace{-0,34cm}\glmptb \gnot{\hspace{-0,2cm} \beta^{-1}}\grmptb \gnl
\gvac{2} \gcl{1} \glmptb \gnot{\hspace{-0,2cm} \zeta^{-1}}\grmp \gnl
\gvac{2} \gcl{1} \gcu{1} \gnl
\gvac{2} \gob{1}{A}
\gend} \stackrel{nat.}{\stackrel{\beta}{=}}
\scalebox{0.9}[0.9]{
\gbeg{4}{9}
\got{1}{} \got{2}{H} \gnl
\gu{1} \gcmu \gnl
\gcl{1} \gbmp{\zeta} \gcl{2} \gnl
\glmptb \gnot{\hspace{-0,28cm} \beta}\grmp \gnl
\gcn{1}{1}{1}{2} \gcu{1} \gcn{1}{1}{1}{0} \gnl
\gvac{1} \hspace{-0,34cm}\glmptb \gnot{\hspace{-0,2cm} \beta^{-1}}\grmptb \gnl
\gvac{1} \gcl{1} \glmptb \gnot{\hspace{-0,2cm} \zeta^{-1}}\grmp \gnl
\gvac{1} \gcl{1} \gcu{1} \gnl
\gvac{1} \gob{1}{A}
\gend}= \scalebox{0.9}[0.9]{
\gbeg{4}{10}
\got{1}{} \got{1}{H} \gnl
\gu{1} \gbmp{\zeta} \gnl
\glmptb \gnot{\hspace{-0,26cm} \beta}\grmptb \gnl
\gcl{1} \gcn{1}{1}{1}{2} \gnl
\gcl{1} \gcmu \gnl
\gcn{1}{1}{1}{2} \gcu{1} \gcn{1}{1}{1}{0} \gnl
\gvac{1} \hspace{-0,34cm}\glmptb \gnot{\hspace{-0,2cm} \beta^{-1}}\grmptb \gnl
\gvac{1} \gcl{1} \glmptb \gnot{\hspace{-0,2cm} \zeta^{-1}}\grmp \gnl
\gvac{1} \gcl{1} \gcu{1} \gnl
\gvac{1} \gob{1}{A}
\gend}= \scalebox{0.9}[0.9]{
\gbeg{3}{6}
\got{1}{H} \gnl
\gcl{1} \gnl
\gcu{1} \gnl
\gu{1} \gnl
\gcl{1} \gnl
\gob{1}{A}
\gend}
\end{array}$$
where in the fourth diagram we applied that $\zeta$ and $\beta$ are
right $H$-colinear. We have:
$$\begin{array}{ll}
\scalebox{0.9}[0.9]{
\gbeg{2}{4}
\got{1}{H} \got{1}{A} \gnl
\glm \gnl
\gvac{1} \gcl{1} \gnl
\gvac{1} \gob{1}{A}
\gend} & = \scalebox{0.9}[0.9]{
\gbeg{5}{7}
\gvac{2} \got{1}{H} \gvac{1} \got{1}{A} \gnl
\gvac{1} \gwcm{3} \gcl{1} \gnl
\gwcm{3} \glm \gnl
\gbmp{v} \gvac{1} \gbmp{u} \gvac{1} \gcl{2} \gnl
\gwmu{3} \gnl
\gvac{1} \gwmu{4} \gnl
\gvac{2} \gob{2}{A}
\gend} \stackrel{coass.}{\stackrel{u}{=}}
\scalebox{0.9}[0.9]{
\gbeg{5}{8}
\gvac{1} \got{1}{H} \gvac{2} \got{1}{A} \gnl
\gwcm{3} \gvac{1} \gcl{2} \gnl
\gbmp{v} \gwcm{3} \gnl
\gcl{2} \gbmp{\zeta} \gvac{1} \glm \gnl
\gvac{1} \glmptb \gnot{\hspace{-0,28cm} j_{A\#H}}\grmp \gvac{1} \gcl{1} \gnl
\gmu \gcu{1} \gcn{1}{1}{3}{2} \gnl
\gvac{1} \hspace{-0,36cm}\gwmu{4} \gnl
\gvac{2} \gob{2}{A}
\gend} \stackrel{\zeta, j}{\stackrel{H\x colin.}{=}}
\scalebox{0.9}[0.9]{
\gbeg{6}{8}
\got{3}{H} \gvac{2} \got{1}{A} \gnl
\gwcm{3} \gvac{2} \gcl{3} \gnl
\gbmp{v} \gvac{1} \gbmp{\zeta} \gnl
\gcl{1} \gvac{1} \glmptb \gnot{\hspace{-0,28cm} j_{A\#H}}\grmp \gnl
\gwmu{3} \hspace{-0,24cm} \gcmu \hspace{-0,18cm} \gcn{1}{1}{3}{2} \gnl
\gvac{2} \gcn{1}{1}{1}{2} \gvac{1} \hspace{-0,38cm} \gcu{1} \glm \gnl
\gvac{4} \hspace{-0,4cm}\gwmu{4} \gnl
\gvac{5} \gob{2}{A}
\gend} \stackrel{nat.}{\stackrel{ass.}{=}}
\scalebox{0.9}[0.9]{
\gbeg{5}{12}
\gvac{1} \got{1}{H} \gvac{1} \got{1}{A} \gnl
\gwcm{3} \gcl{1} \gnl
\gbmp{v} \gvac{1} \gibr \gnl
\gcl{7} \gcn{1}{1}{3}{1} \gvac{1} \gbmp{\zeta} \gnl
\gvac{1} \gcl{1} \glmpb \gnot{\hspace{-0,28cm} j_{A\#H}}\grmptb \gnl
\gvac{1} \gbr \gcn{1}{1}{1}{2} \gnl
\gvac{1} \gcl{3} \gcl{1} \gcmu \gnl
\gvac{2} \gbr \gcl{2} \gnl
\gvac{2} \glm \gnl
\gvac{1} \gwmu{3} \gcu{1} \gnl
\gwmu{3} \gnl
\gvac{1} \gob{1}{A}
\gend} \vspace{-5mm} \\
 & \stackrel{Prop. \ref{braid-lin}}{\stackrel{(\ref{j-Phi})}{=}}
\scalebox{0.9}[0.9]{
\gbeg{5}{8}
\gvac{1} \got{1}{H} \gvac{1} \got{1}{A} \gnl
\gwcm{3} \gcl{1} \gnl
\gbmp{v} \gvac{1} \gbr \gnl
\gcl{2} \gcn{1}{1}{3}{1} \gvac{1} \gbmp{\zeta} \gnl
\gvac{1} \gcl{1} \glmpb \gnot{\hspace{-0,28cm} j_{A\#H}}\grmptb \gnl
\gcn{1}{1}{1}{2} \gmu \gcu{1} \gnl
\gvac{1} \hspace{-0,34cm} \gmu \gnl
\gvac{1} \gob{2}{A}
\gend} \stackrel{u}{=}
\scalebox{0.9}[0.9]{
\gbeg{5}{7}
\gvac{1} \got{1}{H} \gvac{1} \got{1}{A} \gnl
\gwcm{3} \gcl{1} \gnl
\gbmp{v} \gvac{1} \gbr \gnl
\gcl{1} \gvac{1} \gcl{1} \gbmp{u} \gnl
\gcn{1}{1}{1}{2} \gvac{1} \gmu \gnl
\gvac{1} \hspace{-0,22cm} \gwmu{3} \gnl
\gvac{2} \gob{1}{A}
\gend}
\end{array}$$
This means that the $H$-action on $A$ is inner. \par \smallskip

We finally show that $\Upsilon':\BM_{inn}(\C;H) \to \Gal_{nb}(\C;H),[A] \mapsto [(A \# H)^A]$ is a group morphism. Observe that the condition $\Phi_{T,X}=\Phi^{-1}_{X,T}$ required in Proposition \ref{Pi multiplicative} for any $H$-Galois object $T$ and $X \in \C$ was only used there for Galois objects of the form $(A \# H)^A$ for an $H$-Azumaya algebra $A$. This condition is fulfilled by any Galois object $T$ with normal basis. This is because $T \cong H$ as comodules, in particular as objects in $\C$, and since the braiding is $H$-linear, Proposition \ref{braid-lin} gives that $\Phi_{H,X}=\Phi^{-1}_{X,H}$ for any $X \in \C$. We just proved that if $A$ is an $H$-Azumaya algebra with inner action, then $(A \# H)^A$ is an $H$-Galois object with normal basis. As in the proof of Proposition \ref{gmorphism}, one proves that $\Upsilon':\BM_{inn}(\C;H) \to \Gal_{nb}(\C;H)$ is a group morphism.
\qed\end{proof}

\subsection{From an $H$-Galois object to an Azumaya algebra}

In this subsection we will prove that every $H$-Galois object gives rise to an Azumaya
algebra. This will be needed to show in the next subsection that the map $\Upsilon:
\BM(\C; H) \to \Gal(\C; H)$ defined before is surjective. We will start by recalling from \cite[Section 2]{Tak2}
some facts about duality of Hopf algebras. The following proposition collects \cite[2.5, 2.14 and 2.16]{Tak2}:

\begin{proposition} \label{alg<->coalg}
Let $\C$ be a closed braided monoidal category.
\begin{enumerate}
\item[(i)] If $H$ is a coalgebra in $\C$, then $[H, I]$ is an algebra.
\item[(ii)] If $H$ is a finite algebra in $\C$, then $H^*=[H, I]$ is a coalgebra.
\item[(iii)] If $H$ is a finite Hopf algebra in $\C$, then so is $H^*=[H, I]$.
\end{enumerate}
\end{proposition}

We recall here the structure morphisms since we will need them later. The multiplication and unit for $H^*$ are defined via the universal property of $(H^*, \ev)$ by
\begin{center}
\begin{tabular}{p{6.8cm}p{1cm}p{6.8cm}}
\begin{eqnarray} \label{multH*}
\scalebox{0.9}[0.9]{
\gbeg{4}{5}
\got{1}{H^*} \got{3}{H^*} \got{1}{\hspace{-0,66cm} H} \gnl
\gwmu{3} \gcl{1} \gnl
\gvac{1} \gcl{1} \gcn{2}{1}{3}{1}  \gnl
\gvac{1} \gev \gvac{1} \gnl
\gvac{2} \gob{2}{}
\gend} = \scalebox{0.9}[0.9]{
\gbeg{4}{5}
\got{1}{H^*} \got{1}{H^*} \got{3}{\hspace{-0,38cm}H} \gnl
\gcl{2} \gcl{1} \gcmu \gnl
\gvac{1} \gbr \gcl{1} \gnl
\gev \gev \gnl
\gvac{2} \gob{1}{}
\gend}
\end{eqnarray} & &
\begin{eqnarray} \label{unitH*}
\scalebox{0.9}[0.9]{
\gbeg{4}{3}
\got{1}{} \got{1}{H} \gnl
\gu{1} \gcl{1} \gnl
\gev \gnl
\gvac{1} \gob{2}{}
\gend} = \scalebox{0.9}[0.9]{
\gbeg{1}{3}
\got{1}{H} \gnl
\gcl{1} \gnl
\gcu{1} \gnl
\gob{2}{}
\gend}
\end{eqnarray}
\end{tabular}
\end{center}

The finiteness condition in (ii), and hence also in (iii), is needed in order to be
able to consider $H^*\ot H^*\iso(H\ot H)^*$. Then one may apply 2.1.9 to define a comultiplication on $H^*$ using the universal property of $([H\ot H, I], \ev)$ and applying the isomorphism induced by (\ref{duals isom delta}).
The comultiplication and counit are given by the following diagrams:
\begin{center}
\begin{tabular}{p{6.8cm}p{1cm}p{6.8cm}}
\begin{eqnarray} \label{codiagH*}
\scalebox{0.9}[0.9]{
\gbeg{4}{5}
\got{2}{H^*} \got{1}{H}\got{1}{H} \gnl
\gcmu \gcl{1} \gcl{2} \gnl
\gcl{1} \gbr \gnl
\gev \gev \gnl
\gvac{2} \gob{1}{}
\gend} = \scalebox{0.9}[0.9]{
\gbeg{4}{5}
\got{1}{H^*} \got{1}{H} \got{3}{H} \gnl
\gcl{1} \gwmu{3} \gnl
\gcn{2}{1}{1}{3} \gcl{1} \gnl
\gvac{1} \gev \gvac{1} \gnl
\gvac{2} \gob{2}{}
\gend}
\end{eqnarray} & &
\begin{eqnarray} \label{counitH*}
\scalebox{0.9}[0.9]{
\gbeg{1}{3}
\got{1}{H^*} \gnl
\gcl{1} \gnl
\gcu{1} \gnl
\gob{1}{}
\gend} = \scalebox{0.9}[0.9]{
\gbeg{4}{3}
\got{1}{H^*} \got{1}{} \gnl
\gcl{1} \gu{1} \gnl
\gev \gnl
\gvac{1} \gob{2}{}
\gend}
\end{eqnarray}
\end{tabular}
\end{center}
Via the universal property of $(H^*, \ev)$ we may define the antipode $S^*$
for $H^*$ as:
\begin{eqnarray} \label{antipodeH*}
\scalebox{0.9}[0.9]{
\gbeg{2}{3}
\got{1}{H^*} \got{1}{H} \gnl
\gnot{\hspace{0,12cm}S^*} \gmp \gnl \gcl{1} \gnl
\gev \gnl
\gob{1}{}
\gend} = \scalebox{0.9}[0.9]{
\gbeg{1}{3}
\got{1}{H^*} \got{1}{H} \gnl
\gcl{1} \gnot{S} \gmp \gnl \gnl
\gev \gnl
\gvac{1} \gob{2}{}
\gend}
\end{eqnarray}

It is easy to see that a finite algebra $A$ in $\C$ is commutative if and only if $A^*$ is a
cocommutative coalgebra. The proof of the following proposition is not difficult. The first statement is proved in
\cite[Proposition 2.7]{Tak2}.

\begin{proposition} \label{Hmod-H*comod}
Let $H\in \C$ be a finite coalgebra. If $M\in\C^H$, then $M\in {}_{H^*}\C$ with the structure morphism given
in (\ref{H*mod-S}). If $N\in {}_{H^*}\C$, then $N\in\C^H$ with the structure morphism given in (\ref{Hcomod}). The categories $\C^H$ and ${}_{H^*}\C$ are isomorphic via the previous assignments.
\begin{center}
\begin{tabular}{p{6.8cm}p{1cm}p{6.8cm}}
\begin{eqnarray} \label{H*mod-S}
\scalebox{0.9}[0.9]{
\gbeg{2}{3}
\got{1}{H^*} \got{1}{\hspace{0,1cm}M} \gnl
\glm \gnl
\gvac{1} \gob{1}{M}
\gend} = \scalebox{0.9}[0.9]{
\gbeg{2}{5}
\got{1}{H^*} \got{1}{M} \gnl
\gcl{1} \grcm \gnl
\gbr \gcl{1} \gnl
\gcl{1} \gev \gnl
\gob{1}{M} \gvac{1} \gob{1}{}
\gend}
\end{eqnarray} & &
\begin{eqnarray}\label{Hcomod}
\scalebox{0.9}[0.9]{
\gbeg{2}{3}
\got{1}{N} \gnl
\grcm \gnl
\gob{1}{N} \gob{1}{H}
\gend} = \scalebox{0.9}[0.9]{
\gbeg{8}{5}
\got{1}{} \got{1}{} \got{1}{N} \gnl
\gdb \gcl{1} \gnl
\gcn{1}{1}{1}{3} \glm \gnl
\gvac{1} \gbr \gnl
\gvac{1} \gob{1}{N} \gob{1}{H}
\gend}
\end{eqnarray}
\end{tabular}
\end{center}
Moreover, if $T$ is a right $H$-comodule algebra, then $T$ is a left $H^*$-module algebra with the structure (\ref{H*mod-S}).
\end{proposition}

\begin{lemma} \label{module S}
Let $H$ be a finite Hopf algebra and $T$ a right $H$-comodule algebra in $\C$.
Then $T$ is a left $T\# H^*$-module via:
$$\scalebox{0.9}[0.9]{
\gbeg{2}{3}
\got{1}{\hspace{-0,2cm}T\# H^*} \got{1}{\hspace{0,1cm}T} \gnl
\glm \gnl
\gvac{1} \gob{1}{T}
\gend} = \scalebox{0.9}[0.9]{
\gbeg{2}{5}
\got{1}{T} \got{1}{H^*} \got{1}{T} \gnl
\gcl{1} \gcl{1} \grcm \gnl
\gcl{1} \gbr \gcl{1} \gnl
\gmu \gev \gnl
\gob{2}{T} \gvac{1} \gob{1}{}
\gend}
$$
\end{lemma}

\begin{proof}
As a right $H$-comodule, $T$ has the structure of a left $H^*$-module given by Diagram
\ref{H*mod-S}. We compute:
$$\begin{array}{ll}
\scalebox{0.85}[0.85]{
\gbeg{5}{6}
\got{1}{\hspace{-0,24cm}T\# H^*}\got{3}{T\# H^*}\got{1}{T} \gnl
\gwmu{3} \gvac{1} \gcl{3} \gnl
\gvac{1} \gcn{1}{2}{1}{5} \gnl \gnl
\gvac{3} \glm \gnl
\gvac{4} \gob{1}{T}
\gend} & = \scalebox{0.85}[0.85]{
\gbeg{7}{8}
\got{3}{T\# H^*} \got{3}{\hspace{-0,24cm}T\# H^*} \got{1}{\hspace{-0,6cm}T} \gnl
\gcl{1} \gcmu \gcl{1} \gcl{2} \gcl{3} \gnl
\gcl{1} \gcl{1} \gbr \gnl
\gcl{1} \glm \gmu \gnl
\gwmu{3} \gvac{1} \hspace{-0,36cm} \gcn{1}{1}{1}{2} \gvac{1} \hspace{-0,34cm} \grcm \gnl
\gvac{2} \gcn{1}{1}{1}{3} \gvac{2} \gbr \gcl{1} \gnl
\gvac{3} \gwmu{3} \gev \gnl
\gvac{4} \gob{1}{T}
\gend} \stackrel{(\ref{H*mod-S})}{=} \scalebox{0.85}[0.85]{
\gbeg{7}{9}
\got{3}{T\# H^*} \got{3}{\hspace{-0,24cm}T\# H^*} \got{1}{\hspace{-0,6cm}T} \gnl
\gcl{4} \gcmu \gcl{1} \gcl{2} \gcl{3} \gnl
\gvac{1} \gcl{2} \gbr \gnl
\gvac{2} \grcm \hspace{-0,42cm} \gmu \gnl
\gvac{2} \gbr \gcl{1}  \gcn{1}{1}{0}{1}  \grcm \gnl
\gvac{1} \gmu \gev \gbr \gcl{1} \gnl
\gvac{1} \gcn{1}{1}{2}{3} \gvac{3} \gcl{1} \gev \gnl
\gvac{2} \gwmu{4} \gnl
\gvac{3} \gob{2}{T}
\gend} \stackrel{(\ref{multH*})}{\stackrel{nat.}{=}}
\scalebox{0.85}[0.85]{\gbeg{9}{8}
\got{3}{T\# H^*} \got{3}{\hspace{-0,24cm}T\# H^*} \got{1}{\hspace{-0,6cm}T} \gnl
\gcl{1} \gcmu \gcl{1} \gcl{1} \grcm \gnl
\gcl{1} \gcl{1} \gbr \gbr \gcn{1}{1}{1}{3} \gnl
\gcl{1} \gcl{1} \grcm \hspace{-0,42cm} \gcn{1}{1}{1}{2} \gcn{1}{1}{1}{2}
    \gcn{1}{1}{1}{2} \gvac{1} \hspace{-0,34cm} \gcmu \gnl
\gvac{2} \hspace{-0,22cm} \gcl{1} \gbr \gcl{1} \hspace{-0,22cm} \gbr \gbr \gcl{1} \gnl
\gvac{3} \hspace{-0,34cm} \gmu \gev \hspace{-0,2cm} \gcl{1} \gev \gev \gnl
\gvac{4} \gwmu{4} \gnl
\gvac{5} \gob{2}{T}
\gend} \vspace{3mm} \\
& \stackrel{(\ref{codiagH*})}{\stackrel{nat.}{=}}
\scalebox{0.85}[0.85]{
\gbeg{9}{8}
\got{2}{T\# H^*} \got{1}{T} \got{1}{} \got{1}{H^*} \got{1}{T} \gnl
\gcl{1} \gcl{1} \grcm \gcl{1} \grcm \gnl
\gcl{1} \gbr \gcl{1} \gbr \gcn{1}{1}{1}{3} \gnl
\gmu \gcl{1} \gbr \gcn{1}{1}{1}{2} \gvac{1} \hspace{-0,34cm} \gcmu \gnl
\gvac{1} \gcn{1}{1}{1}{2} \gvac{1} \hspace{-0,34cm} \gbr
    \gcn{1}{1}{1}{2} \gvac{1} \hspace{-0,34cm}\gbr \gcl{1} \gnl
\gvac{3} \hspace{-0,34cm} \gmu \gcn{1}{1}{1}{3} \hspace{-0,2cm} \gvac{1} \gmu \gev \gnl
\gvac{4} \gcl{1} \gvac{2} \hspace{-0,34cm} \gev \gnl
\gvac{4} \gob{2}{T}
\gend}  \stackrel{comod.}{\stackrel{nat.}{\stackrel{ass.}{=}}}
\scalebox{0.85}[0.85]{
\gbeg{7}{9}
\got{2}{T\# H^*} \got{1}{T} \got{1}{} \got{1}{H^*} \got{1}{T} \gnl
\gcl{4} \gcl{4} \grcm \gcl{1} \grcm \gnl
\gvac{2} \gcl{3} \gcl{2} \gbr \gcl{1} \gnl
\gvac{4} \grcm \hspace{-0,42cm} \gev \gnl
\gvac{4} \gbr \gcl{1} \gnl
\gvac{1} \gcn{1}{2}{1}{2} \gcn{1}{1}{1}{2} \gmu \gmu \gnl
\gvac{3} \hspace{-0,22cm} \gbr \gcn{1}{1}{3}{1} \gnl
\gvac{2} \gmu  \gev \gnl
\gvac{2} \gob{2}{T}
\gend} \stackrel{nat.}{\stackrel{comod.\x alg.}{=}}
\scalebox{0.85}[0.85]{
\gbeg{7}{8}
\got{2}{T\# H^*} \got{1}{T} \got{1}{} \got{1}{H^*} \got{1}{T} \gnl
\gcl{4} \gcl{3} \gcl{2} \gvac{1} \gcl{1} \grcm \gnl
\gvac{2} \gvac{2} \gbr \gcl{1} \gnl
\gvac{2} \gwmu{3} \gev \gnl
\gvac{1} \gcn{1}{1}{1}{3} \gvac{1} \grcm \gnl
\gcl{1} \gvac{1} \gbr \gcl{1} \gnl
\gwmu{3}  \gev \gnl
\gvac{1} \gob{1}{T}
\gend} \vspace{3mm} \\
& \hspace{10pt} = \scalebox{0.85}[0.85]{
\gbeg{5}{7}
\got{1}{\hspace{-0,24cm}T\# H^*}\got{3}{T\# H^*}\got{1}{T} \gnl
\gcn{2}{2}{1}{5} \gcn{2}{1}{1}{3} \gcl{1} \gnl
\gvac{3} \glm \gnl
\gvac{2} \gcn{2}{1}{1}{3} \gcl{1} \gnl
\gvac{3} \glm \gnl
\gvac{4} \gcl{1} \gnl
\gvac{4} \gob{1}{T}
\gend}
\end{array}$$
We also have:
$$\scalebox{0.85}[0.85]{
\gbeg{4}{5}
\got{1}{} \got{1}{} \got{1}{T} \gnl
\gu{1} \gu{1} \grcm \gnl
\gcl{1} \gbr \gcl{1} \gnl
\gmu \gev \gnl
\gob{2}{T}
\gend} \stackrel{(\ref{unitH*})}{=}
\scalebox{0.85}[0.85]{
\gbeg{2}{4}
\got{1}{T} \gnl
\grcm \gnl
\gcl{1} \gcu{1} \gnl
\gob{1}{T}
\gend}=\id_T.
$$
\qed\end{proof}

\begin{proposition}\label{ChaseSweedler Thm}
Let $H$ be a finite commutative Hopf algebra and $T$ an $H$-Galois object. Assume that the braiding of $\C$ is $H$-linear and that $\Phi_{T,T}=\Phi_{T,T}^{-1}$. Then $T$ is faithfully projective and $T\# H^*$ is an Azumaya algebra in $\C$.
\end{proposition}

\begin{proof}
To prove that $T$ is faithfully projective we need several intermediate results whose proofs are omitted because they are not difficult and require long diagramatic computations. Consider $H^*$ as a left $H$-module via
$$\scalebox{0.9}[0.9]{
\gbeg{3}{4}
\got{1}{H} \got{4}{\hspace{-0,3cm}H^*} \gnl
\gcn{2}{1}{1}{3} \gcl{1} \gnl
\gvac{1} \glm \gvac{1} \gnl
\gvac{2} \gob{2}{\hspace{-0,3cm}H^*}
\gend} = \scalebox{0.9}[0.9]{
\gbeg{3}{6}
\got{1}{H} \got{3}{\hspace{-0,3cm}H^*} \gnl
\gcl{1} \gcmu \gnl
\gibr \gcl{1} \gnl
\gcl{2} \gibr \gnl
\gvac{1} \gev \gnl
\gob{2}{\hspace{-0,3cm}H^*}
\gend} \stackrel{nat.}{=} \scalebox{0.9}[0.9]{
\gbeg{4}{6}
\got{3}{H} \got{1}{\hspace{-0,3cm}H^*} \gnl
\gvac{1} \gibr \gnl
\gcn{1}{1}{3}{2} \gvac{1} \gcl{1} \gnl
\gcmu \gcl{1} \gnl
\gcl{1} \gev \gnl
\gob{2}{\hspace{-0,3cm}H^*}\gob{2}{}
\gend}$$
With the above structure $H^*$ becomes a left $H$-module algebra. One needs for the proof that
$\Phi_{H, H^*}=\Phi_{H^*, H}^{-1}$. This is assured because we are assuming that the braiding is $H$-linear, Proposition \ref{braid-lin}. Let $A$ be a right $H$-comodule algebra. We now consider the smash product $\crta{A} \star H^*$ that is defined as follows: $\crta{A} \star H^*=A \otimes H$ as an object in $\C$, the multiplication is given by
$$\scalebox{0.9}[0.9]{
\gbeg{6}{4}
\got{3}{\crta A\star H^*} \gvac{1} \got{1}{\crta A\star H^*} \gnl
\gvac{1} \gwmu{4} \gnl
\gvac{2} \gcn{1}{1}{2}{2} \gnl
\gvac{2} \gob{2}{\crta A\star H^*}
\gend} = \scalebox{0.9}[0.9]{
\gbeg{4}{6}
\got{1}{A} \got{1}{H^*} \got{1}{A} \got{3}{H^*} \gnl
\gcl{1} \gbr \gvac{1} \gcl{3} \gnl
\gcl{1} \grcm \gcn{1}{1}{-1}{1} \gnl
\gbr \glm \gnl
\gmu \gvac{1} \gmu \gnl
\gob{2}{A} \gvac{1} \gob{2}{H^*}
\gend}$$
and the unit is $\eta_{A} \otimes \eta_{H^*}$. Notice that this smash product is a priori different from the one we have
used before. The reader may check the associativity and unit property of this new multiplication ($\Phi_{H,A}=\Phi_{A, H}^{-1}$ and $\Phi_{H, H^*}=\Phi_{H^*, H}^{-1}$ will be needed). This smash product allows to recognize the category of relative Hopf modules $\C_{A}^H$ as the category of left modules over ${\crta A} \star H^*$. The isomorphism between these categories is established as follows: $M\in\C^H_A$ becomes a left $\crta A\star H^*$-module via (\ref{mod1}) and $N\in{}_{\crta A\star H^*}\C$ is made into an object in $\C^H_A$ by the left $H^*$-action and right $A$-action given in (\ref{mod2}):
\begin{center}
\begin{tabular}{p{4.7cm}p{0.8cm}p{9.4cm}}
\begin{eqnarray}\label{mod1}
\hspace{-1cm}
\scalebox{0.9}[0.9]{
\gbeg{2}{3}
\got{1}{\hspace{-0,2cm}\crta A\star H^*} \got{1}{\hspace{0,1cm}M} \gnl
\glm \gnl
\gvac{1} \gob{1}{M}
\gend}\hspace{-5pt} := \hspace{-7pt} \scalebox{0.9}[0.9]{
\gbeg{2}{6}
\got{1}{\crta A} \got{1}{H^*} \got{1}{M} \gnl
\gcl{1} \gcl{1} \grcm \gnl
\gcl{1} \gbr \gcl{1} \gnl
\gbr \gev \gnl
\grm \gnl
\gob{1}{M} \gvac{2} \gob{1}{}
\gend}
\end{eqnarray} & &
\begin{eqnarray}\label{mod2}
\scalebox{0.9}[0.9]{
\gbeg{2}{4}
\got{1}{H^*} \got{1}{N} \gnl
\glm \gnl
\gvac{1} \gcl{1} \gnl
\gvac{1} \gob{1}{N}
\gend} \hspace{-2pt} = \hspace{-2pt} \scalebox{0.9}[0.9]{
\gbeg{4}{5}
\got{1}{ } \got{3}{H^*} \got{1}{\hspace{-0,6cm}N} \gnl
\gu{1} \gvac{1} \gcl{1} \gcl{2} \gnl
\glmpt \gcmp \gnot{\hspace{-0,7cm}\crta A\star H^*} \grmptb \gnl
\gvac{2} \glm \gnl
\gvac{2} \gob{3}{N}
\gend} \qquad
\scalebox{0.9}[0.9]{
\gbeg{2}{4}
\got{1}{N} \got{1}{A} \gnl
\grm \gnl
\gcl{1} \gnl
\gob{1}{N}
\gend} \hspace{-4pt} = \hspace{-2pt} \scalebox{0.9}[0.9]{
\gbeg{4}{6}
\got{2}{ } \got{1}{N} \got{1}{A} \gnl
\gvac{2} \gibr \gnl
\gvac{1} \gcn{1}{1}{3}{-1} \gu{1} \gcl{2} \gnl
\glmpt \gcmp \gnot{\hspace{-0,7cm}\crta A\star H^*} \grmptb \gnl
\gvac{2} \glm \gnl
\gvac{2} \gob{3}{N}
\gend}
\end{eqnarray}
\end{tabular}
\end{center}
Now we view $N$ as a right $H$-comodule by Proposition \ref{Hmod-H*comod}. These two assignments define functors by acting as the identity on morphisms and they are inverse to each other. Again for the proof the hypothesis $\Phi_{H,A}=\Phi_{A,H}^{-1}$ and $\Phi_{H, H^*}=\Phi_{H^*,H}^{-1}$ is needed. \par \medskip

By Theorem \ref{HopfThm} the functors
$\bfig
\putmorphism(0,30)(1,0)[-\ot T: \C` \C^H_T: (-)^{coH}`]{800}1a
\putmorphism(-40,-10)(1,0)[\phantom{\C^H: (-)^{coH}}`\phantom{-\ot T: \C}` ]{760}{-1}b
\efig
$
establish a $\C$-equivalence of categories. On the other hand, consider the functors
$\G:\C^H_T\to {}_{\crta T\star H^*}\C$ and $\F:{}_{\crta T\star H^*}\C\to\C^H_T$ defining the above isomorphism of categories. They are also $\C$-functors. Then $\G(-\ot T): \C\to {}_{\crta T \star H^*}\C$ is a $\C$-equivalence of categories. By Morita Theorem II (2.4.5), there is a strict Morita context $(I, \crta T\star H^*, P, Q, f, g)$ with $Q:=\G(T)$.
From Morita Theorem I (2.4.4), $Q=T$ is faithfully projective (over $I$). \par \smallskip

By the right version of Morita Theorem II, $[T,T] \cong \crta T\star H^*$ as algebras. Since $H$ is commutative, $\crta{T}$ is a right $H$-comodule algebra when endowed with the comodule structure of $T$. Then $\crta{T}$ is a left $H^*$-module algebra and we can consider the smash product $\crta{T}\# H^*$. The reader may easily verify that the multiplications in $\crta{T}\star H^*$ and $\crta{T}\# H^*$ are the same (cocommutativity of $H^*$ is needed). Hence $\crta{T}\# H^*$ is an Azumaya algebra in $\C$ for any $H$-Galois object $T$. Finally, observe that $\crta{T}$, equipped with the comodule structure of $T$, is an $H$-Galois object and that $\crta{\crta{T}}=T$ as algebras (because $\Phi_{T,T}=\Phi_{T,T}^{-1}$).
\qed\end{proof}

\subsection{Surjectivity of $\Upsilon$}

In what follows we will equip $T\# H^*$ with an $H$-module structure
so that it becomes an $H$-module algebra. This will make it an $H$-Azumaya algebra
in view of Proposition \ref{ChaseSweedler Thm}. Our goal then will be to prove that $[T\# H^*]$ is a
preimage in $\BM(\C; H)$ of $[T]\in\Gal(\C; H)$ through $\Upsilon$. Thus the map
$\Upsilon: \BM(\C; H) \to \Gal(\C; H)$ will be surjective.

\begin{lemma} \label{Ssm.H* left H-mod.alg.}
Let $H\in\C$ be a finite commutative Hopf algebra and $T\in\C$ a right $H$-comodule algebra.
The object $T\# H^*$ is a left $H$-module with the structure:
$$\scalebox{0.9}[0.9]{
\gbeg{3}{4}
\got{1}{H} \got{3}{T\# H^*} \gnl
\gcn{2}{1}{1}{3} \gcl{1} \gnl
\gvac{1} \glm \gvac{1} \gnl
\gvac{2} \gob{1}{T\# H^*}
\gend} = \scalebox{0.9}[0.9]{
\gbeg{4}{6}
\got{1}{H} \got{1}{T} \got{3}{\hspace{-0,3cm}H^*} \gnl
\gcl{1} \gcl{1} \gcmu \gnl
\gcl{1} \gibr \gcl{1} \gnl
\gibr \gcl{2} \gcl{2} \gnl
\gev \gnl
\gvac{1} \gob{3}{T} \gob{1}{\hspace{-0,6cm}H^*}
\gend}
$$
If furthermore $\Phi$ is $H$-linear, then the above makes $T\# H^*$ into a left $H$-module algebra.
\end{lemma}

\begin{proof}
We first have to check that the following diagrams are equal:
$$L:=\scalebox{0.9}[0.9]{
\gbeg{5}{5}
\got{1}{H} \got{1}{H} \got{3}{T\#H^*} \gnl
\gmu \gvac{1} \gcl{2} \gnl
\gcn{1}{1}{2}{5} \gnl
\gvac{2} \glm \gnl
\gvac{2} \gob{3}{T\#H^*}
\gend} = \scalebox{0.9}[0.9]{
\gbeg{5}{7}
\got{1}{H} \got{1}{H} \got{1}{T} \got{2}{\hspace{-0,14cm}H^*} \gnl
\gmu \gcl{2} \gcn{1}{1}{1}{2} \gnl
\gcn{1}{1}{2}{3} \gvac{1} \gcl{1} \gcmu \gnl
\gvac{1} \gcl{1} \gibr \gcl{1} \gnl
\gvac{1} \gibr \gcl{2} \gcl{2} \gnl
\gvac{1} \gev \gnl
\gvac{2} \gob{3}{T} \gob{1}{\hspace{-0,6cm}H^*}
\gend} \quad\textnormal{and}\quad \scalebox{0.9}[0.9]{
\gbeg{6}{8}
\got{1}{H} \got{1}{H} \got{1}{T} \got{3}{\hspace{-0,3cm}H^*} \gnl
\gcl{4} \gcl{2} \gcl{1} \gcmu \gnl
\gvac{2} \gibr \gcn{1}{1}{1}{2} \gnl
\gvac{1} \gibr \gcl{1} \gcmu \gnl
\gvac{1} \gev \gibr \gcl{3} \gnl
\gcn{1}{1}{1}{3} \gvac{1} \gcn{1}{1}{3}{1} \gvac{1} \gcl{2} \gnl
\gvac{1} \gibr \gnl
\gvac{1} \gev \gnl
\gvac{3} \gob{3}{T} \gob{1}{\hspace{-0,6cm}H^*}
\gend} = \scalebox{0.9}[0.9]{
\gbeg{5}{6}
\got{1}{H} \got{1}{H} \got{3}{T\#H^*} \gnl
\gcl{1} \gcn{1}{1}{1}{3} \gvac{1} \gcl{1} \gnl
\gcn{1}{2}{1}{5} \gvac{1} \glm \gnl
\gvac{3} \gcl{1} \gnl
\gvac{2} \glm \gnl
\gvac{2} \gob{3}{T\#H^*}
\gend} =:R
$$
Starting by applying commutativity of $H$ and naturality, we develop $L$ as
follows:
$$L=\scalebox{0.85}[0.85]{
\gbeg{5}{8}
\got{1}{H} \got{1}{H} \got{1}{T} \got{3}{\hspace{-0,3cm}H^*} \gnl
\gcl{3} \gcl{2} \gcl{1} \gcmu \gnl
\gvac{2} \gibr \gcl{5}\gnl
\gvac{1} \gibr \gcl{4} \gnl
\gibr \gcl{1} \gvac{1} \gnl
\gcl{1} \gibr \gnl
\gcn{1}{1}{1}{2} \gmu \gnl
\gvac{1} \hspace{-0,22cm} \gev \gnl
\gvac{3} \gob{2}{T} \gob{1}{\hspace{-0,2cm}H^*}
\gend} \hspace{-0,14cm}\stackrel{(\ref{codiagH*})}{=}\hspace{-0,2cm}
\scalebox{0.9}[0.9]{
\gbeg{6}{10}
\got{1}{} \got{1}{H} \got{1}{H} \got{1}{T} \got{3}{\hspace{-0,3cm}H^*} \gnl
\gvac{1} \gcl{3} \gcl{2} \gcl{1} \gcmu \gnl
\gvac{3} \gibr \gcl{7}\gnl
\gvac{2} \gibr \gcl{6} \gnl
\gvac{1} \gibr \gcl{2} \gnl
\gvac{1} \gcn{1}{1}{1}{0} \gcl{1} \gnl
\gcmu \gibr \gnl
\gcl{1} \gbr \gcl{1} \gnl
\gev \gev \gnl
\gvac{4} \gob{1}{T} \gob{1}{\hspace{0,2cm}H^*}
\gend} \hspace{-0,18cm}\stackrel{2\times nat.}{=}\hspace{-0,12cm}
\scalebox{0.9}[0.9]{
\gbeg{6}{11}
\got{1}{H} \got{1}{H} \got{1}{} \got{1}{T} \got{3}{\hspace{-0,3cm}H^*} \gnl
\gcl{3} \gcl{3} \gvac{1} \gcl{1} \gcmu \gnl
\gvac{3} \gibr \gcl{8} \gnl
\gvac{2} \gcn{1}{1}{3}{2} \gvac{1} \gcl{7} \gnl
\gibr \gcmu \gnl
\gcl{1} \gibr \gcl{1} \gnl
\gibr \gibr \gnl
\gcl{1} \gibr \gcl{1} \gnl
\gcl{1} \gbr \gcl{1} \gnl
\gev \gev \gnl
\gvac{4} \gob{1}{T} \gob{1}{\hspace{0,2cm}H^*}
\gend}\hspace{-0,18cm}\stackrel{nat.}{=}\hspace{-0,12cm}
\scalebox{0.9}[0.9]{
\gbeg{6}{11}
\got{1}{H} \got{1}{H} \got{1}{} \got{1}{T} \got{3}{\hspace{-0,3cm}H^*} \gnl
\gcl{6} \gcl{4} \gvac{1} \gcl{1} \gcmu \gnl
\gvac{3} \gibr \gcl{8} \gnl
\gvac{2} \gcn{1}{1}{3}{2} \gvac{1} \gcl{7} \gnl
\gvac{2} \gcmu \gnl
\gvac{1} \gibr \gcl{2} \gnl
\gvac{1} \gev \gnl
\gcn{1}{1}{1}{3} \gvac{1} \gcn{1}{1}{3}{1} \gnl
\gvac{1} \gibr \gnl
\gvac{1} \gev \gnl
\gvac{3} \gob{3}{T} \gob{1}{\hspace{-0,6cm}H^*}
\gend} \hspace{-0,22cm}\stackrel{nat.}{\stackrel{coass.}{=}}\hspace{-0,14cm}R
$$
The compatibility with unit is also satisfied:
$$\scalebox{0.85}[0.85]{
\gbeg{3}{5}
\got{1}{} \got{2}{T\#H^*} \gnl
\gu{1} \gcl{1} \gnl
\glm \gnl
\gvac{1} \gcl{1} \gnl
\gvac{1} \gob{2}{T\#H^*}
\gend} = \scalebox{0.85}[0.85]{
\gbeg{4}{6}
\got{1}{} \got{1}{T} \got{3}{\hspace{-0,3cm}H^*} \gnl
\gu{1} \gcl{1} \gcmu \gnl
\gcl{1} \gibr \gcl{1} \gnl
\gibr \gcl{2} \gcl{2} \gnl
\gev \gnl
\gvac{1} \gob{3}{T} \got{1}{\vspace{-0,15cm}\hspace{-0,6cm}H^*}
\gend} \stackrel{(\ref{counitH*})}{=} \scalebox{0.85}[0.85]{
\gbeg{4}{5}
\got{1}{T} \got{3}{\hspace{-0,3cm}H^*} \gnl
\gcl{1} \gcmu \gnl
\gibr \gcl{2} \gnl
\gcu{1} \gcl{1} \gnl
\gvac{1} \gob{1}{T} \gob{1}{\hspace{0,26cm}H^*}
\gend} = \scalebox{0.85}[0.85]{
\gbeg{4}{4}
\got{2}{T\#H^*} \gnl
\gcn{1}{2}{2}{2}  \gnl
\gob{2}{T\#H^*}
\gend}
$$
This $H$-module structure will be compatible with the algebra structure of
$T\# H^*$. To prove this we should show first that
$$\scalebox{0.85}[0.85]{
\gbeg{5}{5}
\got{1}{H} \got{1}{} \got{1}{T\# H^*} \got{4}{\hspace{-0,2cm}T\# H^*} \gnl
\gcn{2}{2}{2}{5} \gwmu{3} \gnl
\gvac{3} \gcl{1} \gnl
\gvac{2} \glm \gvac{1} \gnl
\gvac{3} \gob{1}{T\# H^*}
\gend} = \scalebox{0.85}[0.85]{
\gbeg{4}{6}
\got{2}{H} \got{2}{\hspace{-0,1cm}T\# H^*} \got{3}{\hspace{-0,2cm}T\# H^*} \gnl
\gcmu \gcl{1} \gcn{1}{2}{3}{1} \gnl
\gcl{1} \gbr \gnl
\glm \glm \gnl
\gvac{1} \gwmu{3} \gnl
\gvac{2} \gob{1}{T\# H^*}
\gend}
$$
As a right $H$-comodule, $T$ is a left $H^*$-module with the structure given in Diagram \ref{H*mod-S}. Together with the above proved $H$-module structure on $T\# H^*$ the preceding question transforms to
$$\scalebox{0.85}[0.85]{
\gbeg{7}{11}
\got{1}{H} \got{1}{T} \got{2}{H^*} \got{1}{T} \got{1}{} \got{1}{H^*} \gnl
\gcl{5} \gcl{4} \gcmu \gcl{1} \gvac{1} \gcl{3} \gnl
\gvac{2} \gcl{2} \gbr \gnl
\gvac{3} \grcm  \gcn{1}{1}{-1}{1} \gnl
\gvac{2} \gbr \gcl{1} \gmu \gnl
\gvac{1} \gmu \gev \gcn{1}{1}{2}{-1}  \gnl
\gcn{1}{2}{1}{3} \gcn{1}{1}{2}{3} \gvac{1} \gwcm{3} \gnl
\gvac{2} \gibr \gvac{1} \gcl{3} \gnl
\gvac{1} \gibr \gcl{2} \gnl
\gvac{1} \gev \gnl
\gvac{3} \gob{1}{T} \gvac{1} \gob{2}{\hspace{-0,2cm}H^*}
\gend} = \scalebox{0.85}[0.85]{
\gbeg{8}{12}
\got{1}{} \got{2}{H} \got{1}{T} \got{1}{\hspace{0,14cm}H^*} \got{1}{T} \got{1}{\hspace{0,2cm}H^*} \gnl
\gvac{1} \gcmu \gcl{1} \gcl{1} \gcl{1} \gcn{1}{1}{1}{2} \gnl
\gvac{1} \gcl{2} \gbr \gcl{1} \gcl{1} \gcmu \gnl
\gvac{2} \gcl{2} \gbr \gibr \gcl{6} \gnl
\gvac{1} \gcn{1}{1}{1}{0} \gvac{1} \gcl{1} \gibr \gcl{1} \gnl
\gvac{1} \hspace{-0,34cm} \gcl{2} \gcn{1}{1}{2}{1} \gcmu \hspace{-0,22cm} \gev \gcl{1} \gnl
\gvac{3} \hspace{-0,34cm} \gibr \hspace{-0,22cm} \gcmu \gcn{1}{1}{3}{1} \gnl
\gvac{3} \hspace{-0,34cm} \gibr \gcn{1}{1}{1}{0} \hspace{-0,22cm} \gcl{2} \gbr \gnl
\gvac{4} \hspace{-0,2cm} \gev \hspace{-0,22cm} \gcl{2} \gvac{1} \grcm \gcn{1}{1}{-1}{1} \gnl
\gvac{7} \gbr \gcl{1} \gmu \gnl
\gvac{6} \gmu \gev \gcn{1}{1}{2}{2} \gnl
\gvac{6} \gob{2}{T} \gvac{3} \gob{1}{H^*}
\gend}
$$
Let $\Sigma$ and $\Omega$ denote the left and right-hand side diagrams, respectively.
We develop $\Sigma$ as follows:
$$\begin{array}{ll}
\Sigma & \stackrel{bialg.}{=} \scalebox{0.85}[0.75]{
\gbeg{9}{12}
\got{1}{H} \got{1}{T} \got{2}{H^*} \got{1}{T} \got{2}{} \got{1}{H^*} \gnl
\gcl{5} \gcl{4} \gcmu \gcl{1} \gvac{2} \gcl{2} \gnl
\gvac{2} \gcl{2} \gbr \gnl
\gvac{3} \grcm  \gcn{1}{1}{-1}{2} \gcn{1}{1}{3}{4} \gnl
\gvac{2} \gbr \gcl{1} \gcmu \gcmu \gnl
\gvac{1} \gmu \gev \gcl{1} \gbr \gcl{1} \gnl
\gcl{1} \gcn{1}{1}{2}{5} \gvac{1} \gvac{2} \gmu \gmu \gnl
\gcn{1}{2}{1}{5} \gvac{2} \gcl{1} \gcn{2}{1}{4}{1} \gvac{1} \gcn{1}{1}{2}{1} \gnl
\gvac{3} \gibr \gvac{2} \gcl{3} \gnl
\gvac{2} \gibr \gcl{2} \gnl
\gvac{2} \gev \gnl
\gvac{4} \gob{1}{T} \gvac{2} \gob{2}{\hspace{-0,2cm}H^*}
\gend} \stackrel{nat.}{\stackrel{(\ref{multH*})}{=}}
\scalebox{0.85}[0.75]{
\gbeg{9}{14}
\got{1}{H} \got{1}{T} \got{2}{H^*} \got{1}{T} \got{2}{} \got{1}{H^*} \gnl
\gcl{4} \gcl{4} \gcmu \gcl{1} \gvac{2} \gcl{2} \gnl
\gvac{2} \gcl{2} \gbr \gnl
\gvac{3} \grcm  \gcn{1}{1}{-1}{2} \gcn{1}{1}{3}{4} \gnl
\gvac{2} \gbr \gcl{1} \gcmu \gcmu \gnl
\gcn{1}{1}{1}{2} \gmu \gev \gcl{1} \gbr \gcl{1} \gnl
\gcmu \gcn{1}{1}{0}{2} \gvac{1} \gcn{2}{1}{3}{0} \gcn{1}{1}{1}{-2} \gmu \gnl
\gcn{1}{1}{1}{2} \gcn{1}{1}{1}{2} \gvac{1} \hspace{-0,34cm} \gibr \gcl{1} \gvac{2} \gcl{6} \gnl
\gvac{1} \gcl{1} \gibr \gibr \gvac{2} \gnl
\gvac{1} \gibr \gibr \gcl{1} \gnl
\gvac{1} \gcl{2} \gibr \gcl{2} \gcl{3} \gnl
\gvac{2} \gbr \gnl
\gvac{1} \gev \gev \gnl
\gvac{5} \gob{1}{T} \gvac{2} \gob{2}{\hspace{-0,2cm}H^*}
\gend} \stackrel{nat.}{\stackrel{coass.}{=}}
\scalebox{0.85}[0.75]{
\gbeg{10}{15}
\got{2}{H} \got{1}{T} \got{1}{} \got{2}{H^*} \got{1}{} \got{1}{T} \got{1}{H^*} \gnl
\gcmu \gcl{4} \gvac{1} \gcmu \gvac{1} \gcl{1} \gcl{2} \gnl
\gcl{7} \gcl{4} \gvac{2} \hspace{-0,34cm} \gcmu \gcn{1}{1}{0}{1} \gcn{1}{1}{2}{1} \gnl
\gvac{4} \gcl{3} \gcl{1} \gbr \gcmu \gnl
\gvac{5} \gbr \gbr \gcl{1} \gnl
\gvac{2} \gcn{1}{1}{2}{3} \gvac{2} \grcm \hspace{-0,42cm} \gcn{1}{1}{1}{2} \gcn{1}{1}{1}{2} \gmu \gnl
\gvac{2} \gcn{1}{1}{2}{4} \gvac{1} \gcl{1} \gbr \gcl{1} \hspace{-0,21cm} \gcl{2} \gcl{5} \gcl{8} \gnl
\gvac{4} \gcl{2} \hspace{-0,2cm} \gmu \gev \gnl
\gvac{6} \hspace{-0,22cm} \gcn{1}{1}{1}{3} \gvac{1} \gcn{1}{1}{3}{1} \gnl
\gvac{3} \gcn{1}{2}{1}{5} \gvac{1} \gcn{1}{1}{1}{3} \gvac{1} \gibr \gvac{1} \gnl
\gvac{6} \gibr \gcl{1} \gnl
\gvac{5} \gibr \gbr \gcn{1}{1}{3}{1} \gnl
\gvac{5} \gev \gcl{2} \gibr \gnl
\gvac{8} \gev \gnl
\gvac{7} \gob{1}{T} \gvac{3} \gob{2}{\hspace{-0,2cm}H^*}
\gend} \vspace{2mm} \\
 & \stackrel{Prop. \ref{braid-lin}}{\stackrel{\Phi_{H^*, H^*}}{\stackrel{Prop.\ref{braid-lin}}
{\stackrel{\Phi_{H, H^*}}{=}}}}
\scalebox{0.85}[0.75]{
\gbeg{9}{15}
\got{2}{H} \got{1}{T} \got{1}{} \got{2}{H^*} \got{1}{} \got{1}{T} \got{1}{H^*} \gnl
\gcmu \gcl{4} \gvac{1} \gcmu \gvac{1} \gcl{1} \gcl{2} \gnl
\gcl{7} \gcl{4} \gvac{2} \hspace{-0,34cm} \gcmu \gcn{1}{1}{0}{1} \gcn{1}{1}{2}{1} \gnl
\gvac{4} \gcl{3} \gcl{1} \gbr \gcmu \gnl
\gvac{5} \gbr \gibr \gcl{1} \gnl
\gvac{2} \gcn{1}{1}{2}{3} \gvac{2} \grcm \hspace{-0,42cm} \gcn{1}{1}{1}{2} \gcn{1}{1}{1}{2} \gmu \gnl
\gvac{2} \gcn{1}{1}{2}{4} \gvac{1} \gcl{1} \gbr \gcl{1} \hspace{-0,21cm} \gcl{2} \gcl{5} \gcl{8} \gnl
\gvac{4} \gcl{2} \hspace{-0,2cm} \gmu \gev \gnl
\gvac{6} \hspace{-0,22cm} \gcn{1}{1}{1}{3} \gvac{1} \gcn{1}{1}{3}{1} \gnl
\gvac{3} \gcn{1}{2}{1}{5} \gvac{1} \gcn{1}{1}{1}{3} \gvac{1} \gibr \gvac{1} \gnl
\gvac{6} \gbr \gcl{1} \gnl
\gvac{5} \gibr \gbr \gcn{1}{1}{3}{1} \gnl
\gvac{5} \gev \gcl{2} \gibr \gnl
\gvac{8} \gev \gnl
\gvac{7} \gob{1}{T} \gvac{3} \gob{2}{\hspace{-0,2cm}H^*}
\gend} \stackrel{nat.}{=}
\scalebox{0.85}[0.75]{
\gbeg{8}{17}
\got{1}{} \got{1}{H} \got{1}{} \got{1}{T} \got{1}{\hspace{0,14cm}H^*} \got{1}{T} \got{1}{\hspace{0,2cm}H^*} \gnl
\gwcm{3} \gcl{1} \gcl{1} \gcl{1} \gcn{1}{1}{1}{2} \gnl
\gcl{10} \gvac{1} \gbr \gcl{1} \gcl{1} \gcmu \gnl
\gvac{2} \gcl{2} \gbr \gibr \gcl{5} \gnl
\gvac{3} \gcl{1} \gibr \gcl{1} \gnl
\gvac{1} \gcn{1}{1}{3}{1} \gvac{1} \hspace{-0,34cm} \gcmu \hspace{-0,22cm} \gev \gcl{1} \gnl
\gvac{2} \gcl{5} \gcmu \gcn{1}{1}{0}{1} \gcn{1}{1}{3}{1} \gnl
\gvac{3} \gcl{3} \gcl{1} \gbr \gnl
\gvac{4} \gbr \gwmu{3} \gnl
\gvac{4} \grcm \gcn{1}{1}{-1}{1} \gcl{7} \gnl
\gvac{3} \gbr \gcl{1} \gcl{1} \gnl
\gvac{2} \gmu \gev \gcn{1}{1}{1}{0} \gnl
\gvac{1} \gcn{1}{2}{1}{4} \gcn{1}{1}{2}{4} \gvac{1} \gcn{1}{1}{4}{2} \gnl
\gvac{4} \hspace{-0,34cm} \gibr \gnl
\gvac{3} \gibr \gcl{2} \gnl
\gvac{3} \gev \gnl
\gvac{5} \gob{1}{T} \gvac{2} \gob{1}{\hspace{-0,2cm}H^*}
\gend} \stackrel{nat.}{=}
\scalebox{0.85}[0.75]{
\gbeg{8}{13}
\got{1}{} \got{1}{H} \got{1}{} \got{1}{T} \got{1}{\hspace{0,14cm}H^*} \got{1}{T} \got{1}{\hspace{0,2cm}H^*} \gnl
\gwcm{3} \gcl{1} \gcl{1} \gcl{1} \gcn{1}{1}{1}{2} \gnl
\gcl{7} \gvac{1} \gbr \gcl{1} \gcl{1} \gcmu \gnl
\gvac{2} \gcl{2} \gbr \gibr \gcl{5} \gnl
\gvac{3} \gcl{1} \gibr \gcl{1} \gnl
\gvac{1} \gcn{1}{1}{3}{1} \gvac{1} \hspace{-0,34cm} \gcmu \hspace{-0,22cm} \gev \gcl{1} \gnl
\gvac{2} \gcl{2} \gcmu \gcn{1}{1}{0}{1} \gcn{1}{1}{3}{1} \gnl
\gvac{3} \gibr \gbr \gnl
\gvac{2} \gibr \gcl{2} \gcl{1} \gwmu{3} \gnl
\gvac{1} \gibr \gcl{2} \gvac{1} \grcm \gcl{3} \gnl
\gvac{1} \gev \gvac{1} \gbr \gcl{1} \gcl{1} \gnl
\gvac{3} \gmu \gev \gnl
\gvac{3} \gob{2}{T} \gvac{2} \gob{2}{\hspace{-0,2cm}H^*}
\gend} \vspace{1mm} \\
& \stackrel{H^*}{\stackrel{coc.}{=}}
\scalebox{0.85}[0.75]{
\gbeg{8}{12}
\got{1}{} \got{1}{H} \got{1}{} \got{1}{T} \got{1}{\hspace{0,14cm}H^*} \got{1}{T} \got{1}{\hspace{0,2cm}H^*} \gnl
\gwcm{3} \gcl{1} \gcl{1} \gcl{1} \gcn{1}{1}{1}{2} \gnl
\gcl{6} \gvac{1} \gbr \gcl{1} \gcl{1} \gcmu \gnl
\gvac{2} \gcl{2} \gbr \gibr \gcl{6} \gnl
\gvac{3} \gcl{1} \gibr \gcl{1} \gnl
\gvac{1} \gcn{1}{1}{3}{1} \gvac{1} \hspace{-0,34cm} \gcmu \hspace{-0,22cm} \gev \gcl{1} \gnl
\gvac{2} \gcl{1} \gcmu \gcn{1}{1}{0}{1} \gcn{1}{1}{3}{1} \gnl
\gvac{2} \gibr \gcl{2} \gbr \gnl
\gvac{1} \gibr \gcl{2} \gvac{1} \grcm \gcn{1}{1}{-1}{1} \gnl
\gvac{1} \gev \gvac{1} \gbr  \gcl{1} \gmu \gnl
\gvac{3} \gmu \gev \gvac{1} \hspace{-0,2cm} \gcl{1} \gnl
\gvac{4} \gob{1}{T} \gvac{3} \gob{1}{\hspace{0,2cm}H^*}
\gend} \stackrel{coass.}{=}\Omega.
\end{array}$$
In the fourth and penultimate equality we have used that $H^*$ is cocommutative since $H$ is finite and commutative by assumption. Finally, the $H$-module structure of $T\# H^*$ is compatible with the unit:
$$\scalebox{0.85}[0.85]{
\gbeg{2}{5}
\got{1}{H} \got{1}{} \gnl
\gcl{1} \gu{1} \gnl
\glm \gnl
\gvac{1} \gcl{1} \gnl
\gvac{1} \gob{1}{T\# H^*}
\gend} = \scalebox{0.85}[0.85]{
\gbeg{4}{7}
\got{1}{H} \got{1}{} \got{1}{} \gnl
\gcl{3} \gu{1} \gvac{1} \hspace{-0,34cm} \gu{1} \gnl
\gvac{2} \hspace{-0,36cm} \gcl{1} \gcmu \gnl
\gvac{2} \gibr \gcl{1} \gnl
\gvac{1} \gibr \gcl{2} \gcl{2} \gnl
\gvac{1} \gev \gnl
\gvac{3} \gob{1}{T} \gob{1}{\hspace{0,2cm}H^*}
\gend} = \scalebox{0.85}[0.85]{
\gbeg{4}{6}
\got{1}{H} \got{1}{} \got{1}{} \gnl
\gcl{1} \gu{1} \gu{1} \gu{1} \gnl
\gcl{1} \gibr \gcl{3} \gnl
\gibr \gcl{2} \gnl
\gev \gnl
\gvac{2} \gob{1}{T} \gob{1}{\hspace{0,2cm}H^*}
\gend} \stackrel{(\ref{unitH*})}{=}
\scalebox{0.85}[0.85]{
\gbeg{3}{4}
\got{1}{H} \got{1}{} \got{1}{} \gnl
\gcl{1} \gu{1} \gu{1} \gnl
\gcu{1} \gcl{1} \gcl{1} \gnl
\gvac{1} \gob{1}{T} \gob{1}{\hspace{0,2cm}H^*}
\gend} =
\scalebox{0.85}[0.85]{
\gbeg{2}{4}
\got{1}{H} \gnl
\gcu{1} \gnl
\gu{1} \gnl
\gob{1}{\hspace{0,3cm}T\# H^*}
\gend} $$
\qed\end{proof}

As announced we now prove that $\Upsilon$ is surjective by showing that $\Upsilon([T\# H^*])=[T]$.
Let $\gamma: [(T\# H^*)\# H]^{T\# H^*} \to T$ be the morphism given by $\gamma=(T\ot\Epsilon_{H^*}\ot\Epsilon_H)j.$ We prove that $\gamma$ is an isomorphism of $H$-Galois objects by proving that it is an $H$-comodule algebra morphism
(Proposition \ref{comod-alg-iso}). Before this let us deduce several identities that hold on $[(T\# H^*)\# H]^{T\# H^*}$. To its equalizer property, expressed in Diagram \ref{Asmash etaH-prop.} with
$A:=T\# H^*$, we will apply $\sigma:=T\ot\Epsilon_{H^*}\ot\Epsilon_H$. Recalling the algebra structure of $A=T\# H^*$, which we already dealt with in Lemma \ref{Ssm.H* left H-mod.alg.} (taking into account the left $H^*$-module
structure of the right $H$-comodule $T$), we obtain:
$$\scalebox{0.85}[0.85]{
\gbeg{9}{8}
\got{1}{T} \got{2}{H^*} \got{2}{} \got{2}{[(T\# H^*)\# H]^{T\# H^*}}  \gnl
\gcl{4} \gcmu \glmpb \gnot{j}\gcmptb\grmpb \gnl
\gvac{1} \gcl{3} \gbr \gcl{1} \gcl{3}   \gnl
\gvac{2} \grcm \hspace{-0,42cm} \gmu   \gnl
\gvac{3} \gibr \gcn{1}{1}{0}{1} \gnl
\gvac{1} \gcn{1}{1}{1}{2} \gev \gcn{1}{1}{1}{0} \gcu{1} \gcu{1} \gnl
\gvac{2} \hspace{-0,35cm} \gwmu{3} \gnl
\gvac{3} \gob{1}{T}
\gend} = \scalebox{0.85}[0.75]{
\gbeg{6}{21}
\got{1}{T} \got{1}{H^*} \got{2}{} \got{2}{[(T\# H^*)\# H]^{T\# H^*}}  \gnl
\gcl{1} \gcl{1} \gnot{\hspace{1,3cm}j} \glmp \gcmptb \gcmpb \grmpb \gnl
\gcn{1}{1}{1}{3} \gcn{1}{1}{1}{3} \gvac{1} \gcl{1} \gcl{2} \gcl{3} \gnl
\gvac{1} \gcl{1} \gbr  \gnl
\gvac{1} \gbr \gbr  \gnl
\gvac{1} \gcl{1} \gbr \gbr \gnl
\gvac{1} \gcl{1} \gcl{1} \gbr \gcl{1} \gnl
\gcn{1}{1}{3}{1} \gcn{1}{1}{3}{1} \gvac{1} \hspace{-0,36cm} \gcmu \gcn{1}{1}{0}{1}
    \gcn{1}{1}{0}{1} \gnl
\gvac{1} \hspace{-0,34cm} \gcl{10} \gcl{5} \gcn{1}{1}{2}{1} \gvac{1} \hspace{-0,36cm}
   \gbr \gcl{1} \gnl
\gvac{4} \hspace{-0,34cm} \gcl{3} \gcn{1}{1}{2}{1} \gvac{1} \hspace{-0,36cm} \gbr \gnl
\gvac{6} \hspace{-0,22cm} \gcl{1} \gcmu \gcn{1}{1}{0}{1} \gnl
\gvac{6} \gbr \gcl{5} \gcl{8} \gnl
\gvac{5} \gibr \gcl{2} \gnl
\gvac{4} \gcn{1}{1}{1}{2} \gev \gnl
\gvac{4} \gcmu \gcn{1}{1}{3}{1} \gnl
\gvac{4} \gcl{3} \gbr \gnl
\gvac{5} \grcm \hspace{-0,4cm} \gwmu{3} \gnl
\gvac{6} \gibr \gcl{2} \gnl
\gvac{4} \gcn{1}{1}{1}{2} \gev \gcn{1}{1}{1}{0} \gnl
\gvac{5} \hspace{-0,34cm} \gwmu{3} \gvac{1} \hspace{-0,36cm} \gcu{1} \gvac{1}
   \hspace{-0,16cm} \gcu{1} \gnl
\gvac{6} \gob{2}{T} \gvac{3} \gob{2}{}
\gend}
$$

By the corresponding properties of the counits this simplifies to
\begin{eqnarray} \label{surjec-equal}
\scalebox{0.85}[0.85]{
\gbeg{8}{7}
\got{1}{T} \got{1}{H^*} \got{3}{} \got{1}{[(T\# H^*)\# H]^{T\# H^*}}  \gnl
\gcl{3} \gcl{3} \glmpb \gnot{j}\gcmptb\grmpb \gnl
\gvac{2} \grcm \gcn{1}{1}{-1}{1} \gcn{1}{1}{-1}{1} \gnl
\gvac{2} \gibr \gcu{1} \gcu{1} \gnl
\gcn{1}{1}{1}{2} \gev \gcn{1}{1}{1}{0} \gnl
\gvac{1} \hspace{-0,35cm} \gwmu{3} \gnl
\gvac{2} \gob{1}{T}
\gend} = \scalebox{0.85}[0.75]{
\gbeg{6}{14}
\got{1}{T} \got{1}{H^*} \got{3}{} \got{1}{[(T\# H^*)\# H]^{T\# H^*}}  \gnl
\gcl{1} \gcl{1} \gnot{\hspace{1,3cm}j} \glmp \gcmptb \gcmp \grmp \gnl
\gcn{1}{1}{1}{3} \gcn{1}{1}{1}{3} \gvac{1} \gcl{1} \gcl{2} \gcl{3} \gnl
\gvac{1} \gcl{1} \gbr  \gnl
\gvac{1} \gbr \gbr  \gnl
\gvac{1} \gcl{5} \gbr \gbr \gnl
\gvac{2} \gcl{4} \gbr \gcl{1} \gnl
\gvac{3} \gcl{1} \gibr \gnl
\gvac{3} \gibr \grcm \gnl
\gvac{3} \gev \gibr \gnl
\gvac{1} \gcn{1}{2}{1}{4} \gcn{1}{1}{1}{3} \gvac{1} \gcn{1}{1}{3}{1} \gcn{1}{2}{3}{0} \gnl
\gvac{3} \gev \gnl
\gvac{3} \hspace{-0,2cm}\gwmu{3} \gnl
\gvac{4} \gob{1}{T}
\gend}
\end{eqnarray}
Neutralizing the braiding $\Phi_{T\# H^*, [(T\# H^*)\# H]^{T\# H^*}}$ on the
right-hand side, we compose the whole equation from above (in the braided diagrams
orientation) with $\Phi_{[(T\# H^*)\# H]^{T\# H^*},T\# H^*}^{-1}$. Then the
above expression takes the form:
\begin{eqnarray} \label{surjec-equal-1}
\scalebox{0.85}[0.85]{
\gbeg{7}{11}
\got{2}{[(T\# H^*)\# H]^{T\# H^*}} \got{1}{} \got{3}{T} \got{1}{\hspace{-0,56cm}H^*} \gnl
\glmptb \gnot{j} \gcmpb\grmpb \gcn{1}{3}{3}{-1} \gcn{1}{4}{3}{-1} \gnl
\gcn{1}{2}{1}{3} \gcu{1} \gcu{1} \gnl \gnl
\gvac{1} \gibr \gnl
\gvac{1} \gcl{3} \gibr \gnl
\gvac{2} \gcl{2} \grcm \gnl
\gvac{3} \gibr \gnl
\gvac{1} \gcn{1}{1}{1}{2} \gev \gcn{1}{1}{1}{0} \gnl
\gvac{2} \hspace{-0,35cm} \gwmu{3} \gnl
\gvac{3} \gob{1}{T}
\gend} = \scalebox{0.85}[0.85]{
\gbeg{6}{10}
\gvac{1} \got{2}{[(T\# H^*)\# H]^{T\# H^*}} \got{2}{} \got{1}{T} \got{1}{H^*} \gnl
\gnot{\hspace{1,3cm}j} \glmp \gcmptb \gcmpb \grmpb \gcn{1}{2}{3}{1} \gcn{1}{2}{3}{1} \gnl
\gvac{1} \gcl{4} \gcl{4} \gcl{2}  \gnl
\gvac{4} \gibr \gnl
\gvac{3} \gibr \grcm \gnl
\gvac{3} \gev \gibr \gnl
\gvac{1} \gcn{1}{2}{1}{4} \gcn{1}{1}{1}{3} \gvac{1} \gcn{1}{1}{3}{1} \gcn{1}{2}{3}{0} \gnl
\gvac{3} \gev \gnl
\gvac{3} \hspace{-0,2cm}\gwmu{3} \gnl
\gvac{4} \gob{1}{T}
\gend}
\end{eqnarray}

On the other hand, from (\ref{surjec-equal}) we can obtain further equations,
\begin{eqnarray} \label{surjec-equal-2}
\scalebox{0.85}[0.85]{
\gbeg{6}{4}
\got{1}{T} \got{1}{} \got{4}{[(T\# H^*)\# H]^{T\# H^*}}  \gnl
\gcl{1} \glmpb \gnot{j}\gcmptb\grmpb \gnl
\gmu \gcu{1} \gcu{1} \gnl
\gob{2}{T}
\gend} \hspace{-0,1cm}=\hspace{-0,14cm} \scalebox{0.85}[0.85]{
\gbeg{7}{7}
\got{1}{T} \got{2}{} \got{4}{[(T\# H^*)\# H]^{T\# H^*}}  \gnl
\gcl{3} \gu{1} \glmpb \gnot{j}\gcmptb\grmpb \gnl
\gvac{1} \gcl{2} \grcm \gcn{1}{1}{-1}{1} \gcn{1}{1}{-1}{1} \gnl
\gvac{2} \gibr \gcu{1} \gcu{1} \gnl
\gcn{1}{1}{1}{2} \gev \gcn{1}{1}{1}{0} \gnl
\gvac{1} \hspace{-0,34cm} \gwmu{3} \gnl
\gvac{2} \gob{1}{T}
\gend} \hspace{-0,38cm}\stackrel{(\ref{surjec-equal})}{=}\hspace{-0,38cm}
\scalebox{0.85}[0.85]{
\gbeg{7}{14}
\got{1}{T} \got{2}{} \got{4}{[(T\# H^*)\# H]^{T\# H^*}}  \gnl
\gcl{1} \gu{1} \gvac{1} \glmpb \gnot{j}\gcmptb\grmpb \gnl
\gcn{1}{1}{1}{3} \gcn{1}{1}{1}{3} \gvac{1} \gcl{1} \gcl{2} \gcl{3} \gnl
\gvac{1} \gcl{1} \gbr  \gnl
\gvac{1} \gbr \gbr  \gnl
\gvac{1} \gcl{5} \gbr \gbr \gnl
\gvac{2} \gcl{4} \gbr \gcl{1} \gnl
\gvac{3} \gcl{1} \gibr \gnl
\gvac{3} \gibr \grcm \gnl
\gvac{3} \gev \gibr \gnl
\gvac{1} \gcn{1}{2}{1}{4} \gcn{1}{1}{1}{3} \gvac{1} \gcn{1}{1}{3}{1} \gcn{1}{2}{3}{0} \gnl
\gvac{3} \gev \gnl
\gvac{3} \hspace{-0,2cm}\gwmu{3} \gnl
\gvac{4} \gob{1}{T}
\gend} \hspace{-0,36cm}\stackrel{(\ref{unitH*})}{=}\hspace{-0,2cm}
\scalebox{0.85}[0.85]{
\gbeg{7}{9}
\got{1}{T} \got{1}{} \got{4}{[(T\# H^*)\# H]^{T\# H^*}}  \gnl
\gcl{1} \glmpb \gnot{j}\gcmptb\grmpb \gnl
\gbr \gcl{1} \gcu{1} \gnl
\gcl{3} \gbr  \gnl
\gvac{1} \gcl{2} \grcm  \gnl
\gvac{2} \gibr \gnl
\gcn{1}{1}{1}{2} \gev \gcn{1}{1}{1}{0} \gnl
\gvac{1} \hspace{-0,34cm} \gwmu{3} \gnl
\gvac{2} \gob{1}{T}
\gend}
\end{eqnarray}
and similarly
$$\scalebox{0.85}[0.85]{
\gbeg{7}{7}
\got{1}{} \got{1}{H^*} \got{1}{} \got{4}{[(T\# H^*)\# H]^{T\# H^*}}  \gnl
\gu{1} \gcl{3} \glmpb \gnot{j}\gcmptb\grmpb \gnl
\gcl{2} \gcl{2} \grcm \gcn{1}{1}{-1}{1} \gcn{1}{1}{-1}{1} \gnl
\gvac{2} \gibr \gcu{1} \gcu{1} \gnl
\gcn{1}{1}{1}{2} \gev \gcn{1}{1}{1}{0} \gnl
\gvac{1} \hspace{-0,34cm} \gwmu{3} \gnl
\gvac{2} \gob{1}{T}
\gend} \stackrel{(\ref{surjec-equal})}{=} \scalebox{0.85}[0.75]{
\gbeg{7}{14}
\got{1}{} \got{1}{H^*} \got{1}{} \got{4}{[(T\# H^*)\# H]^{T\# H^*}}  \gnl
\gu{1} \gcl{1} \gvac{1} \glmpb \gnot{j}\gcmptb\grmpb \gnl
\gcn{1}{1}{1}{3} \gcn{1}{1}{1}{3} \gvac{1} \gcl{1} \gcl{2} \gcl{3} \gnl
\gvac{1} \gcl{1} \gbr  \gnl
\gvac{1} \gbr \gbr  \gnl
\gvac{1} \gcl{5} \gbr \gbr \gnl
\gvac{2} \gcl{4} \gbr \gcl{1} \gnl
\gvac{3} \gcl{1} \gibr \gnl
\gvac{3} \gibr \grcm \gnl
\gvac{3} \gev \gibr \gnl
\gvac{1} \gcn{1}{2}{1}{4} \gcn{1}{1}{1}{3} \gvac{1} \gcn{1}{1}{3}{1} \gcn{1}{2}{3}{0} \gnl
\gvac{3} \gev \gnl
\gvac{3} \hspace{-0,2cm}\gwmu{3} \gnl
\gvac{4} \gob{1}{T}
\gend}
$$
which is equivalent to
\begin{eqnarray} \label{surjec-equal-3}
\scalebox{0.9}[0.9]{
\gbeg{6}{6}
\got{1}{H^*} \got{1}{} \got{4}{[(T\# H^*)\# H]^{T\# H^*}}  \gnl
\gcl{3} \glmpb \gnot{j}\gcmptb\grmpb \gnl
\gvac{1} \grcm \gcn{1}{1}{-1}{1} \gcn{1}{1}{-1}{1} \gnl
\gvac{1} \gibr \gcu{1} \gcu{1} \gnl
\gev \gcl{1} \gnl
\gvac{2} \gob{1}{T}
\gend} = \scalebox{0.9}[0.9]{
\gbeg{5}{8}
\got{1}{H^*} \got{1}{} \got{4}{[(T\# H^*)\# H]^{T\# H^*}}  \gnl
\gcl{1} \glmpb \gnot{j}\gcmptb\grmpb \gnl
\gbr \gcu{1} \gcl{2} \gnl
\gcl{4} \gcn{1}{1}{1}{3} \gnl
\gvac{2} \gbr \gnl
\gvac{2} \gibr \gnl
\gvac{2} \gev \gnl
\gob{1}{T}
\gend} \stackrel{nat.}{=}
\scalebox{0.85}[0.85]{
\gbeg{6}{6}
\got{1}{H^*} \got{1}{} \got{4}{[(T\# H^*)\# H]^{T\# H^*}}  \gnl
\gcl{1} \glmpb \gnot{j}\gcmptb\grmpb \gnl
\gbr \gcu{1} \gcl{2} \gnl
\gcl{2} \gcn{1}{1}{1}{3} \gnl
\gvac{2} \gev \gnl
\gob{1}{T} \gvac{2} \gob{1}{}
\gend}
\end{eqnarray}

Now we are going to prove that $\gamma: [(T\# H^*)\# H]^{T\# H^*} \to T$ is
an algebra morphism. With the same notation as above, $\sigma=
T\ot\Epsilon_{H^*}\ot\Epsilon_H$ and $A=T\# H^*$, we have that
$\gamma$ will be multiplicative if
$$\scalebox{0.9}[0.9]{
\gbeg{5}{5}
\got{1}{(A\# H)^A} \got{1}{} \gvac{1} \got{1}{(A\# H)^A}  \gnl
\gwmu{3} \gnl
\gvac{1} \gbmp{j} \gnl
\gvac{1} \gbmp{\sigma} \gnl
\gvac{1} \gob{1}{T}
\gend} = \scalebox{0.85}[0.85]{
\gbeg{4}{5}
\gvac{1} \got{1}{(A\# H)^A} \got{1}{} \gvac{1} \got{1}{(A\# H)^A}  \gnl
\gvac{1} \gbmp{j} \gvac{1} \gbmp{j} \gnl
\gvac{1} \gbmp{\sigma} \gvac{1} \gbmp{\sigma} \gnl
\gvac{1} \gwmu{3} \gnl
\gvac{2} \gob{1}{T}
\gend}
$$
Since $j$ is an algebra morphism (Lemma \ref{(AsmashH)A-alg}), this amounts to
$$\scalebox{0.85}[0.85]{
\gbeg{5}{5}
\got{1}{(A\# H)^A} \got{1}{} \gvac{1} \got{1}{(A\# H)^A}  \gnl
\gbmp{j} \gvac{1} \gbmp{j} \gnl
\gwmu{3} \gnl
\gvac{1} \gbmp{\sigma} \gnl
\gvac{1} \gob{1}{T}
\gend} = \scalebox{0.85}[0.85]{
\gbeg{4}{5}
\gvac{1} \got{1}{(A\# H)^A} \got{1}{} \gvac{1} \got{1}{(A\# H)^A}  \gnl
\gvac{1} \gbmp{j} \gvac{1} \gbmp{j} \gnl
\gvac{1} \gbmp{\sigma} \gvac{1} \gbmp{\sigma} \gnl
\gvac{1} \gwmu{3} \gnl
\gvac{2} \gob{1}{T}
\gend}
$$
It would be satisfied if $\sigma$ were multiplicative. However, this does
not seem to be the case. This is why we make a more elaborate computation.
We substitute back $A=T\# H^*$, then the last equation that is to prove becomes
$$\scalebox{0.85}[0.75]{
\gbeg{8}{17}
\got{4}{[(T\# H^*)\# H]^{T\# H^*}} \gvac{2} \got{4}{[(T\# H^*)\# H]^{T\# H^*}} \gnl
\glmpb \gnot{j} \gcmptb\grmpb \gvac{2} \glmpb \gnot{j}\gcmptb\grmpb \gnl
\gcl{6} \gcl{5} \gcn{1}{1}{1}{3} \gvac{2} \gcl{2} \gcl{3} \gcl{4} \gnl
\gvac{2} \gwcm{3} \gnl
\gvac{2} \gcl{2} \gvac{1} \gbr  \gnl
\gvac{4} \gcl{1} \gbr  \gnl
\gvac{2} \gcl{1} \gcn{1}{1}{3}{2} \gvac{1} \hspace{-0,34cm} \gcmu \hspace{-0,22cm} \gmu \gnl
\gvac{2} \hspace{-0,2cm} \gcmu \gcn{1}{1}{0}{1} \gibr \gcl{5} \gcu{1} \gnl
\gvac{2} \hspace{-0,34cm} \gcl{3} \hspace{-0,22cm} \gcl{2} \gcl{1} \gibr \gcl{1} \gnl
\gvac{4} \gcn{1}{2}{1}{3} \gev \gcn{1}{2}{1}{-1} \gnl
\gvac{3}  \gcn{1}{2}{1}{3} \gnl
\gvac{3} \gcn{1}{1}{0}{1} \gvac{1} \gbr \gnl
\gvac{3} \gcl{2} \gcl{2} \grcm \hspace{-0,4cm} \gwmu{3} \gnl
\gvac{6} \gibr \gcu{1} \gnl
\gvac{4} \gcn{1}{1}{1}{2} \gev \gcn{1}{1}{1}{0} \gnl
\gvac{5} \hspace{-0,34cm} \gwmu{3} \gnl
\gvac{6} \gob{1}{T}
\gend} = \scalebox{0.85}[0.85]{
\gbeg{9}{6}
\gvac{1} \got{4}{[(T\# H^*)\# H]^{T\# H^*}} \gvac{2} \got{4}{[(T\# H^*)\# H]^{T\# H^*}} \gnl
\gvac{1} \glmpb \gnot{j} \gcmptb\grmpb \gvac{2} \glmpb \gnot{j}\gcmptb\grmpb \gnl
\gvac{1} \gcn{1}{2}{1}{3} \gcu{1} \gcu{1} \gvac{1} \gcn{1}{2}{3}{1} \gvac{1} \gcu{1} \gcu{1} \gnl
\gnl
\gvac{2} \gwmu{4} \gnl
\gvac{3} \gob{2}{T}
\gend}$$
We denote the left-hand side by $\Sigma$ and the right one by $\Omega$. Using the properties of the counits,
we can rewrite and further develop $\Sigma$ as follows:
$$\Sigma=\scalebox{0.9}[0.9]{
\gbeg{10}{9}
\got{4}{[(T\# H^*)\# H]^{T\# H^*}} \gvac{2} \got{4}{[(T\# H^*)\# H]^{T\# H^*}} \gnl
\glmpb \gnot{j} \gcmptb\grmpb \gvac{1} \glmpb \gnot{j}\gcmptb\grmpb \gnl
\gcl{2} \gcl{1} \gcn{1}{1}{1}{3} \gvac{1} \gibr \gcn{1}{1}{1}{3} \gnl
\gvac{1} \gcn{1}{2}{1}{3} \gvac{1} \gibr \grcm \gcu{1} \gnl
\gcn{1}{1}{1}{3} \gvac{2} \gev \gibr \gnl
\gvac{1} \gcn{1}{2}{1}{4} \gcn{1}{1}{1}{3} \gvac{1} \gcn{1}{1}{3}{1} \gcn{1}{2}{3}{0} \gnl
\gvac{3} \gev \gnl
\gvac{3} \hspace{-0,2cm} \gwmu{3} \gnl
\gvac{4} \gob{1}{T}
\gend} \stackrel{(\ref{surjec-equal-1})}{=} \scalebox{0.85}[0.75]{
\gbeg{8}{11}
\gvac{1} \got{4}{[(T\# H^*)\# H]^{T\# H^*}} \gvac{2} \got{4}{[(T\# H^*)\# H]^{T\# H^*}} \gnl
\gvac{1} \glmpb \gnot{j} \gcmptb\grmpb \gvac{1} \glmpb \gnot{j}\gcmptb\grmpb \gnl
\gvac{1} \gcn{1}{2}{1}{3} \gcu{1} \gcu{1} \gcn{1}{2}{3}{-1} \gcn{1}{3}{3}{-1} \gvac{1} \gcu{1} \gnl
\gnl \gvac{2} \gibr \gnl
\gvac{2} \gcl{3} \gibr \gnl
\gvac{3} \gcl{2} \grcm \gnl
\gvac{4} \gibr \gnl
\gvac{2} \gcn{1}{1}{1}{2} \gev \gcn{1}{1}{1}{0} \gnl
\gvac{3} \hspace{-0,35cm} \gwmu{3} \gnl
\gvac{4} \gob{1}{T}
\gend} \stackrel{(\ref{surjec-equal-2})}{=} \Omega.
$$
In the last equation we applied that $\Phi_{T,X}=\Phi^{-1}_{X,T}$ for any $H$-Galois object $T$ and any
$X \in \C$. The morphism $\gamma$ will be compatible with unit if $\gamma\eta_{M^A}=\sigma j\eta_{M^A}=\eta_M$. But from  Lemma \ref{(AsmashH)A-alg} we know that $j\eta_{M^A}=\eta_M$. With $M=(T\# H^*)\# H$ it is clear that $\sigma\eta_M=\eta_T$. With this we have proved that $\gamma$ is an algebra morphism. \par \medskip

We prove that $\gamma$ is right $H$-colinear by showing that it is left $H^*$-linear
(recall Proposition \ref{Hmod-H*comod}). We consider $T$ as a left $H^*$-module by Diagram \ref{H*mod-S}.
Recall that $A\# H$ is a right $H$-comodule via $A\ot\Delta_H$ where $A=T\# H^*$. Then it is a left $H^*$-module with the structure given in Diagram \ref{H*mod-S}. Analogously as in Lemma \ref{M^A-comod}, $(A\# H)^A$ inherits its left
$H^*$-module structure from $A\# H$, which satisfies
\begin{equation} \label{H*mod-Asm.H}
\scalebox{0.9}[0.9]{
\gbeg{3}{5}
\got{1}{H^*} \got{1}{} \got{1}{(A\# H)^A} \gnl
\gcn{2}{1}{1}{3} \gcl{1} \gnl
\gvac{1} \glm \gvac{1} \gnl
\gvac{2} \gbmp{j} \gnl
\gvac{2} \gob{1}{A\# H}
\gend} = \scalebox{0.9}[0.9]{
\gbeg{4}{7}
\got{1}{H^*} \got{1}{} \got{2}{(A\# H)^A} \gnl
\gcl{4} \gvac{1} \hspace{-0,36cm} \glmpb \gnot{\hspace{-0,4cm}j} \grmptb \gnl
\gvac{2} \hspace{-0,21cm} \gcn{1}{1}{2}{1} \gcmu \gnl
\gvac{2} \gcl{1} \gibr \gnl
\gvac{2} \gibr \gcl{2} \gnl
\gvac{1} \gev \gcl{1} \gnl
\gvac{3} \gob{1}{A} \gob{1}{H}
\gend}
\end{equation}
Putting back $A=T\# H^*$ and using $\sigma=T\ot\Epsilon_{H^*}\ot\Epsilon_H$,
then $\gamma$ will be left $H^*$-linear if we show
$$\scalebox{0.9}[0.9]{
\gbeg{6}{6}
\got{1}{H^*} \got{1}{} \got{4}{[(T\# H^*)\# H]^{T\# H^*}}  \gnl
\gcn{2}{1}{1}{5} \gvac{1} \gcl{1} \gnl
\gvac{2} \glm \gvac{1} \gnl
\gvac{3} \gbmp{j} \gnl
\gvac{3} \gbmp{\sigma} \gnl
\gvac{3} \gob{1}{T}
\gend} = \scalebox{0.9}[0.9]{
\gbeg{4}{6}
\got{1}{H^*} \got{1}{} \got{4}{[(T\# H^*)\# H]^{T\# H^*}}  \gnl
\gcl{1} \gvac{2} \gbmp{j} \gnl
\gcn{2}{2}{1}{5} \gvac{1} \gbmp{\sigma} \gnl
\gvac{3} \gcl{1} \gnl
\gvac{2} \glm \gnl
\gvac{3} \gob{1}{T}
\gend}
$$
Again, $\sigma$ is not left $H^*$-linear, so we have to do a more involved computation.
Applying Diagram \ref{H*mod-Asm.H}, the definition of $\sigma$ and the
$H^*$-module structure of $T$, we get that the above question becomes
$$\scalebox{0.9}[0.9]{
\gbeg{5}{8}
\got{1}{H^*} \got{1}{} \got{4}{[(T\# H^*)\# H]^{T\# H^*}} \gnl
\gcl{5} \gvac{1} \hspace{-0,34cm} \glmpb \gnot{j} \gcmptb \grmpb \gnl
\gvac{2} \hspace{-0,21cm} \gcn{1}{1}{2}{1} \gcn{1}{1}{2}{1} \gcmu \gnl
\gvac{2} \gcl{2} \gcl{1} \gibr \gnl
\gvac{3} \gibr \gcl{2} \gnl
\gvac{2} \gibr \gcl{1} \gnl
\gvac{1} \gev \gcl{1} \gcu{1} \gcu{1} \gnl
\gvac{3} \gob{1}{T}
\gend} = \scalebox{0.9}[0.9]{
\gbeg{4}{7}
\got{1}{H^*} \got{1}{} \got{4}{[(T\# H^*)\# H]^{T\# H^*}}  \gnl
\gcl{1} \gvac{1} \glmpb \gnot{j} \gcmptb \grmpb \gnl
\gcn{1}{1}{1}{3} \gvac{1} \gcl{1} \gcu{1} \gcu{1} \gnl
\gvac{1} \gcl{2} \grcm \gnl
\gvac{2} \gibr \gnl
\gvac{1} \gev \gcl{1} \gnl
\gvac{3} \gob{1}{T}
\gend}
$$
But this is fulfilled because of (\ref{surjec-equal-3}). Hence $\gamma:
[(T\# H^*)\# H]^{T\# H^*}\to T$ is right $H$-colinear and summing up it
is an isomorphism of $H$-Galois objects. Recall from 2.1.8 that since $H$ is finite, it is also flat. Then we have finally established:

\begin{proposition}\label{surjective}
Let $\C$ be a closed braided monoidal category with equalizers and coequalizers.
Let $H$ be a finite and commutative Hopf algebra. Suppose that the braiding $\Phi$
is $H$-linear and that $\Phi_{T,X}=\Phi^{-1}_{X,T}$ for any $H$-Galois object $T$ and any
$X \in \C$.  Then the map $\Upsilon: \BM(\C; H)\to\Gal(\C; H), [A] \mapsto [(A\# H)^A]$
is surjective.
\end{proposition}

\subsection{The split exact sequence}

Consider the sequence
$$
\bfig
\putmorphism(0, 0)(1, 0)[1`\Br(\C)`]{400}1a
\putmorphism(400, 0)(1, 0)[\phantom{\Br(\C)}`\BM(\C;H)` q]{600}1a
\putmorphism(410, -50)(1, 0)[\phantom{\Br(\C)}`\phantom{\Br(\C)}`p]{470}{-1}b
\putmorphism(1000, 0)(1, 0)[\phantom{\BM(\C;H)}`\Gal(\C;H)`\Upsilon]{700}1a
\putmorphism(1700, 0)(1, 0)[\phantom{\Gal(\C;H)}`1.`]{400}1a
\efig
$$
Recall that the map $q$ sends $[A] \in \Br(\C)$ to $[A]$ in $\BM(\C;H)$, where $A$ is equipped with the trivial $H$-module structure. The map $p: \BM(\C;H) \to \Br(\C)$, induced by forgetting the $H$-module structure of an $H$-Azumaya algebra, satisfies $pq=id_{\Br(\C)}$. Thus the sequence is split and we know from the previous subsection that $\Upsilon$ is surjective. We prove exactness at $\BM(\C;H)$. \par \medskip

We show that $\Bi(q) \subseteq \Ker(\Upsilon)$. For an Azumaya algebra $A\in\C$ with trivial $H$-module structure it is $A\# H=A\ot H$ as algebras. By the equivalence of categories given by the pair of functors
$(A\ot -, (-)^A)$ we get that $(A\ot H)^A\iso H$ as objects in $\C$.
In order to prove that this is an isomorphism of $H$-Galois objects we prove that
it is a right $H$-comodule algebra morphism, Proposition \ref{comod-alg-iso}. Observe Diagram \ref{counitazu}
of Lemma \ref{Az-counit-inverse}, defining the unit $\zeta:H\to (A\ot H)^A$ of the above
adjunction, with $M=H$. The morphism $\eta_A\ot H$ is clearly a right $H$-comodule
algebra one, as so is $j_{A\ot H}$, by Corollary \ref{(Asm.H)A+j com.alg.}. Now by 2.2.3(a),
$\zeta$ is such a morphism as well. Thus $\Upsilon([A])=[(A\# H)^A]=[(A\ot H)^A]=[H]$, so $[A]\in \Ker(\Upsilon)$.
\par \smallskip

Suppose that $\Bi(q) \subsetneq \Ker(\Upsilon)$. Since $p(\Ker(\Upsilon))\subseteq \Br(\C)$, there exists an $H$-Azumaya algebra $A$ which determines two different classes $[A_1]$ and $[A_2]$ in $\Ker(\Upsilon)$. This algebra has two different $H$-module structures in $\BM(\C;H)$, a non-trivial one and a trivial one. As the underlying Azumaya algebra is the same, it is $p([A_1])=p([A_2])$. We will show that $p$ is injective when restricted to $\Ker(\Upsilon)$, reaching a contradiction with $[A_1]\not=[A_2]$ in $\Ker(\Upsilon)$. This will prove that $\Bi(q)=\Ker(\Upsilon)$. \par \smallskip

We proceed to prove that the morphism $p: \BM(\C;H) \to \Br(\C)$ restricted to $\Ker(\Upsilon)$ is injective.
Suppose that $p([A])$ is trivial in $\Br(\C)$, for $[A] \in \Ker(\Upsilon)$. Then there is an algebra isomorphism $\delta: A\to [P,P]$, for some faithfully projective object $P\in\C$. On the other hand, we have an $H$-comodule algebra isomorphism $\omega : H \to (A\# H)^A$. Consider the equalizer algebra morphism $j: (A\# H)^A \to A\# H$.
The composition $(A\# \Epsilon_H)j:(A\# H)^A \to A$ is denoted by $\epsilon$.
Then $\epsilon\omega: H\to A$ is an algebra morphism. We are going to define an $H$-module structure on $P$. This will induce an $H$-module algebra structure on $[P,P]$ by Remark \ref{Hmod-innerhom}. Then we will
prove that $\delta: A\to [P,P]$ is a morphism of $H$-module algebras, thus $A$ will become trivial in $\BM(\C;H)$ and we will have the claim.

\begin{lemma}
Let $\xi:=\delta\epsilon\omega: H\to [P, P]$. Then the following morphism defines a left $H$-module structure on $P$:
$$\scalebox{0.9}[0.9]{
\gbeg{2}{3}
\got{1}{H} \got{1}{P} \gnl
\glm \gnl
\gvac{1} \gob{1}{P}
\gend} := \scalebox{0.9}[0.9]{
\gbeg{3}{6}
\got{1}{H} \gvac{1} \got{1}{P} \gnl
\gmp{S} \gvac{1} \gcl{2}\gnl
\gbmp{\xi} \gnl
\gcn{1}{1}{1}{3}\gcn{1}{1}{3}{1} \gnl
\gvac{1} \gmp{\ev} \gnl
\gvac{1} \gob{1}{}
\gend}
$$
\end{lemma}

\begin{proof} The compatibility with the multiplication in $H$ is easily proved using that $S$ is an algebra antimorphism, $H$ is commutative, $\xi$ is an algebra morphism and the definition of the multiplication in $[P, P]$. For the compatibility with the unit use furthermore the definition of the unit of $[P, P]$ (2.1.3).
\qed\end{proof}

Let $\Fi: H\to [P, P]$ denote the algebra morphism from Lemma \ref{mod vs algEnd}(i) stemming
from the left $H$-module structure of $P$ shown in the above lemma. Since $H$ is (co)commutative, $S^2=id_H$. This implies further:
$$\scalebox{0.85}[0.85]{
\gbeg{3}{6}
\got{1}{H} \gvac{1} \got{1}{P} \gnl
\gmp{S} \gvac{1} \gcl{2}\gnl
\gbmp{\Fi} \gnl
\gcn{1}{1}{1}{3}\gcn{1}{1}{3}{1} \gnl
\gvac{1} \gmp{\ev} \gnl
\gvac{1} \gob{1}{P}
\gend} \stackrel{Lem. \ref{mod vs algEnd}}{=}
\scalebox{0.85}[0.85]{
\gbeg{2}{4}
\got{1}{H} \got{1}{P} \gnl
\gmp{S} \gcl{1} \gnl
\glm \gnl
\gvac{1} \gob{1}{P}
\gend} = \scalebox{0.85}[0.85]{
\gbeg{3}{7}
\got{1}{H} \gvac{1} \got{1}{P} \gnl
\gmp{S} \gvac{1} \gcl{3} \gnl
\gmp{S} \gnl
\gbmp{\xi} \gnl
\gcn{1}{1}{1}{3}\gcn{1}{1}{3}{1} \gnl
\gvac{1} \gmp{\ev} \gnl
\gvac{1} \gob{1}{P}
\gend} = \scalebox{0.85}[0.85]{
\gbeg{3}{5}
\got{1}{H} \gvac{1} \got{1}{P} \gnl
\gbmp{\xi} \gvac{1} \gcl{1} \gnl
\gcn{1}{1}{1}{3}\gcn{1}{1}{3}{1} \gnl
\gvac{1} \gmp{\ev} \gnl
\gvac{1} \gob{1}{P}
\gend}
$$
This yields
\begin{eqnarray}\label{Fi-comp-S}
\Fi S=\xi.
\end{eqnarray}
Recall from Lemma \ref{mod vs algEnd}(ii) that the $H$-module structure on $P$ induces an
$H$-module structure on $[P, P]$ making it into an $H$-module algebra. This structure
was given by:
$$\scalebox{0.85}[0.85]{
\gbeg{4}{4}
\got{1}{H} \got{3}{[P, P]} \gnl
\gcn{1}{1}{1}{3} \gvac{1} \gcl{1} \gnl
\gvac{1} \glm \gnl
\gvac{2}\gob{1}{[P, P]}
\gend} = \scalebox{0.85}[0.85]{
\gbeg{6}{8}
\gvac{1} \got{2}{H}\got{1}{[P, P]} \gnl
\gvac{1} \gcmu  \gcl{2} \gvac{1} \gnl
\gvac{1} \gcl{1} \gmp{S} \gnl
\gvac{1} \gcl{1} \gbr  \gnl
\gvac{1} \gbmp{\Fi} \gcl{1} \gbmp{\Fi}  \gnl
\gvac{1} \gcn{1}{1}{1}{0} \gmu \gnl
\gvac{1} \hspace{-0,2cm}\gwmu{3} \gnl
\gvac{2} \gob{1}{[P, P]}
\gend}
$$

\begin{lemma}
With the above $H$-module structure on $[P,P]$, the morphism $\delta: A\to [P,P]$ is left $H$-linear.
\end{lemma}

\begin{proof}
We know from Corollary \ref{(Asm.H)A+j com.alg.} that the morphism $j:(A\# H)^A\to A\# H$
is right $H$-colinear. Using the equalizer property of $((A\# H)^A, j)$ and writing out
$\epsilon=(A\#\Epsilon_H)j$, we deduce
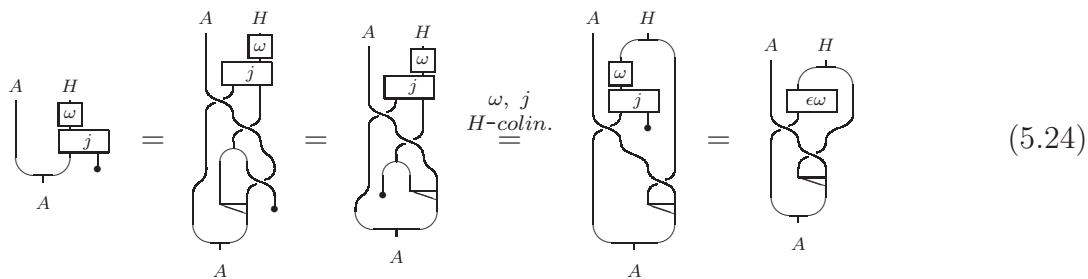
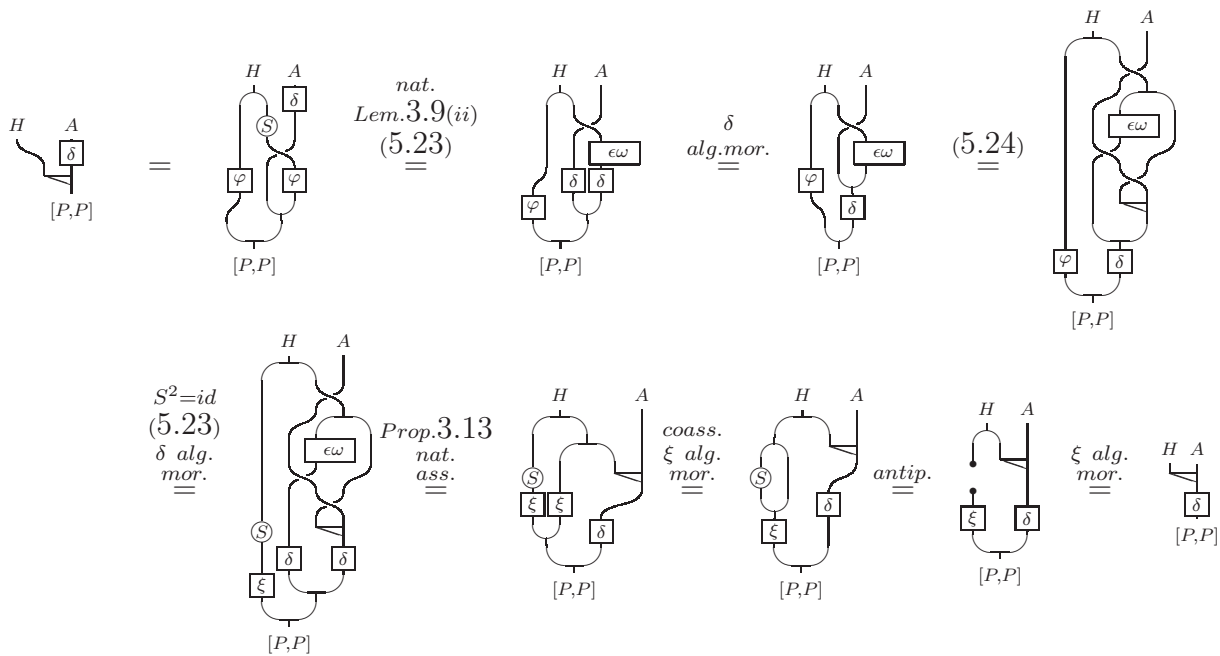
\begin{eqnarray} \label{split-omega}
\hspace{-0,3cm}
\scalebox{0.85}[0.85]{
\gbeg{5}{5}
\got{1}{A} \gvac{1} \got{1}{H} \gnl
\gcl{2} \gvac{1} \gbmp{\omega} \gnl
\gvac{2} \glmptb\gnot{\hspace{-0,3cm}j}\grmpb  \gnl
\gvac{1} \hspace{-0,42cm} \gwmu{3} \gcu{1} \gnl
\gvac{2} \gob{1}{A}
\gend} \hspace{-0,2cm}= \scalebox{0.85}[0.85]{
\gbeg{4}{10}
\got{1}{A} \gvac{1} \got{1}{H} \gnl
\gcl{2} \gvac{1} \gbmp{\omega} \gnl
\gvac{1} \glmpb\gnot{\hspace{-0,3cm}j}\grmptb  \gnl
\gbr \gcl{1} \gnl
\gcl{2} \gbr  \gnl
\gvac{1} \hspace{-0,2cm} \gcmu \gcn{1}{1}{0}{1} \gnl
\gcn{1}{1}{2}{1} \gcl{1} \gbr \gnl
\gcl{1}\glm \gcu{1} \gnl
\gwmu{3} \gnl
\gvac{1} \gob{1}{A}
\gend} \hspace{-0,2cm} =
\scalebox{0.85}[0.85]{
\gbeg{4}{9}
\got{1}{A} \gvac{1} \got{1}{H} \gnl
\gcl{2} \gvac{1} \gbmp{\omega} \gnl
\gvac{1} \glmpb\gnot{\hspace{-0,3cm}j}\grmptb  \gnl
\gbr \gcl{1} \gnl
\gcl{2} \gbr \gnl
\gvac{1} \hspace{-0,2cm} \gcmu \gcn{1}{1}{0}{1} \gnl
\gcn{1}{1}{2}{1} \gcu{1} \glm \gnl
\gwmu{4} \gnl
\gvac{1} \gob{2}{A}
\gend} \hspace{-0,2cm}\stackrel{\omega, \ j}{\stackrel{H\x colin.}{=}}
\scalebox{0.85}[0.85]{
\gbeg{4}{10}
\got{1}{A} \gvac{1} \got{1}{H} \gnl
\gcl{3} \gwcm{3} \gnl
\gvac{1} \gbmp{\omega} \gvac{1} \gcl{4} \gnl
\gvac{1} \glmptb\gnot{\hspace{-0,3cm}j}\grmpb  \gnl
\gbr \gcu{1} \gnl
\gcl{3} \gcn{1}{1}{1}{3} \gnl
\gvac{2} \gbr \gnl
\gvac{2} \glm \gnl
\gwmu{4} \gnl
\gvac{1} \gob{2}{A}
\gend} = \scalebox{0.85}[0.85]{
\gbeg{4}{8}
\got{1}{A} \gvac{1} \got{1}{H} \gnl
\gcl{2} \gwcm{3} \gnl
\gvac{1} \glmptb\gnot{\hspace{-0,3cm}\epsilon\omega}\grmp \gcl{1} \gnl
\gbr \gcn{1}{1}{3}{1} \gnl
\gcl{2} \gbr \gnl
\gvac{1}  \glm \gnl
\gwmu{3} \gnl
\gvac{1} \gob{1}{A}
\gend} \hspace{0,3cm}
\end{eqnarray}
We have:
$$\begin{array}{ll}
\scalebox{0.85}[0.85]{
\gbeg{4}{4}
\got{1}{H} \gvac{1} \got{1}{A} \gnl
\gcn{1}{1}{1}{3} \gvac{1} \gbmp{\delta} \gnl
\gvac{1} \glm \gnl
\gvac{2}\gob{1}{[P, P]}
\gend} & =\scalebox{0.85}[0.85]{
\gbeg{5}{8}
\gvac{1} \got{2}{H}\got{1}{A} \gnl
\gvac{1} \gcmu \gbmp{\delta} \gnl
\gvac{1} \gcl{1} \gmp{S} \gcl{1} \gnl
\gvac{1} \gcl{1} \gbr  \gnl
\gvac{1} \gbmp{\Fi} \gcl{1} \gbmp{\Fi}  \gnl
\gvac{1} \gcn{1}{1}{1}{0} \gmu \gnl
\gvac{1} \hspace{-0,2cm}\gwmu{3} \gnl
\gvac{2} \gob{1}{[P, P]}
\gend} \stackrel{nat.}{\stackrel{Lem. \ref{mod vs algEnd}(ii)}{\stackrel{(\ref{Fi-comp-S})}{=}}}
\scalebox{0.85}[0.85]{
\gbeg{6}{8}
\gvac{1} \got{2}{H}\got{1}{A} \gnl
\gvac{1} \gcmu \gcl{1} \gnl
\gvac{1} \gcl{2} \gbr \gnl
\gvac{2} \gcl{1} \glmptb\gnot{\hspace{-0,3cm}\epsilon\omega}\grmp \gnl
\gvac{1} \gcn{1}{1}{1}{0} \gbmp{\delta} \gbmp{\delta} \gnl
\gvac{1} \hspace{-0,34cm} \gbmp{\Fi} \gvac{1} \hspace{-0,36cm} \gmu \gnl
\gvac{2} \hspace{-0,2cm} \gwmu{3} \gnl
\gvac{3} \gob{1}{[P, P]}
\gend} \stackrel{\delta}{\stackrel{alg. mor.}{=}}
\scalebox{0.85}[0.85]{
\gbeg{5}{8}
\got{2}{H}\got{1}{A} \gnl
\gcmu \gcl{1} \gnl
\gcl{2} \gbr \gnl
\gvac{1} \gcl{1} \glmptb\gnot{\hspace{-0,3cm}\epsilon\omega}\grmp \gnl
\gbmp{\Fi} \gmu \gnl
\gcn{1}{1}{1}{2} \gvac{1} \hspace{-0,34cm} \gbmp{\delta}  \gnl
\gvac{1} \gmu \gnl
\gvac{1} \gob{2}{[P, P]}
\gend} \stackrel{(\ref{split-omega})}{=}
\scalebox{0.85}[0.85]{
\gbeg{5}{12}
\gvac{1} \got{1}{H} \gvac{1} \got{1}{A} \gnl
\gwcm{3} \gcl{1} \gnl
\gcl{7} \gvac{1} \gbr \gnl
\gvac{1} \gcn{1}{1}{3}{1} \gwcm{3} \gnl
\gvac{1} \gcl{1} \glmptb\gnot{\hspace{-0,3cm}\epsilon\omega}\grmp \gcl{1} \gnl
\gvac{1} \gbr \gcn{1}{1}{3}{1} \gnl
\gvac{1} \gcl{2} \gbr \gnl
\gvac{2} \glm \gnl
\gvac{1} \gwmu{3} \gnl
\gbmp{\Fi} \gvac{1} \gbmp{\delta} \gnl
\gwmu{3} \gnl
\gvac{1} \gob{1}{[P, P]}
\gend} \\
 & \stackrel{S^2=id}{\stackrel{(\ref{Fi-comp-S})}{\stackrel{\delta\ alg.}{\stackrel{mor.}{=}}}}
\scalebox{0.85}[0.85]{
\gbeg{5}{12}
\gvac{1} \got{1}{H} \gvac{1} \got{1}{A} \gnl
\gwcm{3} \gcl{1} \gnl
\gcl{5} \gvac{1} \gbr \gnl
\gvac{1} \gcn{1}{1}{3}{1} \gwcm{3} \gnl
\gvac{1} \gcl{1} \glmptb\gnot{\hspace{-0,3cm}\epsilon\omega}\grmp \gcl{1} \gnl
\gvac{1} \gbr \gcn{1}{1}{3}{1} \gnl
\gvac{1} \gcl{2} \gbr \gnl
\gmp{S} \gvac{1} \glm \gnl
\gcl{1} \gbmp{\delta} \gvac{1} \gbmp{\delta} \gnl
\gbmp{\xi} \gwmu{3} \gnl
\gwmu{3} \gnl
\gvac{1} \gob{1}{[P, P]}
\gend}
\hspace{-0,3cm}\stackrel{Prop. \ref{braid-lin}}{\stackrel{nat.}{\stackrel{ass.}{=}}}
\scalebox{0.85}[0.85]{
\gbeg{5}{8}
\gvac{1} \got{1}{H} \gvac{2} \got{1}{A} \gnl
\gwcm{3} \gvac{1} \gcl{2} \gnl
\gcl{1} \gwcm{3} \gnl
\gmp{S} \gcl{1} \gvac{1} \glm \gnl
\gbmp{\xi} \gbmp{\xi} \gvac{1} \hspace{-0,14cm} \gcn{2}{1}{3}{0} \gnl
\gmu \gvac{1} \hspace{-0,2cm} \gbmp{\delta} \gnl
\gvac{1} \hspace{-0,14cm} \gwmu{3} \gnl
\gvac{2} \gob{1}{[P, P]}
\gend} \hspace{-0,14cm}\stackrel{coass.}{\stackrel{\xi \ alg.}{\stackrel{mor.}{=}}}\hspace{-0,2cm}
\scalebox{0.85}[0.85]{
\gbeg{5}{8}
\gvac{2} \got{1}{H} \gvac{1} \got{1}{A} \gnl
\gvac{1} \gwcm{3} \gcl{1} \gnl
\gvac{1} \hspace{-0,36cm} \gcmu \gvac{1} \hspace{-0,2cm} \glm \gnl
\gvac{2} \hspace{-0,22cm} \gmp{S} \gcl{1} \gvac{1} \gcn{1}{1}{2}{0} \gnl
\gvac{2} \gmu \gvac{1} \hspace{-0,22cm} \gbmp{\delta} \gnl
\gvac{3} \gbmp{\xi} \gvac{1} \gcl{1} \gnl
\gvac{3} \hspace{-0,14cm} \gwmu{3} \gnl
\gvac{4} \gob{1}{[P, P]}
\gend} \hspace{-0,2cm}\stackrel{antip.}{=}
\scalebox{0.85}[0.85]{
\gbeg{4}{7}
\got{2}{H} \got{1}{A} \gnl
\gcmu \gcl{1} \gnl
\gcu{1} \glm \gnl
\gu{1}  \gvac{1} \hspace{-0,14cm} \gcl{1} \gnl
\gbmp{\xi} \gvac{1} \gbmp{\delta} \gnl
\hspace{-0,14cm} \gwmu{3} \gnl
\gvac{1} \gob{1}{[P, P]}
\gend} \hspace{-0,2cm}\stackrel{\xi \ alg.}{\stackrel{mor.}{=}}
\scalebox{0.85}[0.85]{
\gbeg{4}{4}
\got{1}{H} \got{1}{A} \gnl
\glm \gnl
\gvac{1} \gbmp{\delta} \gnl
\gvac{1} \gob{1}{[P, P]}
\gend}
\end{array}
$$
\qed\end{proof}

We have done practically all the work to prove the main theorem of this paper.

\begin{theorem}\label{Beattie}
Let $\C$ be a closed braided monoidal category with equalizers and coequalizers.
Let $H$ be a finite and commutative Hopf algebra. Suppose that the braiding is
$H$-linear and that $\Phi_{T,X}=\Phi^{-1}_{X,T}$ for every $H$-Galois object $T$ and
$X \in \C$. Then there is a split exact sequence
$$
\bfig
\putmorphism(0, 0)(1, 0)[1`\Br(\C)`]{400}1a
\putmorphism(400, 0)(1, 0)[\phantom{\Br(\C)}`\BM(\C;H)`q]{600}1a
\putmorphism(410, -50)(1, 0)[\phantom{\Br(\C)}`\phantom{\Br(\C)}`p]{470}{-1}b
\putmorphism(1000, 0)(1, 0)[\phantom{\BM(\C;H)}`\Gal(\C;H)`\Upsilon]{700}1a
\putmorphism(1700, 0)(1, 0)[\phantom{\Gal(\C;H)}`1.`]{400}1a
\efig
$$
Furthermore, $\BM(\C;H)\iso \Br(\C)\times \Gal(\C;H).$
\end{theorem}

\begin{proof}
It just remains to prove the last statement. By the definition of the morphism
$q:\Br(\C)\to\BM(\C;H)$ and surjectivity of $\Upsilon:\BM(\C;H) \to \Gal(\C;H)$ we have the following situation
$$
\bfig
\putmorphism(0, 0)(1, 0)[\Br(\C)`\BM(\C;H)`q]{700}1a
\putmorphism(0, -140)(1, 0)[[A]`[A]`]{700}1b
\putmorphism(700, 0)(1, 0)[\phantom{\BM(\C;H)}`\Gal(\C;H)`\Upsilon]{800}1a
\putmorphism(700, -280)(1, 0)[[T\# H^*]`[T]`]{800}{-1}a
\efig
$$
We are going to prove that $[A]$ and $[T\# H^*]$ commute in $\BM(\C;H)$.
This will be done by showing that the braiding acting between $T\#H^*$ and $A$ is a left
$H$-module algebra morphism. Since $\Phi$ is left $H$-linear, $\Phi$ is right $H$-colinear,
by Proposition \ref{braid-lin}(iii). In view of Proposition \ref{Hmod-H*comod} this means that
$\Phi$ is left $H^*$-linear. In particular, $\Phi_{H^*, A}=
\Phi_{A, H^*}^{-1}$, because of Proposition \ref{braid-lin}. On the other hand, by hypothesis,
if $T$ is an $H$-Galois object, then $\Phi_{T, A}=\Phi_{A, T}^{-1}$. We will show that
$\Phi_{T\# H^*, A}$ is a left $H$-module algebra morphism.

Consider $(T\# H^*)\ot A$ and $A\ot (T\# H^*)$ as left $H$-modules by the codiagonal
structures. However, since $[A]\in\Br(\C)$, we have that $A$ is a trivial $H$-module,
so the above two respective $H$-module structures will be induced by the one of $T\#H^*$.
Recalling the left $H$-module structure of $T\#H^*$ from Lemma \ref{Ssm.H* left H-mod.alg.}
we find:
$$\scalebox{0.9}[0.9]{
\gbeg{4}{5}
\got{1}{H} \got{2}{T\ot H^*} \got{1}{A} \gnl
\gcn{2}{1}{1}{3} \gcl{1} \gcl{2} \gnl
\gvac{1} \glm \gvac{1} \gnl
\gvac{2} \gbr \gnl
\gvac{2} \gob{1}{A} \gob{2}{T\ot H^*}
\gend} = \scalebox{0.9}[0.9]{
\gbeg{5}{6}
\got{1}{H} \got{1}{T} \got{3}{\hspace{-0,3cm}H^*} \got{1}{\hspace{-0,8cm}A}\gnl
\gcl{1} \gcl{1} \gcmu \gcl{2} \gnl
\gcl{1} \gibr \gcl{1} \gnl
\gibr \gcl{1} \gbr\gnl
\gev \gbr \gcl{1} \gnl
\gvac{1} \gob{3}{A} \gob{1}{\hspace{-0,7cm}T} \gob{1}{\hspace{-0,5cm}H^*}
\gend} \stackrel{nat.}{=} \scalebox{0.9}[0.9]{
\gbeg{5}{8}
\got{1}{H} \got{1}{T} \got{3}{\hspace{-0,6cm}H^*} \got{1}{\hspace{-1,6cm}A}\gnl
\gcl{2} \gcl{1} \gbr \gnl
\gvac{1} \gbr \gcn{1}{1}{1}{2} \gnl
\gbr \gcl{1} \gcmu \gnl
\gcl{3} \gcl{1} \gibr \gcl{3} \gnl
\gvac{1} \gibr \gcl{2} \gnl
\gvac{1} \gev \gnl
\gob{1}{A} \gvac{2} \gob{1}{T} \gob{1}{H^*}
\gend} = \scalebox{0.9}[0.9]{
\gbeg{4}{5}
\got{1}{\hspace{-0,5cm}H} \got{1}{T\ot H^*} \got{2}{\hspace{-0,2cm}A} \gnl
\gcn{1}{1}{0}{1} \gbr \gnl
\gbr \gcl{1} \gnl
\gcl{1} \glm \gnl
\gob{1}{A} \gob{3}{T\ot H^*}
\gend}
$$
This proves that $\Phi_{T\# H^*, A}$ is left $H$-linear. That $\Phi_{T\# H^*, A}$ is compatible with multiplication follows by naturality:
$$\begin{array}{ll}
\scalebox{0.85}[0.85]{
\gbeg{9}{6}
\gvac{1} \got{1}{T\#H^*} \got{2}{\hspace{-0,2cm}A} \got{1}{T\#H^*} \got{2}{\hspace{-0,2cm}A} \gnl
\gvac{1} \gbr \gvac{1} \gbr \gnl
\gcn{1}{1}{3}{1} \gvac{1} \gcl{1} \gvac{1} \gcl{1} \gcn{1}{1}{1}{3} \gnl
\glmpt\gnot{\hspace{0,5cm}A\ot (T\#H^*)}\gcmpb\gcmpt\grmp
   \glmpt\gnot{\hspace{0,5cm}A\ot (T\#H^*)}\gcmpb\gcmpt\grmp\gnl
\gvac{1} \gwmu{5}  \gnl
\gvac{2} \gob{3}{A\ot (T\#H^*)}
\gend} & = \scalebox{0.85}[0.85]{
\gbeg{7}{10}
\got{1}{T} \got{1}{H^*} \got{1}{A} \got{1}{T} \got{1}{H^*} \got{1}{A} \gnl
\gcl{1} \gbr \gcl{1} \gbr \gnl
\gbr \gcl{1} \gbr \gcl{1} \gnl
\gcl{2} \gcl{1} \gbr \gcn{1}{2}{1}{3} \gcn{1}{2}{1}{3} \gnl
\gvac{1} \gbr \gcn{1}{1}{1}{2} \gnl
\gmu \gcl{1} \gcmu \gcl{1} \gcl{2} \gnl
\gcn{1}{1}{2}{3} \gvac{1} \gcl{1} \gcl{1} \gbr \gnl
\gvac{1} \gcl{2} \gcl{1} \glm \gmu \gnl
\gvac{2} \gwmu{3} \gcn{1}{1}{2}{1}  \gnl
\gob{3}{A} \gob{1}{T} \gvac{1} \gob{1}{H^*}
\gend} \stackrel{nat.}{=}
\scalebox{0.85}[0.85]{
\gbeg{7}{8}
\got{1}{T} \got{2}{H^*} \got{1}{A} \got{1}{T} \got{1}{H^*} \got{1}{A} \gnl
\gcl{1} \gcmu \gibr \gcl{1} \gcl{2} \gnl
\gcl{1} \gcl{1} \gbr \gibr \gnl
\gcl{1} \glm \gmu \gmu \gnl
\gwmu{3} \gcn{1}{1}{2}{3} \gvac{1} \gcn{1}{1}{2}{1} \gnl
\gvac{1} \gcn{2}{1}{1}{5} \gvac{1} \gbr \gnl
\gvac{3} \gbr \gcl{1} \gnl
\gvac{3} \gob{1}{A} \gob{1}{T} \gob{1}{H^*}
\gend} \stackrel{cond.}{\stackrel{\Phi_{T, A}}{\stackrel{\Phi_{H^*, A}}{=}}}
\scalebox{0.85}[0.85]{
\gbeg{7}{8}
\got{1}{T} \got{2}{H^*} \got{1}{A} \got{1}{T} \got{1}{H^*} \got{1}{A} \gnl
\gcl{1} \gcmu \gbr \gcl{1} \gcl{2} \gnl
\gcl{1} \gcl{1} \gbr \gbr \gnl
\gcl{1} \glm \gmu \gmu \gnl
\gwmu{3} \gcn{1}{1}{2}{3} \gvac{1} \gcn{1}{1}{2}{1} \gnl
\gvac{1} \gcn{2}{1}{1}{5} \gvac{1} \gbr \gnl
\gvac{3} \gbr \gcl{1} \gnl
\gvac{3} \gob{1}{A} \gob{1}{T} \gob{1}{H^*}
\gend} \vspace{3mm} \\
 & = \scalebox{0.85}[0.85]{
\gbeg{9}{6}
\gvac{1} \got{1}{(T\#H^*)\ot A} \gvac{3} \got{1}{(T\#H^*)\ot A} \gnl
\gvac{1} \hspace{-0,26cm} \gwmu{5}  \gnl
\gvac{2} \hspace{-0,3cm} \glmp\gnot{\hspace{0,5cm}(T\#H^*)\ot A}\gcmpb\gcmp\grmpb\gnl
\gvac{3} \gcl{1} \gcn{1}{1}{3}{1} \gnl
\gvac{3} \gbr \gnl
\gvac{3} \gob{1}{A} \gob{2}{T\#H^*}
\gend}
\end{array}$$
Obviously $\Phi_{T\# H^*, A}$ is compatible with unit, so it is an algebra morphism.
Thus we have proved that $\BM(\C;H)\iso\Br(\C)\times\Gal(\C;H)$.
\qed\end{proof}

The following result is a byproduct of the proof of Theorem \ref{Beattie}.

\begin{theorem}\label{innerbeat}
Let $\C$ be a closed braided monoidal category with equalizers and coequalizers. Let $H \in \C$ be a finite and commutative Hopf algebra. Assume that the braiding is $H$-linear. Then there is a split exact sequence
\begin{equation}\label{2ndseq}
\bfig
\putmorphism(0, 0)(1, 0)[1`\Br(\C)`]{400}1a
\putmorphism(400, 10)(1, 0)[\phantom{\Br(\C)}`\BM_{inn}(\C;H)`q]{680}1a
\putmorphism(305, -35)(-1, 0)[\phantom{\BM(R;H)}`\phantom{\Br(R)}`p]{600}{-1}b
\putmorphism(1080, 0)(1, 0)[\phantom{\BM_{inn}(\C;H)}`\Gal_{nb}(\C;H)`\Upsilon']{800}1a
\putmorphism(1900, 0)(1, 0)[\phantom{\Gal_{nb}(\C;H)}`1.`]{500}1a
\efig
\end{equation}
Moreover, $\BM_{inn}(\C;H)\iso \Br(\C)\times \Hc(\C;H,I).$
\end{theorem}

\begin{proof}
Observe that the condition $\Phi_{T,X}=\Phi^{-1}_{X,T}$ for any $H$-Galois object $T$ and any
$X \in \C$ used in the proof of Theorem \ref{Beattie} is not required here. This is because when dealing with an $H$-Galois object with normal basis $T$ this condition is automatically fulfilled. Since $T \cong H$ as right $H$-comodules, $T \cong H$ as objects in $\C$. The $H$-linearity of the braiding assures that $\Phi_{H,X}=\Phi^{-1}_{X,H}$ for any $X \in \C$. So this conditions also holds for any $T$ isomorphic to $H$. \par \smallskip

By Proposition \ref{Gal-nb-BM-inn}, the map $\Upsilon':\BM_{inn}(\C;H) \to \Gal_{nb}(\C;H),[A] \mapsto [(A \# H)^A]$ is a group morphism. For an $H$-Galois object $T$, the algebra $T\# H^*$ is an $H$-Azumaya algebra in virtue of Proposition \ref{ChaseSweedler Thm} and Lemma \ref{Ssm.H* left H-mod.alg.}. In order to establish the isomorphism $((T\# H^*)\# H)^{T\# H^*} \cong T$ in Proposition \ref{surjective} we have used the aformentioned symmetricity condition. This isomorphism then holds for any $H$-Galois object with normal basis property. The $H$-Azumaya algebra $T\# H^*$ has inner action when $T$ has a normal basis, Proposition \ref{Gal-nb-BM-inn}. Then $\Upsilon'$ is surjective and we have the split sequence (\ref{2ndseq}). \par \smallskip
Arguing as in the proof of Theorem \ref{Beattie} we may show that this sequence is exact. When proving that $[A] \in \Br(\C)$ and $[T \# H^*]$ commute in $\BM(\C; H)$ we have used that $\Phi_{T,X}=\Phi^{-1}_{X,T}$. For $T$ having normal basis property this holds, showing that $\BM_{inn}(\C;H)\iso \Br(\C)\times \Gal_{nb}(\C;H).$ Finally, apply Theorem \ref{cohpical} stating that $\Gal_{nb}(\C;H) \iso \Hc(\C;H,I).$
\qed\end{proof}

\begin{corollary}\label{bmeqbminn}
Let $\C$ be a closed braided monoidal category with equalizers and coequalizers. Let $H \in \C$ be a finite and commutative Hopf algebra. Assume that the braiding is $H$-linear and any $H$-Galois object has normal basis property.
Then $\BM_{inn}(\C;H)=\BM(\C;H).$
\end{corollary}

\begin{proof}
The hypothesis $\Phi_{T,X}=\Phi^{-1}_{X,T}$ for any $H$-Galois object $T$ in Theorem \ref{Beattie} holds because we are assuming that any $H$-Galois object has normal basis property. From this theorem, any class [A] in $\BM(\C;H)$ is of the form $[A]=[B][T \# H^*]$ where $[B] \in \Br(\C)$ and $T$ is an $H$-Galois object. By hypothesis, $T$ has normal basis property and by Proposition \ref{Gal-nb-BM-inn}, $[T \# H^*] \in \BM_{inn}(\C;H).$ Since $\Br(\C)\subset \BM_{inn}(\C;H),$ we are done.
\qed \end{proof}

We can now derive the results of Alonso \'Alvarez and Fern\'andez Vilaboa  \cite[Proposition 4.2 and Theorem 4.5]{AV2} for the decomposition of $\BM(\C;H)$ and $\BM_{inn}(\C;H)$ in case $\C$ is a symmetric monoidal category and $H \in \C$ a finite commutative and cocommutative Hopf algebra. In this situation, the braiding is automatically $H$-linear, Proposition \ref{braid-lin}. They require that every $H$-Galois object is, in our terminology, faithfully projective \cite[Definition 2.3]{AV2} (instead of faithfully flat as we do). Both definitions coincide when $H$ is finite. From Proposition \ref{ChaseSweedler Thm}, an $H$-Galois object is faithfully projective. On the other hand, faithfully projective implies faithfully flat (2.5.1).

\begin{corollary}\label{seqsymm}
Let $\C$ be a closed symmetric monoidal category with equalizers and coequalizers and $H \in \C$ a finite commutative and cocommutative Hopf algebra. Then $\BM(\C;H)\iso \Br(\C)\times \Gal(\C;H)$ and $\BM_{inn}(\C;H)\iso \Br(\C)\times \Hc(\C;H,I).$
\end{corollary}

For $\C$ symmetric notice that $\BM_{inn}(\C;H)=\BM(\C;H)$ if and only if every $H$-Galois object has normal basis property.
\medskip

\begin{problem}
For a Hopf algebra $H$ in a closed braided monoidal category $\C$ such that the braiding $\Phi$ is $H$-linear we have a split exact sequence
$$
\bfig
\putmorphism(0, 0)(1, 0)[1`\Br(\C)`]{400}1a
\putmorphism(400, 0)(1, 0)[\phantom{\Br(\C)}`\BM(\C;H)`q]{600}1a
\putmorphism(410, -50)(1, 0)[\phantom{\Br(\C)}`\phantom{\Br(\C)}`p]{470}{-1}b
\putmorphism(1000, 0)(1, 0)[\phantom{\BM(\C;H)}`Coker(q)`\Pi]{700}1a
\putmorphism(1700, 0)(1, 0)[\phantom{Coker(q)}`1.`]{400}1a
\efig
$$
Under the additional hypothesis that $\C$ has equalizers and coequalizers, that $H$ is finite and commutative, and that $\Phi_{T,T'}\Phi_{T',T}=id_{T' \otimes T}$ for every two $H$-Galois objects $T$ and  $T'$ we have proved that $Coker(q) \cong \Gal(\C;H)$. Would it be possible to describe $Coker(q)$ without the symmetricity assumption on the braiding for $H$-Galois objects?
\end{problem}

\section{Applications}
\setcounter{equation}{0}

In this final section we will apply our two main results to a certain family of Radford biproduct (indeed bosonization) Hopf algebras including relevant examples of Hopf algebras like Sweedler Hopf algebra, Radford Hopf algebra, Nichols Hopf algebra and modified supergroup algebras. It is important to mention that the latter exhausts, up to Drinfel'd twist, all triangular Hopf algebras when the base field is algebraically closed of characteristic zero \cite{EG}. We will show that Beattie's exact sequence underlies in the computation of the Brauer group of these Hopf algebras, allowing thus to attain a complete understanding of it. A family of examples generalizing Radford Hopf algebra and Nichols Hopf algebra will be treated in detail to illustrate the different topics discussed in this paper. New computations of Brauer groups of quasi-triangular Hopf algebras are presented.

\subsection{Radford biproducts and Majid's bosonization}

In this subsection we mainly recollect known facts on the Radford biproduct, with the exception of Proposition \ref{Radf extends} and Corollary \ref{braid monoid iso Radf}, which may be of independent interest. \par \medskip

If $(H, \R)$ is a quasi-triangular Hopf algebra with bijective antipode over a field $K$, then a braiding for the category ${}_H\M$
is given by $\Phi_{\R}:=\tau\R$, where $\tau$ is the usual flip map. Conversely, if $({}_H\M, \Phi)$ is a  braided monoidal category, then $\R:=\tau\Phi(1_H\ot 1_H)$ defines a quasi-triangular structure on $H$. Concretely, if $\R=\R^{(1)}\ot\R^{(2)} \in H \otimes H$ is a quasi-triangular structure, the braiding $\Phi_{\R}$ for $M, N\in  {}_H\M$ and its inverse are given by
\begin{eqnarray} \label{R-braiding}
\Phi_{\R}(m\ot n)=\R^{(2)}\pH n\ot\R^{(1)}\pH m \vspace{3mm} \\
\Phi_{\R}^{-1}(n\ot m)=\R^{(1)}\pH m\ot S^{-1}(\R^{(2)})\pH n,
\end{eqnarray}
for $m\in M, n\in N$, where $\hspace{0,2cm}\pH\hspace{0,1cm}$ denotes the action of $H$ on $M,N.$ \par\medskip

Let $B$ be an algebra in ${}_H \M$ and a coalgebra in ${}^H\M$. We use the abbreviated form of Sweedler's notation for the comultiplication $\Delta_B$ and the coaction $\rho_B$ of $B$: given $b \in B$ we write $\Delta_B(b)=b_{(1)} \otimes b_{(2)} \in B \otimes B$ and $\rho_B(b)=b_{[-1]} \otimes b_{[0]} \in H \otimes B.$ We denote by $B\times H$ the space $B\ot H$ and the element $b\ot h$ in $B\times H$ is written as $b\times h$. The action of $H$ on  $B$ is denoted by $\tr$. We equip $B\times H$ with the following operations:
$$\begin{array}{rl}
 \textnormal{smash product:\index{smash product}} \label{Radf-bp}
\hskip-1em&(b\tm h)(b'\tm h')= b(h_{(1)}\tr b')\tm h_{(2)}h' \\
\textnormal{smash coproduct:\index{smash coproduct}} \quad &\Delta(b\tm h)=
(b_{(1)}\tm b_{{(2)}_{[-1]}} h_{(1)})\ot
 (b_{{(2)}_{[0]}}\tm h_{(2)}) \\
\textnormal{unit:} \quad &1_{B\times H}= 1_B\tm 1_H \\
\textnormal{counit:} \quad &\Epsilon_{B\times H}(b\tm h)=\Epsilon_B(b)\Epsilon_H(h)
\end{array}$$
for $b, b'\in B$ and $h, h'\in H$. Radford biproduct Theorem \cite[Theorem 2.1 and Proposition 2]{Rad1} characterizes when $B\ot H$ is a bialgebra and a Hopf algebra with the above operations. Majid observed that Radford biproduct construction is better understood in the framework of a certain braided monoidal category. With this observation, Radford biproduct Theorem states that $B\otimes H$ is a bialgebra with the above operations if and only if $B$ is a bialgebra in the braided monoidal category ${\ }^H_H\YD$ of left Yetter-Drinfel'd $H$-modules. This category is braided monoidal with braiding
$\Psi:M\ot N\to N\ot M$ and its inverse given by
$$\Psi(m\ot n)=m_{[-1]}\pH n\ot m_{[0]},\qquad \Psi^{-1}(n\ot m)=m_{[0]}\ot S^{-1}(m_{[-1]})\pH n,$$
for $m\in M, n\in N$ and $M, N\in{\ }^H_H\YD$. Moreover, if $B$ is a Hopf algebra in ${\ }^H_H\YD$, then $B \times H$ becomes a Hopf algebra with antipode
$$S(b\times h)=(1\times S_H(b^{[-1)]}h))(S_B(b^{[0]})\times 1).$$
\par\smallskip

For a quasi-triangular Hopf algebra $(H, \R)$ every left $H$-module $M$ belongs to
${\ }^H_H\YD$ with coaction $\lambda:M\to H\ot M$ given by
\begin{eqnarray} \label{leftcomMaj}
\lambda(m):=\R^{(2)}\ot \R^{(1)}\pH m, \quad m\in M,
\end{eqnarray}
and  $({}_H\M, \Phi_R)$ can be seen as a braided monoidal subcategory of $({\ }^H_H\YD, \Psi)$. As a particular case of Radford's Theorem, if $B$ is a Hopf algebra in ${}_H\M$, then $B\times H$ is a Hopf algebra. The process of obtaining an ordinary Hopf algebra $B\times H$ out of a Hopf algebra $B$ in ${}_H\M$ as above is called {\em bosonization} by Majid. The following proposition is a simplified version of \cite[Theorem 4.2]{Maj3}.

\begin{proposition} \label{cat-iso BxH}
Let $H$ be a quasi-triangular Hopf algebra and $B$ a Hopf algebra in ${}_H\M$.
Consider the category ${}_B({}_H\M)$
of $B$-modules in ${}_H\M$. Given $M\in {}_B({}_H\M)$ the compatibility condition is
$$h\pH(b\pB m)=(h_{(1)}\tr b)\pB(h_{(2)}\pH m)$$
for $h\in H, b\in B$ and $m\in M$. Then there is an isomorphism of monoidal
categories ${}_B({}_H\M)\iso {}_{B\times H}\M$.
\end{proposition}

This isomorphism is the identity on morphisms and on objects is defined as follows. We make a module $L$ in ${}_B({}_H\M)$ into a $B\times H$-module via $$(b\times h) \cdot l:=b\pB(h\pH l), \quad l \in L.$$
Conversely, on a $B\times H$-module $M$ we define a $B$- and an $H$-action via
$$b\pB m:= (b\times 1_H)\cdot m\quad\textnormal{and}\quad h\pH m:= (1_B\times h)\cdot m$$
respectively, for $b\in B, h\in H$ and $m\in M$. \par \medskip

Let $\iota:H \rightarrow B \times H, h \mapsto 1_B \times h$ denote the inclusion map, which is a Hopf algebra map.

\begin{proposition} \label{Radf extends}
Let $(H, \R)$ be a quasi-triangular Hopf algebra and $B$ a Hopf algebra in ${}_H\M$. Consider the Radford biproduct Hopf algebra $B\times H$. Then $\crta\R:=(\iota\ot\iota)(\R)$ is a quasi-triangular structure on $B\times H$ if and only if the braiding $\Phi_{\R}$ is $B$-linear in ${}_H\M$.
\end{proposition}

\begin{proof}
Assume that $\crta\R$ is a quasi-triangular structure for $B\times H$. Then the category $({}_{B\times H}\M, \Phi_{\crta\R})$ is  braided monoidal. In particular, $\Phi_{\crta\R}$ is $B\times H$-linear. Due to Proposition \ref{cat-iso BxH}
its corresponding map $\Phi_{\R}$ is $B$-linear in ${}_H\M$. Indeed, for $M, N\in{}_B({}_H\M)$ one has that the maps $\Phi_{\crta\R}:M\ot N\to N\ot M$ in ${}_{B\times H}\M$ and $\Phi_{\R}:M\ot N\to N\ot M$ in ${}_B({}_H\M)$ are equal. For, observing that $\crta\R=(1_B\times \R^{(1)})\ot(1_B\times \R^{(2)})$, we find that
$$\Phi_{\crta\R}(m\ot n)=(1_B\times \R^{(2)}) \cdot n\ot(1_B\times \R^{(1)}) \cdot m=\R^{(2)} \pH n\ot
\R^{(1)} \pH m=\Phi_{\R}(m\ot n)$$ for $m\in M$ and $n\in N$.

Conversely, in $({}_H\M, \Phi_{\R})$ suposse $\Phi_{\R}$ is $B$-linear. Similarly as above, by Proposition \ref{cat-iso BxH}, $\Phi_{\crta\R}$ is an isomorphism in ${}_{B\times H}\M$. Moreover, $\Phi_{\R}$ satisfies the two hexagon axioms for a braided monoidal category in ${}_H\M$ and it is $B$-linear. Then $\Phi_{\crta\R}$ satisfies the two hexagon axioms in ${}_{B\times H}\M$. Thus $({}_{B\times H}\M, \Phi_{\crta\R})$ is a braided monoidal category and then $\crta\R$ is a quasi-triangular structure  on $B\times H$. \qed
\end{proof}

\begin{corollary} \label{braid monoid iso Radf}
Let $(H, \R)$ be a quasi-triangular Hopf algebra so that $\R$ extends to a quasi-triangular
structure $\crta\R$ on the Radford biproduct $B\times H$. Then the braided monoidal categories
$({}_B({}_H\M), \Phi_{\R})$ and $({}_{B\times H}\M, \Phi_{\crta{\R}})$ are isomorphic.
\end{corollary}

\begin{proof}
The monoidal category ${}_{B\times H}\M$ is braided with braiding $\Phi_{\crta{\R}}.$  The other one is braided because of Proposition \ref{Radf extends} and Proposition \ref{H-C closed}. That the two braidings are equal we saw in the previous proposition.
\qed
\end{proof}

\begin{remark}{\em Under the hypothesis of the previous corollary, we know from Proposition \ref{braid-lin} that $B$ is cocommutative in ${}_H\M$. According to \cite[Corollary 5]{Sch}, there is a triangular Hopf algebra $(\tilde{H},\tilde{\R})$ and a surjective map $g:(H,\R) \rightarrow (\tilde{H},\tilde{\R})$ such that $B$ is a cocommutative Hopf algebra in ${}_{\tilde{H}}\M$ and $h\cdot b=g(h)b$ for all $h \in H, b \in B$. This makes more precise the description of the family of Radford biproducts we are regarding.}
\end{remark}

\subsection{Beattie's sequence as the root of the known computations}

For the quasi-triangular Hopf algebra $(H, \R)$ we will denote as usual the Brauer group of the category $({}_H\M,\Phi_{\R})$
by $\BM(K,H,\R).$ The next two results are consequences of our main theorems applied to the particular type of Radford biproduct Hopf algebras described before.

\begin{theorem}\label{appth1}
Let $(H, \R)$ be a quasi-triangular Hopf algebra and $B \in {}_H\M$ a Hopf algebra. Suppose that $\overline{\R}=(\iota\otimes \iota)(\R)$ is a quasi-triangular structure for $B \times H$. Assume that the braiding $\Phi_{\R}$ is symmetric on $T \otimes X$ for every $B$-Galois object $T \in {}_H\M$ and $X \in {}_H\M$. Then $\BM(K, B \times H ,\overline{\R}) \cong \BM(K,H,\R) \times \Gal(B;{}_H\M).$
\end{theorem}

After the group of lazy 2-cocycles was studied in \cite{BiCar}, and knowing the computations of Brauer groups done in \cite{VZ2}, \cite{CC1}, \cite{CC2} and \cite{Car2}, it was suspected that the group of lazy 2-cocycles would embed in the Brauer group, as mentioned in the introduction of \cite{BiCar}. In view of the following result, there is an embeding of the second braided cohomology group of the braided Hopf algebra into the Brauer group of the corresponding Radford biproduct. This cohomology group embeds in (and coincides in some cases with) the second lazy cohomology group for the examples analyzed, as verified in \cite{CF}, that clarifies this suspicion.

\begin{theorem}\label{appth2}
Let $(H, \R)$ be a quasi-triangular Hopf algebra and $B \in {}_H\M$ a Hopf algebra. Suppose that $\overline{\R}=(\iota\otimes \iota)(\R)$ is a quasi-triangular structure for $B \times H$. Then $\BM(K,H,\R) \times H^2({}_H\M; B,K)$ is a subgroup of $\BM(K,B \times H ,\overline{\R})$.
\end{theorem}

In the sequel we will show that Beattie's sequence underlies in the computations of the Brauer group of Sweedler Hopf algebra, Radford Hopf algebra, Nichols Hopf algebra and modified supergroup algebras carried out in the aforementioned papers. All these Hopf algebras are examples of the sort of Radford biproducts we are considering.  The strategy in the computation of the Brauer group of these Hopf algebras is to consider the group homomorphism $\iota^*:\BM(K,B \times H ,\overline{\R}) \rightarrow \BM(K,H,\R)$ induced by the inclusion map $\iota:H \rightarrow B \times H$. Then we have a split exact sequence
$$
\bfig
\putmorphism(-490, 0)(1, 0)[1`Ker(\iota^*)`]{450}1a
\putmorphism(-40, 0)(1, 0)[\phantom{Ker(\iota^*)}`\BM(K,B \times H ,\overline{\R})`]{940}1a
\putmorphism(900, 25)(1, 0)[\phantom{\BM(K,B \times H ,\overline{\R})}`\phantom{\BM(K,H,\R)}`\iota^*]{1100}1a
\putmorphism(1180, -25)(1, 0)[\phantom{\Br(\C)}`\phantom{\Br(\C)}`\pi^*]{630}{-1}b
\putmorphism(2000, 0)(1, 0)[\BM(K,H,\R)`1.`]{600}1a
\efig
$$
Being $B \times H$  a Radford biproduct, the map $\pi:B \times H \rightarrow H, b \times h \mapsto \varepsilon_B(b)h$ is a Hopf algebra map such that $\pi \iota=id_H$. Since $(\pi\otimes \pi)(\overline{\R})=\R$, the map $\pi$ induces a group homomorphism $\pi^*:\BM(K,H,\R) \rightarrow \BM(K,B \times H ,\overline{\R})$ such that $\iota^*\pi^*=id_{\BM(K,H,\R)}$. One next computes $Ker(\iota^*)$ by attaching to each class a certain invariant constructed from the fact that the action of $B\times H$ on a representative of the class, which is a matrix algebra, is inner (Skolem-Noether Theorem for Hopf algebras).

\begin{theorem}\label{Beattieroot}
With hypothesis as in the above paragraph and assuming that the braiding $\Phi_{\R}$ is symmetric on $T \otimes X$ for every $B$-Galois object $T \in {}_H\M$ and $X \in {}_H\M$, we have $Ker(\iota^*) \cong \Gal({}_H\M;B)$.
\end{theorem}

\begin{proof}
Put $\C={}_H\M$ and consider the isomorphism of categories $_B\C \cong {}_{B \times H}\M$ of Corollary \ref{braid monoid iso Radf}. Pick $M \in \C$ and consider it as a right $B \times H$-module via $\pi$. View it as an object in $_B\C$. The corresponding $B$-module structure of $M$ in $\C$ is the trivial one. Pick now $N \in {}_{B \times H}\M$ and equip it with the $H$-module structure via $\iota$. This is the same as forgetting the $B$-module structure of $N$ (as an object in $_B\C$). This shows that under the identification of categories $_B\C \cong {}_{B \times H}\M$ the morphisms $p$ and $q$ of Theorem \ref{Beattie} coincide with $\iota^*$ and $\pi^*$ respectively, that is, the following diagram is commutative:
$$
\bfig
\putmorphism(0,420)(0,-1)[``\iso]{400}1r
\putmorphism(-30,425)(1,0)[\Br(\C;B)`\Br(\C)`p]{1130}1a
\putmorphism(-100,0)(1,0)[\BM(K,B \times H ,\overline{\R})` \BM(K,H,\R)`
\iota^*]{1210}1a
\putmorphism(1100,420)(0,-1)[`` \iso]{400}1r
\putmorphism(1100,425)(1,0)[\phantom{\Br(\C)}`1`]{650}1a
\putmorphism(1100,0)(1,0)[\phantom{\BM(K,H,\R)}`1`]{650}1a
\efig
$$
From here it follows that $Ker(\iota^*)  \cong Ker(p)$ and $Ker(p)=\Gal({}_H\M;B)$ by Theorem \ref{Beattie}. \qed
\end{proof}

\subsection{Examples}

We finish this paper by analyzing a family of examples of Radford biproduct Hopf algebras that will illustrate our main theorems. They  will provide us with an example of a non-symmetric braided monoidal category and a braided Hopf algebra where the symmetricity assumption on the braiding for a Galois object and any other object is fulfilled. They will also give us an example of braided Hopf algebra possessing Galois objects without the normal basis property, and consequently an example where the subgroup of Azumaya algebras with inner action is proper. Using our sequence we will compute the Brauer group for each of these Hopf algebras,
recovering the computation for Radford Hopf algebra $H_{\nu}$ and Nichols Hopf algebras $E(n)$ carried out in \cite{CC1} and \cite{CC2} respectively. \par \bigskip

The following family of Hopf algebras appears in \cite[Section 4]{N}. Let $n,m$ be natural numbers, $K$ a field such that $char(K) \nmid 2m$ and $\omega$ a $2m$-th primitive root of unity. For $i=1,...,n$ pick $1 \leq d_i < 2m$ an odd number and set $d=(d_1,...,d_n)$. Consider the Hopf algebra
$$H(n,d)=K\langle g,x_1,...,x_n\vert g^{2m}=1, x_i^2=0, gx_i=\omega^{d_i} x_ig, x_ix_j=-x_jx_i \rangle,$$
where $g$ is group-like whereas $x_i$ is a $(g^m, 1)$-primitive element, that is, $\Delta(x_i)=1\ot x_i+x_i\ot g^m$. The antipode is given by $S(g)=g^{-1}$ and $S(x_i)=-x_ig^{m}$. View $K\Ent_{2m}$ as a Hopf subalgebra of $H(n,d)$. The quasi-triangular structures of $K\Ent_{2m}$ were classified in \cite[Page 219]{Rad3} and are of the form:
\begin{eqnarray} \label{qtr-Hnu}
\R_s=\frac{1}{2m}\Big(\sum_{j,t=0}^{2m-1}\omega^{-jt}g^j \otimes g^{st}\Big)
\end{eqnarray}
for $0\leq s< 2m$. It is easy to check that $\R_s$ is a quasi-triangular structure of $H(n,d)$ if and only if $sd_i \equiv m \ (mod.\ 2m)$ for every $i=1,...,n$. Moreover, $\R_s$ is triangular if and only if $s=m$. As particular instances of this family we obtain Radford Hopf algebra \cite[Page 252]{Rad3} ($n=1, m$ odd, $d_1=m$) and Nichols Hopf algebra ($m=1$). \par \smallskip

Set $d^{\leq n-1}=(d_1,...,d_{n-1})$. We will describe $H(n,d)$ as a Radford biproduct of $L=H(n-1,d^{\leq n-1})$ and the exterior algebra $B=K[x_n]/(x_n^2)$, which is an $L$-module algebra with action $g\cdot x_n=\omega^{d_n}x_n$ and $x_i \cdot x_n=0$ for $i=1,...,n-1$. It becomes a commutative and cocommutative Hopf algebra in $_L\M$ with coalgebra structure and antipode given as follows:
$$\Delta(x_n)=1\ot x_n+x_n\ot 1,\ \Epsilon(x_n)=0, \ S(x_n)=-x_n.$$
We now can consider the Hopf algebra $B \times L$ obtained by Majid's bosonization. It is not difficult to show that the map
$$\Psi:H(n,d)\to B \times H(n-1,d^{\leq n-1}), G \mapsto 1\times g, X_i \mapsto 1 \times x_i, X_n \mapsto x_n \times g^m$$ is a Hopf algebra isomorphism. Here we denote the generators of $H(n,d)$ by $G$ and $X_i$ instead of $g$ and $x_i$. That $\R_s$ is a quasi-triangular structure both on $H(n-1,d^{\leq n-1})$ and $H(n,d)$ means that $\R_s$, as a quasi-triangular structure on $H(n-1,d^{\leq n-1})$, extends to the quasi-triangular structure of $H(n,d)$. The extension $\crta{\R_s}=(\iota\ot \iota) (\R_s)$ lying in $(B\times L)\ot (B\times L)$ corresponds to $(\Psi^{-1}\ot\Psi^{-1})(\iota\ot \iota)(\R_s)=\R_s$ in $H(n,d)\ot H(n,d)$. By Theorem \ref{appth1} we will obtain the direct sum decomposition
$$\BM(K,H(n,d),\R_s)\iso \BM(K,L,\R_s) \times \Gal({}_L\M;B),$$
once we check that the symmetricity condition on the braiding is fulfilled. Our goal now is to describe the group of $B$-Galois objects in ${}_L\M$. For this we must distinguish two cases: $d_n=m$ and $d_n \neq m$. \par \smallskip

\begin{proposition}
Suppose that $d_n \neq m$ and set $I_n=\{1 \leq i <n : d_i \equiv -d_n \ (mod.\ 2m)\}.$ Then $\Gal({}_L\M;B) \cong
(K, +)^{\vert I_n \vert}.$
\end{proposition}

\begin{proof}
Let $\alpha=(\alpha_1,...,\alpha_{n-1}) \in K^{n-1}$.  Consider the algebra $C(\alpha)=K\langle w : w^2=0 \rangle$ with $L$-action and right $B$-comodule structure given by:
\begin{equation}\label{lbstruct}
\begin{array}{l}
g\cdot 1= 1 \hspace{1.65cm} \rho(1)=1 \otimes 1 \vspace{3pt} \\
g\cdot w=\omega^{d_n} w \hspace{0.85cm} \rho(w)=1 \otimes x_n +w \otimes 1 \vspace{3pt} \\
x_i \cdot w= \alpha_i
\end{array}
\end{equation}
Observe that $\omega^{d_n+d_i}\alpha_i=\omega^{d_n+d_i}x_i \cdot w = (\omega^{d_i}x_ig)\cdot w =(gx_i)\cdot w=\alpha_i.$
Then $C(\alpha)$ is an $L$-module if and only if $\alpha_i=0$ for $i \notin I_n$. Under this condition, it is easy to check that $C(\alpha)$ is a $B$-Galois object in ${}_L\M$. \par \medskip

We will next prove that any $B$-Galois object $A$ in ${}_L\M$ is of this form. Denote its comodule structure map by $\rho:A\to A\ot B$. From the isomorphism $can:A \otimes A \rightarrow A \otimes B$ we obtain $dim_K(A)=dim_K(B)=2$. We may choose a basis $\{1,y\}$ of $A$ such that $g\cdot 1=1$ and $g\cdot y=\omega^k y$ for $0 \leq k \leq 2m-1$. Write
$$\rho(y)=a 1 \otimes 1 +b 1 \otimes x_n +c y \otimes 1 +ey \otimes x_n$$
where $a,b,c,e \in K$. Since $A$ is a $B$-comodule, the counit axiom entails $a=0$ and $c=1$. Then,
\begin{equation}\label{rhoy}
\rho(y)=b 1 \otimes x_n +y \otimes 1 +e y \otimes x_n
\end{equation}
Using that $\rho$ is $L$-linear, from $\rho(g \cdot y)=g \cdot \rho(y)$ we get $b(\omega^{d_n}-\omega^k)=0$ and $e=0$. If $b=0$, then $\rho(y)=y \otimes 1$ and consequently $A^{co(B)}=A$, contradicting the fact $A^{co(B)}=K.$ Thus $b \neq 0$ and $d_n=k$. Write $z=b^{-1}y$. Then $\rho(z)=1 \otimes x_n+z \otimes 1$. We check that $z^2=0$. As $g\cdot y^2=\omega^{2k} y^2$, if $y^2 \neq 0$, then either $\omega^{2k}=1$ or $\omega^{2k}=\omega^k$. This would imply $d_n=0$ or $d_n=m$, a contradiction. We claim that, for $i=1,...,n-1$, it is $x_i \cdot z=\alpha_i$ for some $\alpha_i \in K$. If $x_i \cdot z=0$, it is clear. Suppose $x_i \cdot z \neq 0$. From the equality
$$\omega^{d_n+d_i}x_i \cdot z = (\omega^{d_i}x_ig)\cdot z =(gx_i)\cdot z=g \cdot (x_i \cdot z).$$
it follows $\omega^{d_n+d_i}=1$ or $\omega^{d_n+d_i}=\omega^{d_n}$. The latter case is not possible because $d_i$ is odd. Therefore $x_i \cdot z=\alpha_i$. Since $A$ is an $L$-module, $\alpha_i=0$ for $i \notin I_n$. Put $\alpha=(\alpha_1,...,\alpha_{n-1})$. Then $A$ is isomorphic to $C(\alpha)$ by mapping $1$ to $1$ and $z$ to $w$. \par \medskip

We have thus described the set $\Gal({}_L\M;B).$ To find its group structure we first verify the symmetricity condition of the braiding $\Phi$ on $C(\alpha) \otimes X$ for every $X \in {}_L\M.$ We may find a basis $\{v_l: l \in \Gamma\}$ of $X$ such that $g \cdot v_l=\omega^{\theta_l}v_l$ for $1 \leq \theta_l \leq 2m$. It is easy to check that
$$\begin{array}{lcl}
\Phi_{C(\alpha),X}(1 \otimes v_l)=v_l \otimes 1 & & \Phi_{X,C(\alpha)}(v_l \otimes 1)=1 \otimes v_l \vspace{3pt} \\
\Phi_{C(\alpha),X}(z \otimes v_l)= (-1)^{\theta_l}v_l \otimes z & & \Phi_{X,C(\alpha)}(v_l \otimes z)=(-1)^{\theta_l}z \otimes v_l
\end{array}$$
Hence $\Phi_{X,C(\alpha)}\Phi_{C(\alpha),X}=id_{C(\alpha) \otimes X}$. \par \medskip

Set $r=\vert I_n \vert$ and $I_n=\{i_1,...,i_r\}$ with $i_1<...<i_r$. For $\alpha$ as before, define  $\alpha_{I_n}=(\alpha_{i_1},...,\alpha_{i_r})$. We claim that the map  $\Omega:\Gal({}_L\M;B) \rightarrow (K, +)^r, [C(\alpha)] \mapsto \alpha_{I_n}$ is a group isomorphism. That $\Omega$ is well-defined and injective follows from the fact $C(\alpha) \cong C(\beta)$ as $B$-comodule algebras in ${}_L\M$ if and only if $\alpha=\beta$. For, let $\vartheta:C(\alpha)\to C(\beta)$ be such an isomorphism. Denote the algebra generators of $C(\alpha)$ and $C(\beta)$ by $y$ and $z$ respectively. Then $\vartheta(1)=1$ and $\vartheta(y)=\kappa z$ for some $\kappa\in K$ because $\vartheta$ is $g$-linear. Since $\vartheta$ is right $B$-colinear,
$(\vartheta \ot \Id_B)\rho_{\alpha}(y)=\rho_{\beta}\vartheta(y)$, where $\rho_{\alpha}$ and $\rho_{\beta}$ are the comodule structure maps of $C(\alpha)$ and $C(\beta)$ respectively. This yields $\kappa=1$. On the other hand, $\beta_i=x_i \cdot \vartheta(y)= \vartheta(x_i \cdot y)=\alpha_i$ for every $i=1,...,n-1$. For the injectivity, recall that $\alpha_i=\beta_i=0$ for $i \notin I_n$. \par \smallskip

For the surjectivity, take $\bar{\alpha}=(\bar{\alpha}_1,...,\bar{\alpha}_r) \in K^r$ and construct $\alpha \in K^{n-1}$ as follows: $\alpha_i=0$ if $i \notin I_n$ and $\alpha_i=\bar{\alpha}_j$ if $i=i_j$. Then $\Omega([C(\alpha)])=\bar{\alpha}.$ To verify  compatibility with the product, consider $C(\alpha)$ and $C(\beta)$. Let $y,z,w$ be the algebra generators of $C(\alpha), C(\beta)$ and $C(\alpha+\beta)$ respectively. We turn $C(\beta)$ into a left $B$-co\-mo\-du\-le via the braiding. If $\lambda_{\beta}$ denotes the left $B$-comodule structure morphism, then $\lambda_{\beta}(1)=1\ot 1$ and $\lambda_{\beta}(z)=1\otimes z+x_n\otimes 1$. It is easily checked that $\{1 \otimes 1, 1\ot z+y\ot 1\}$ is a $K$-basis of $C(\alpha)\Box_B C(\beta).$  Recall that we view $C(\alpha)\Box_B C(\beta)$ as a right $B$-comodule via $C(\alpha)\Box_B \rho_{\beta}$. It is also easy to show that $\theta:C(\alpha+\beta) \to C(\alpha)\Box_B C(\beta),$ defined by $\theta(1)=1\ot 1$ and $\theta(w)=1\ot z+y\ot 1,$ is a
$B$-comodule algebra morphism in ${}_L\M.$ Then $\Omega([C(\alpha)][C(\beta)])=\Omega([C(\alpha+\beta)])=(\alpha+\beta)_{I_n}=\alpha_{I_n}+\beta_{I_n}=
\Omega([C(\alpha)])+\Omega([C(\beta)])$. \qed
\end{proof}

\begin{proposition}
Suppose that $d_n=m$ and set $I_n=\{1 \leq i <n : d_i \equiv m \ (mod.\ 2m)\}.$ Then $\Gal({}_L\M;B) \cong
(K, +)^{\vert I_n \vert+1}.$
\end{proposition}

\begin{proof}
This proof follows the lines of the above one and we will just point out where the differences are. Let $a \in K$ and $\alpha=(\alpha_1,...,\alpha_{n-1}) \in K^{n-1}$ such that $\alpha_i=0$ for $i \notin I_n$. Consider the algebra
$C(a;\alpha)=K\langle w : w^2=a \rangle$ with $L$-action and right $B$-comodule structure given by (\ref{lbstruct}). One may easily check that $C(a;\alpha)$ is a $B$-Galois object in ${}_L\M$. Unlike the case $d_n \neq m$, notice that $g \cdot w^2=(g \cdot w)^2$ even for $a \neq 0$. To prove that any $B$-Galois object $A$ in ${}_L\M$ is of this form argue as in the previous proof but take into account that now can occur $y^2=a$, for some not necessarily zero $a \in K$, since $g \cdot y^2=(g \cdot y)^2=y^2$. \par \smallskip

The map $\Omega:\Gal({}_L\M;B) \rightarrow (K, +)^{r+1}$ is now defined by $[C(a;\alpha)] \mapsto (a,\alpha_{I_n})$. Note that
$C(a;\alpha) \cong C(b;\beta)$ as $B$-comodule algebras in ${}_L\M$ if and only if $a=b$ and $\alpha=\beta$. With notation as above, since $\vartheta(y)=z$, then $a=\vartheta(y^2)=z^2=b$. Proving that $\alpha=\beta$ is similar. Surjectivity is clear. For the compatibility with the product, observe that $(1\ot z+y\ot 1)^2=1 \otimes z^2-y \otimes z +y \otimes z+y^2 \otimes 1=(a+b)$. The part concerning $\alpha_{I_n}$ and $\beta_{I_n}$ is as before. \qed
\end{proof}

\begin{remark}{\em For $d_n=m$ any normalized $2$-cocycle $\sigma: B \otimes B \rightarrow K$ in ${}_L\M$ is of the form
$\sigma(1 \otimes 1)=1, \sigma(1 \otimes x_n)=\sigma(x_n \otimes 1)=0$ and $\sigma(x_n \otimes x_n)=a$ for some $a \in K$. Then $B_{\sigma}=C(a;0)$. Observe that every $B$-Galois object with the normal basis property is of the form $C(a;0).$ It is easy to see that $\Gal_{nb}({}_L\M;B) \cong (K,+)$ and hence $\BM_{inn}({}_L\M;B) \cong \Br({}_L\M)\times (K,+).$ If $n>1$, then $C(a;\alpha)$ has not the normal basis property for $\alpha \neq 0$. If $n=1$, then every $B$-Galois object has the normal basis property. \par

However, if $d_n \neq m$, there are no non-trivial normalized $2$-cocycles because $\sigma(x_n \otimes x_n)=g \cdot \sigma(x_n \otimes x_n)=\sigma(g \cdot (x_n \otimes x_n))=\omega^{2d_n}\sigma(x_n \otimes x_n).$ In this case the group of $B$-Galois objects with the normal basis property is trivial. This explains the differences between the two cases.}
\end{remark}
\medskip

Applying the previous propositions and recursion we are able to compute the Brauer group of $H(n,d)$ with respect to the quasi-triangular structure $\R_s$.

\begin{theorem}\label{bmnd}
For each $j=1,...,n$ set $I_j=\{1 \leq i < j: d_i \equiv -d_j \ (mod. \ 2m)\}$. Let $r_j=\vert I_j \vert$ if $d_j \neq m$ and $r_j=\vert I_j \vert+1$ otherwise. Then
$$\BM(K,H(n,d),\R_s) \cong \BM(K,K\Ent_{2m},\R_s) \times (K,+)^{r_1+...+r_n}.$$
\end{theorem}

Since $K\Ent_{2m}$ is self-dual, $\R_s$ may be seen as the bicharacter on $\Ent_{2m}$ given by $\theta_s(x,y)=\omega^{sxy}$ for all $x,y \in \Ent_{2m}$ and the Brauer group $\BM(K,K\Ent_{2m},\R_s)$ is that studied by Childs, Garfinkel and Orzech in \cite{CGO} and denoted by $\Bcgo_{\theta_s}(K,\Ent_{2m})$. For $m=1$, the group $\Bcgo_{\theta_1}(K,\Ent_{2})$ is nothing but the Brauer-Wall group $\BW(k)$. So this part in the above decomposition belongs to the Brauer-Long group theory and there are tools to compute it. \par \medskip

From our previous theorem we can recover the computation of the Brauer group of Radford Hopf algebra
$H_{m}$ ($n=1, m$ odd, $d_1=m$ in our notation) and Nichols Hopf algebra $E(n)$ ($m=1$) done in \cite{CC1} and
\cite{CC2} respectively. In the first case, $d_1=m$ and then $\BM(K,H_m,\R_s) \cong \Bcgo_{\theta_s}(K,\Ent_{2m}) \times (K,+).$ In the second case, $d_i=m$ and $r_i=i$ for $i=1,...,n$. Then $\BM(K,E(n),\R_1) \cong \BW(K) \times (K,+)^{\frac{n(n+1)}{2}}.$

\end{document}